\documentclass[12pt]{book}
\usepackage{geometry}  
\usepackage{amssymb, amsthm, mathrsfs, pstricks, pst-node, amsmath, longtable, fancyhdr, sectsty, dsfont, relsize, framed, import, appendix, mathabx, verbatim}
\usepackage{graphicx}
\usepackage[normalem]{ulem}
\usepackage[all]{xy}
\xyoption{color} 
\UseCrayolaColors 
\definecolor{dgreen}{rgb}{0,0.4,0}
\usepackage[titles]{tocloft}

\chapternumberfont{\scshape \bfseries} 
\chaptertitlefont{\bfseries}
\sectionfont{\bfseries} 
\subsectionfont{\bfseries} 

\makeatletter
\let\cleardouble@page\cleardoublepage
\AtBeginDocument{%
  \ifx\cleardouble@page\cleardoublepage
    \def\cleardoublepage{\clearpage
      {\pagestyle{empty}\cleardouble@page}}%
  \fi
} \makeatother 

\begin{document}

\newtheoremstyle{prp}
{}
{}
{\itshape}
{}
{\bfseries\upshape}
{.}
{5pt}
{}

\newtheoremstyle{algor}
{}
{}
{\upshape\sffamily}
{}
{\bfseries\upshape}
{}
{5pt}
{}

\newtheoremstyle{def}
{}
{}
{\upshape}
{}
{\bfseries\upshape}
{.}
{5pt}
{}

\newtheoremstyle{bew}
{}
{}
{\upshape}
{}
{\itshape}
{}
{5pt}
{}

\newtheoremstyle{rmk}
{}
{}
{\upshape}
{-5pt}
{\bfseries \scshape}
{.}
{5pt}
{}

\newtheoremstyle{sctext} 
{}
{}
{\upshape}
{0pt}
{\bfseries \scshape} 
{.} 
{5pt} 
{} 

\theoremstyle{prp}
\newtheorem{thm}{Theorem}[section]

\theoremstyle{prp}
\newtheorem{propn}[thm]{Proposition}

\theoremstyle{prp}
\newtheorem{corl}[thm]{Corollary}

\theoremstyle{prp}
\newtheorem{lemm}[thm]{Lemma}

\theoremstyle{def}
\newtheorem{defn}[thm]{Definition}

\theoremstyle{def}
\newtheorem{remind}[thm]{Reminder}

\theoremstyle{def}
\newtheorem{condition}[thm]{Condition}

\theoremstyle{def}
\newtheorem{constr}[thm]{Construction}

\newtheorem{stz}{Satz}[section]

\theoremstyle{def}
\newtheorem{beispiel}[stz]{Beispiel}

\theoremstyle{def}
\newtheorem{ex}[thm]{Example}

\theoremstyle{def}
\newtheorem{notat}[thm]{Notation}

\theoremstyle{def}
\newtheorem{conven}[thm]{Convention}

\theoremstyle{def}
\newtheorem{konstruktion}[stz]{Konstruktion}

\theoremstyle{def}
\newtheorem{bemerkung}[stz]{Bemerkung}

\theoremstyle{bew}
\newtheorem*{beweisa}{Beweis.}

\theoremstyle{def}
\newtheorem{remk}[thm]{Remark}

\theoremstyle{sctext}
\newtheorem*{prf}{Proof}

\newtheorem{korl}[stz]{Korollar}
\newtheorem{verm}[stz]{Vermutung}
\newtheorem{beisp}[stz]{Beispiel}
\newtheorem{eg}[thm]{Example}

\newenvironment{satz}{\begin{stz}}{\end{stz}}
\newenvironment{theorem}{\begin{thm}}{\end{thm}}
\newenvironment{proposition}{\begin{propn}}{\end{propn}}
\newenvironment{lemma}{\begin{lemm}}{\end{lemm}}
\newenvironment{korollar}{\begin{korl}}{\end{korl}}
\newenvironment{corollary}{\begin{corl}}{\end{corl}}
\newenvironment{vermutung}{\begin{verm}}{\end{verm}}
\newenvironment{example}{\begin{ex}}{\end{ex}}
\newenvironment{notation}{\begin{notat}}{\mbtu \end{notat}}
\newenvironment{convention}{\begin{conven}}{\mbtu \end{conven}}
\newenvironment{beweis}{\begin{beweisa}}{\qed \end{beweisa} }

\newenvironment{prooff}{\begin{prf}}{\qed \end{prf}}

\newenvironment{proofof}[1]{\leavevmode \newline \textsc{Proof}{ (of #1)}\textsc{.}}{\qed }
\newenvironment{proofoftwo}[1]{\leavevmode \newline \emph{Proof}{ (of #1)}\textsc{.}}{\qed }
\newenvironment{construction}{\begin{constr}}{\mbtu \end{constr}}
\newenvironment{definition}{\begin{defn}}{\mbtu \end{defn}}
\newenvironment{remark}{\begin{remk}}{\mbtu \end{remk}}

\def\mbtu {{
\parfillskip=0pt 
\widowpenalty=10000 
\displaywidowpenalty=10000 
\finalhyphendemerits=0 
%
\leavevmode 
\unskip 
\nobreak 
\hfil 
\penalty50 
\hskip.2em 
\null 
\hfill 
$\tu$
%
\par}} 

\def\mbtd {{
\parfillskip=0pt 
\widowpenalty=10000 
\displaywidowpenalty=10000 
\finalhyphendemerits=0 
%
\leavevmode 
\unskip 
\nobreak 
\hfil 
\penalty50 
\hskip.2em 
\null 
\hfill 
$\btd$
%
\par}} 

\newcommand{\tick}{\ding{52}} 
\newcommand{\sseq}{\subseteq}
\newcommand{\ot}{\otimes}
\newcommand{\ts}{\textsection}
\newcommand{\tts}{\textsection \textsection}
\newcommand{\Gr}{\textnormal{Gr}}
\newcommand{\sqk}{\tn{sq}_{k}}
\newcommand{\cosqk}{\tn{cosq}_{k}}
\newcommand{\sq}{\tn{sq}}
\newcommand{\cosq}{\tn{cosq}}
\newcommand{\qaq}{\qquad \tn{and} \qquad}
\newcommand{\qbq}{\qquad \tn{but} \qquad}
\newcommand{\qfq}{\quad \tn{for all} \quad}
\newcommand{\quq}{\qquad \tn{und} \qquad}
\newcommand{\nt}{\natural}
\newcommand{\shp}{\sharp}
\newcommand{\Biext}{\mathbf{Biext}}
\newcommand{\Biexz}{\tn{Biext}^{0}}
\newcommand{\Biexo}{\tn{Biext}^{1}}
\newcommand{\Biex}{\tn{Biext}}
\newcommand{\bBiex}{\ul{\tn{Biext}}}
\newcommand{\cl}{\tn{cl}}
\newcommand{\lc}{\tn{lc}}
\newcommand{\cts}{\tn{cts}}
\newcommand{\lcXA}{\tn{lc}_{X, A}}
\newcommand{\ctsXA}{\tn{cts}_{X, A}}
\newcommand{\x}{\times}
\newcommand{\Cont}{\tn{Cont}}
\newcommand{\Gb}{\mathbb{G}}
\newcommand{\Ga}{\Gb_{a}}
\newcommand{\Gm}{\Gb_{m}}
\newcommand{\wtGm}{\wt{\Gb}_{m}}
\newcommand{\Gmk}{\mathbb{G}_{m, k}}
\newcommand{\Gmkb}{\mathbb{G}_{m, \cj{k}}}
\newcommand{\GmX}{\mathbb{G}_{m, X}}
\newcommand{\GmA}{\mathbb{G}_{m, A}}
\newcommand{\Gmak}{\mathbb{G}_{m, \cj{k}}}
\newcommand{\uGm}{\ul{\Gm}}
\newcommand{\Gmr}{\mathbb{G}_{m}^{r}}
\newcommand{\GmAr}{\mathbb{G}_{m, A}^{r}}
\newcommand{\Gmkr}{\mathbb{G}_{m, k}^{r}}
\newcommand{\Gmakr}{\mathbb{G}_{m, \cj{k}}^{r}}
\newcommand{\Ao}{\A^{1}}
\newcommand{\Aor}{(\A^{1})^{r}}
\newcommand{\Aon}{(\A^{1})^{n}}

\newcommand{\Aok}{\A^{1}_{k}}
\newcommand{\cosec}{\tn{cosec} \,}
\newcommand{\DS}{\displaystyle}
\newcommand{\Gal}{\textnormal{Gal}}
\newcommand{\Ell}{\textnormal{Ell}}
\newcommand{\Emb}{\textnormal{Emb}}
\newcommand{\Opt}{\textnormal{Opt}}
\newcommand{\Nrm}{\textnormal{N}}
\newcommand{\Tr}{\textnormal{Tr}}
\newcommand{\adj}{\textnormal{adj}}
\newcommand{\cof}{\tn{cof}}
\newcommand{\Ta}{\textnormal{T}}
\newcommand{\Hpq}{H^{p,q}}
\newcommand{\Fr}{\textnormal{Fr}}
\newcommand{\Xb}{X_{\bl}}
\newcommand{\Yb}{Y_{\bl}}
\newcommand{\tr}{\textnormal{tr}}
\newcommand{\Hilb}{\textnormal{Hilb}}
\newcommand{\Ideals}{\textnormal{Ideals}}
\newcommand{\Fix}{\textnormal{Fix}}
\newcommand{\Ann}{\textnormal{Ann}}
\newcommand{\Nil}{\textnormal{Nil}}
\newcommand{\nil}{\textnormal{nil}}
\newcommand{\Disc}{\textnormal{Disc}}
\newcommand{\Orb}{\textnormal{Orb}}
\newcommand{\Mat}{\textnormal{Mat}}
\newcommand{\Mn}{\textnormal{M}_{n}}
\newcommand{\Mm}{\textnormal{M}_{m}}
\newcommand{\Mmn}{\textnormal{M}_{m \x n}}
\newcommand{\MmR}{\Mm(\R)}
\newcommand{\MmnR}{\Mmn(\R)}
\newcommand{\MnR}{\textnormal{M}_{n}(\R)}
\newcommand{\MnK}{\textnormal{M}_{n}(K)}
\newcommand{\GL}{\textnormal{GL}}
\newcommand{\PGL}{\tn{PGL}}
\newcommand{\PSL}{\tn{PSL}}
\newcommand{\PGLt}{\tn{PGL}_{2}}
\newcommand{\PSLt}{\tn{PSL}_{2}}
\newcommand{\GLn}{\GL_{n}}
\newcommand{\GLnK}{\GL_{n}(K)}
\newcommand{\SLnK}{\SL_{n}(K)}
\newcommand{\GLt}{\GL_{2}}
\newcommand{\SL}{\textnormal{SL}}
\newcommand{\StL}{\tn{S}^{*}\tn{L}}
\newcommand{\PStL}{\tn{PS}^{*}\tn{L}}
\newcommand{\StLt}{\StL_{2}}
\newcommand{\PStLt}{\PStL_{2}}
\newcommand{\SLn}{\textnormal{SL}_{n}}
\newcommand{\SLt}{\SL_{2}}
\newcommand{\Un}{\mathbb{U}_{n}}
\newcommand{\Tt}{\mathbb{T}}
\newcommand{\Tn}{\mathbb{T}_{n}}
\newcommand{\Sp}{\textnormal{Sp}}
\newcommand{\GSp}{\textnormal{GSp}}
\newcommand{\SU}{\textnormal{SU}}
\newcommand{\SO}{\textnormal{SO}}
\newcommand{\Lf}{\mathbb{L}}
\newcommand{\Lfd}{\Lf^{d}}
\newcommand{\Lfm}{\Lf^{m}}
\newcommand{\Rd}{\mathbb{R}}
\newcommand{\dg}{\dagger}
\newcommand{\id}{\mathrm{id}}
\newcommand{\mult}{\tn{mult}}
\newcommand{\chr}{\mathrm{char}}
\newcommand{\lcm}{\mathrm{lcm}}
\newcommand{\ggt}{\mathrm{ggt}}
\newcommand{\Lie}{\tn{Lie}}

\newcommand{\Frac}{\mathrm{Frac \:}}
\newcommand{\Aut}{\textnormal{Aut}}
\newcommand{\End}{\textnormal{End}}
\newcommand{\Sym}{\textnormal{Sym}}
\newcommand{\SymX}{\Sym_{X}}
\newcommand{\AltX}{\Alt_{X}}
\newcommand{\MX}{M_{X}}
\newcommand{\MtX}{\wt{M}_{X}}
\newcommand{\Alt}{\tn{Alt}}
\newcommand{\Altn}{\tn{Alt}^{n}}
\newcommand{\Symn}{\Sym^{n}}
\newcommand{\Symnp}{\Sym^{n}_{\tn{pre}}}
\newcommand{\Altd}{\tn{Alt}^{d}}
\newcommand{\Symd}{\Sym^{d}}
\newcommand{\Symi}{\Sym^{i}}
\newcommand{\Altk}{\Alt^{k}}
\newcommand{\Alti}{\Alt^{i}}
\newcommand{\Symk}{\Sym^{k}}

\newcommand{\Symm}{\Sym^{m}}
\newcommand{\Drg}{\textnormal{Drg}}
\newcommand{\sgn}{\textnormal{sgn}}
\newcommand{\Jac}{\tn{Jac}}
\newcommand{\JacC}{\Jac \, C}
\newcommand{\uJac}{\ul{\tn{Jac}}}
\newcommand{\Alb}{\tn{Alb}}
\newcommand{\AlbqS}{\tn{Alb}_{?/S}}
\newcommand{\Albz}{\tn{Alb}^{\circ}}
\newcommand{\Albo}{\tn{Alb}^{1}}
\newcommand{\uAlb}{\ul{\tn{Alb}}}
\newcommand{\uAlbz}{\uAlb^{\circ}}

\newcommand{\wAlb}{\wt{\tn{Alb}}}
\newcommand{\wAlbz}{\wt{\tn{Alb}^{\circ}}}

\newcommand{\uA}{\ul{A}}
\newcommand{\Albm}{\Alb^{-}}
\newcommand{\Albp}{\Alb^{+}}
\newcommand{\rd}{\tn{red}}
\newcommand{\Funct}{\textnormal{Funct}}
\newcommand{\Fun}{\tn{\tbf{Fun}}}
\newcommand{\Nat}{\textnormal{Nat}}
\newcommand{\Pic}{\tn{Pic}}
\newcommand{\Picp}{\tn{Pic}^{+}}
\newcommand{\Picm}{\tn{Pic}^{-}}
\newcommand{\PicXS}{\tn{Pic}_{X/S}}
\newcommand{\PicoXS}{\tn{Pic}^{\circ}_{X/S}}
\newcommand{\PicoCS}{\tn{Pic}^{\circ}_{C/S}}
\newcommand{\PicoCk}{\tn{Pic}^{\circ}_{C/k}}
\newcommand{\PicoAk}{\Pico_{A/k}}
\newcommand{\PicAk}{\Pic_{A/k}}

\newcommand{\PicXs}{\tn{Pic}_{X_{s}/k_{s}}}
\newcommand{\PicoXs}{\tn{Pic}^{\circ}_{X_{s}/k_{s}}}
\newcommand{\AlbXS}{\tn{Alb}_{X/S}}
\newcommand{\AlbYS}{\tn{Alb}_{Y/S}}
\newcommand{\AlbzXS}{\tn{Alb}^{\circ}_{X/S}}
\newcommand{\AlbzYS}{\tn{Alb}^{\circ}_{Y/S}}
\newcommand{\AlbzfS}{\tn{Alb}^{\circ}_{f/S}}

\newcommand{\AlbXk}{\tn{Alb}_{X/k}}
\newcommand{\AlbAk}{\tn{Alb}_{A/k}}
\newcommand{\AlbGk}{\tn{Alb}_{G/k}}
\newcommand{\AlbYk}{\tn{Alb}_{Y/k}}
\newcommand{\AlbzXk}{\tn{Alb}^{\circ}_{X/k}}
\newcommand{\AlbzAk}{\tn{Alb}^{\circ}_{A/k}}
\newcommand{\AlboXk}{\tn{Alb}^{1}_{X/k}}
\newcommand{\AlbnXk}{\tn{Alb}^{n}_{X/k}}
\newcommand{\AlbCk}{\tn{Alb}_{C/k}}
\newcommand{\AlbzCk}{\tn{Alb}^{\circ}_{C/k}}

\newcommand{\AlboXS}{\tn{Alb}^{1}_{X/S}}
\newcommand{\AlbiXS}{\tn{Alb}^{i}_{X/S}}
\newcommand{\AlbXs}{\tn{Alb}_{X_{s}/k_{s}}}
\newcommand{\AlbzXs}{\tn{Alb}^{\circ}_{X_{s}/k_{s}}}
\newcommand{\AzXS}{A^{\circ}_{X/S}}

\newcommand{\fPicXS}{\tn{\tbf{Pic}}_{X/S}}
\newcommand{\fPicoXS}{\tn{\tbf{Pic}}_{X/S}}
\newcommand{\fAlbXS}{\tn{\tbf{Alb}}_{X/S}}
\newcommand{\fAlbzXS}{\tn{\tbf{Alb}}^{\circ}_{X/S}}
\newcommand{\fAlbXs}{\tn{\tbf{Alb}}_{X_{s}/k_{s}}}
\newcommand{\fAlbzXs}{\tn{\tbf{Alb}}^{\circ}_{X_{s}/k_{s}}}
\newcommand{\fPicXSt}{\tn{\tbf{Pic}}_{(X/S) \; (\tau)}}
\newcommand{\fPicXSs}{\tn{\tbf{Pic}}_{(X/S) \; (\sigma)}}
\newcommand{\fPicXSfp}{\tn{\tbf{Pic}}_{(X/S) \; (\fppf)}}
\newcommand{\fPicXSfl}{\tn{\tbf{Pic}}_{(X/S) \; (\fl)}}
\newcommand{\fPicXSet}{\tn{\tbf{Pic}}_{(X/S) \; (\et)}}
\newcommand{\fPicXSZ}{\tn{\tbf{Pic}}_{(X/S) \; (\Zar)}}

\newcommand{\uPicXS}{\ul{\tn{Pic}}_{X/S}}
\newcommand{\uPicoXS}{\ul{\tn{Pic}}^{\circ}_{X/S}}
\newcommand{\uAlbXS}{\ul{\tn{Alb}}_{X/S}}
\newcommand{\uAlbXk}{\ul{\tn{Alb}}_{X/k}}
\newcommand{\uAlbAk}{\ul{\tn{Alb}}_{A/k}}
\newcommand{\uAlbCk}{\ul{\tn{Alb}}_{C/k}}
\newcommand{\uAlbYS}{\ul{\tn{Alb}}_{Y/S}}
\newcommand{\uAlbYk}{\ul{\tn{Alb}}_{Y/k}}
\newcommand{\uAlbzXS}{\ul{\tn{Alb}}^{\circ}_{X/S}}
\newcommand{\uAlbzXk}{\ul{\tn{Alb}}^{\circ}_{X/k}}
\newcommand{\uAlbzCk}{\ul{\tn{Alb}}^{\circ}_{C/k}}
\newcommand{\uAlbzYS}{\ul{\tn{Alb}}^{\circ}_{Y/S}}
\newcommand{\wAlbXS}{\wAlb_{X/S}}
\newcommand{\wAlbXk}{\wAlb_{X/k}}
\newcommand{\wAlbAk}{\wAlb_{A/k}}
\newcommand{\wAlbGk}{\wAlb_{G/k}}
\newcommand{\wAlbCk}{\wAlb_{C/k}}
\newcommand{\wAlbYS}{\wAlb_{Y/S}}
\newcommand{\wAlbYk}{\wAlb_{Y/k}}
\newcommand{\wAlbzXS}{\wAlbz_{X/S}}
\newcommand{\wAlbzXk}{\wAlbz_{X/k}}
\newcommand{\wAlbzCk}{\wAlbz_{C/k}}
\newcommand{\wAlbzYS}{\wAlbz_{Y/S}}

\newcommand{\bC}{b_{C}}
\newcommand{\bCo}{b_{C}^{\circ}}

\newcommand{\uAlbXs}{\ul{\tn{Alb}}_{X_{s}/k_{s}}}
\newcommand{\uAlbzXs}{\ul{\tn{Alb}}^{\circ}_{X_{s}/k_{s}}}
\newcommand{\uPicXSt}{\ul{\tn{Pic}}_{(X/S) \; (\tau)}}
\newcommand{\uPicXSs}{\ul{\tn{Pic}}_{(X/S) \; (\sigma)}}
\newcommand{\uPicXSfp}{\ul{\tn{Pic}}_{(X/S) \; (\fppf)}}
\newcommand{\uPicXSfl}{\ul{\tn{Pic}}_{(X/S) \; (\fl)}}
\newcommand{\uPicXSet}{\ul{\tn{Pic}}_{(X/S) \; (\et)}}
\newcommand{\uPicXSZ}{\ul{\tn{Pic}}_{(X/S) \; (\Zar)}}

\newcommand{\bPic}{\tn{\tbf{Pic}}} 
\newcommand{\bPicX}{\bPic_{X}} 
\newcommand{\bPicp}{\tn{\tbf{Pic}}^{+}}
\newcommand{\bPicm}{\tn{\tbf{Pic}}^{-}}
\newcommand{\bPico}{\tn{\tbf{Pic}}^{\circ}}
\newcommand{\uPic}{\ul{\tn{Pic}}}
\newcommand{\uPico}{\ul{\tn{Pic}}^{\circ}}
\newcommand{\Pico}{\tn{Pic}^{\circ}}
\newcommand{\wPic}{\wt{\tn{Pic}}}
\newcommand{\wPico}{\wt{\tn{Pic}^{\circ}}}
\newcommand{\Ca}{\textnormal{Ca}}
\newcommand{\CaCl}{\textnormal{CaCl}}
\newcommand{\Cl}{\textnormal{Cl}}
\newcommand{\Div}{\tn{Div}}
\newcommand{\wDiv}{\wt{\Div}}
\newcommand{\DivXk}{\Div_{X/k}}
\newcommand{\DivCk}{\Div_{C/k}}
\newcommand{\wDivXk}{\wDiv_{X/k}}
\newcommand{\wDivCk}{\wDiv_{C/k}}
\newcommand{\Divo}{\tn{Div}^{\circ}}
\newcommand{\Princ}{\textnormal{Princ}}
\newcommand{\divr}{\textnormal{div}}
\newcommand{\NS}{\tn{NS}}
\newcommand{\Th}{\textnormal{Th}}
\newcommand{\ThTriv}{\textnormal{Th}_{\textnormal{triv}}}
\newcommand{\Rie}{\textnormal{Rie}}
\newcommand{\Herm}{\textnormal{Herm}}
\newcommand{\inv}{\textnormal{inv}}
\newcommand{\diag}{\textnormal{diag}}
\newcommand{\grcd}{\textnormal{gcd}}
\newcommand{\eval}{\textnormal{eval}}
\newcommand{\Supp}{\textnormal{Supp}}
\newcommand{\N}{\mathbb{N}}
\newcommand{\Z}{\mathbb{Z}}
\newcommand{\Zm}{\mathbb{Z}_{m}}
\newcommand{\Zmx}{\mathbb{Z}_{m}^{\x}}
\newcommand{\uZ}{\ul{\Z}}
\newcommand{\Znn}{\mathbb{Z}_{\ge 0}}
\newcommand{\Zpn}{\mathbb{Z}^{n}}
\newcommand{\Zr}{\mathbb{Z}^{r}}
\newcommand{\Zrn}{(\mathbb{Z}^{r})^{n}}
\newcommand{\Zcs}{\cs{\Z}}
\newcommand{\wtZcs}{\wt{\cs{\Z}}}
\newcommand{\Q}{\mathbb{Q}}
\newcommand{\Ql}{\mathbb{Q}_{\ell}}
\newcommand{\R}{\mathbb{R}}
\newcommand{\C}{\mathbb{C}}
\newcommand{\K}{\mathbb{K}}
\newcommand{\F}{\mathbb{F}}
\newcommand{\Uy}{\mathbb{U}}
\newcommand{\Rn}{\R^{n}}
\newcommand{\Rfour}{\R^{4}}
\newcommand{\Rthree}{\R^{3}}
\newcommand{\Rtwo}{\R^{2}}
\newcommand{\Cfour}{\C^{4}}
\newcommand{\Cthree}{\C^{3}}
\newcommand{\Ctwo}{\C^{2}}
\newcommand{\RN}{\R^{\N}}
\newcommand{\Rm}{\R^{m}}
\newcommand{\Cn}{\C^{n}}
\newcommand{\Kn}{\K^{n}}
\newcommand{\Fn}{\F^{n}}
\newcommand{\Fib}{\tn{Fib}}

\newcommand{\uX}{\ul{X}}
\newcommand{\wtX}{\wt{X}}
\newcommand{\wtG}{\wt{G}}
\newcommand{\wtH}{\wt{H}}
\newcommand{\wtf}{\wt{f}}
\newcommand{\uC}{\ul{C}}
\newcommand{\uY}{\ul{Y}}
\newcommand{\uGt}{\ul{G}}
\newcommand{\uH}{\ul{H}}
\newcommand{\poss}{\tn{pos}_{*}}

\newcommand{\Hyn}{\Hy^{n}}
\newcommand{\Hyz}{\Hy_{0}}
\newcommand{\HynR}{\Hy^{n}_{R}}
\newcommand{\Bn}{\B^{n}}
\newcommand{\BnR}{\B^{n}_{R}}
\newcommand{\Uyn}{\Uy^{n}}
\newcommand{\UynR}{\Uy^{n}_{R}}

\newcommand{\Hyt}{\Hy^{2}}
\newcommand{\HytR}{\Hy^{2}_{R}}
\newcommand{\Bt}{\B^{2}}
\newcommand{\BtR}{\B^{2}_{R}}
\newcommand{\BytR}{\BtR}
\newcommand{\Uyt}{\Uy^{2}}
\newcommand{\UytR}{\Uy^{2}_{R}}
\newcommand{\Sn}{\tn{S}^{n}}
\newcommand{\Sd}{\tn{S}^{d}}
\newcommand{\Se}{\tn{S}^{e}}
\newcommand{\Sde}{\tn{S}^{de}}
\newcommand{\Sinf}{\tn{S}^{\infty}}
\newcommand{\kba}{\mathbf{k}}
\newcommand{\X}{\mathbb{X}}
\newcommand{\Fp}{\mathbb{F}_{p}}
\newcommand{\Fq}{\mathbb{F}_{q}}
\newcommand{\Fpx}{\mathbb{F}_{p}^{\times}}
\newcommand{\Fqx}{\mathbb{F}_{q}^{\times}}
\newcommand{\A}{\mathbb{A}}
\newcommand{\An}{\mathbb{A}^{n}}
\newcommand{\Ano}{\mathbb{A}^{n+1}}
\newcommand{\Hy}{\mathbb{H}}
\newcommand{\wg}{\wedge}
\newcommand{\Zyk}{\textnormal{Zyk}}
\newcommand{\catc}{\shf{C}}
\newcommand{\catd}{\shf{D}}
\newcommand{\cata}{\shf{A}}
\newcommand{\catb}{\shf{B}}
\newcommand{\catt}{\shf{T}}
\newcommand{\cats}{\shf{S}}
\newcommand{\catq}{\shf{Q}}
\newcommand{\shfF}{\shf{F}}
\newcommand{\shfFx}{\shf{F}_{x}}
\newcommand{\shfFxg}{\shf{F}_{\cj{x}}}
\newcommand{\xg}{\cj{x}}
\newcommand{\shfO}{\shf{O}}
\newcommand{\shfOx}{\shf{O}^{\x}}
\newcommand{\shfL}{\shf{L}}
\newcommand{\shfV}{\shf{V}}
\newcommand{\shfR}{\shf{R}}
\newcommand{\shfLi}{\shf{L}_{i}}
\newcommand{\shfLix}{\shf{L}_{i}^{\x}}
\newcommand{\shfLj}{\shf{L}_{j}}
\newcommand{\shfLjx}{\shf{L}_{j}^{\x}}

\newcommand{\shfM}{\shf{M}}
\newcommand{\shfN}{\shf{N}}
\newcommand{\clyOD}{\cly{O}\cly{D}}

\newcommand{\shfA}{\shf{A}}
\newcommand{\shfB}{\shf{B}}
\newcommand{\shfC}{\shf{C}}
\newcommand{\shfD}{\shf{D}}
\newcommand{\shfE}{\shf{E}}
\newcommand{\shfG}{\shf{G}}
\newcommand{\shfGL}{\shfG(\shfL)}
\newcommand{\shfKL}{\shfK(\shfL)}
\newcommand{\shfH}{\shf{H}}
\newcommand{\shfK}{\shf{K}}
\newcommand{\shfP}{\shf{P}}
\newcommand{\shfQ}{\shf{Q}}
\newcommand{\shfT}{\shf{T}}
\newcommand{\shfS}{\shf{S}}
\newcommand{\shfX}{\shf{X}}
\newcommand{\shfSk}{\shf{S}^{k}}
\newcommand{\sSk}{\idl{S}^{k}}
\newcommand{\idlp}{\idl{p}}
\newcommand{\idlq}{\idl{q}}
\newcommand{\idla}{\idl{a}}
\newcommand{\idlb}{\idl{b}}

\newcommand{\shfI}{\shf{I}}
\newcommand{\sOX}{\shf{O}_{X}}
\newcommand{\sOY}{\shf{O}_{Y}}
\newcommand{\Xet}{X_{\et}}

\newcommand{\Deltan}{\Delta^{n}}
\newcommand{\Deltaz}{\Delta^{0}}
\newcommand{\Deltao}{\Delta^{1}}
\newcommand{\lieg}{\idl{g}}
\newcommand{\Vect}{\tbf{Vec}}
\newcommand{\Vectk}{\tbf{Vec}_{k}}
\newcommand{\VectQ}{\tbf{Vec}_{\Q}}
\newcommand{\VectQi}{\tbf{Vec}_{\Q}^{\infty}}
\newcommand{\VecR}{\tbf{Vec}_{R}}
\newcommand{\RepkG}{\tbf{Rep}_{k}(G)}
\newcommand{\Rep}{ \tbf{Rep}}
\newcommand{\Repn}{ \tbf{Rep}^{\tn{nil}}}
\newcommand{\ModR}{ \tbf{Mod}_{R}}
\newcommand{\Mod}{ \tbf{Mod}}
\newcommand{\fggrRmod}{ \tbf{fggrRmod}}
\newcommand{\ProjR}{\tbf{Proj}_{R}}
\newcommand{\ds}{\displaystyle}
\newcommand{\Ha}{\mathbf{H}}
\newcommand{\HA}{\mathbf{H}_{\shf{A}}}
\newcommand{\HaB}{\mathbf{H}_{\shf{B}}}
\newcommand{\Htorp}{\mathbf{H}_{\Torp}}
\newcommand{\Hdr}{H_{\tn{dR}}}
\newcommand{\Het}{H_{\tn{\'et}}}
\newcommand{\HB}{H_{\tn{B}}}
\newcommand{\Hcr}{H_{\tn{crys}}}
\newcommand{\Hs}{H_{\tn{sing}}}
\newcommand{\anis}{a_{\tn{Nis}}}

\newcommand{\one}{\mathds{1}}
\newcommand{\oneX}{\one_{X}}
\newcommand{\oneY}{\one_{Y}}
\newcommand{\onec}{\one_{\tn{\tiny{Ch}}}}
\newcommand{\onev}{\one_{\tn{\tiny{Vo}}}}
\newcommand{\onedd}{\one_{\delta d}}
\newcommand{\oners}{\one(r)[s]}
\newcommand{\tate}{\one(1)}
\newcommand{\tato}{\one(1)[1]}
\newcommand{\tatt}{\one(1)[2]}
\newcommand{\tatn}{\one(n)[n]}
\newcommand{\tatta}{\one(a)[2a]}
\newcommand{\tatk}{\one(k)[k]}
\newcommand{\Zb}{\mathbf{Z}}
\newcommand{\TO}{\mathbf{Z}(1)}
\newcommand{\TN}{\mathbf{Z}(n)}
\newcommand{\Zq}{\mathbf{Z}(q)}
\newcommand{\Zn}{\mathbf{Z}(n)}
\newcommand{\HQ}{\mathbb{H}}
\newcommand{\Aff}{\mathbb{A}}
\newcommand{\aff}{\tn{aff}}
\newcommand{\Pj}{\mathbb{P}}
\newcommand{\Pjo}{\Pj^{1}}
\newcommand{\Pjor}{(\Pjo)^{r}}
\newcommand{\Pn}{\mathbb{P}^{n}}
\newcommand{\Proj}{\textnormal{Proj}}
\newcommand{\cjG}{\cj{G}}
\newcommand{\cjT}{\cj{T}}
\newcommand{\Latt}{\curly{L}att}
\newcommand{\LS}{\mathcal{L}}
\newcommand{\rk}{{\textnormal{rank}}}
\newcommand{\coker}{\textnormal{coker} \,}
\newcommand{\im}{\textnormal{im} \,}
\newcommand{\imess}{\tn{im}_{\tn{ess}} \,}
\newcommand{\Hom}{\textnormal{Hom}}
\newcommand{\iHom}{\ul{\tn{Hom}}}
\newcommand{\iHomet}{\ul{\tn{Hom}}_{\et}}
\newcommand{\iHomn}{\ul{\tn{Hom}}_{\tn{Nis}}}
\newcommand{\Hm}{\textnormal{H}}
\newcommand{\Ob}{\textnormal{Ob}}
\newcommand{\Abc}{\shf{A}}
\newcommand{\Bbc}{\shf{B}}
\newcommand{\prmn}{P^{R}_{M, N}}
\newcommand{\frmn}{F^{R}_{M, N}}
\newcommand{\Mor}{\textnormal{Mor}}
\newcommand{\op}{\tn{op}}
\newcommand{\Ext}{\tn{Ext}}
\newcommand{\Exto}{\tn{Ext}^{1}}
\newcommand{\YExt}{\tn{YnExt}}
\newcommand{\Extz}{\tn{Ext}^{0}}
\newcommand{\Exti}{\tn{Ext}^{i}}
\newcommand{\EExt}{\mathbb{E}\tn{xt}}
\newcommand{\hit}{{}^{\tn{t}}}
\newcommand{\lot}{{}_{\tn{t}}}
\newcommand{\Lin}{\tn{Lin}}
\newcommand{\Bil}{\tn{Bil}}
\newcommand{\tu}{\triangle} 
\newcommand{\btd}{\bigtriangledown}
\newcommand{\dom}{\textnormal{dom}}
\newcommand{\codom}{\textnormal{codom}}
\newcommand{\range}{\textnormal{range}}
\newcommand{\codim}{\textnormal{codim}}
\newcommand{\infim}{\textnormal{inf}}
\newcommand{\Isom}{\textnormal{Isom}}
\newcommand{\Spec}{\textnormal{Spec}}
\newcommand{\Speck}{\Spec(k)}
\newcommand{\rad}{\textnormal{rad}}
\newcommand{\spc}{\textnormal{sp}}
\newcommand{\Tag}{\tn{T}}
\newcommand{\nrm}{\textnormal{norm}}
\newcommand{\Imp}{\textnormal{Im}}
\newcommand{\Imag}{\textnormal{Im}}
\newcommand{\Rl}{\tn{Re}}
\newcommand{\Pu}{\tn{Pu}}
\newcommand{\vol}{\textnormal{vol}}
\newcommand{\Res}{\textnormal{Res}}
\newcommand{\Ress}{\textnormal{Res}\limits}
\newcommand{\expb}{\mathbf{e}}
\newcommand{\ua}{\uparrow}
\newcommand{\ra}{\rightarrow}
\newcommand{\la}{\leftarrow}
\newcommand{\Ra}{\Rightarrow}
\newcommand{\La}{\Leftarrow}
\newcommand{\Lefa}{L_{\alpha}}
\newcommand{\Lra}{\Leftrightarrow}
\newcommand{\lLra}{\Longleftrightarrow}
\newcommand{\llra}{\longleftrightarrow}
\newcommand{\lra}{\longrightarrow}
\newcommand{\lla}{\longleftarrow}
\newcommand{\hra}{\hookrightarrow}
\newcommand{\thra}{\twoheadrightarrow}
\newcommand{\Li}{\textnormal{Li}}
\newcommand{\rprod}{\prod}
\newcommand{\vp}{\varphi}
\newcommand{\vpm}{\vp^{-1}}
\newcommand{\erf}{\tn{erf}}
\newcommand{\sech}{\tn{sech}}
\newcommand{\cosech}{\tn{cosech}}
\newcommand{\Curl}{\tn{Curl}\,}
\newcommand{\sse}{\subseteq}
\newcommand{\suse}{\supseteq}
\newcommand{\susn}{\supsetneq}
\newcommand{\ssn}{\subsetneq}
\newcommand{\nsuse}{\nsupseteq}
\newcommand{\nsse}{\nsubseteq}
\newcommand{\sst}{\subset}
\newcommand{\Alg}{\tn{\tbf{Alg}}}
\newcommand{\Algk}{\tn{\tbf{Alg}}_{k}}
\newcommand{\AlgA}{\tn{\tbf{Alg}}_{A}}
\newcommand{\Gp}{\tn{\tbf{Gp}}}
\newcommand{\eff}{\tn{eff}}
\newcommand{\GpSch}{\tn{\tbf{GpSch}}}
\newcommand{\LgGpSch}{\tn{\tbf{LgGpSch}}}
\newcommand{\GpSchS}{\tn{\tbf{GpSch}}/S}
\newcommand{\LgGpSchS}{\tn{\tbf{LgGpSch}}/S}
\newcommand{\GpSchk}{\tn{\tbf{GpSch}}/k}
\newcommand{\CGSFTk}{\tn{\tbf{CmGSch}}/k}
\newcommand{\Ab}{\tn{\tbf{Ab}}}
\newcommand{\SAb}{\tn{\tbf{SAb}}}
\newcommand{\SAbk}{\SAb(k)}
\newcommand{\SAbS}{\SAb(S)}
\newcommand{\AbS}{\tn{\tbf{AbSch}}}
\newcommand{\AbV}{\tn{\tbf{AbVar}}}
\newcommand{\AbVk}{\tn{\tbf{AbVar}}/k}
\newcommand{\AbVkQ}{\tn{\tbf{AbVar}}_{\Q}/k}
\newcommand{\AVI}{\tn{\tbf{AbVarIsog}}}
\newcommand{\AVIk}{\tn{\tbf{AbVarIsog}}/k}
\newcommand{\AVIkQ}{(\tn{\tbf{AbVarIsog}}/k)^{\Q}}
\newcommand{\tAbS}{{}^{\tn{t}}\tn{\tbf{AbS}}}
\newcommand{\CGS}{\tn{\tbf{CommGpSch}}}
\newcommand{\CGSS}{\tn{\tbf{CommGpSch}}_{S}}
\newcommand{\fgAb}{\tn{\tbf{fgAb}}}
\newcommand{\Open}{\tn{\tbf{Open}}}
\newcommand{\Soe}{\tn{\tbf{S}}^{\tn{eff}}_{1}}
\newcommand{\So}{\tn{\tbf{S}}_{1}}
\newcommand{\res}{\tn{res}}
\newcommand{\an}{\tn{an}}
\newcommand{\Prim}{\tn{Prim}}
\newcommand{\rsa}{\rightsquigarrow}
\newcommand{\lbk}{\lbrack}
\newcommand{\rbk}{\rbrack}

\newcommand{\aXYZ}{\alpha_{X, Y}^{Z}}
\newcommand{\aXYU}{\alpha_{X, Y}^{U}}

\newcommand{\tXYZ}{\theta_{X, Y}^{Z}}
\newcommand{\tXYU}{\theta_{X, Y}^{U}}
\newcommand{\tXYV}{\theta_{X, Y}^{V}}
\newcommand{\tXYW}{\theta_{X, Y}^{W}}

\newcommand{\Sets}{\tn{\tbf{Sets}}}
\newcommand{\SFL}{S_{\tn{FL}}}
\newcommand{\kFL}{\Speck_{\tn{FL}}}
\newcommand{\ket}{\Speck_{\et}}
\newcommand{\D}{\tn{D}}
\newcommand{\DD}{\mathbb{D}}
\newcommand{\GoS}{\mathbb{G}/S}
\newcommand{\uG}{\ul{\mathbb{G}}}
\newcommand{\uGS}{\ul{\mathbb{G}}/S}
\newcommand{\Dm}{\tn{D}^{-}}
\newcommand{\Dp}{\tn{D}^{+}}
\newcommand{\Db}{\tn{D}^{\tn{b}}}
\newcommand{\bDm}{\tn{\tbf{D}}^{-}}
\newcommand{\bDmc}{\tn{\tbf{D}}^{-}_{\tn{Cor}}}
\newcommand{\bDms}{\tn{\tbf{D}}^{-}_{\tn{St}}}

\newcommand{\bS}{\tn{\tbf{S}}}
\newcommand{\bSc}{\tn{\tbf{S}}_{\tn{Cor}}}
\newcommand{\bSt}{\tn{\tbf{S}}_{\tn{St}}}

\newcommand{\qcoh}{\tn{qcoh}}
\newcommand{\Dbqcoh}{\Db_{\qcoh}}

\newcommand{\Cb}{\cly{C}^{\tn{b}}}
\newcommand{\Kb}{\cly{K}^{\tn{b}}}

\newcommand*{\longhookrightarrow}{\ensuremath{\lhook\joinrel\relbar\joinrel\rightarrow}}
\newcommand{\lhra}{\longhookrightarrow}

\newcommand{\gr}{\tn{gr}}
\newcommand{\grwn}{\tn{gr}^{W}_{n}}
\newcommand{\grwp}{\tn{gr}^{W}_{p}}
\newcommand{\grwo}{\tn{gr}^{W}_{1}}
\newcommand{\grwi}{\tn{gr}^{W}_{i}}
\newcommand{\grw}{\tn{gr}^{W}}
\newcommand{\et}{\tn{\'{e}t}}
\newcommand{\fl}{\tn{fl}}
\newcommand{\fppf}{\tn{fppf}}
\newcommand{\nis}{\tn{Nis}}
\newcommand{\eh}{\tn{\'{e}h}}
\newcommand{\Dlg}{\tn{Dlg}}
\newcommand{\Tori}{\tn{\tbf{Tori}}}
\newcommand{\fgfrAb}{\tn{\tbf{fgfrAb}}}
\newcommand{\fgf}{\tn{\tbf{fgfrAb}}}
\newcommand{\cch}{\chi_{\bl}}
\newcommand{\Fnlm}{F^{\nu}_{\lambda, \mu}}
\newcommand{\etl}{\'{e}tale}
\newcommand{\Etl}{\'{E}tale}
\newcommand{\Vo}{Voevodsky}
\newcommand{\Vov}{\tn{Vo}}
\newcommand{\Dgl}{D\'{e}glise}
\newcommand{\SeAV}{Serre-Albanese variety}
\newcommand{\SeAT}{Serre-Albanese torsor}
\newcommand{\ASc}{Albanese scheme}
\newcommand{\sesq}{short exact sequence}
\newcommand{\Cech}{\check{C}ech}
\newcommand{\cech}{\check{c}ech}
\newcommand{\sng}{\tn{sing}}
\newcommand{\co}{c_{1}}
\newcommand{\cosng}{c_{1}^{\sng}}
\newcommand{\comot}{c_{1}^{\mot}}
\newcommand{\mot}{\tn{mot}}

\newcommand{\ner}{N\'{e}ron}
\newcommand{\ners}{N\'{e}ron-Severi}
\newcommand{\nsv}{Nisnevich}
\newcommand{\pnc}{Poincar\'{e}}
\newcommand{\gc}{Grothendieck-Chow}
\newcommand{\gcm}{Grothendieck-Chow motive}
\newcommand{\gcms}{Grothendieck-Chow motives}
\newcommand{\fa}{Fourier analysis}
\newcommand{\rie}{Riemann}
\newcommand{\rien}{Riemannian}
\newcommand{\rienme}{Riemannian metric}
\newcommand{\rienmes}{Riemannian metrics}
\newcommand{\rienma}{Riemannian manifold}
\newcommand{\rienmas}{Riemannian manifolds}

\newcommand{\Coh}{\tn{\tbf{Coh}}}
\newcommand{\Cop}{\tn{\tbf{Coh}}(\Pj^{1})}
\newcommand{\Torp}{\tn{\tbf{Tor}}(\Pj^{1})}
\newcommand{\Tor}{\tn{\tbf{Tor}}}
\newcommand{\slt}{\idl{sl}_{2}}

\newcommand{\De}{\tn{DM}^{\tn{eff}}_{\,\tn{\pcl{X}}}}
\newcommand{\Dg}{\tn{DM}_{\tn{gm}}}
\newcommand{\Deg}{\tn{DM}^{\tn{eff}}_{\tn{gm}}}
\newcommand{\Chow}{\tn{\tbf{Chow}}} 
\newcommand{\Chowk}{\Chow(k)}
\newcommand{\tChowk}{{}^{\tn{t}}\Chow(k)}
\newcommand{\Chowuk}{\tn{\tbf{Chow}}_{\tn{un}}(k)}
\newcommand{\ChowS}{\tn{\tbf{Chow}}(S)}
\newcommand{\Chowe}{\tn{\tbf{Chow}}^{\tn{eff}}}
\newcommand{\ChoweS}{\tn{\tbf{Chow}}^{\tn{eff}}(S)}
\newcommand{\Chowek}{\Chow^{\tn{eff}}(k)}
\newcommand{\ChowZek}{\Chow^{\tn{eff}}_{\Z}(k)}
\newcommand{\tChowek}{{}^{\tn{t}}\Chow^{\tn{eff}}(k)}
\newcommand{\ChowkQ}{\tn{\tbf{Chow}}(k)_{\Q}}
\newcommand{\ChowekQ}{\tn{\tbf{Chow}}^{\tn{eff}}(k)_{\Q}}
\newcommand{\CHp}{\tn{CH}^{p}}
\newcommand{\CHo}{\tn{CH}_{0}}
\newcommand{\tCHo}{\wt{\tn{CH}}_{0}}
\newcommand{\req}{\equiv_{\tn{rat}}}
\newcommand{\DMe}{\DM^{\tn{eff}}} 
\newcommand{\DTM}{\tn{DTM}} 
\newcommand{\DTMn}{\tn{DTM}_{n}} 
\newcommand{\DTMet}{\tn{DTM}_{\et}} 
\newcommand{\DTMeet}{\tn{DTM}^{\tn{eff}}_{\,\tn{\pcl{X}}, \et}}
\newcommand{\DTMetk}{\tn{DTM}_{\et}(k)} 
\newcommand{\DTMeetk}{\tn{DTM}^{\tn{eff}}_{\,\tn{\pcl{X}}, \et}(k)}
\newcommand{\DMeet}{\DMe_{\,\tn{\pcl{X}}, \et}}
\newcommand{\DMen}{\DMe_{\,\tn{\pcl{X}}, \nis}}
\newcommand{\DMett}{\DMe_{\,\tn{\pcl{X}}, \tau}}
\newcommand{\DMettk}{\DMett(k)}
\newcommand{\DMenk}{\DMen(k)}
\newcommand{\DMenkQ}{\DMen(k, \Q)}
\newcommand{\DMenkR}{\DMen(k, R)}

\newcommand{\DMettkR}{\DMett(k, R)}
\newcommand{\DMettkZ}{\DMett(k, \Z)}
\newcommand{\DMettkQ}{\DMett(k, \Q)}
\newcommand{\DMeetk}{\DMeet(k)}
\newcommand{\DMeetkQ}{\DMeet(k, \Q)}
\newcommand{\DMeetkZ}{\DMeet(k, \Z)}
\newcommand{\DMeetkR}{\DMeet(k, R)}
\newcommand{\DMeetA}{\DMeet(A)}
\newcommand{\DMeetX}{\DMeet(X)}
\newcommand{\DMX}{\DM(X)}
\newcommand{\DMY}{\DM(Y)}
\newcommand{\DMgmX}{\DMgm(X)}
\newcommand{\DMgmY}{\DMgm(Y)}

\newcommand{\DMetk}{\tn{\tbf{DM}}_{\,\tn{\pcl{X}}, \et}(k)}
\newcommand{\DMetkQ}{\tn{\tbf{DM}}_{\,\tn{\pcl{X}}, \et}(k)_{\Q}}
\newcommand{\DMeget}{\tn{\tbf{DM}}^{\tn{eff}}_{\tn{gm}, \et}}
\newcommand{\DMegetk}{\tn{\tbf{DM}}^{\tn{eff}}_{\tn{gm}, \et}(k)}
\newcommand{\DMgetk}{\tn{\tbf{DM}}_{\tn{gm}, \et}(k)}

\newcommand{\DMentkQ}{\tn{\tbf{DM}}^{\tn{eff}}_{\,\tn{\pcl{X}}, \nis}(k)_{\Q}}
\newcommand{\DMegetQ}{\tn{DM}^{\tn{eff}}_{\tn{gm}, \et} \ox \Q}
\newcommand{\DMegetkQ}{\tn{DM}^{\tn{eff}}_{\tn{gm}, \et}(k) \ox \Q}
\newcommand{\DMg}{\tn{DM}_{\tn{gm}}}
\newcommand{\DMeg}{\tbf{DM}^{\tn{eff}}_{\tn{gm}}}
\newcommand{\DMgS}{\tn{\tbf{DM}}_{\tn{gm}}(S)}
\newcommand{\DMeS}{\tn{\tbf{DM}}^{\tn{eff}}(S)}
\newcommand{\DMegS}{\tn{\tbf{DM}}^{\tn{eff}}_{\tn{gm}}(S)}

\newcommand{\DMegQ}{\DMeg \ox \Q}
\newcommand{\dDMeg}{d_{\le 1}\tn{DM}^{\tn{eff}}_{\tn{gm}}}
\newcommand{\dDMeget}{d_{\le 1}\tn{DM}^{\tn{eff}}_{\tn{gm}, \et}}
\newcommand{\dDMegetk}{d_{\le 1}\tn{DM}^{\tn{eff}}_{\tn{gm}, \et}(k)}
\newcommand{\as}{\alpha_{*}}

\newcommand{\Mot}{\tn{Mot}}
\newcommand{\Motk}{\Mot/k}

\newcommand{\DbM}{\Db(\Mok)}
\newcommand{\DbMp}{\Db(\Mok[1/p])}
\newcommand{\DbMQ}{\Db(\Mok \ox \Q)}
\newcommand{\CbMp}{\Cb(\Mok[1/p])}
\newcommand{\CbMQ}{\Cb(\Mok \ox \Q)}
\newcommand{\dlo}{d_{\le 1}} 
\newcommand{\Dlo}{\tn{D}_{\le 1}}
\newcommand{\Deto}{\tn{D}^{\et}_{\le 1}}
\newcommand{\Dnto}{\tn{D}^{\tn{Nis}}_{\le 1}}
\newcommand{\Mo}{\cly{M}_{1}}
\newcommand{\Mof}{M_{1}}
\newcommand{\hof}{h^{1}}
\newcommand{\Do}{D_{1}}
\newcommand{\shtG}{\sht{G}}
\newcommand{\shtH}{\sht{H}}
\newcommand{\Mog}{\mathbf{M}_{1}}
\newcommand{\Mg}{\mathbf{M}}
\newcommand{\Mogr}{\mathbf{M}_{1}^{r}}
\newcommand{\Mgr}{\mathbf{M}^{r}}
\newcommand{\Mea}{\cly{M}_{\tn{anc}}^{\tn{eff}}}
\newcommand{\tMeo}{{}^{t}\cly{M}_{1}^{\tn{eff}}}
\newcommand{\ctMeo}{{}_{t}\cly{M}_{1}^{\tn{eff}}}
\newcommand{\tMo}{{}^{t}\cly{M}_{1}}
\newcommand{\Mok}{\cly{M}_{1}(k)}
\newcommand{\MoS}{\cly{M}_{1}(S)}
\newcommand{\Meak}{\cly{M}_{\tn{anc}}^{\tn{eff}}(k)}
\newcommand{\tMeok}{{}^{t}\cly{M}_{1}^{\tn{eff}}(k)}
\newcommand{\tMok}{{}^{t}\cly{M}_{1}(k)}
\newcommand{\ctMo}{{}_{t}\cly{M}_{1}}
\newcommand{\ctMok}{{}_{t}\cly{M}_{1}(k)}

\newcommand{\tH}{{}^{t} H}
\newcommand{\tHn}{{}^{t} H^{n}}
\newcommand{\tHln}{{}^{t} H_{n}}
\newcommand{\Mz}{\cly{M}_{0}}
\newcommand{\Mzk}{\cly{M}_{0}(k)}
\newcommand{\MzS}{\cly{M}_{0}(S)}
\newcommand{\tMz}{{}^{t}\cly{M}_{0}}
\newcommand{\tMzk}{{}^{t}\cly{M}_{0}(k)}

\newcommand{\MS}{\cly{M}(S)}
\newcommand{\MChS}{\cly{M}^{0}(S)}
\newcommand{\MChpS}{\cly{M}^{0}_{+}(S)}

\newcommand{\LAlb}{\tn{LAlb}}
\newcommand{\LAlbc}{\tn{LAlb}^{\tn{c}}}
\newcommand{\LAlbcs}{\tn{LAlb}^{\tn{c}} \;}
\newcommand{\LiAlb}{\tn{L}_{i}\tn{Alb}}
\newcommand{\LoAlb}{\tn{L}_{1}\tn{Alb}}
\newcommand{\LtAlb}{\tn{L}_{2}\tn{Alb}}
\newcommand{\LAlbi}[1]{\tn{L}_{#1}\tn{Alb}}
\newcommand{\LiAlbc}{\tn{L}_{i}\tn{Alb}^{\tn{c}}}
\newcommand{\LzAlb}{\tn{L}_{0}\tn{Alb}}
\newcommand{\LAlbQ}{\tn{LAlb}^{\Q}}
\newcommand{\LiAlbQ}{\LiAlb^{\Q}}

\newcommand{\RPic}{\tn{RPic}}
\newcommand{\RPicc}{\tn{RPic}^{\tn{c}}}
\newcommand{\RiPic}{\tn{R}_{i}\tn{Pic}}

\newcommand{\RPici}[1]{\tn{R}_{#1}\tn{Pic}}

\newcommand{\RiPicc}{\tn{R}_{i}\tn{Pic}^{\tn{c}}}
\newcommand{\LAlbs}{\tn{LAlb} \;}
\newcommand{\RPics}{\tn{RPic} \;}
\newcommand{\AXk}{\cly{A}_{X/k}}
\newcommand{\AGk}{\cly{A}_{G/k}}
\newcommand{\AUk}{\cly{A}_{U/k}}
\newcommand{\AzXk}{\cly{A}^{0}_{X/k}}
\newcommand{\AsXk}{\cly{A}^{*}_{X/k}}
\newcommand{\uAXk}{\ul{\cly{A}}_{X/k}}
\newcommand{\uAsXk}{\ul{\cly{A}}^{*}_{X/k}}
\newcommand{\PicXk}{\tn{Pic}_{X/k}}
\newcommand{\PicCk}{\tn{Pic}_{C/k}}
\newcommand{\uPicXk}{\ul{\tn{Pic}}_{X/k}}
\newcommand{\uPicCk}{\ul{\tn{Pic}}_{C/k}}
\newcommand{\wPicXk}{\wPic_{X/k}}
\newcommand{\wPicCk}{\wPic_{C/k}}

\newcommand{\PicoXk}{\tn{Pic}^{\circ}_{X/k}}
\newcommand{\uPicoXk}{\ul{\tn{Pic}}^{\circ}_{X/k}}

\newcommand{\NSXk}{\tn{NS}_{X/k}}
\newcommand{\NSZX}{\tn{NS}_{Z}(X)}
\newcommand{\NSsXk}{\tn{NS}_{X/k}^{*}}

\newcommand{\DM}{\tn{\tbf{DM}}} 
\newcommand{\DMgm}{\DM_{\tn{gm}}}
\newcommand{\DT}{\tn{\tbf{DT}}}
\newcommand{\Sm}{\tn{\tbf{Sm}}}
\newcommand{\Smc}{\tn{\tbf{SmCor}}}
\newcommand{\Sch}{\tn{\tbf{Sch}}}
\newcommand{\Schsf}{\Sch_{\tn{sf}}}
\newcommand{\SchsfX}{\Schsf/X}
\newcommand{\SchS}{\tn{\tbf{Sch}}/S}
\newcommand{\Var}{\tn{\tbf{Var}}}
\newcommand{\SmPrVar}{\tn{\tbf{SmPrVar}}}
\newcommand{\GrVekt}{\tn{\tbf{GrVekt}}}
\newcommand{\Vark}{\tn{\tbf{Var}}/k}
\newcommand{\Mgm}{M_{\tn{gm}}}
\newcommand{\tMgm}{\wt{\tn{M}}_{\tn{gm}}}
\newcommand{\Mgmc}{\tn{M}_{\tn{gm}}^{\tn{c}}}
\newcommand{\Mgmi}{\tn{M}_{\tn{gm}}^{i}}
\newcommand{\Mgmt}{\widetilde{\tn{M}}_{\tn{gm}}}

\newcommand{\DMek}{\DMe(k)}

\newcommand{\psul}[2]{{#1}_{\rlap{--}{\phantom{#2}}}}
\newcommand{\pl}{\partial}
\newcommand{\bl}{\bullet}
\newcommand{\sbl}{\mathsmaller{\bullet}}
\newcommand{\ssbl}{\mathsmaller{\mathsmaller{\bullet}}}
\newcommand{\w}{\wedge}
\newcommand{\Dgk}{\tn{\tbf{DM}}_{\tn{gm}}(k)}
\newcommand{\Degk}{\tn{\tbf{DM}}^{\tn{eff}}_{\tn{gm}}(k)}
\newcommand{\DMk}{\tn{\tbf{DM}}(k)}
\newcommand{\DMS}{\tn{\tbf{DM}}(S)}
\newcommand{\DMgk}{\tn{\tbf{DM}}_{\tn{gm}}(k)}
\newcommand{\DMgkQ}{\tn{\tbf{DM}}_{\tn{gm}}(k)_{\Q}}

\newcommand{\DMegk}{\tn{\tbf{DM}}^{\tn{eff}}_{\tn{gm}}(k)}
\newcommand{\DMegkQ}{\tn{\tbf{DM}}^{\tn{eff}}_{\tn{gm}}(k)_{\Q}}
\newcommand{\Dgr}{\tn{Deg}}
\newcommand{\DTegk}{\tn{DT}^{\tn{eff}}_{\tn{gm}}(k)}
\newcommand{\DTk}{\tn{DT}(k)}
\newcommand{\DTS}{\tn{DT}(S)}
\newcommand{\CX}{\cat{C}/X}
\newcommand{\Cx}{\ul{C}_{*}}
\newcommand{\RCx}{\tn{\tbf{R}}\ul{C}_{*}}
\newcommand{\Cnx}{\ul{C}_{n}}
\newcommand{\rat}{\tn{\tiny{rat}}}
\newcommand{\xo}{x_{0}}
\newcommand{\yo}{y_{0}}
\newcommand{\xoT}{\xo^{T}}
\newcommand{\yoT}{\yo^{T}}

\newcommand{\zxoo}{z_{\xo}^{1}}
\newcommand{\zxon}{z_{\xo}^{n}}

\newcommand{\simeqcan}{\st{\tn{\tiny{can.}}}{\simeq}}

\newcommand{\AXz}{A_{X}^{0}}
\newcommand{\AXo}{A_{X}^{1}}
\newcommand{\poC}{p_{1}^{C}}
\newcommand{\fpo}{f_{P_{0}}}
\newcommand{\fxo}{f_{x_{0}}}
\newcommand{\jxo}{j_{x_{0}}}
\newcommand{\jxon}{j_{x_{0}}^{n}}
\newcommand{\jxoo}{j_{x_{0}}^{1}}
\newcommand{\jyo}{j_{y_{0}}}
\newcommand{\ixo}{\iota_{\xo}}
\newcommand{\ixon}{\ixo^{n}}
\newcommand{\iyo}{\iota_{\yo}}
\newcommand{\pxo}{p_{x_{0}}}
\newcommand{\axo}{a_{x_{0}}}
\newcommand{\cxo}{c_{x_{0}}}
\newcommand{\vxo}{v_{x_{0}}}
\newcommand{\upxo}{\ul{p}_{x_{0}}}
\newcommand{\uaxo}{\ul{a}_{x_{0}}}
\newcommand{\qX}{q_{X}}
\newcommand{\qG}{q_{G}}
\newcommand{\Go}{G^{\circ}}

\newcommand{\vnM}{\tn{V}_{n}(M)}
\newcommand{\vmM}{\tn{V}_{m}(M)}
\newcommand{\vnN}{\tn{V}_{n}(N)}
\newcommand{\vmN}{\tn{V}_{m}(N)}
\newcommand{\vn}{\tn{V}_{n}}
\newcommand{\vm}{\tn{V}_{m}}
\newcommand{\Zar}{\tn{Zar}}
\newcommand{\Cov}{\tn{Cov}}
\newcommand{\Et}{\tn{\tbf{Et}}}
\newcommand{\EtX}{\Et/X}
\newcommand{\Etk}{\tn{\tbf{Et}}/k}
\newcommand{\Hetz}{H_{\tn{{\'e}t}}^{0}}
\newcommand{\Heto}{H_{\tn{{\'e}t}}^{1}}
\newcommand{\Mc}{M^{\tn{c}}}
\newcommand{\Chi}{\cly{X}}
\newcommand{\toh}{\textonehalf}
\newcommand{\Lc}{L^{\tn{c}}}
\newcommand{\Smk}{\tn{\tbf{Sm}}/k}
\newcommand{\Smck}{\tn{\tbf{SmCor}}/k}
\newcommand{\SmckQ}{\tn{\tbf{SmCor}}_{\Q}/k}
\newcommand{\SmcS}{\tn{\tbf{SmCor}}/S}

\newcommand{\scd}{[\, \cdot \, ]_{\tn{sc}}}
\newcommand{\scf}[1]{[#1]_{\tn{sc}}}

\newcommand{\Smcko}{\tn{\tbf{SmCor}}\; k}
\newcommand{\Schk}{\tn{\tbf{Sch}}/k}

\newcommand{\SchX}{\Sch/X}

\newcommand{\pio}{\pi_{0}}
\newcommand{\nat}{\natural}

\newcommand{\pioX}{\pi_{0}(X)}
\newcommand{\pioG}{\pi_{0}(G)}

\newcommand{\cZpioX}{\cly{Z}_{\pi_{0}(X)}}
\newcommand{\cZX}{\cly{Z}_{X}}
\newcommand{\cAX}{\cly{A}_{X}}
\newcommand{\AX}{A_{X}}
\newcommand{\pAX}{\tn{p}A_{X}}
\newcommand{\AoX}{A^{\circ}_{X}}
\newcommand{\cAoX}{\cly{A}^{\circ}_{X}}
\newcommand{\SmProj}{\tn{\tbf{SmProj}}}
\newcommand{\CSmProj}{\tn{\tbf{CSmProj}}}
\newcommand{\SmProjk}{\SmProj/k}
\newcommand{\CSmProjk}{\CSmProj/k}
\newcommand{\CSmProjek}{\CSmProj^{\tn{eff}}(k)}
\newcommand{\Sh}{\tn{\tbf{Sh}}}
\newcommand{\Shv}{\tn{\tbf{Shv}}}
\newcommand{\Shvs}{\tn{\tbf{Shv}}^{\tn{s}}}
\newcommand{\PSh}{\tn{\tbf{PShv}}}
\newcommand{\Ch}{\tn{\tbf{Ch}}}
\newcommand{\Chp}{\tn{\tbf{Ch}}^{+}}
\newcommand{\Chm}{\tn{\tbf{Ch}}^{-}}
\newcommand{\CoCh}{\tn{\tbf{CoCh}}}
\newcommand{\CoChp}{\tn{\tbf{CoCh}}^{+}}
\newcommand{\CoChm}{\tn{\tbf{CoCh}}^{-}}
\newcommand{\SN}{\tn{\tbf{Sh}}_{\tn{Nis}}}
\newcommand{\SE}{\tn{\tbf{Sh}}_{\tn{\'et}}}
\newcommand{\SNk}{\tn{\tbf{Sh}}_{\tn{Nis}}(\Smck)}
\newcommand{\SEk}{\tn{\tbf{Sh}}_{\et}(\Smck)}
\newcommand{\SEkQ}{\SE(k, \Q)}
\newcommand{\fgt}{\tn{\tbf{fgt}}}

\newcommand{\SEsmk}{\SE(\Smk)}
\newcommand{\SEsmck}{\SE(\Sm/\cj{k})}

\newcommand{\ST}{\tn{\tbf{ShT}}}
\newcommand{\STt}{\ST_{\tau}}
\newcommand{\STE}{\ST_{\et}}
\newcommand{\STN}{\ST_{\nis}}
\newcommand{\STtk}{\STt(k)}
\newcommand{\STEk}{\STE(k)}
\newcommand{\STNk}{\STN(k)}
\newcommand{\STEkR}{\STE(k, R)}
\newcommand{\STNkR}{\STN(k, R)}
\newcommand{\STEkZ}{\STE(k, \Z)}
\newcommand{\STNkZ}{\STN(k, \Z)}
\newcommand{\STEkQ}{\STE(k, \Q)}
\newcommand{\STNkQ}{\STN(k, \Q)}

\newcommand{\STtkR}{\STt(k, R)}
\newcommand{\STtkZ}{\STt(k, \Z)}
\newcommand{\STtkQ}{\STt(k, \Q)}

\newcommand{\PST}{\tn{\tbf{PShT}}}
\newcommand{\PSTk}{\PST(k)}
\newcommand{\PSTkR}{\PST(k, R)}
\newcommand{\PSTkZ}{\PST(k, \Z)}
\newcommand{\PSTkQ}{\PST(k, \Q)}

\newcommand{\SEsk}{\tn{\tbf{Sh}}_{\et}(\Smk)}
\newcommand{\SNkQ}{\tn{\tbf{Sh}}_{\tn{Nis}}(\SmckQ)}
\newcommand{\SNks}{\tn{\tbf{Sh}}_{\tn{Nis}}(k)}
\newcommand{\DmSNk}{\Dm \; \tn{\tbf{Sh}}_{\tn{Nis}}(\Smcko)}
\newcommand{\DmSN}{\Dm \;  \tn{\tbf{Sh}}_{\tn{Nis}}}
\newcommand{\MM}{\cly{MM}}
\newcommand{\MMk}{\cly{MM}(k)}
\newcommand{\Mnum}{\cly{M}^{\tn{eff}}_{\tn{num}}}
\newcommand{\Mnumk}{\cly{M}^{\tn{eff}}_{\tn{num}}(k)}
\newcommand{\Mh}{\cly{M}^{\tn{eff}}_{\tn{hom}}}
\newcommand{\Mhk}{\cly{M}^{\tn{eff}}_{\tn{hom}}(k)}
\newcommand{\Ma}{\cly{M}^{\tn{eff}}_{\tn{alg}}}
\newcommand{\Mak}{\cly{M}^{\tn{eff}}_{\tn{alg}}(k)}
\newcommand{\Mer}{\cly{M}^{\tn{eff}}_{\tn{rat}}}
\newcommand{\Merk}{\cly{M}^{\tn{eff}}_{\tn{rat}}(k)}
\newcommand{\Mek}{\cly{M}^{\tn{eff}}_{\sim}}
\newcommand{\Meek}{\cly{M}^{\tn{eff}}_{\sim}(k)}
\newcommand{\RMp}{\cly{M}_{\sim}}
\newcommand{\RMkp}{\cly{M}_{\sim}(k)}
\newcommand{\Mnump}{\cly{M}_{\tn{num}}}
\newcommand{\Mnumkp}{\cly{M}_{\tn{num}}(k)}
\newcommand{\HI}{\tn{HI}}
\newcommand{\Hb}{\cly{H}^{\tn{b}}}
\newcommand{\Hmi}{H^{-i}}
\newcommand{\Hpi}{H^{i}}
\newcommand{\Hmz}{H^{0}}
\newcommand{\Hmo}{H^{-1}}

\newcommand{\Heti}{\Hpi_{\et}}
\newcommand{\Hnisi}{\Hpi_{\nis}}
\newcommand{\Hmoti}{\Hpi_{\mot}}

\newcommand{\HIet}{\tn{HI}_{\tn{\'et}}}
\newcommand{\HIset}{\tn{HI}_{\tn{\'et}}^{\tn{s}}}
\newcommand{\cH}{\tn{H}}
\newcommand{\cxH}{\widecheck{H}}
\newcommand{\cxHo}{\widecheck{H}^{1}}
\newcommand{\HIk}{\tn{HI}(k)}
\newcommand{\HIkR}{\tn{HI}(k, R)}
\newcommand{\EA}{\shf{E}_{\A^{1}}}
\newcommand{\Ztr}{\Z_{\tn{tr}}}
\newcommand{\Rtr}{R_{\tn{tr}}}
\newcommand{\ZtrX}{\Z_{\tn{tr}}(X)}
\newcommand{\Cbl}{C_{\bullet}}
\newcommand{\Cobl}{C^{\bullet}}
\newcommand{\CblF}{C_{\bullet}\shf{F}}
\newcommand{\Hot}{\tn{\tbf{Hot}}}
\newcommand{\Hotb}{\tn{Hot}^{\tn{b}}}
\newcommand{\Hotp}{\tn{Hot}^{+}}
\newcommand{\Hotm}{\tn{Hot}^{-}}
\newcommand{\Tot}{\tn{Tot}}
\newcommand{\TotQ}{\tn{Tot}^{\Q}}
\newcommand{\Tots}{\tn{Tot} \;}
\newcommand{\CH}{\tn{CH}}
\newcommand{\CHx}{\tn{CH}^{*}}

\newcommand{\CHMo}{\tn{CH}_{\Mo}}

\newcommand{\uCH}{\ul{\tn{CH}}}
\newcommand{\Rat}{\tn{Rat}}
\newcommand{\Ratr}{\tn{Rat}_{r}}
\newcommand{\Cor}{\tn{Cor}}
\newcommand{\Zc}{\tn{Z}}
\newcommand{\Zfin}{\tn{Z}_{\tn{fin}}}
\newcommand{\Cors}{\tn{Cor}_{\sim}}
\newcommand{\Corsk}{\tn{Cor}_{\sim}(k)}
\newcommand{\Cork}{\tn{Cor}_{k}}
\newcommand{\CCork}{\tn{Cor}_{\tn{Ch}/k}}
\newcommand{\CCorkZ}{\tn{Cor}_{\tn{Ch}/k, \Z}}

\newcommand{\VCork}{\tn{Cor}_{\tn{Vo}/k}}
\newcommand{\VCorke}{\tn{Cor}_{\tn{Vo}/k}^{\eff}}
\newcommand{\CCor}{\tn{Cor}_{\tn{Ch}}}
\newcommand{\VCor}{\tn{Cor}_{\tn{Vo}}}
\newcommand{\Sgn}{\Sigma_{n}}
\newcommand{\SmCor}{\tn{\tbf{SmCor}}}
\newcommand{\SchCor}{\tn{\tbf{SchCor}}}
\newcommand{\SchCork}{\tn{\tbf{SchCor}}_{k}}
\newcommand{\SmCork}{\tn{\tbf{SmCor}}_{k}}

\newcommand{\SH}{\tn{\tbf{SH}}}
\newcommand{\SHX}{\tn{\tbf{SH}}(X)}
\newcommand{\SHA}{\tn{\tbf{SH}}(A)}

\newcommand{\pr}{\tn{pr}}
\newcommand{\cone}{\tn{cone}}
\newcommand{\GS}[1]{\Gamma(X, #1)}
\newcommand{\GSX}{\Gamma(X, \pul{U})}
\newcommand{\GammaX}{\Gamma^{X}}
\newcommand{\muonn}{\nu^{\ge n}}
\newcommand{\mlnn}{\nu_{< n}}
\newcommand{\mlenn}{\nu_{\le n}}
\newcommand{\munn}{\nu_{n}}
\newcommand{\cn}{c_{n}}
\newcommand{\ox}{\otimes}
\newcommand{\oxp}{\ox_{\tn{pre}}}
\newcommand{\oxpn}{\oxp^{n}}
\newcommand{\oxtr}{\ox_{\tn{tr}}}
\newcommand{\oxtrD}{\ox_{\tn{tr}}^{\tn{\tiny{doub}}}}
\newcommand{\oxL}{\ox^{\Lf}}
\newcommand{\oxtrLet}{\oxtr^{\Lf, \et}}
\newcommand{\bx}{\bigotimes}
\newcommand{\os}{\oplus}
\newcommand{\od}{\odot}
\newcommand{\bos}{\bigoplus}
\newcommand{\oxo}{\ox_{1}}
\newcommand{\sx}{\boxtimes}
\newcommand{\ost}{\circledast}
\newcommand{\Part}{\tn{Part}}
\newcommand{\supp}{\tn{supp}}
\newcommand{\Tp}{\tn{T}_{p}}
\newcommand{\TpM}{\tn{T}_{p}M}
\newcommand{\TpN}{\tn{T}_{p}N}
\newcommand{\TpV}{\tn{T}_{p}V}
\newcommand{\TLpW}{\tn{T}_{L(p)}W}

\newcommand{\TFpM}{\tn{T}_{F(p)}M}
\newcommand{\TFpN}{\tn{T}_{F(p)}N}
\newcommand{\TM}{\tn{T} M}
\newcommand{\TpX}{\tn{T}_{p}X}
\newcommand{\TX}{\tn{T} X}
\newcommand{\TqM}{\tn{T}_{q}M}

\newcommand{\hinis}{\ul{h}_{i}^{\tn{Nis}}}
\newcommand{\honis}{\ul{h}_{1}^{\tn{Nis}}}
\newcommand{\hznis}{\ul{h}_{0}^{\tn{Nis}}}

\newcommand{\hitau}{\ul{h}_{i}^{\tau}}
\newcommand{\hotau}{\ul{h}_{1}^{\tau}}
\newcommand{\hztau}{\ul{h}_{0}^{\tau}}

\newcommand{\hiet}{\ul{h}_{i}^{\et}}
\newcommand{\hoet}{\ul{h}_{1}^{\et}}
\newcommand{\hzet}{\ul{h}_{0}^{\et}}

\newcommand{\Hoet}{H^{1}_{\et}}

\newcommand{\cHoet}{\widecheck{H}^{1}_{\et}}

\newcommand{\hz}{h^{0}}
\newcommand{\ho}{h^{1}}
\newcommand{\hn}{h^{n}}
\newcommand{\htw}{h^{2}}
\newcommand{\htd}{h^{2d}}
\newcommand{\hgo}{h^{\ge 1}}
\newcommand{\hi}{h^{i}}

\newcommand{\chinis}{\ul{h}^{i}_{\nis}}
\newcommand{\chonis}{\ul{h}^{1}_{\nis}}
\newcommand{\chznis}{\ul{h}^{0}_{\nis}}
\newcommand{\chmnis}{\ul{h}^{-1}_{\nis}}
\newcommand{\chnis}{\ul{h}_{\nis}}

\newcommand{\cs}[1]{{#1}^{\tn{cs}}}

\newcommand{\Slm}{S_{\lambda}}
\newcommand{\Tk}{\tn{T}^{k}}
\newcommand{\Ak}{\Lambda^{k}}
\newcommand{\Sk}{\Sigma^{k}}
\newcommand{\TkM}{\tn{T}^{k} M}
\newcommand{\TkN}{\tn{T}^{k} N}
\newcommand{\AkM}{\Lambda^{k} M}
\newcommand{\SkM}{\Sigma^{k} M}
\newcommand{\TkV}{\tn{T}^{k} (V)}
\newcommand{\AkV}{\Lambda^{k} (V)}
\newcommand{\SkV}{\Sigma^{k} (V)}
\newcommand{\TkW}{\tn{T}^{k} (W)}
\newcommand{\AkW}{\Lambda^{k} (W)}
\newcommand{\SkW}{\Sigma^{k} (W)}

\newcommand{\TaukM}{\tau^{k}(M)}
\newcommand{\TaukN}{\tau^{k}(N)}

\newcommand{\Td}{\tn{T}^{\bl}}
\newcommand{\Ad}{\Lambda^{\bl}}
\newcommand{\TdM}{\tn{T}^{\bl} M}
\newcommand{\AdM}{\Lambda^{\bl} M}
\newcommand{\SdM}{\Sigma^{\bl} M}
\newcommand{\TdV}{\tn{T}^{\bl} (V)}
\newcommand{\AdV}{\Lambda^{\bl} (V)}
\newcommand{\SdV}{\Sigma^{\bl} (V)}
\newcommand{\TdW}{\tn{T}^{\bl} (W)}
\newcommand{\AdW}{\Lambda^{\bl} (W)}
\newcommand{\SdW}{\Sigma^{\bl} (W)}

\newcommand{\Tuk}{\tn{T}_{k}}
\newcommand{\Auk}{\Lambda_{k}}
\newcommand{\Suk}{\Sigma_{k}}
\newcommand{\TukM}{\tn{T}_{k} M}
\newcommand{\AukM}{\Lambda_{k} M}
\newcommand{\SukM}{\Sigma_{k} M}
\newcommand{\TukV}{\tn{T}_{k} (V)}
\newcommand{\AukV}{\Lambda_{k} (V)}
\newcommand{\SukV}{\Sigma_{k} (V)}
\newcommand{\TukW}{\tn{T}_{k} (W)}
\newcommand{\AukW}{\Lambda_{k} (W)}
\newcommand{\SukW}{\Sigma_{k} (W)}

\newcommand{\TenkV}{\tn{Ten}^{k}(V)}
\newcommand{\Tenk}{\tn{Ten}^{k}}
\newcommand{\TenukV}{\tn{Ten}_{k}(V)}
\newcommand{\Tenuk}{\tn{Ten}_{k}}

\newcommand{\TenV}{\tn{Ten}(V)}
\newcommand{\Ten}{\tn{Ten}}
\newcommand{\TendV}{\tn{Ten}^{\bl}(V)}
\newcommand{\Tend}{\tn{Ten}^{\bl}}
\newcommand{\TenuV}{\tn{Ten}(V)}
\newcommand{\Tenu}{\tn{Ten}}
\newcommand{\shF}{\shf{F}}
\newcommand{\shV}{\shf{V}}
\newcommand{\shG}{\shf{G}}
\newcommand{\shH}{\shf{H}}
\newcommand{\shT}{\shf{T}}
\newcommand{\shO}{\shf{O}}
\newcommand{\shOX}{\shf{O}_{X}}
\newcommand{\shOn}{\shf{O}(n)}
\newcommand{\shOm}{\shf{O}(m)}
\newcommand{\shOd}{\shf{O}(d)}
\newcommand{\shOdi}{\shf{O}(d_{i})}
\newcommand{\shL}{\shf{L}}
\newcommand{\ZX}{\Z^{X}}
\newcommand{\ZY}{\Z^{Y}}
\newcommand{\ZC}{\Z^{C}}
\newcommand{\DX}{D_{X}}
\newcommand{\ZpioX}{\Z^{\pioX}}
\newcommand{\ZqX}{\Z^{\qX}}
\newcommand{\rhoX}{\rho_{X}}
\newcommand{\etaX}{\eta_{X}}
\newcommand{\thetaX}{\theta_{X}}
\newcommand{\psiX}{\psi_{X}}
\newcommand{\muX}{\mu_{X}}
\newcommand{\nuX}{\nu_{X}}
\newcommand{\muXz}{\mu_{X}^{0}}
\newcommand{\nuXz}{\nu_{X}^{0}}
\newcommand{\muXo}{\mu_{X}^{1}}
\newcommand{\muC}{\mu_{C}}
\newcommand{\muCz}{\mu_{C}^{0}}
\newcommand{\muCo}{\mu_{C}^{1}}
\newcommand{\muY}{\mu_{Y}}
\newcommand{\muYo}{\mu_{Y}^{0}}
\newcommand{\unuX}{\ul{\nu}_{X}}
\newcommand{\gammaX}{\gamma_{X}}
\newcommand{\veX}{\ve_{X}}
\newcommand{\deltaX}{\delta_{X}}
\newcommand{\degX}{\deg_{X}}
\newcommand{\alphaX}{\alpha_{X}}
\newcommand{\aX}{a_{X}}
\newcommand{\betaX}{\beta_{X}}
\newcommand{\iotaX}{\iota_{X}}

\newcommand{\utri}{\bigtriangleup}
\newcommand{\TNm}{\tn{T} N}
\newcommand{\cTp}{\tn{T}_{p}^{\vee}}
\newcommand{\cTpM}{\tn{T}_{p}^{\vee}M}
\newcommand{\cTM}{\tn{T}^{\vee} M}
\newcommand{\cTpN}{\tn{T}_{p}^{\vee}N}
\newcommand{\cTN}{\tn{T}^{\vee} N}
\newcommand{\CiM}{C^{\infty}(M)}
\newcommand{\CiN}{C^{\infty}(N)}
\newcommand{\Ci}{C^{\infty}}
\newcommand{\ev}{\tn{ev}}
\newcommand{\trf}{\st{\sim}{\hra}}

\newcommand{\lrai}{\st{\sim}{\lra}}
\newcommand{\rai}{\st{\sim}{\ra}}
\newcommand{\trc}{\st{\sim}{\twoheadrightarrow}}

\newcommand{\dlV}{V^{\vee}}
\newcommand{\dlW}{W^{\vee}}
\newcommand{\dlf}{f^{\vee}}
\newcommand{\dl}[1]{{#1}^{\vee}}
\newcommand{\dldl}[1]{{#1}^{\vee \vee}}
\newcommand{\dldlW}{\dldl{W}}
\newcommand{\dldlV}{\dldl{V}}
\newcommand{\dldlf}{\dldl{f}}
\newcommand{\dlA}{{A}^{\vee}}
\newcommand{\wt}[1]{\widetilde{#1}}

\newcommand{\muon}[1]{\nu^{\ge {#1}}}
\newcommand{\mln}[1]{\nu_{< {#1}}}
\newcommand{\mlen}[1]{\nu_{\le {#1}}}
\newcommand{\mun}[1]{\nu_{{#1}}}
\newcommand{\dpqr}{d^{p,q}_{r}}
\newcommand{\dpqro}{d^{p+r,q-r+1}_{r}}
\newcommand{\dpqrm}{d^{p-r,q+r-1}_{r}}
\newcommand{\Epq}{E^{p,q}}
\newcommand{\Epqr}{E^{p,q}_{r}}
\newcommand{\Epqro}{E^{p+r,q-r+1}_{r}}
\newcommand{\Zpqr}{Z^{p,q}_{r}}
\newcommand{\Bpqr}{B^{p,q}_{r}}
\newcommand{\apqr}{\alpha^{p,q}_{r}}
\newcommand{\bpqr}{\beta^{p,q}_{r}}
\newcommand{\Epqi}{E^{p,q}_{\infty}}

\newcommand{\Hr}{H^{\textnormal{r}}}
\newcommand{\Hl}{H^{\textnormal{l}}}
\newcommand{\Autr}{\textnormal{Aut}_{\textnormal{r}}}
\newcommand{\Autl}{\textnormal{Aut}_{\textnormal{l}}}
\newcommand{\CMt}{K/K_{0}/\Q}

\newcommand{\phn}[1]{\phantom{#1}}
\newcommand{\pul}[1]{\ul{\phantom{#1}}}
\newcommand{\pcl}[1]{\sout{\phantom{#1}}}
\newcommand{\sds}{\, \cdot \,}
\newcommand{\tc}[2]{\textcolor{#1}{#2}}
\newcommand{\pol}[1]{\ol{\phantom{#1}}}
\newcommand{\st}[2]{\stackrel{#1}{#2}}

\makeatletter
\def\equalsfill{$\m@th\mathord=\mkern-7mu
\cleaders\hbox{$\!\mathord=\!$}\hfill
\mkern-7mu\mathord=$}
\makeatother

\newcommand{\lgeqt}[1]{\st{\tn{#1}}{\hbox{\equalsfill}}}
\newcommand{\lgeq}[1]{\st{#1}{\hbox{\equalsfill}}}
\newcommand{\defeq}{\lgeq{\tn{\tiny{def}}}}

\newcommand{\xra}{\xrightarrow}
\newcommand{\xla}{\xleftarrow}
\newcommand{\sxra}[1]{\xra{#1}}
\newcommand{\sxla}[1]{\xla{#1}}

\newcommand{\sra}[1]{\stackrel{#1}{\ra}}
\newcommand{\sla}[1]{\stackrel{#1}{\la}}
\newcommand{\slra}[1]{\stackrel{#1}{\lra}}
\newcommand{\sllra}[1]{\stackrel{#1}{\llra}}
\newcommand{\slla}[1]{\stackrel{#1}{\lla}}

\newcommand{\ira}{\stackrel{\simeq}{\ra}}
\newcommand{\ila}{\stackrel{\simeq}{\la}}
\newcommand{\ilra}{\stackrel{\simeq}{\lra}}
\newcommand{\illa}{\stackrel{\simeq}{\lla}}

\newcommand{\wh}[1]{\widehat{#1}}
\newcommand{\fn}[1]{\footnote{#1}}
\newcommand{\mbf}[1]{\mathbf{#1}}
\newcommand{\modulo}[1]{\; \left( \textnormal{mod} \; {#1} \right)}
\newcommand{\ArtSym}[3]{\left( \frac{{#1} / {#2}}{{#3}} \right)}
\newcommand{\SqrArt}[3]{\left[ \frac{{#1} / {#2}}{{#3}} \right]}
\newcommand{\ArtMap}[2]{\left( \frac{{#1} / {#2}}{\cdot} \right)}
\newcommand{\recmap}[2]{[\, \cdot \,, {#1}/{#2}]}
\newcommand{\rec}[3]{[{#1}, {#2}/{#3}]}
\newcommand{\conj}[1]{\overline{#1}}
\newcommand{\cj}[1]{\overline{#1}}
\newcommand{\leg}[2]{\left( \frac{#1}{#2} \right)}
\newcommand{\sep}[1]{{#1}^{\textnormal{sep}}}
\newcommand{\alg}{\tn{alg}}
\newcommand{\ab}[1]{{#1}^{\textnormal{ab}}}
\newcommand{\abs}[1]{\left|{#1} \right|}
\newcommand{\er}[2]{\textnormal{End}_{#1}(#2)}
\newcommand{\ea}[2]{\textnormal{End}^{0}_{#1}(#2)}
\newcommand{\plh}{( \cdot )}
\newcommand{\absmap}{| \cdot |}
\newcommand{\paar}{\langle \cdot, \cdot \rangle}
\newcommand{\nrmmap}{\rVert \, . \, \rVert}
\newcommand{\nmap}[1]{\rVert {#1} \rVert}
\newcommand{\bnmap}[2]{\langle {#1}, {#2} \rangle}
\newcommand{\sqmap}{[ \, \cdot \, ]}
\newcommand{\smprod}{\; \Pi \;}
\newcommand{\car}{\curvearrowright}
\newcommand{\limm}{\lim\limits}
\newcommand{\pair}{\langle \, , \, \rangle}
\newcommand{\Top}{\tn{\tbf{Top}}}
\newcommand{\TopX}{\tn{\tbf{Top}}(X)}
\newcommand{\nfn}[1]{\rVert #1 \rVert}
\newcommand{\sheaf}[1]{\mathscr{#1}}
\newcommand{\shf}[1]{\mathscr{#1}}
\newcommand{\sh}[1]{\wt{G}}
\newcommand{\sht}[1]{\wt{#1}^{\tr}}
\newcommand{\goth}[1]{\mathfrak{#1}}
\newcommand{\cat}[1]{\textnormal{\textbf{#1}}}
\newcommand{\ideal}[1]{\mathfrak{#1}}
\newcommand{\idl}[1]{\mathfrak{#1}}
\newcommand{\functor}[1]{\mathcal{#1}}
\newcommand{\curly}[1]{\mathcal{#1}}
\newcommand{\cly}[1]{\mathcal{#1}}
\newcommand{\und}[1]{\underline{#1}}
\newcommand{\floor}[1]{\lfloor #1 \rfloor}
\newcommand{\ceil}[1]{\lceil #1 \rceil}
\newcommand{\dirlim}[1]{\lim_{\stackrel{\longrightarrow}{#1}}}
\newcommand{\homol}[3]{H_{#1}(#2, #3)}
\newcommand{\cohom}[3]{H^{#1}(#2, #3)}
\newcommand{\cohomzero}[2]{H^{0}(#1, #2)}
\newcommand{\cohomone}[2]{H^{1}(#1, #2)}
\newcommand{\labeleq}[1]{\stackrel{\textnormal{\tiny{#1}}}{=}}
\newcommand{\diff}[2]{\frac{d{#1}}{d{#2}}}
\newcommand{\mdiff}[3]{\frac{d^{#3}{#1}}{d{#2^{#3}}}}
\newcommand{\pderiv}[2]{\frac{\partial{#1}}{\partial{#2}}}
\newcommand{\pd}{\partial}
\newcommand{\nf}[2]{\nicefrac{#1}{#2}}
\newcommand{\nfh}{\nicefrac{1}{2}}
\newcommand{\nfq}{\nicefrac{1}{4}}
\newcommand{\sth}{{\small\textonehalf}}
\newcommand{\stq}{{\small\textonequarter}}

\newcommand{\pdx}{\frac{\pd}{\pd x}}
\newcommand{\pdxi}{\frac{\pd}{\pd x_{i}}}
\newcommand{\pdxj}{\frac{\pd}{\pd x_{j}}}
\newcommand{\pdxn}[1]{\frac{\pd}{\pd x_{#1}}}

\newcommand{\ord}[2]{\textnormal{ord}_{#2}(#1)}
\newcommand{\mpderiv}[3]{\frac{\partial^{#3}{#1}}{\partial{#2^{#3}}}}
\newcommand{\mpderivb}[3]{\frac{\partial^{#3}{#1}}{\partial{#2}}}
\newcommand{\dotp}[2]{\langle #1, #2 \rangle}
\newcommand{\dotpa}[2]{\langle #1, #2 \rangle_{\tn{a}}}
\newcommand{\dotpm}[2]{\langle #1, #2 \rangle_{\tn{m}}}
\newcommand{\euls}[2]{(#1,  #2)_{\tn{a}}}
\newcommand{\eulsym}{(\pul{x},  \pul{x})}
\newcommand{\eulsyma}{(\pul{x},  \pul{x})_{\tn{a}}}
\newcommand{\dotpsym}{\langle \, \cdot \, \rangle}
\newcommand{\norm}[1]{\Vert #1 \Vert}
\newcommand{\bfm}{(\, \cdot \, , \, \cdot \,)}
\newcommand{\bfmr}[1]{(\, \cdot \, , #1)}
\newcommand{\bfml}[1]{(#1, \, \cdot \,)}
\newcommand{\bfma}{\langle \, \cdot \, , \, \cdot \, \rangle}
\newcommand{\normm}[1]{\rVert #1 \rVert}
\newcommand{\refl}[1]{{#1}^{\textnormal{r}}}
\newcommand{\unit}[1]{{#1}^{\times}}
\newcommand{\tcr}[1]{\textcolor{red}{#1}}
\newcommand{\tcg}[1]{\textcolor{grey}{#1}}
\newcommand{\tbr}[1]{\textcolor{red}{\tbf{#1}{}}}
\newcommand{\tbb}[1]{\textcolor{blue}{\tbf{#1}{}}}
\newcommand{\tbg}[1]{\textcolor{dgreen}{\tbf{#1}{}}}
\newcommand{\remph}[1]{\emph{\tcr{#1}}}
\newcommand{\tp}[1]{{#1}^{\tn{t}}}
\newcommand{\pt}[1]{{}^{\tn{t}}{#1}}
\newcommand{\fns}[1]{\footnotesize{#1}}
\newcommand{\scs}[1]{\scriptsize{#1}}
\newcommand{\tpt}[1]{{}^{\tn{t}}{#1}}
\newcommand{\gbx}[1]{\tcg{\fbox{\tn{\small{#1}}}}}
\newcommand{\gbl}{\tcg{\fbox{\tn{\small{\phn{blah}}}}}}

\newcommand{\tcb}[1]{\textcolor{blue}{#1}}
\newcommand{\NrmF}[2]{\textnormal{N}^{#1}_{#2}}
\newcommand{\prm}{\prime}
\newcommand{\p}{\prime}
\newcommand{\bs}{\backslash}
\newcommand{\ve}{\varepsilon}
\newcommand{\vth}{\vartheta}
\newcommand{\vsig}{\varsigma}
\newcommand{\tn}[1]{\textnormal{#1}}
\newcommand{\tbf}[1]{\textbf{#1}}
\newcommand{\ul}[1]{\underline{#1}}
\newcommand{\ol}[1]{\overline{#1}}
\newcommand{\ub}[1]{\underbrace{#1}}
\newcommand{\ubu}[2]{\underbrace{#1}_{#2}}
\newcommand{\ob}[1]{\overbrace{#1}}
\newcommand{\obo}[2]{\overbrace{#1}^{#2}}
\newcommand{\csf}[1]{\stackrel{\smallfrown}{#1}}
\newcommand{\tcl}[1]{\textcolor{#1}}
\newcommand{\ora}[1]{\overrightarrow{#1}}

\newcommand{\explicitbijection}[4] 
{\begin{array}{ccccccc} 
\rnode{left}{#1} & & & & & & \rnode{right}{#2}  \\\\\\
\end{array} 
\psset{nodesep=3pt, offsetA=3pt, offsetB=3pt} 
\everypsbox{\scriptstyle} 
\ncline{->}{left}{right}\Aput{#3} 
\ncline{->}{right}{left}\Aput{#4}}

\newcommand{\tricommdiag}[6]
{{\xymatrix{
{#1} \ar@{->}[d]_{#6} \ar@{->}[r]^{#4} & {#2}\\
{#3} \ar@{->}[ur]_{#5} &}}}

\newcommand{\hightricommdiag}[6]
{{\xymatrix{
{#1} \ar@{->}[dd]_{#6} \ar@{->}[rr]^{#4} & & {#2}\\
& & \\
{#3} \ar@{->}[uurr]_{#5} & & }}}

\newcommand{\tricommdiagtwo}[6]
{\begin{array}{ccccc} 
\rnode{left}{#1} & & & & \rnode{right}{#2}  \\\\\\
& & \rnode{bottom}{#3} & &
\end{array} 
\psset{nodesep=3pt} 
\everypsbox{\scriptstyle} 
\ncline{<-}{left}{right}\Aput{#4} 
\ncline{<-}{left}{bottom}\Aput{#5}
\ncline{<-}{bottom}{right}\Bput{#6}}

\newcommand{\tricommdiagthree}[7]
{\begin{array}{ccccc} 
& & & & \rnode{top}{#1}  \\\\\\\\
\rnode{left}{#2} & & & & \rnode{right}{#3}
\end{array} 
\psset{nodesep=3pt} 
\everypsbox{\scriptstyle} 
\ncline{->}{top}{left}\Bput{#4} 
\ncline{->}{right}{left}\Bput{#5}
\ncarc[linestyle=dashed,ncurv=0.7,arcangleB=10,arcangleA=10,nodesep=2pt]{<-}{top}{right} \Aput{{#6}}
\ncarc[linestyle=dashed,ncurv=0.7,arcangleB=10,arcangleA=10,nodesep=2pt]{<-}{right}{top} \Aput{{#7}}}

\newcommand{\amultdiag}[6]
{\begin{array}{ccccccc} 
& & & & & & \rnode{top}{#1^{#2}}  \\\\\\
\rnode{midleft}{#1} & & & & & & \rnode{midright}{#1} \\\\\\
& & & & & & \rnode{bottom}{#3}
\end{array} 
\psset{nodesep=3pt} 
\everypsbox{\scriptstyle} 
\ncline{->}{midleft}{top}\Aput{#4^{#2}} 
\ncline{->}{midleft}{midright}\Aput{\iota(#5)}
\ncline{->}{midleft}{bottom}\Bput{#6}
\ncline[linestyle=dashed]{->}{top}{midright}\Bput{f_{#5}}
\ncline[linestyle=dashed]{->}{bottom}{midright}\Aput{g_{#5}}
\ncarc[linestyle=dashed,ncurv=0.7,arcangleB=30,arcangleA=30,nodesep=2pt]{<-}{top}{bottom}\Aput{h_{#5}}}

\newcommand{\hightricommdiagtwo}[6]
{\begin{array}{ccccc} 
\rnode{left}{#1} & & & & \rnode{right}{#2}  \\\\\\\\
& & \rnode{bottom}{#3} & &
\end{array} 
\psset{nodesep=3pt} 
\everypsbox{\scriptstyle} 
\ncline{<-}{left}{right}\Aput{#4} 
\ncline{<-}{left}{bottom}\Aput{#5}
\ncline{<-}{bottom}{right}\Bput{#6}}

\newcommand{\diamcommdiag}[8]
{\begin{array}{ccccc} 
& & \rnode{t}{#1} & & \\\\\\
\rnode{l}{#4} & & & & \rnode{r}{#2}  \\\\\\
& & \rnode{b}{#3} & &
\end{array} 
\psset{nodesep=3pt} 
\everypsbox{\scriptstyle} 
\ncline{->}{t}{r}\Aput{#5} 
\ncline{->}{r}{b}\Aput{#6}
\ncline{->}{l}{b}\Bput{#7}
\ncline{->}{t}{l}\Bput{#8}}

\newcommand{\fibrediag}[7] 
{\begin{array}{ccccc} 
& & \rnode{t}{#3} & & \\\\\\\\
& & \rnode{m}{#1 \times_{#4} #2} & & \\\\\\
\rnode{l}{#1} & & & & \rnode{r}{#2}  \\\\\\
& & \rnode{b}{#4} & &
\end{array} 
\psset{nodesep=3pt} 
\everypsbox{\scriptstyle} 
\ncline[linestyle=dashed]{->}{t}{m}\Aput{#7} 
\ncline{->}{r}{b}
\ncline{->}{l}{b}
\ncline{->}{m}{r}\Bput{p_{2}}
\ncline{->}{m}{l}\Aput{p_{1}}
\ncarc{<-}{l}{t}\Aput{#5}
\ncarc{->}{t}{r}\Aput{#6}}

\newcommand{\newfibrediag}[8] 
{\begin{array}{ccccc} 
& & \rnode{t}{#3} & & \\\\\\\\
& & \rnode{m}{#4} & & \\\\\\
\rnode{l}{#1} & & & & \rnode{r}{#2}  \\\\\\
& & \rnode{b}{#5} & &
\end{array} 
\psset{nodesep=3pt} 
\everypsbox{\scriptstyle} 
\ncline[linestyle=dashed]{->}{t}{m}\Aput{#8} 
\ncline{->}{r}{b}
\ncline{->}{l}{b}
\ncline{->}{m}{r}\Bput{p_{2}}
\ncline{->}{m}{l}\Aput{p_{1}}
\ncarc{<-}{l}{t}\Aput{#6}
\ncarc{->}{t}{r}\Aput{#7}}

\newcommand{\specfibrediag}[7] 
{\begin{array}{ccccc} 
& & \rnode{t}{#4} & & \\\\\\\\
& & \rnode{m}{\Spec(#1 \otimes_{#3} #2)} & & \\\\\\
\rnode{l}{\Spec(#1)} & & & & \rnode{r}{\Spec(#2)}  \\\\\\
& & \rnode{b}{\Spec(#3)} & &
\end{array} 
\psset{nodesep=3pt} 
\everypsbox{\scriptstyle} 
\ncline[linestyle=dashed]{->}{t}{m}\Aput{#7} 
\ncline{->}{r}{b}
\ncline{->}{l}{b}
\ncline{->}{m}{r}\Bput{p_{2}}
\ncline{->}{m}{l}\Aput{p_{1}}
\ncarc{<-}{l}{t}\Aput{#5}
\ncarc{->}{t}{r}\Aput{#6}}

\newcommand{\contrafibrediag}[7] 
{\begin{array}{ccccc} 
& & \rnode{t}{\Gamma(#4, \sheaf{O}_{#4})} & & \\\\\\\\
& & \rnode{m}{#1 \otimes_{#3} #2} & & \\\\\\
\rnode{l}{#1} & & & & \rnode{r}{#2}  \\\\\\
& & \rnode{b}{#3} & &
\end{array} 
\psset{nodesep=3pt} 
\everypsbox{\scriptstyle} 
\ncline[linestyle=dashed]{<-}{t}{m}\Aput{#7} 
\ncline{<-}{r}{b}
\ncline{<-}{l}{b}
\ncline{<-}{m}{r}\Bput{p_{2}^{*}}
\ncline{<-}{m}{l}\Aput{p_{1}^{*}}
\ncarc{->}{l}{t}\Aput{#5^{*}}
\ncarc{<-}{t}{r}\Aput{#6^{*}}}

\newcommand{\revtricommdiag}[6]
{\begin{array}{ccccc} 
\rnode{i}{#1} & & & & \rnode{j}{#2}  \\\\\\
& & \rnode{g}{#3} & &
\end{array} 
\psset{nodesep=3pt} 
\everypsbox{\scriptstyle} 
\ncline{<-}{i}{j}\Aput{#4} 
\ncline{<-}{j}{g}\Aput{#5}
\ncline{<-}{i}{g}\Bput{#6}}

\newcommand{\catequivcommdiag}[6] 
{\begin{array}{cccccc} 
& & & \rnode{tu}{#1(#3)} & & \\\\\\
\rnode{tl}{#3} & & & & & \rnode{tr}{(#2 \circ #1)(#3)} \\\\\\\\
& & & \rnode{bu}{#1(#4)} & & \\\\\\
\rnode{bl}{#4} & & & & & \rnode{br}{(#2 \circ #1)(#4)}  
\end{array} 
\psset{nodesep=3pt} 
\everypsbox{\scriptstyle} 
\ncline{->}{tl}{tu}\Aput{#1} 
\ncline{->}{tu}{tr}\Aput{#2}
\ncline{->}{tl}{tr}\Aput{#6(#3)}
\ncline{->}{bl}{bu}\Aput{#1} 
\ncline{->}{bu}{br}\Aput{#2}
\ncline{->}{bl}{br}\Bput{#6(#4)}
\ncline{->}{tl}{bl}\Bput{#5} 
\ncline[linestyle=dashed]{->}{tu}{bu}\Aput{#1(#5)}
\ncline{->}{tr}{br}\Aput{(#2 \circ #1)(#5)}}

\newcommand{\rectcommdiag}[8]
{\xymatrix{
{#1} \ar@{->}[d]_{#8} \ar@{->}[r]^{#5} & {#2} \ar@{->}[d]^{#6} \\
{#3} \ar@{->}[r]_{#7} & {#4}}}

\newcommand{\highrectcommdiag}[8]
{\xymatrix{
{#1} \ar@{->}[dd]_{#8} \ar@{->}[rr]^{#5} & & {#2} \ar@{->}[dd]^{#6} \\
& & \\
{#3} \ar@{->}[rr]_{#7} & & {#4}}}

\newcommand{\sheafrectcommdiag}[5] 
{\begin{array}{ccccc} 
\rnode{tl}{#1(#3)} & & & & \rnode{tr}{#2(#3)}  \\\\\\
\rnode{bl}{#1(#4)} & & & & \rnode{br}{#2(#4)}
\end{array} 
\psset{nodesep=3pt} 
\everypsbox{\scriptstyle} 
\ncline{->}{tl}{tr}\Aput{#5(#3)} 
\ncline{->}{tr}{br}\Aput{\rho^{#2}_{#3 #4}}
\ncline{->}{bl}{br}\Aput{#5(#4)}
\ncline{->}{tl}{bl}\Aput{\rho^{#1}_{#3 #4}}}

\newcommand{\rectcommdashdiag}[9]
{\begin{array}{ccccc} 
\rnode{tl}{#1} & & & & \rnode{tr}{#2}  \\\\\\
\rnode{bl}{#4} & & & & \rnode{br}{#3}
\end{array} 
\psset{nodesep=3pt} 
\everypsbox{\scriptstyle} 
\ncline{->}{tl}{tr}\Aput{#5} 
\ncline{->}{tr}{br}\Aput{#6}
\ncline{->}{bl}{br}\Aput{#7}
\ncline{->}{tl}{bl}\Bput{#8}
\ncline[linestyle=dashed]{->}{bl}{tr}\Aput{#9}}

\newcommand{\contrahighrectcommdiag}[8]
{\begin{array}{ccccc} 
\rnode{tl}{#1} & & & & \rnode{tr}{#2}  \\\\\\\\
\rnode{bl}{#4} & & & & \rnode{br}{#3}
\end{array} 
\psset{nodesep=3pt} 
\everypsbox{\scriptstyle} 
\ncline{->}{tl}{tr}\Aput{#5} 
\ncline{<-}{tr}{br}\Aput{#6}
\ncline{->}{bl}{br}\Aput{#7}
\ncline{<-}{tl}{bl}\Aput{#8}}

\newcommand{\sheafhighrectcommdiag}[5] 
{\begin{array}{ccccc} 
\rnode{tl}{#1(#3)} & & & & \rnode{tr}{#2(#3)}  \\\\\\\\
\rnode{bl}{#1(#4)} & & & & \rnode{br}{#2(#4)}
\end{array} 
\psset{nodesep=3pt} 
\everypsbox{\scriptstyle} 
\ncline{->}{tl}{tr}\Aput{#5(#3)} 
\ncline{->}{tr}{br}\Aput{\rho^{#2}_{#3 #4}}
\ncline{->}{bl}{br}\Aput{#5(#4)}
\ncline{->}{tl}{bl}\Aput{\rho^{#1}_{#3 #4}}}

\newcommand{\dirlimdiag}[6] 
{\begin{array}{ccccccc} 
\rnode{tl}{#1_{#4}} & & & & & & \rnode{tr}{#1_{#5}}  \\\\\\\\
& & & \rnode{m}{#1} & & & \\\\\\\\
& & & \rnode{b}{#2} & & &
\end{array} 
\psset{nodesep=3pt} 
\everypsbox{\scriptstyle} 
\ncline{->}{tl}{tr}\Aput{#3_{#4 #5}} 
\ncline{->}{tr}{m}\Bput{#3^{#1}_{#5}}
\ncline{->}{tl}{m}\Aput{#3^{#1}_{#4}}
\ncarc{<-}{b}{tl}\Aput{#3^{#2}_{#4}}
\ncline{->}{m}{b}\Aput{#6}
\ncarc{->}{tr}{b}\Aput{#3^{#2}_{#5}}}

\newcommand{\dirlimdiagdash}[6] 
{\begin{array}{ccccccc} 
\rnode{tl}{#1_{#4}} & & & & & & \rnode{tr}{#1_{#5}}  \\\\\\\\
& & & \rnode{m}{#1} & & & \\\\\\\\
& & & \rnode{b}{#2} & & &
\end{array} 
\psset{nodesep=3pt} 
\everypsbox{\scriptstyle} 
\ncline{->}{tl}{tr}\Aput{#3_{#4 #5}} 
\ncline{->}{tr}{m}\Bput{#3^{#1}_{#5}}
\ncline{->}{tl}{m}\Aput{#3^{#1}_{#4}}
\ncarc{<-}{b}{tl}\Aput{#3^{#2}_{#4}}
\ncline[linestyle=dashed]{->}{m}{b}\Aput{#6}
\ncarc{->}{tr}{b}\Aput{#3^{#2}_{#5}}}

\newcommand{\twomatrix}[4]
{\left(  \begin{array}{cc} 
\rnode{tl}{#1} & \rnode{tr}{#2}  \\
\rnode{bl}{#3} & \rnode{br}{#4}
\end{array} \right)}

\newcommand{\bigmatrix}[9]
{\left(  \begin{array}{cccc} 
\rnode{tl}{#1} & \rnode{tm}{#2}  & \cdots & \rnode{tr}{#3} \\
\rnode{ml}{#4} & \rnode{mm}{#5}  & \cdots & \rnode{mr}{#6} \\
\vdots & \vdots  & \ddots & \vdots \\
\rnode{bl}{#7} & \rnode{bm}{#8}  & \cdots & \rnode{br}{#9}
\end{array} \right)}

\newcommand{\bigmatrixb}[9]
{\left(  \begin{array}{cccc} 
\rnode{tl}{#1} & \rnode{tm}{#2}  & \cdots & \rnode{tr}{#3} \\\\
\rnode{ml}{#4} & \rnode{mm}{#5}  & \cdots & \rnode{mr}{#6} \\\\
\vdots & \vdots  & \ddots & \vdots \\\\
\rnode{bl}{#7} & \rnode{bm}{#8}  & \cdots & \rnode{br}{#9}
\end{array} \right)}

\newcommand{\mathsmap}[5] 
{\begin{array}{cccc}
#1 \, : & #2 & \longrightarrow & #3\\
& #4 & \longmapsto & #5
\end{array}}

\newcommand{\mathsmapbi}[7] 
{\begin{array}{cccc}
#1 \, : & #2 & \longrightarrow & #3 \\
& #4 & \longmapsto & #5 \\
& #6 & \longmapsfrom & #7 \\
\end{array}}

\newcommand{\mathsmapxy}[5] 
{\xymatrix{
{#1}: & {#2} \ar@{->}[r] & {#3} \\
& {#4} \ar@{|->}[r] & {#5}}}

\newcommand{\matrixtwo}[4] 
{\left( \begin{array}{ccc}
{#1} & {#2} \\
{#3} & {#4} 
\end{array}\right)}

\newcommand{\mathsmaptwo}[4] 
{\begin{array}{ccc}
#1 & \longrightarrow & #2\\
#3 & \longmapsto & #4
\end{array}}

\newcommand{\mathsmaptwox}[4] 
{\begin{eqnarray*}
\xymatrix{
#1 \ar@{->}[rr] && #2 \\
#3 \ar@{|->}[rr] && #4}
\end{eqnarray*}}

\newcommand{\mathsmapvert}[5] 
{\xymatrix{
{#2} \ar@{->}[dd]_{{#1}} & \ni & {#4} \ar@{|->}[dd]\\
& & \\
{#3} & \ni & {#5}}}

\newcommand{\doublemathsmap}[7] 
{\begin{array}{cccc}
#1: & #2 & \longrightarrow & #3\\
& #4 & \longmapsto & #5 \\
& #6 & \longmapsto & #7
\end{array}}

\newcommand{\ses}[7]
{\xymatrix{#1 \ar[r] & #2 \ar[r]^{#6} & #3 \ar[r]^{#7} & #4 \ar[r] & #5}}

\newcommand{\sestwo}[5]
{\xymatrix{#1 \ar[r] & #2 \ar[r] & #3 \ar[r] & #4 \ar[r] & #5}}

\newcommand{\sesthree}[5]
{\xymatrix{0 \ar[r] & #1 \ar[r]^{#4} & #2 \ar[r]^{#5} & #3 \ar[r] & 0}}

\newcommand{\splitses}[8]
{\xymatrix{#1 \ar[r] & #2 \ar[r]^{#6} & #3 \ar@/^/[r]^{#7} & #4 \ar@/^/[l]^{#8} \ar[r] & #5}}

\begin{titlepage}

\centering

\leavevmode\\
\leavevmode\\
\leavevmode\\

\huge{\tbf{The Voevodsky motive of a rank one semiabelian variety}}\\[1.2cm]
\large{\scshape Dissertation zur Erlangung des Doktorgrades\\ der Fakult{\"a}t f{\"u}r Mathematik und Physik\\ der Albert-Ludwigs-Universit{\"a}t \\Freiburg im Breisgau}\\[7cm]
\Large{vorgelegt von}\\
\Large{Stephen Enright-Ward}\\[1.2cm]
\Large{November 2012}

\end{titlepage}

\begin{flushleft}

\phn{a}

\begin{tabular}{lll}
& & \\
& & \\
& & \\
& & \\
& & \\
& & \\
& & \\
& & \\
& & \\
& & \\
{\scshape Dekan:} & & Prof. Dr. Michael Rů\v{z}i\v{c}ka \\
& & \\
{\scshape 1. Gutachter:} & & Prof. Dr. Annette Huber-Klawitter \\
& & \\
{\scshape 2. Gutachter:} & & Prof. Dr. Uwe Jannsen \\
& & \\
{\scshape  Datum der Promotion:} & & 17.April 2013 \\
\end{tabular}

\end{flushleft}

\thispagestyle{empty}

\newpage 

\phn{a}
\vspace{\stretch{2}}

\begin{center}
{\huge To Joan Enright (1915--2011)} 
\end{center}
\vspace{\stretch{4}}

\thispagestyle{empty}

\frontmatter

\tableofcontents

\include{introduction_FIN} 

\mainmatter

\chapter{Categories of motives}

Let $k$ be a perfect field. In this chapter we will recall the definitions, main properties, and relationships between various categories of motives over $\Speck$ used in the sequel. We shall see that all categories are idempotent complete and have a tensor structure, which allows us to define symmetric powers $\Symn$. The symmetric powers on $\Chowek$, $\DMegk \sx \Q$ and $\DMeetkQ$ are compatible in the sense that there are tensor functors:
\begin{align*}
& \Psi: \Chowek \lra \DMegk \sx \Q \qaq \\
& \Phi: \Chowek \lra \DMeetkQ
\end{align*}
commuting with $\Symn$ for all non-negative integers $n$. We stress that all material presented here is revision: \tbf{none of the concepts or results in this chapter is due to the author}. The material on symmetric powers in general additive idempotent complete $\Q$-linear categories is a superficial generalisation of the classical construction of symmetric powers in the category of $R$-modules, as treated for example in \cite{BAlg}. The definition of the category of Chow motives is originally due to Grothendieck, and details of its construction can be found for example in the article \cite{Scholl} or the book \cite{Andre}. The construction and results concerning the categories $\DMegk$ and $\DMettkR$ are all due to Voevodsky, and the main reference is \cite{TMF}. 

Let us now describe more precisely the content of each section. It is convenient to introduce the following notation.

\begin{notation} \label{sxdef}
For any additive category $\catc$, we denote by $\catc \ox \Q$ the category whose objects are identical to those of $\catc$, and whose morphisms are given by:
\begin{align*}
\Mor_{\catc \ox \Q}(X, Y) := \Mor_{\catc} (X, Y) \ox \Q,
\end{align*}
and borrow the notation  $\catc \sx \Q$ from Kahn and Barbieri-Viale \cite[Definition 1.1.3]{BVK} to denote the idempotent completion of $\catc \ox \Q$.
\end{notation}

\tbf{(i)} In the first section, we recall the definitions of the symmetric powers $\Symn$ for an idempotent complete, $\Q$-linear tensor category, following Biglari \cite[\S 2.5]{Big}. This general construction will be used to define symmetric powers in the ($\Q$-linear versions of the) categories of motives discussed in subsequent sections.

\tbf{(ii)} In the second section, we recall the constructions and main properties of the categories $\Chowek$ and $\Chowk$ of Chow motives. Note that there are two conflicting sets of conventions on the definition of Chow motives, which differ in two points. Classically, categories of Chow motives are (i) $\Q$-linear, and defined in such a way that (ii) the natural functor $h: \SmProjk \ra \Chowk$ is contravariant. This is the case for example in the papers \cite{Scholl} of Scholl, \cite{Ku1} of K{\"u}nnemann, and \cite{DeMu} of Deninger and Murre, upon which our exposition is closely modelled. However, under the modern convention, used by Voevodsky \cite[Preamble to Proposition 2.1.4]{TMF}, Chow motives are (i) $\Z$-linear, and defined in such a way that (ii) the natural functor $h: \SmProjk \ra \Chowk$ is covariant. We are ultimately interested in Voevodsky motives with $\Q$-coefficients, and the category of Chow motives is a stepping stone. Hence we mix the conventions, defining Chow motives to be covariant, but $\Q$-linear.

\tbf{(iii)} In the third section, we recall the definition of the category $\DMegk$ of effective geometric Voevodsky motives, as well as that of the category $\DMettkR$ of Voevodsky motivic complexes, where $\tau$ stands for either the {\etl} or Nisnevich topologies on $\Smk$, and $R$ is a ring containing $\Z$. These categories were first defined by Voevodsky in \cite{TMF}, and are both tensor triangulated and idempotent complete. Consequently, the $\Q$-linear categories:
\begin{align*} 
\DMettkQ \qaq \DMegk \sx \Q
\end{align*}
each possess symmetric powers as defined in the first section, where the symbol ``$\DMegk \sx \Q$'' denotes the idempotent completion of $\DMegk \ox \Q$. We also recall the construction of an exact tensor functor:
\begin{align*}
j: \DMegk \sx \Q \lra \DMeetkQ.
\end{align*}

\tbf{(iv)} In the fourth section, we will see that there is a natural tensor functor $\Psi: \Chowek \ra \DMegk \sx \Q$ which commutes with symmetric powers; composing with $j: \DMegk \sx \Q \ra \DMeetkQ$ then yields a tensor functor: 
\begin{align*}
\Phi: \Chowek \lra \DMeetkQ,
\end{align*}
which will function as a bridge between the worlds of Chow and Voevodsky motives in the fourth chapter. 

\section{Alternating and symmetric powers}

Let $\catc$ be an additive, $\Q$-linear, and idempotent complete tensor category. The goal of this section is to define 
certain endofunctors $\Symn$, for each non-negative $n$, which we refer to as \emph{symmetric powers}. They generalise the familiar constructions of the symmetric powers associated for example to a $\Q$-vector space. In particular, we will see below that these constructions may be applied in the categories of Chow motives $\Chowek$ and $\Chowk$, as well as for the categories $\DMegk \sx \Q$ and $\DMeetkQ$ of Voevodsky motives. The exposition follows that of Biglari \cite[\S 2.5]{Big} in large part, although we note that the categories we work in are not required to be triangulated.

\subsection{The symmetric and alternating functors}

\begin{definition}
We refer to a \emph{symmetric monoidal category} as defined in \cite[\S 1.1]{Ke1} as a \emph{tensor category}. 
\end{definition}

\begin{notation} 
For a positive integer $n$, we denote by $\Sym(n)$ the group of permutations on $n$ letters.
\end{notation}

\begin{notation} \label{tensym}
Let $(\catc, \ox, \one)$ be a $\Q$-linear tensor category and $n$ a positive integer. Let us agree to the convention that the tensor product of three or more objects $X_{1} \ox X_{2} \ox \cdots \ox X_{n}$ means by definition the object:
\begin{align*}
(( \ldots ((X_{1} \ox X_{2}) \ox X_{3}) \ox \cdots \ox X_{n}) \qquad \tn{(left bracketing).}
\end{align*}
When all $X_{i}$ are the same object $X$, we denote this by $X^{\ox n}$. We employ the convention that zeroth tensor powers of all objects and morphisms are equal to $\one$ and $1_{\one}$, respectively. 
\end{notation}

\begin{notation} \label{sigmanot}
As the permutation group $\Sym(n)$ is generated by transpositions, the transposition and associativity isomorphisms: 
\begin{align*}
t_{XY}: X \ox Y \ilra Y \ox X \qaq a_{X,Y,Z}: X \ox (Y \ox Z) \ilra (X \ox Y) \ox Z
\end{align*} 
which are built into the definition of $\catc$ induce well-defined isomorphisms:
\begin{align*}
X_{1} \ox X_{2} \ox \cdots \ox X_{n} \ilra X_{\sigma(1)} \ox X_{\sigma(2)} \ox \cdots \ox X_{\sigma(n)}
\end{align*} 
for any $n$ objects $X_{1}, \ldots, X_{n}$ in $\catc$, which we think of as reordering of factors. This gives a representation: 
\begin{align*}
\Sym(n) \lra \End(X^{\ox n}).
\end{align*}
We denote by $\sigma_{X}$ the endomorphism of $X^{\ox n}$ associated to $\sigma \in \Sym(n)$ by this representation. We recall that since $\catc$ is $\Q$-linear, $\End(X^{\ox n})$ is a $\Q$-algebra.
\end{notation} 

\begin{notation} \label{ansn}
Let $\catc$ be a $\Q$-linear tensor category, let $X$ be an object in $\catc$, and let $n$ denote a non-negative integer. We define\tn{:}
\begin{align*}
s_{X}^{n} & := \left\{ \begin{array}{cc}
1_{\one}, & \tn{ if } n = 0, \\
\frac{1}{n!} \left( \sum_{\sigma \in \Sym(n)} \sigma_ {X} \right), & \tn{ if } n > 0.
\end{array} \right.
\end{align*}
We will sometimes omit the subscript ``$X$'', where this is clear from the context. 
\end{notation}

\begin{proposition}[\tbf{Classical}] \label{imsn}
Let $\catc$ be a $\Q$-linear tensor category, $X$ an object in $\catc$, and let $n$ denote a positive integer. Then the morphism\tn{:}
\begin{align*}
s_{X}^{n}: X^{\ox n} \lra X^{\ox n} 
\end{align*}
of Notation \ref{ansn} is idempotent.
\end{proposition}
\begin{prooff}
This is well-known for example in the category of $\Q$-vector spaces. The proof in this more general context is formally similar. Indeed, for any $\sigma \in \Sym(n)$, we have:
\begin{align*}
s^{n}_{X} \circ \sigma_{X} = \frac{1}{n!} \sum_{\tau \in \Sym(n)} (\tau_{X} \circ \sigma_{X}) = 
\frac{1}{n!} \sum_{\tau \in \Sym(n)} \tau_{X} = s^{n}_{X},
\end{align*}
since multiplying all the $\tau_{X}$ by $\sigma_{X}$ merely rearranges the order of the summands. It follows that:
\begin{align*}
(s^{n}_{X})^{2} & = \frac{1}{n!} \sum_{\sigma \in \Sym(n)} s^{n}_{X} \circ \sigma_{X} 
= \frac{1}{n!} \sum_{\sigma \in \Sym(n)} s^{n}_{X} \\
& = \frac{1}{n!} \Big( n! \cdot s^{n}_{X} \Big) = s^{n}_{X},
\end{align*}
showing that $s^{n}_{X}$ is idempotent.
\end{prooff}

\begin{definition}[\tbf{Classical}]  \label{symaltn}
Let $\catc$ be an idempotent complete $\Q$-linear tensor category, $X$ an object in $\catc$, and let $n$ denote a positive integer. We define\tn{:}
\begin{align*}
\Sym^{n}(X) := \im \Big( s_{X}^{n} \Big). 
\end{align*}
This image exists, due to Proposition \ref{imsn} and idempotent completeness of the category $\catc$. Thus by definition the endomorphism $s_{X}^{n}$ of $X^{\ox n}$ can be thought of as a morphism from $X^{\ox n}$ to $\Symn(X)$, which we denote by:
\begin{align*}
\pi_{X}^{n}: X^{\ox n} \lra \Sym^{n}(X).
\end{align*}
We refer to this as the \emph{canonical projection associated to $\Sym^{n}(X)$}. By definition of a categorical image, we also have morphism:
\begin{align*}
\iota_{X}^{n}: \Sym^{n}(X) \lra X^{\ox n},
\end{align*}
which we refer to as the \emph{canonical embedding associated to $\Sym^{n}(X)$}. By construction, we have 
$\pi_{X}^{n} \circ \iota_{X}^{n} = 1_{\Symn(X)}$.
\end{definition} 

\begin{remark} \label{kerpq}
Let $\catc$ be an idempotent complete $\Q$-linear tensor category, $X$ an object in $\catc$, and let $n$ denote a positive integer. One can alternatively define $\Symn(X)$ to be $\ker(p_{X}^{n})$, where: $p_{X}^{n} := 1_{X} - s_{X}^{n}$. To see that this is equivalent to Definition \ref{symaltn}, we first show that $p^{n} := p_{X}^{n}$ is an idempotent in $\End(X)$. Indeed, by bilinearity of the composition in the additive category $\catc$, we have in $\End(X)$ the relation: 
\begin{align*}
(p^{n}_{X})^{2} & = \left( 1_{X} - s^{n}_{X} \right)^{2} = 1_{X} - 2 s^{n}_{X} + (s^{n}_{X})^{2} \\
& = 1_{X} - s^{n}_{X},
\end{align*}
where the last equality follows from idempotency of $s^{n}_{X}$ (Proposition \ref{imsn}). Thus $s_{n}$ and $p_{n}$ are mutually orthogonal idempotents, from which it follows that $\ker(p_{X}^{n}) = \im(s_{X}^{n})$.
\end{remark}

\begin{proposition}[\tbf{Classical}]  \label{symunprop}
Let $X$ and $Y$ be objects of an idempotent complete $\Q$-linear tensor category $\catc$, and let $n$ be a non-negative integer. For any morphism $f: X^{\ox n} \ra Y$ satisfying $f \circ \sigma_{X} = f$ for all $\sigma \in \Sym(n)$, and for any morphism $g: Y \ra X^{\ox n}$ satisfying $\sigma_{X} \circ g = g$ for all $\sigma \in \Sym(n)$, there exist unique morphisms $\wt{f}: \Symn(X) \ra Y$ and
$\wt{g}: Y \ra \Symn(X)$ such that the diagrams\tn{:}
\begin{align*}
\xymatrix{
X^{\ox n} \ar@{->}[rr]^(0.45){\pi^{n}_{X}} \ar@{->}[ddrr]_(0.5){f} & & 
\Symn(X) \ar@{->}[dd]^(0.5){\wt{f}} \\
& & \\
& & Y} \qaq
\xymatrix{
Y \ar@{->}[dd]_(0.5){\wt{g}} \ar@{->}[ddrr]^(0.5){g} & & \\
& & \\
\Symn(X) \ar@{->}[rr]_(0.55){\iota^{n}_{X}} & & X^{\ox n}}
\end{align*}
commute.
\end{proposition}
\begin{prooff}
This follows from the universal property of the categorical image $\Symn(X) = \im(s_{X}^{n})$.
\end{prooff}

\begin{proposition}[\tbf{Classical}]  \label{symnfunc}
Let $X$ and $Y$ be objects of an idempotent complete $\Q$-linear tensor category $\catc$, and let $n$ be a positive integer. For any morphism $f: X \ra Y$, there exists a unique morphism\tn{:} 
\begin{align*}
\Symn(f) : \Symn(X) \ra \Symn(Y)
\end{align*}
such that the diagram\tn{:}
\begin{align*}
\xymatrix{
X^{\ox n} \ar@{->}[rr]^(0.45){\pi^{n}_{X}} \ar@{->}[dd]_(0.5){f^{\ox n}} & & 
\Symn(X) \ar@{->}[dd]^(0.5){\Symn(f)} \ar@{->}[rr]^(0.55){\iota^{n}_{X}} & & X^{\ox n} \ar@{->}[dd]^(0.5){f^{\ox n}} \\
& & & & \\
Y^{\ox n} \ar@{->}[rr]_(0.45){\pi^{n}_{Y}} & & \Symn(Y) \ar@{->}[rr]_(0.55){\iota^{n}_{Y}} & & Y^{\ox n}}
\end{align*}
commutes. Consequently, $\Symn$ is an endofunctor $\catc \ra \catc$.
\end{proposition}
\begin{prooff}
Evidently $\pi_{Y}^{n} \circ f^{\ox n} \circ \sigma_{X} = \pi_{Y}^{n} \circ f^{\ox n}$, and so by Proposition 
\ref{symunprop} there is a unique morphism $\Symn(f) : \Symn(X) \ra \Symn(Y)$ such that the left-hand square 
commutes. The proof that the right-hand square commutes is similar.
\end{prooff}

\begin{notation} \label{symnot}
For any object $X$ of an idempotent complete $\Q$-linear tensor category $\catc$, we use the notation:
\begin{align*}
\Sym(X) := \bos_{n \ge 0} \Symn(X),
\end{align*}
provided this direct sum exists. By the universal property of the direct sum, this is extends to an endofunctor of $\Sym: \catc \ra \catc$.
\end{notation}

\begin{remark} 
From now on, we will always assume that the direct sum $\os$ in an idempotent complete $\Q$-linear tensor category $\catc$ distributes over the tensor product $\ox$. 
\end{remark}

\begin{notation} \label{psimndef}
Let $Z$ be an object of an idempotent complete $\Q$-linear tensor category $\catc$, and let $m$ and $n$ be non-negative integers. We denote by:
\begin{align*}
\psi_{Z}^{m,n}: \Symm(Z) \ox \Symn(Z) \lra \Sym^{m+n}(Z)
\end{align*}
the unique morphism of Proposition \ref{symunprop} such that the diagram:
\begin{align*}
\xymatrix{
Z^{\ox m} \ox Z^{\ox n} \ar@{=}[dd] \ar@{->}[rr]^(0.43){\pi^{m}_{Z} \ox \pi^{n}_{Z}} & & 
\Symm(Z) \ox \Symn(Z) \ar@{->}[dd]^(0.5){\psi_{Z}^{m,n}} \\
& & \\
Z^{\ox N} \ar@{->}[rr]_(0.45){\pi^{m+n}_{Z}} & & \Sym^{m+n}(Z)} 
\end{align*}
commutes. By distributivity of $\ox$ over $\os$, the tensor product $\Sym(Z) \ox \Sym(Z)$ is canonically isomorphic to the direct sum of the $\Symm(Z) \ox \Symn(Z)$; let $i_{Z}^{m,n}: \Symm(Z) \ox \Symn(Z) \ra \Sym(Z) \ox \Sym(Z)$ be the canonical inclusion morphism. We denote now by: 
\begin{align*}
\psi_{Z}: \Sym(Z) \ox \Sym(Z) \lra \Sym(Z)
\end{align*}
the unique morphism induced by the universal property of the direct sum which makes the diagram:
\begin{align*}
\xymatrix{
\Symm(Z) \ox \Symn(Z) \ar@{->}[rr]^(0.52){i_{Z}^{m,n}} \ar@{->}[dd]_(0.5){\psi_{Z}^{m,n}} & & 
\Sym(Z) \ox \Sym(Z) \ar@{->}[dd]^(0.5){\psi_{Z}} \\
& & \\
\Sym^{N}(Z) \ar@{->}[rr]_(0.53){i_{Z}^{m+n}} & & \Sym(Z)}
\end{align*}
commute, for all non-negative $m$ and $n$.
\end{notation}

\begin{notation} \label{psinot}
Let $X$ and $Y$ be objects of an idempotent complete $\Q$-linear tensor category $\catc$, and let $j_{X}: X \ra X \os Y$ and $j_{Y}: Y \ra X \os Y$ be the natural inclusion morphisms. Let us denote by:
\begin{align*}
\psi_{X, Y}: \Sym(X) \ox \Sym(Y) \ra \Sym(X \os Y) 
\end{align*}
the composite morphism:
\begin{align*}
\Sym(X) \ox \Sym(Y) & \sxra{\Sym(j_{X}) \ox \Sym(j_{Y})} \Sym(X \os Y) \ox \Sym(X \os Y) \\
& \sxra{\psi_{X \os Y}} \Sym(X \os Y).
\end{align*}
Moreover, for all non-negative $m$ and $n$, we denote by:
\begin{align*}
\psi_{X, Y}^{m,n}: \Symm(X) \ox \Symn(Y) \ra \Sym^{m+n}(X \os Y)
\end{align*}
the composite morphism: 
\begin{align*}
\Symm(X) \ox \Symn(Y) & \sxra{\Symm(j_{X}) \ox \Symn(j_{Y})} \Symm(X \os Y) \ox \Symn(X \os Y) \\
& \sxra{\psi^{m, n}_{X \os Y}} \Sym^{m+n}(X \os Y)
\end{align*}
and further, by:
\begin{align*}
\psi_{X, Y}^{N}: \bos_{m+n = N} \Symm(X) \ox \Symn(Y) \lra \Sym^{N}(X \os Y)
\end{align*}
the unique morphism, induced by the universal property of the direct sum, such that the diagram: 
\begin{align*}
\xymatrix{
\Symm(X) \ox \Symn(Y) \ar@{->}[rr]^(0.45){i_{X, Y}^{m,n}} \ar@{->}[ddrr]_(0.5){\psi_{X, Y}^{m,n}} & & 
\DS{\bos_{m+n = N}} \Symm(X) \ox \Symn(Y) \ar@{->}[dd]^(0.5){\psi_{X, Y}^{N}} \\
& & \\
& & \Sym^{N}(X \os Y)}
\end{align*}
commutes. 
\end{notation}

\begin{notation} \label{ximndef}
Let $Z$ be an object of an idempotent complete $\Q$-linear tensor category $\catc$, and let $m$ and $n$ be non-negative integers of sum $N$. We denote by:
\begin{align*}
\xi_{Z}^{m,n}: \Sym^{N}(Z) \lra \Symm(Z) \ox \Symn(Z) 
\end{align*}
the unique morphism of Proposition \ref{symunprop} such that the diagram:
\begin{align*}
\xymatrix{
Z^{\ox (m+n)} \ar@{->}[rrrr]^(0.5){\pi^{m+n}_{Z}} \ar@{->}[dd]_(0.5){s^{m+n}_{Z}} & & & & 
\Sym^{m+n}(Z) \ar@{->}[dd]^(0.5){\xi_{Z}^{m,n}} \\
& & & & \\
Z^{\ox (m+n)} \ar@{=}[rr] & & Z^{\ox m} \ox Z^{\ox n} \ar@{->}[rr]_(0.43){\pi^{m}_{Z} \ox \pi^{n}_{Z}} & & \Symm(Z) \ox \Symn(Z)}
\end{align*}
commutes. Letting again $i_{Z}^{m,n}: \Symm(Z) \ox \Symn(Z) \ra \Sym(Z) \ox \Sym(Z)$ be the canonical inclusion morphism, we write: 
\begin{align*}
\xi_{Z}: \Sym(Z) \lra \Sym(Z) \ox \Sym(Z)
\end{align*}
for the unique morphism induced by the universal property of the direct sum $\Sym(Z) = \bos \Symn(Z)$ such that the diagram:
\begin{align*}
\xymatrix{
\Sym^{m+n}(Z) \ar@{->}[rr]^(0.52){i_{Z}^{m+n}} \ar@{->}[dd]_(0.5){\xi_{Z}^{m,n}} & & \Sym(Z) \ar@{->}[dd]^(0.5){\xi_{Z}} \\
& & \\
\Symm(Z) \ox \Symn(Z) \ar@{->}[rr]_(0.53){i_{Z}^{m,n}} & & \Sym(Z) \ox \Sym(Z)}
\end{align*}
is commutative. 
\end{notation}

\begin{notation} \label{xinot}
Let $X$ and $Y$ be objects of an idempotent complete $\Q$-linear tensor category $\catc$, and let $p_{X}: X \os Y \ra X$ and $p_{Y}: X \os Y \ra Y$ be the natural projection morphisms. Let us denote by:
\begin{align*}
\xi_{X, Y}: \Sym(X \os Y) \lra \Sym(X) \ox \Sym(Y)
\end{align*}
the composite morphism:
\begin{align*}
\Sym(X \os Y) & \sxra{\xi_{X \os Y}} \Sym(X \os Y) \ox \Sym(X \os Y) \\
& \sxra{\Sym(p_{X}) \ox \Sym(p_{Y})} \Sym(X) \ox \Sym(Y)\tn{;}
\end{align*}
moreover, for all non-negative $m$ and $n$ with sum $N$, we write: 
\begin{align*}
\xi_{X, Y}^{m,n}: \Sym^{N}(X \os Y) \lra \Symm(X) \ox \Symn(Y)
\end{align*}
for the composite morphism: 
\begin{align*}
\Sym^{N}(X \os Y) & \sxra{\xi^{m, n}_{X \os Y}} \Sym^{N}(X \os Y) \ox \Sym^{N}(X \os Y) \\
& \sxra{\Symm(p_{X}) \ox \Symn(p_{Y})} \Symm(X) \ox \Symn(Y).
\end{align*}
Finally, we denote by 
\begin{align*}
\xi_{X, Y}^{N}: \Sym^{N}(X \os Y) \lra \bos_{m+n = N} \Sym^{m}(X) \ox \Sym^{n}(Y)
\end{align*}
the sum of $\left( i_{X, Y}^{m, n} \circ \xi_{X, Y}^{m, n} \right):  \Sym^{N}(X \os Y) \ra 
\DS{\bos_{m + n = N}} \Sym^{m}(X) \ox \Sym^{n}(Y)$ over all $m$ and $n$ with sum $N$.
\end{notation}

\begin{remark}
We remark here that the morphism: 
\begin{align*}
\xi_{Z}: \Sym(Z) \lra \Sym(Z) \ox \Sym(Z)
\end{align*}
of Notation \ref{xinot} is \emph{not} the same as the comultiplication map: 
\begin{align*}
\Delta: \Sym(Z) \ra \Sym(Z) \ox \Sym(Z)
\end{align*}
associated with the standard Hopf algebra structure on $\Sym(Z)$.
\end{remark}

\begin{proposition}[\tbf{Classical}]  \label{psixiprop}
Let $X$ and $Y$ be objects of an idempotent complete $\Q$-linear tensor category $\catc$. Then for all non-negative integers $m$ and $n$ with sum $N$, the morphisms\tn{:}
\begin{align*}
& \psi_{X, Y}^{N}: \bos_{m+n = N} \Symm(X) \ox \Symn(Y) \lra \Sym^{N}(X \os Y) \qaq \\
& \xi_{X, Y}^{N}: \Sym^{N}(X \os Y) \lra \bos_{m+n = N} \Sym^{m}(X) \ox \Sym^{n}(Y)
\end{align*}
are mutually inverse isomorphisms. Consequently, 
\begin{align*}
& \psi_{X, Y}: \Symm(X) \ox \Symn(Y) \lra \Sym(X \os Y) \qaq \\
& \xi_{X, Y}: \Sym(X \os Y) \lra \Sym(X) \ox \Sym(Y)
\end{align*}
are also mutually inverse isomorphisms. 
\end{proposition}
\begin{prooff}
A version of this proposition is proved in \cite[Chapter III, \S 6.6]{BAlg} over the category of modules over a commutative ring. The proof in our situation is formally identical. 
\end{prooff}

\begin{definition}[\tbf{Classical}]
Let $(\catc, \ox_{\catc})$ and $(\catd, \ox_{\catd})$ be tensor categories. Then a functor $F: \catc \lra \catd$ is called a \emph{tensor functor} if the diagram\tn{:}
\begin{align*}
\xymatrix{
\catc \x \catc \ar@{->}[rr]^(0.55){- \ox_{\catc} -} \ar@{->}[dd]_(0.5){F \x F} & & \catc \ar@{->}[dd]^(0.5){F} \\
& & \\
\catd \x \catd \ar@{->}[rr]_(0.55){- \ox_{\catd} -} & & \catd}
\end{align*}
of categories commutes up to natural isomorphism of functors.
\end{definition}

\begin{proposition}[\tbf{Classical}]  \label{symaltcomm}
Let $(\catc, \ox_{\catc})$ and $(\catd, \ox_{\catd})$ be idempotent complete $\Q$-linear tensor categories, and let\tn{:}
\begin{align*}
F: \catc \lra \catd 
\end{align*}
be a tensor functor. Then $F$ commutes with $\Symn$; i.e. the diagram\tn{:}
\begin{align*}
\xymatrix{
\catc \ar@{->}[rr]^(0.5){\Symn} \ar@{->}[dd]_(0.5){F} & & \catc \ar@{->}[dd]^(0.5){F} \\
& & \\
\catd \ar@{->}[rr]_(0.5){\Symn} & & \catd} 
\end{align*} 
commutes up to natural isomorphism of functors. Consequently, $F$ commutes also with $\Sym$.
\end{proposition} 
\begin{prooff}
A proof is given in \cite[Proposition 2.5.9]{Big}. Note that although the formulation of this result requires that the source category $\catc$ be a triangulated category, the proof goes through without this assumption. 
\end{prooff}

\section{Categories of Chow motives}

Let $k$ be a perfect field. In this section we recall the definitions and main properties of the categories $\Chowk$ and $\Chowek$ of Chow motives over $k$ and effective Chow motives over $k$, respectively. 

\begin{notation}
We write $\Schk$ for the category whose objects are schemes of finite type over $k$, and whose morphisms are ordinary morphisms of schemes. The full subcategory of $\Schk$ whose objects are smooth schemes over $k$ we denote by $\Smk$, and the the full subcategory of $\Schk$ whose objects are smooth projective $k$-schemes we denote by $\SmProjk$. 
\end{notation}

To define Chow motives over $k$, the basic idea is to turn the category $\SmProjk$ into an additive category by replacing usual morphisms of schemes with \emph{correspondences}, which are (equivalence classes of) formal sums of closed subschemes. Our exposition is closely modelled on those of Scholl \cite{Scholl}, K{\"u}nnemann \cite{Ku1}, and Deninger and Murre \cite{DeMu}, but with one important caveat:

\begin{convention} \label{covcon}
The categories $\Chowk$ and $\Chowek$ will be defined to be $\Q$-linear, and such that the functor:
\begin{align*}
h: \SmProjk \lra \Chowk
\end{align*}
is covariant (see Definition \ref{hfunc}). The covariance of $h$ agrees with Voevodsky's definition of Chow motives \cite[preamble to Proposition 2.1.4]{TMF}, but conflicts with the classical literature. On the other hand, $\Q$-linearity belongs to the classical definition, whereas Voevodsky's category of Chow motives is merely $\Z$-linear.
\end{convention}

\noindent Aside from this convention, all of the definitions and results here are classical, and can be found in the above references. Throughout, $X$ will denote a smooth projective $k$-scheme.

\subsection{Category $\CSmProjk$ of Chow Correspondences}

Let $k$ be a perfect field. We review the construction of the categories $\Chowk$ of \emph{Chow motives} over $k$, and the subcategory $\Chowek$ of \emph{effective Chow motives}. 

\begin{definition}[\tbf{Classical}] 
Let $X$ be smooth projective scheme over $k$ and let $r$ be a non-negative integer. We recall that a \emph{cycle in $X$} means a formal $\Z$-linear combination $Z = \sum_{i} n_{i} [Z_{i}]$ of integral closed subschemes. The \emph{support} of $Z$ is the set of $Z_{i}$ for which 
$n_{i} \neq 0$. If all elements of the support are $r$-dimensional, then we say that $Z$ is an $r$-dimensional cycle. We denote by $\CH_{r}(X)$ the \emph{Chow group of cycles of dimension $r$ in $X$}, which is the group of $r$-dimensional cycles in $X$ modulo rational equivalence.
\end{definition}

\begin{notation}[\tbf{Scholl}]  \label{CCornot}
Let $X$ and $Y$ be smooth projective schemes over $k$ and let $r$ be a non-negative integer. Write $X_{i}$ for the connected components of $X$, and let $d_{X_{i}} := \dim(X_{i})$. Following \cite[\S 1]{Scholl}, we set:
\begin{align*}
\CCor^{r}(X, Y) & := \bos_{i} \CH_{d_{X_{i}} + r}(X_{i} \x Y).
\end{align*}
This is called the \emph{group of Chow correspondences of degree $r$ from $X$ to $Y$}.
\end{notation}

\begin{definition}[\tbf{Classical}]  \label{chowcorcomp}
Let $X, Y$ and $Z$ be smooth projective $k$-schemes, and assume $X$ and $Y$ are connected. For non-negative integers $r$ and $s$, we define following \cite[\S 1]{Scholl} a composition:
\begin{align*}
\mathsmap{- \: \circ -}{\CH_{s}(Y \x Z) \x \CH_{r}(X \x Y)}{\CH_{r+s}(X \x Z)}{(\beta, \alpha)}
{(p_{XZ})_{*}(p_{XY}^{*}(\alpha) \cdot p_{YZ}^{*}(\beta)),}
\end{align*}
where the morphisms $p_{XY}$, $p_{XZ}$ and $p_{YZ}$ denote the canonical projections from $X \x Y \x Z$ onto the product of the two factors indicated by the subscript, and $p_{XY}^{*}(\alpha) \cdot p_{YZ}^{*}(\beta))$ denotes the intersection product (see \cite[\S 6.1]{Ful}). By definition of the group of Chow correspondences in Notation \ref{CCornot}, composition extends bilinearly to a composition: 
\begin{align*}
\mathsmap{- \: \circ -}{\CCork^{s}(Y \x Z) \x \CCork^{r}(X \x Y)}{\CCork^{r+s}(X \x Z)}{(\beta, \alpha)}
{(p_{XZ})_{*}(p_{XY}^{*}(\alpha) \cdot p_{YZ}^{*}(\beta)),}
\end{align*}
where $X$ and $Y$ are not necessarily connected. 
\end{definition}

\begin{definition}[\tbf{Classical}] \label{chowcondef}
Following \cite[\S 1]{Scholl}, we define the category $\CSmProjk$ of \emph{Chow correspondences} as follows. The objects of $\CSmProjk$ are smooth projective schemes over $k$. The set of morphisms from a smooth projective scheme $X$ to a smooth projective scheme $Y$ is defined to be:
\begin{align*}
\Mor_{\CSmProjk}(X, Y) & = \CCor^{0}(X, Y), 
\end{align*}
and composition of morphisms is given by Definition \ref{chowcorcomp}. Observe that $\CSmProjk$ is an additive category.
\end{definition}

\begin{remark} \label{onecyc}
There is a natural functor $\SmProjk \ra \CSmProjk$ which is the identity on objects, and which sends a morphism 
$f:X \ra Y$ of smooth projective $k$-schemes to (the rational equivalence class of) its graph 
$[\Gamma(f)] \in \CCor^{0}(X, Y)$. Note that the identity morphism $X \ra X$ in $\CSmProjk$ is given by 
$[\Gamma(1_{X})]$, where $1_{X}: X \ra X$ is the identity morphism in $\SmProjk$.
\end{remark}

\begin{definition}[\tbf{Classical}] \label{dsumc}
We define a direct sum on $\CSmProjk$, as follows. On objects, we set $X \os Y := X \amalg Y$. To define the direct sum: 
\begin{align*}
f \os g \in \CCork^{0}(X_{1} \amalg X_{2}, Y_{1} \amalg Y_{2})
\end{align*}
of two morphisms $f \in \CCork^{0}(X_{1}, Y_{1})$ and $g \in \CCork^{0}(X_{2}, Y_{2})$, we observe first that by Definition \ref{chowcondef} of the morphisms for $\CSmProjk$, we may assume that $X_{1}$ and $X_{2}$ are connected. In view of the general relations:
\begin{align*}
\CH_{n}(A \amalg B) & \simeq \CH_{n}(A) \os \CH_{n}(B) \qaq \\
(X_{1} \amalg X_{2}) \x (Y_{1} \amalg Y_{2}) & \simeq 
(X_{1} \x Y_{1}) \amalg (X_{1} \x Y_{2}) \amalg (X_{2} \x Y_{1}) \amalg (X_{2} \x Y_{2}),
\end{align*}
it then follows that:
\begin{align*}
\CCork^{0}(X_{1} \amalg X_{2}, Y_{1} \amalg Y_{2}) \simeq &
\CH_{d_{X_{1}}} \big( X_{1} \x (Y_{1} \amalg Y_{2}) \big) \os 
\CH_{d_{X_{2}}} \big( X_{2} \x (Y_{1} \amalg Y_{2}) \big) \\
\simeq & \CH_{d_{X_{1}}} \big( X_{1} \x Y_{1} \big) \os \CH_{d_{X_{1}}} \big( X_{1} \x Y_{2} \big) \\
& \os \CH_{d_{X_{2}}} \big( X_{2} \x Y_{1} \big) \os \CH_{d_{X_{2}}} \big( X_{2} \x Y_{2} \big).
\end{align*}
Therefore: 
\begin{align*}
\CH_{d_{X_{1}}}(X_{1} \x Y_{1}) \os \CH_{d_{X_{2}}}(X_{2} \x Y_{2}) = 
\CCork^{0}(X_{1}, Y_{1}) \os \CCork^{0}(X_{2}, Y_{2})
\end{align*}
is a natural direct summand of $\CCork^{0}(X_{1} \amalg X_{2}, Y_{1} \amalg Y_{2})$, and so there is a natural map:
\begin{align*}
- \os -: \CCork^{0}(X_{1}, Y_{1}) \x \CCork^{0}(X_{2}, Y_{2}) \lra 
\CCork^{0}(X_{1} \amalg X_{2}, Y_{1} \amalg Y_{2}) 
\end{align*}
defining direct sum of morphisms in $\CSmProjk$.
\end{definition} 

\begin{definition}[\tbf{Classical}] \label{chowprod}
Given smooth projective $k$-schemes $X$ and $Y$, and non-negative integers $r$ and $s$, there is a natural product:
\begin{align*}
\mathsmap{- \x -}{\CH_{r}(X) \x \CH_{s}(Y)}{\CH_{r+s}(X \x Y)}{(\alpha, \beta)}{\alpha \x \beta} 
\end{align*}
induced by the fibre product of a closed subscheme of $X$ with a closed subscheme of $Y$.
\end{definition}

\begin{notation} \label{tpchow} 
Let $X = X_{1} \x Y_{1}$ and $Y = X_{2} \x Y_{2}$ be smooth projective $k$-schemes, with all $X_{i}$ and $Y_{i}$ connected, and let $r$ and $s$ be non-negative integers. Then the product of Definition \ref{chowprod} induces a map:
\begin{align*}
- \ox -: \CH_{k}(X_{1} \x Y_{1}) \x \CH_{\ell}(X_{2} \x Y_{2}) & 
\sxra{- \x -} \CH_{k+\ell}(X_{1} \x Y_{1} \x X_{2} \x Y_{2}) \\
& \slra{\tau} \CH_{k+\ell}(X_{1} \x X_{2} \x Y_{1} \x Y_{2}),
\end{align*}
where $\tau$ is the isomorphism exchanging the middle two factors. Using Definition \ref{chowcondef} of Chow correspondences, we extend this bilinearly to produce a tensor product:
\begin{align*}
-\ox -: \CCork^{r}(X_{1}, Y_{1}) \x \CCork^{s}(X_{2}, Y_{2}) \lra
\CCork^{r+s}(X_{1} \x X_{2}, Y_{1} \x Y_{2}),
\end{align*}
for all integers $r$ and $s$, where the $X_{i}$ and $Y_{i}$ are not necessarily connected. 
\end{notation} 

\begin{definition}[\tbf{Classical}] \label{tensorchowcordef}
We define a tensor product structure on the category $\CSmProjk$ as follows. On objects, we define $X \ox Y := X \x Y$. To define the tensor product on morphisms, let $X_{1}, Y_{1}, X_{2}$ and $Y_{2}$ be smooth projective $k$-schemes.
Setting $r = \dim(X_{1})$ and $s = \dim(X_{2})$ in Notation \ref{tpchow} gives us a tensor product: 
\begin{align*}
-\ox -: & \CCork^{0}(X_{1}, Y_{1}) \x \CCork^{0}(X_{2}, Y_{2}) \lra
\CCork^{0}(X_{1} \x X_{2}, Y_{1} \x Y_{2}).
\end{align*}
The identity object for $\ox$ is $\one_{\CSmProjk} = \Speck$. 
\end{definition}

\subsection{The category of effective Chow motives}

\begin{definition}[\tbf{Classical}]
Let $k$ be a perfect field. The category $\Chowek$ of \emph{effective Chow motives} over $k$ is defined to be the idempotent completion of the category $\CSmProjk \ox \Q$. 
\end{definition}

\begin{remark}
By definition $\Chowek$ is $\Q$-linear, whereas $\CSmProjk$ was merely $\Z$-linear. Since $\Chowek$ is moreover idempotent complete, there exist symmetric powers $\Symn(X)$ in the sense of Definition \ref{symaltn}.
\end{remark}

\begin{remark}
By \cite[Definition 1.2]{BaSch}, $\Chowek$ can be realised as the category whose objects consist of pairs $(X, p)$, where $X$ is a smooth projective scheme over $k$ and $p \in \CCor^{0}(X, X)$ is an idempotent, and the morphisms are given by:
\begin{align*}
\Mor_{\Chowek}((X, p), (Y, q)) & = q \circ \Mor_{\CSmProjk}(X, Y) \circ p \\
& = q \circ \CCork^{0}(X, Y) \circ p.
\end{align*} 
Note that since $p^{2} = p$ and $q^{2} = q$, we have:
\begin{align*}
\Mor_{\Chowek}((X, p), (Y, q))  =  \{ f \in \Mor_{\CSmProjk}(X, Y) \mid q \circ f = f = f \circ p \},
\end{align*} 
for all objects $(X, p)$ and $(Y, q)$ of $\Chowek$.
\end{remark}

\begin{remark} \label{chowim}
Let $(X, p)$ be an effective Chow motive, and suppose we have a morphism $s: (X, p) \ra (X, p)$ in $\Chowek$, which is itself idempotent. Then we have $s \circ p = p \circ s$ and $s^{2} = s$, implying that the diagram:
\begin{align*}
\xymatrix{
X \ar@{->}[rr]^(0.45){s} \ar@{->}[dd]_(0.5){p} & & X \ar@{->}[dd]^(0.5){s \circ p} \\
& & \\
X \ar@{->}[rr]_(0.45){s} & & X}
\end{align*}
commutes. Hence $s$ can also be regarded as a morphism $(X, p) \ra (X, s \circ p)$. One can check that $(X, s \circ p)$ is the categorical image of $s$. 
\end{remark}

\begin{definition}[\tbf{Classical}]
Following \cite[Theorem 1.6]{Scholl}, we note that there is a direct sum structure on $\Chowek$ coming from the one on $\CSmProjk$, defined as follows. On objects one sets:
\begin{align*}
(X, p) \os (Y, q) := (X \amalg Y, p \os q), 
\end{align*} 
where the direct sum of $p$ and $q$ is taken in the sense of the natural $\Q$-linear extension of the Definition (\ref{dsumc}), and the direct sum $f \os g$ of two morphisms $f$ and $g$ in $\CSmProjk$ is also in the sense of Definition (\ref{dsumc}).
\end{definition}

\begin{definition}[\tbf{Classical}]
As noted in \cite[\S 1.9]{Scholl}, the category $\Chowek$ also inherits a tensor product structure from $\CSmProjk$, defined as follows. Given two objects $(X, p), (Y, q) \in \Chowek$, set:
\begin{align*}
(X, p) \ox (Y, q) := (X \x Y, p \ox q),
\end{align*} 
where $p \ox q$ is the natural $\Q$-linear extension of the tensor product of Notation \ref{tpchow}. Given two morphisms $f \in \CCor^{0}(X_{1}, Y_{1})$ and $g \in \CCor^{0}(X_{2}, Y_{2})$, the tensor product $f \ox g$ is also taken in the sense of Notation \ref{tpchow}.
\end{definition}

\subsection{The category of Chow motives}

\begin{definition}[\tbf{Classical}]
Following \cite[\S 1.4]{Scholl}, we recall the definition of the category $\Chowk$ of \emph{Chow motives} over $k$. The objects are triples of the form $\cly{M} := (X, p, m)$, where the pair $(X ,p)$ is an effective Chow motive and $m$ is an integer. The morphisms are given by:
\begin{align*}
\Mor_{\Chowk}((X, p, m), (Y, q, n)) & := \{ f \in \CCork^{n-m}(X, Y) \mid p \circ f = f = f \circ q \},
\end{align*} 
which is interpreted as zero when $n-m$ is negative. This is equal to the set $p \circ \CCork^{n-m}(X, Y) \circ q$, since $p$ and $q$ are projectors. 
\end{definition}

\begin{remark} \label{chowsubcat}
In the special case $m=n=0$, the definition of the morphism set $\Mor_{\Chowk}((X, p, 0), (Y, q, 0))$ agrees with the definition of $\Mor_{\Chowek}((X, p), (Y, q))$. Thus the assignment $(X, p) \mapsto (X, p, 0)$ induces a natural embedding of categories $e: \Chowek \ra \Chowk$.
\end{remark}

\begin{definition}[\tbf{Classical}] \label{hfunc}
We define:
\begin{align*}
h: \SmProjk \lra \Chowk
\end{align*} 
to be the (covariant) functor sending a smooth projective scheme $X$ to the triple $(X, 1_{X}, 0)$ and a morphism 
$f: X \ra Y$ to its graph $\Gamma(f) \in \CCork^{0}(X, Y)$. By its definition, this functor factors through the subcategory $\Chowek$ of effective Chow motives, and so we are free to consider $h$ also as a functor $\SmProjk \ra \Chowek$.
\end{definition}

\begin{definition}[\tbf{Classical}]
Following \cite[\S 1.9]{Scholl}, we define the \emph{Lefschetz motive} in $\Chowk$ to be:
\begin{align*}
\Lf := (\Speck, 1_{\Speck}, 1). 
\end{align*} 
For all integers $n$, we adopt Scholl's convention of writing: 
\begin{align*}
\Lf^{n} := (\Speck, 1_{\Speck}, n),
\end{align*} 
to ease the notation. Our definition of the Lefschetz motive differs by a sign from Scholl's: $\Lf = (\Speck, 1_{\Speck}, -1)$, since his Chow motives are contravariant, whereas ours are covariant. After we have defined the tensor product of Chow motives in Definition \ref{chowtensdef}, it will be apparent that $\Lf^{n} = \Lf^{\ox n}$.
\end{definition}

\begin{construction}[\tbf{Classical}]
Let $X$ be an irreducible smooth projective $k$-scheme of dimension $d$, equipped with a fixed $k$-rational point $x_{0}: \Speck \ra X$. Then both $p_{0}^{X} := [ X \x \{ x_{0} \}]$ and $p_{2d}^{X} := [ \{ x_{0} \x X \}]$ are idempotent elements of $\CH_{d}(X \x X)_{\Q} = \CCork^{0}(X, X)$. Following \cite[\S 1.13]{Scholl}, we define objects:
\begin{align*}
\hz(X) := \left( X, p_{0}^{X}, 0 \right) \qaq \htd(X) := \left( X, p_{2d}^{X}, 0 \right).
\end{align*} 
By Remark \ref{chowsubcat}, they are effective Chow motives. Note that the roles of $p_{0}^{X}$ and $p_{2d}^{X}$, hence of $\hz(X)$ and $\htd(X)$ are reversed compared to Scholl's definition \cite[\S 1.11 and 1.13]{Scholl}, because our Chow motives are covariant and his are contravariant.
\end{construction}

\begin{proposition}[\tbf{Classical}] \label{hdecomp} 
Let $X$ be an irreducible smooth projective $k$-scheme of dimension $d$, equipped with a fixed $k$-rational point $x_{0}: \Speck \ra X$. Then\tn{:}
\begin{itemize}
\item[\tn{(i)}] There are isomorphisms $\onec \simeq \hz(X)$ and $\Lfd \simeq \htd(X)$.
\item[\tn{(ii)}] There is a direct sum decomposition\tn{:} 
\begin{align*}
h(X) & \simeq \hz(X) \os R(X) \os \htd(X) \\
& \simeq \onec \os R(X) \os \Lfd,
\end{align*}
where $R(X) := \left( X, q^{X}, 0 \right)$ and $q^{X} := 1_{X} - p_{0}^{X} - p_{2d}^{X}$.
\item[\tn{(iii)}] If $X = \Pjo$ then $R(X)$ vanishes, and we have\tn{:} 
\begin{align*}
h(\Pjo) & \simeq \onec \os \Lf.
\end{align*} 
\end{itemize}
\end{proposition}
\begin{prooff}
\tbf{(i)} See \cite[\S 1.11 -- 1.13]{Scholl}. More precisely, the fact that the $k$-rational point $x_{0}$ is a section to the structure morphism $\pi:X \ra \Speck$ is used in \cite[\S 1.11]{Scholl} to define $\hz(X)$ and $\htd(X)$ to be subobjects of $h(X)$ which are isomorphic to 
$\onec$ and $\Lfd$, respectively. It is then shown in \cite[\S 1.13]{Scholl} that one also has $\onec \simeq \left( X, p_{0}^{X}, 0 \right)$ and $\Lfd \simeq \left( X, p_{2d}^{X}, 0 \right)$.

\tbf{(ii)} Follows immediately from the fact that $p_{0}^{X}$, $p_{2d}^{X}$ and $q^{X}$ are mutually orthogonal idempotents 
with sum $1_{X}$.

\tbf{(iii)} \cite[\S 1.13]{Scholl}.
\end{prooff}

\begin{remark} \label{hodef}
Note that if $(C, \xo)$ is a smooth projective pointed curve, then we write $\ho(C)$ in place of the notation $R(C)$ in Proposition \ref{hdecomp}, and $p_{1}^{C}$ in place of $q^{C}$. The result then reads:
\begin{align*}
h(C) = \hz(C) \os \ho(C) \os h^{2}(C),
\end{align*} 
where $\ho(C) := \left(C, p_{1}^{C}, 0 \right)$, and $p_{0}^{C}$, $p_{2}^{C}$ and $p_{1}^{C} := 1_{C} - p_{0}^{C} - p_{2}^{C}$ are mutually orthogonal idempotents. 
\end{remark}

\begin{notation} \label{Cjpnot}
Let $(C, \xo)$ be a smooth projective pointed curve and $0 \le i \le 2$. Then since $p_{i}^{C}$ is a projector in 
$\CCork(C, C)$, the diagram:
\begin{align*}
\xymatrix{
C \ar@{->}[dd]_(0.5){p_{i}^{C}} \ar@{->}[rr]^(0.5){p_{i}^{C}} & & C \ar@{->}[dd]^(0.5){1_{C}} \\
& & \\
C \ar@{->}[rr]_(0.5){p_{i}^{C}} & & C}
\end{align*}
commutes, and hence $p_{i}^{C}$ determines a morphism from $\hi(C)$ to $h(C)$, which we denote by: 
\begin{align*}
j_{i}^{C}: \hi(C) \lra h(C).
\end{align*}
It also determines a morphism from $h(C)$ to $\hi(C)$, which we will continue to denote by $p_{i}^{C}: h(C) \ra \hi(C)$.
\end{notation}

\begin{proposition} \label{imkernp}
Let $(C, \xo)$ be a smooth projective pointed curve over $k$ with structure morphism $\pi_{C}$. Define the projector $p: C \ra C$ to be the composite morphism\tn{:} 
\begin{align*}
p: C \sxra{\pi_{C}} \Speck \sxra{\xo} C.
\end{align*}
Then the image and kernel of the idempotent morphism $h(p)$ in the category $\Chowek$ are given by\tn{:}
\begin{align*}
h(\pi_{C}): h(C) \lra \onec \qaq j_{1}^{C} \os j_{2}^{C}: \ho(C) \os h^{2}(C) \lra h(C),
\end{align*}
respectively.
\end{proposition}
\begin{prooff}
In general, if $(X, p, 0)$ is an effective Chow motive, then the image of the idempotent $h(p): h(X) \ra h(X)$ is given by:
\begin{align*}
[\Gamma(p)]: h(X) = (X, 1_{X}, 0) \lra (X, \Gamma(p), 0). 
\end{align*}
Observe now that since $p:C \ra C$ factors through $\Speck$, we have $[\Gamma(p)] = [C \x \{ \xo \}] = p_{0}^{C}$. The image of $h(p)$ is thus:
\begin{align*}
p_{0}^{C}: h(C) \lra (C, p_{0}^{C}, 0) = \hz(C) = \onec, 
\end{align*}
which is equal to $h(\pi_{C}): h(C) \ra h(\Speck) = \onec$. By Remark \ref{hodef}, the morphisms $p_{0}^{C}$, $p_{1}^{C}$ and 
$p_{2}^{C}$ are mutually orthogonal idempotents with sum $1_{C}$. It follows that the kernel of $h(p)$ is given by:
\begin{align*}
j_{1}^{C} \os j_{2}^{C}: \ho(C) \os h^{2}(C) \lra h(C).
\end{align*} 
This completes the proof.
\end{prooff}

\begin{definition}[\tbf{Classical}]
Following \cite[\S 1.14]{Scholl}, we define a direct sum on $\Chowk$ as follows. Suppose first that 
$\cly{M} := \left( X, p, m \right)$ and $\cly{N} := \left( Y, q, n \right)$ are Chow motives. If $m$ and $n$ are equal, set:
\begin{align}
\cly{M} \os \cly{N} := \left( X, p, m \right)  \os \left( Y, q, m \right) := \left( X \amalg Y, p \os q, m \right).
\label{mnequal}
\end{align} 
Suppose now that $n > m$. By Proposition \ref{hdecomp} Part (i), we have:
\begin{align*}
\Lf^{n-m} = \htw(\Pjo)^{\ox (n-m)} = \left((\Pjo)^{n-m}, (p^{\Pjo}_{2})^{\ox (n-m)}, 0 \right).
\end{align*} 
Since $\Lf^{n-m} = \left( \Speck, 1_{\Speck}, n - m \right)$, it follows that:
\begin{align*}
\cly{N} & := \left( Y, q, n \right) = \left( Y, q, m + (n - m) \right) = 
\left( Y, q, m \right) \ox \left( \Speck, 1_{\Speck}, n - m \right) \\ 
& = \left( Y, q, m \right) \ox \Lf^{n-m} = 
\left( Y, q, m \right) \ox  \left((\Pjo)^{n-m}, (p^{\Pjo}_{2})^{\ox (n-m)}, 0 \right) \\
& = \left( Y \x (\Pjo)^{n-m}, q \ox (p^{\Pjo}_{2})^{\ox (n-m)}, m \right).
\end{align*}
By \ref{mnequal}, we therefore have:
\begin{align*}
\cly{M} \os \cly{N} & = \left(X, p, m \right) \os \left( Y \x (\Pjo)^{n-m}, q \ox (p^{\Pjo}_{2})^{\ox (n-m)}, m \right) \\
& = \left( X \amalg \left\{ Y \x (\Pjo)^{n-m} \right\}, \: p \os \left\{ q \ox (p^{\Pjo}_{2})^{\ox (n-m)} \right\}, \: m \right).
\end{align*} 
The direct sum and tensor product in $\Chowk$ are compatible with those defined for $\Chowek$, in the sense that the natural embedding:
\begin{align*}
e: \Chowek \ra \Chowk
\end{align*} 
of Remark \ref{chowsubcat} is an additive tensor functor. 
\end{definition}

\begin{definition}[\tbf{Classical}] \label{chowtensdef}
Following \cite[\S 1.9]{Scholl}, we define a tensor product on the category $\Chowk$ as follows. On objects we set:
\begin{align*}
(X, p, m) \ox (Y, q, n) := (X \x Y, p \ox q, m+n),
\end{align*} 
where the tensor product $p \ox q$ of the projectors $p \in \CCork^{0}(X, X)$ and $q \in \CCork^{0}(Y, Y)$ is understood in the sense of Definition \ref{tensorchowcordef}. For morphisms:
\begin{align*}
f & \in \Mor_{\Chowk}((X_{1}, p_{1}, m_{1}), (Y_{1}, q_{1}, n_{1})) \sseq
\CCork^{n_{1} - m_{1}}(X_{1}, Y_{1}) \qaq \\
g & \in \Mor_{\Chowk}((X_{2}, p_{2}, m_{2}), (Y_{2}, q_{2}, n_{2})) \sseq
\CCork^{n_{2} - m_{2}}(X_{2}, Y_{2})
\end{align*} 
the tensor product: 
\begin{align*} 
f \ox g \in \CCork^{n_{1}+ n_{2} - m_{1} - m_{2}}(X_{1} \x X_{2}, Y_{1} \x Y_{2}),
\end{align*} 
is defined in the sense of in Notation \ref{tpchow}. Clearly the object: 
\begin{align*}
\onec := (\Speck, 1_{\Speck}, 0) = h(\Speck)
\end{align*} 
is the identity for $\ox$, where $1_{\Speck}$ is the symbol defined in Notation \ref{onecyc}. 
\end{definition}

\begin{remark} \label{symnchow}
Since the categories $\Chowek$ and $\Chowk$ are $\Q$-linear, additive and idempotent complete by construction, they possess symmetric and alternating powers in the sense of Definition \ref{symaltn}. Recall the morphism: 
\begin{align*}
s_{X}^{n}: X^{\ox n} \ra X^{\ox n}
\end{align*} 
from Notation \ref{ansn}, and let $(X, p, 0)$ be an effective Chow motive such that the relation 
$s_{X}^{n} \circ p^{\ox n} = p^{\ox n} \circ s_{X}^{n}$ is satisfied. Thus $s_{X}^{n}$ can be regarded as a morphism 
$s_{X}^{n}: (X^{\ox n}, p^{\ox n}) \ra (X^{\ox n}, p^{\ox n})$. We therefore have:
\begin{align*}
\Symn((X, p, 0)) = \im(s_{X}^{n}) = (X^{\ox n}, s_{X}^{n} \circ p^{\ox n}, 0),
\end{align*} 
by the Definition \ref{symaltn} and Remark \ref{chowim}. Moreover, Remark \ref{chowim} also says that the canonical projection mapping:
\begin{align*}
\pi^{n}_{(X, p , 0)}: (X, p, 0)^{\ox n} = (X^{\ox n}, p^{\ox n} , 0) \lra \Symn((X, p, 0)) 
\end{align*} 
is also equal to $s_{X}^{n}$.
\end{remark}

\section{Categories of triangulated Voevodsky motives}

Let $k$ be a perfect field. We recall the definition of the category $\DMegk$ of effective geometric Voevodsky motives, following Voevodsky \cite[\S 2]{TMF}. 

\subsection{Finite correspondences}

\begin{definition}[\tbf{Voevodsky}]
Let $X$ and $Y$ be smooth schemes over $k$. Following \cite[\S 2.1]{TMF}, a \emph{finite correspondence} is defined to be a formal $\Z$-linear combination $\sum n_{i} Z_{i}$ of integral closed subschemes of $X \x Y$, such that each $Z_{i}$ is finite and surjective over a connected component of $X$.  
\end{definition}

\begin{notation}
We denote by $\VCork(X, Y)$ the set of all finite correspondences from $X$ to $Y$. Note that it is an abelian group with respect to formal addition. 
\end{notation}

\begin{remark}
Since a subscheme $Z$ of $X \x Y$ which is finite and surjective over a connected component $X_{i}$ of $X$ has the same dimension as $X_{i}$, there is a natural morphism:
\begin{align*}
\VCork(X, Y) \lra \CCork(X, Y)
\end{align*}
for all smooth projective $k$-schemes $X$ and $Y$, which assigns to every integral closed subscheme $Z$ of $X \x Y$ its rational equivalence class $[Z]$. 
\end{remark}

\begin{remark}
Let $X$ and $Y$ be smooth schemes over $k$. Observe that the fibre product of schemes induces a well-defined product:
\begin{align*}
\mathsmap{- \x -}{\VCork(X) \x \VCork(Y)}{\VCork(X \x Y)}{(\alpha, \beta)}{\alpha \x \beta.}
\end{align*}
because the product of two finite and surjective morphisms of schemes is again finite and surjective.
\end{remark}

\begin{definition}[\tbf{Voevodsky}] \label{vcorcomp}
Let $X, Y$ and $Z$ be smooth $k$-schemes. Following \cite[\S 2.1]{TMF}, we define a composition law of finite correspondences:
\begin{align*}
\mathsmap{- \circ -}{\VCork(Y, Z)  \x \VCork(X, Y)}{\VCork(X, Z)}{(\beta, \alpha)}
{(p_{XZ})_{*}\left( p_{XY}^{*}(\alpha) \cdot p_{YZ}^{*}(\beta) \right).}
\end{align*}
Here $p_{XY}^{*}(\alpha) \cdot p_{YZ}^{*}(\beta)$ denotes the intersection product in $X \x Y \x Z$ 
(see \cite[\S 6.1]{Ful}) and the morphisms $p_{XY}$, $p_{XZ}$ and $p_{YZ}$ denote the three canonical projections from $X \x Y \x Z$ onto the specified two factors.
\end{definition}

\begin{definition}[\tbf{Voevodsky}]
Following \cite[\S 2.1]{TMF}, we define the category $\Smck$ of \emph{smooth correspondences over $k$} as follows. The objects are the smooth schemes of finite type over $k$, and the morphisms between two such schemes $X$ and $Y$ are the finite correspondences. Composition of morphisms is given by Definition \ref{vcorcomp}.
\end{definition}

\begin{remark}
As in the categories $\Chowek$ and $\Chowk$, the role of the identity morphism in $\VCork(X, X)$ is played by the 
$\Gamma(1_{X})$, where $1_{X}: X \ra X$ is the identity morphism in the category $\Schk$.
\end{remark}

\begin{definition}
We define a functor:
\begin{align*}
\scd: \Smk \lra \Smck
\end{align*}
as follows. On objects, we set $\scf{X} := X$. On morphisms, $\scf{f: X \ra Y}$ is defined to be the graph $\Gamma(f)$ of the morphism $f$.
\end{definition}

\begin{definition}[\tbf{Voevodsky}] \label{dsumv}
Following \cite[\S 2.1]{TMF}, we define a direct sum in $\Smck$. On objects, we set 
$\scf{X} \os \scf{Y} := \scf{X \amalg Y}$. To define the direct sum: 
\begin{align*}
f \os g \in \VCork(X_{1} \amalg X_{2}, Y_{1} \amalg Y_{2})
\end{align*}
of two morphisms $f \in \VCork(X_{1}, Y_{1})$ and $g \in \VCork(X_{2}, Y_{2})$, we note the general relations:
\begin{align*}
\VCork(A \amalg B) & \simeq \VCork(A) \os \VCork(B) \qaq \\
(X_{1} \amalg X_{2}) \x (Y_{1} \amalg Y_{2}) & \simeq 
(X_{1} \x Y_{1}) \amalg (X_{1} \x Y_{2}) \amalg (X_{2} \x Y_{1}) \amalg (X_{2} \x Y_{2}),
\end{align*}
analogous to those in Definition \ref{dsumc}, imply that there is a natural inclusion:
\begin{align*}
- \os -: \VCork(X_{1}, Y_{1}) \x \VCork(X_{2}, Y_{2}) \lra \VCork(X_{1} \amalg X_{2}, Y_{1} \amalg Y_{2}), 
\end{align*}
which defines $f \os g$.
\end{definition}

\begin{definition}[\tbf{Voevodsky}] \label{extprod}
Let $X_{1}, Y_{1}, X_{2}$ and $Y_{2}$ be smooth schemes over $k$. Following \cite[\S 2.1]{TMF}, we define a product:
\begin{align*}
\mathsmap{- \ox -}{\VCork(X_{1}, Y_{1})  \x \VCork(X_{2}, Y_{2})}
{\VCork(X_{1} \x Y_{1}, X_{2} \x Y_{2})}{(\alpha \x \beta)}{\tau(\alpha \x \beta),}
\end{align*}
where: 
\begin{align*}
\tau: \VCork(X_{1} \x Y_{1}, X_{2} \x Y_{2}) \lra \VCork(X_{1} \x X_{2}, Y_{1} \x Y_{2})
\end{align*}
denotes the morphism induced by the map:
\begin{align*}
X_{1} \x Y_{1} \x X_{2} \x Y_{2} \lra X_{1} \x X_{2} \x Y_{1} \x Y_{2}
\end{align*}
exchanging the middle two factors. It is known as the \emph{external product} of finite correspondences.
\end{definition}

\begin{definition}[\tbf{Voevodsky}] \label{dprodv}
We define a tensor product in $\Smck$. On objects, the tensor product is given by $\scf{X} \ox \scf{Y} := \scf{X \x Y}$. Tensor product of morphisms is taken in the sense of Definition \ref{extprod}.
\end{definition}

\subsection{Effective geometric Voevodsky motives}

\begin{definition}[\tbf{Voevodsky}]
Let $\Hb(\Smck)$ denote the homotopy category of bounded complexes over $\Smck$. Let $\catt$ be the smallest thick subcategory of $\Hb(\Smck)$ containing objects of the following two kinds:
\begin{itemize}
\item[\tn{(i)}] For any smooth scheme $X$ over $k$, the complex:
\begin{align*}
 \scf{X \x \Ao} \sxra{\pr_{X}} \scf{X}
\end{align*}
\item[\tn{(ii)}] For any smooth scheme $X$ and any open covering $X = U \cup V$ of $X$:
\begin{align*}
 \scf{U \cap V} \sxra{ \scf{j_{U}} \os  \scf{j_{V}}}  \scf{U} \os  \scf{V} 
 \sxra{ \scf{i_{U}} \os (- \scf{i_{V}})}  \scf{X},
\end{align*}
where $j_{U}$ and $j_{V}$ denote the inclusions of $U \cap V$ into $U$ and $V$, respectively, and 
$i_{U}$ and $i_{V}$ denote the embeddings of $U$ and $V$ into $X$, respectively.
\end{itemize}
Following \cite[Definition 2.1.1]{TMF}, we define the category $\DMegk$ of \emph{effective geometric Voevodsky motives} to be the idempotent completion of the localisation $\Hb(\Smck)[\catt^{-1}]$ of $\Hb(\Smck)$ at $\catt$.
\end{definition}

\begin{remark}
The category $\DMegk$ has a direct sum induced by disjoint union of smooth $k$-schemes, and a tensor product structure induced by the tensor product of Definition \ref{extprod} for finite correspondences. 
\end{remark}

\begin{notation} \label{MgmNot}
We denote by $\Mgm: \Smk \lra \DMegk$ the functor which assigns to a smooth scheme $X$ the complex consisting of $X$ concentrated in degree zero. Abusing notation, we shall also denote by $\Mgm$ the analogous functor $\Smck \lra \DMegk$.
\end{notation}

\begin{proposition}[\tbf{Voevodsky}]
The category $\DMegk$ of effective geometric Voevodsky motives is a tensor triangulated category, and the functor\tn{:} 
\begin{align*}
\Mgm: \Smck \lra \DMegk
\end{align*}
is an additive tensor functor.
\end{proposition}
\begin{prooff}
Additivity of $\Mgm$ is clear. The rest is \cite[Proposition 2.1.3]{TMF}.
\end{prooff}

\begin{notation}[\tbf{Voevodsky}]
We denote $\Mgm(\Speck)$ by $\one$, since it is the unit of the tensor product in $\DMegk$. For any smooth $k$-scheme $X$, let $\tMgm(X)$ denote the complex $\scf{X} \ra \scf{\Speck}$, concentrated in degrees zero and one, and set $\one(1) := \tMgm(\Pjo)[-2]$. The object $\one(1)$ is known as the \emph{Tate object}. All this notation is borrowed from 
\cite[Remarks after Proposition 2.1.3]{TMF}, except that Voevodsky writes ``$\Z$'' instead of ``$\one$''.
\end{notation}

\begin{remark}
From the definition of $\one(1)$, it follows that there is a direct sum decomposition $\Mgm(\Pjo) = \one \os \tatt$. 
\end{remark}

\begin{remark}
Recall from Notation \ref{sxdef} that $\DMegk \sx \Q$ denotes the idempotent completion of $\DMegk \ox \Q$. The category $\DMegk \sx \Q$ inherits a structure of tensor triangulated category from $\DMegk$. Moreover, since it is idempotent complete by definition, it also possesses, symmetric powers $\Symn(X)$ in the sense of Definition \ref{symaltn}. Note that here one uses the result \cite{BaSch} of Balmer and Schlichting that the idempotent completion of a triangulated category is again a triangulated category. 
\end{remark}

\subsection{The Voevodsky category of motivic complexes}

We recall Voevodsky's definition \cite[\S 3]{TMF} of the triangulated category of effective motivic complexes. Throughout, the symbol $\tau$ will denote either the {\etl} or {\nsv} topology on $\Smk$, and $R$ will be a ring containing the integers $\Z$.

\begin{definition}[\tbf{Nisnevich}]
Recall that an {\etl} covering $p_{i}: \{ U_{i} \ra X \}_{i \in I}$ is said to be a \emph{{\nsv} covering} if for all $x \in X$, there exists an $i \in I$ and $u \in U_{i}$ such that $p_{i}(u) = x$, and the induced morphism 
$k(x) \ra k(u)$ of residue fields is an isomorphism. 
\end{definition}

\begin{definition}[\tbf{Voevodsky}]
Following \cite[Definition 3.1.1]{TMF}, the category $\PSTkR$ of \emph{($R$-module-valued) presheaves with transfers} is defined to be the category of all contravariant functors from $\Smck$ to the category of $R$-modules.
\end{definition}

\begin{definition}[\tbf{Voevodsky}]
An $R$-module-valued presheaf with transfers is said to be a \emph{$\tau$ sheaf with transfers} if the associated presheaf of $R$-modules on $\Smk$ is a sheaf in the $\tau$ topology. We denote this category by $\STtkR$. Note that $\STtkR$ is an abelian category: this is proved in \cite[Theorem 3.1.4]{TMF} for the $\tau$ Nisnevich case and is similar in the {\etl} case.
\end{definition}

\begin{remark} \label{niset}
Since Nisnevich coverings are special kinds of {\etl} coverings, every {\etl} sheaf is a Nisnevich sheaf. Thus there is a natural functor $\STEkR \lra \STNkR$.
\end{remark}

\begin{proposition}[\tbf{Voevodsky}] \label{transord}
Two $\tau$ sheaves with transfers are isomorphic if and only if the underlying ordinary $\tau$ sheaves are isomorphic.
\end{proposition}
\begin{prooff}
This follows from \cite[Lemma 3.1.5]{TMF} in the case of the Nisnevich topology. The proof for the {\etl} topology is similar.
\end{prooff}

\begin{notation}[\tbf{Voevodsky}] \label{Lfuncnot}
Following \cite[Remark after Definition 3.1.1]{TMF}, we denote by:
\begin{align*}
L: \left( \Smck \right)^{\op} \lra \PSTkR
\end{align*}
the functor assigning to an object $X$ the presheaf $\VCork(\pul{X}, X)$ that it represents. By abuse of notation, we also regard $L$ as a functor from $(\Smk)^{\op}$ to $\PSTkR$.
\end{notation}

\begin{proposition}[\tbf{Voevodsky}]
For any smooth scheme $X$ over $k$, the presheaf $L(X)$ is a sheaf in both the Nisnevich and {\etl} topologies. 
\end{proposition}
\begin{prooff}
This is proved for the Nisnevich topology in \cite[Lemma 3.1.2]{TMF} and for the {\etl} topology in \cite[Lemma 6.2]{MVW}, both in the special case $R = \Z$. The proof for general coefficients $R$ is similar.
\end{prooff}

\begin{definition}[\tbf{Voevodsky}]
Let $X^{\bl}$ and $Y^{\bl}$ be objects of $\Dm(\STtkR)$ with projective resolutions $P^{\bl}$ and $Q^{\bl}$, respectively. Then we define:
\begin{align*}
X^{\bl} \oxL Y^{\bl} := \Tot(P^{\bl} \ox Q^{\bl}),
\end{align*}
where $\Tot(P^{\bl} \ox Q^{\bl})$ denotes the total complex of the bicomplex $P^{\bl} \ox Q^{\bl}$. 
\end{definition}

\begin{proposition}[\tbf{Voevodsky}]
The category $\Dm(\STtkR)$ is tensor triangulated with respect to $\oxL$. 
\end{proposition}
\begin{prooff}
\cite[Proposition 8.8]{MVW}.
\end{prooff}

\begin{definition}[\tbf{Voevodsky}]
For any smooth scheme $X$, let $p_{X}: X \x \Ao \ra X$ denote projection onto the first factor. We say that a presheaf $\shfF \in \PSTkR$ is \emph{homotopy invariant} if for all $X$, the induced morphism:
\begin{align*}
F(p_{X}): F(X) \lra F(X \x \Ao)
\end{align*}
is an isomorphism. We denote by $\HIkR$ the full subcategory of $\STtkR$ consisting of homotopy invariant sheaves. 
\end{definition}

\begin{notation}
Let $X^{\bl}$ be a complex in $\Dm(\STtkR)$. The \emph{$i^{\tn{th}}$ cohomology sheaf $\shfH^{i}(X^{\bl})$} is by definition the quotient sheaf:
\begin{align*}
\shfH^{i}(X^{\bl}) := \ker(d^{i})/\im(d^{i-1}),
\end{align*}
where $d^{i}$ denotes the $i^{\tn{th}}$ differential, for all integers $i$.
\end{notation}

\begin{definition}[\tbf{Voevodsky}]
Following \cite[Remark after Proposition 3.1.12]{TMF}, we define $\DMettkR$ to be the full subcategory of the derived category $\Dm(\STtkR)$ consisting of complexes $X^{\bl}$ such that the $i^{\tn{th}}$ cohomology sheaf $\shfH^{i}(X^{\bl})$ is homotopy invariant for all $i \in \Z$.
\end{definition}

\begin{proposition}[\tbf{Voevodsky}]
The category $\DMettkR$ is a triangulated subcategory of $\Dm(\STtkR)$.
\end{proposition}
\begin{prooff}
This is proved for the {\nsv} topology and $R = \Z$ in \cite[Proposition 3.1.12]{TMF}; the proof for the {\etl} topology and for general coefficients is similar.
\end{prooff}

\begin{notation}
Let $\cj{k}$ be a fixed algebraic closure of the perfect field $k$, and let $\pi: \Spec(\cj{k}) \ra \Spec(k)$ be the induced morphism of schemes. This induces pullback functors:
\begin{align*}
& \pi^{*}: \SEsmk \lra \SEsmck \qaq \\
& \pi^{*}: \DMeetkR \lra \DMeet(\cj{k}, R).
\end{align*}
For all $\shfF \in \SEsmk$ and $\cly{X} \in \DMeetkR$, we write $\shfF_{\cj{k}}$ in place of 
$\pi^{*}(\shfF)$ and $\cly{X}_{\cj{k}}$ in place of $\pi^{*}(\cly{X})$.
\end{notation}

\begin{lemma} \label{etldesshf}
Let $\cj{k}$ be a fixed algebraic closure of the perfect field $k$. A morphism $f: \shfF \lra \shfG$ of {\etl} sheaves $\shfF, \shfG \in \SEsmk$ is an isomorphism provided that its base-extension\tn{:}
\begin{align*}
f_{\cj{k}}: \shfF_{\cj{k}} \lra \shfG_{\cj{k}}
\end{align*}
is an isomorphism in the category $\SEsmck$. 
\end{lemma}
\begin{prooff}
The morphism $f: \shfF \ra \shfG$ of sheaves on the big {\etl} site $(\Sm/k)_{\et}$ is an isomorphism if and only if, for every $k$-scheme $X$ and every geometric point $x: \Spec(K) \ra X$, the induced map $f_{x}: \shfF_{x} \ra \shfG_{x}$ of stalks is an isomorphism. Let us fix an embedding $k \hra \cj{k}$. Now, since $f_{\cj{k}}: \shfF_{\cj{k}} \ra \shfG_{\cj{k}}$ is an isomorphism of sheaves on $(\Sm/\cj{k})_{\et}$, it follows that:
\begin{align*}
(f_{\cj{k}})_{x^{\p}}: (\shfF_{\cj{k}})_{x^{\p}} \lra (\shfG_{\cj{k}})_{x^{\p}}
\end{align*}
is an isomorphism, where $x^{\p}: \Spec(K) \ra X_{\cj{k}}$ is the unique map such that the diagram:
\begin{align*}
\xymatrix{
\Spec(K) \ar@{->}[ddrr]^(0.5){x} \ar@{->}[dd]_(0.5){x^{\p}} & & \\
& & \\
X_{\cj{k}} \ar@{->}[rr]_(0.5){\pi} & & X}
\end{align*}
commutes, which is induced by the universal property of the fibre product $X_{\cj{k}} = X \x_{k} \cj{k}$. Here $\pi$ denotes the projection onto $X$. But there are canonical isomorphisms of stalks:
\begin{align*}
\shfF_{x} \ilra (\shfF_{\cj{k}})_{x^{\p}} \qaq \shfG_{x} \ilra (\shfG_{\cj{k}})_{x^{\p}},
\end{align*}
and so the fact that $(f_{\cj{k}})_{x^{\p}}$ is an isomorphism implies that $f_{x}$ is also.
\end{prooff}

\begin{proposition} \label{etldes}
Let $\cj{k}$ be a fixed algebraic closure of the perfect field $k$. A morphism $f: \cly{X} \ra \cly{Y}$ of {\etl} motives $\cly{X}, \cly{Y} \in \DMeetkR$ is an isomorphism provided that its base-extension\tn{:}
\begin{align*}
f_{\cj{k}}: \cly{X}_{\cj{k}} \lra \cly{Y}_{\cj{k}}
\end{align*}
is an isomorphism in the category $\DMeet(\cj{k}, R)$. 
\end{proposition}
\begin{prooff}
Suppose $f_{\cj{k}}: \cly{X}_{\cj{k}} \ra \cly{Y}_{\cj{k}}$ is an isomorphism in $\DMeet(\cj{k}, R)$. Since $\DMeet(\cj{k}, R)$ is a full subcategory of the derived category $\Dm(\STE(\cj{k}, R))$, this implies that:
\begin{align*}
H^{n}(f_{\cj{k}}) : H^{n}(\cly{X}_{\cj{k}}) \lra H^{n}(\cly{Y}_{\cj{k}})
\end{align*}
is an isomorphism of sheaves on the big {\etl} site $(\Sm/\cj{k})_{\et}$, for all integers $n$. By Lemma \ref{etldesshf}, it follows that:
\begin{align*}
H^{n}(f) : H^{n}(\cly{X}) \lra H^{n}(\cly{Y})
\end{align*}
is an isomorphism of sheaves on the big {\etl} site $(\Smk)_{\et}$, for all integers $n$, which means that 
$f: \cly{X} \ra \cly{Y}$ is an isomorphism in $\DMeetkR$.
\end{prooff}

\begin{definition}[\tbf{Voevodsky}] \label{Cxdef}
For any non-negative integer $n$, let $A^{n}$ be the $k$-algebra:
\begin{align*}
k[T_{0}, \ldots, T_{n}]\Big/ \left\langle \left( \sum_{i=0}^{n} t_{i} \right) - 1 \right\rangle
\end{align*}
and let $\Delta^{n}$ be the affine scheme $\Spec(A^{n})$. We set for all $0 \le i \le n$:
\begin{align*}
& \mathsmap{\delta^{n}_{i}}{A^{n}}{A^{n-1}}{t_{j}}
{\left\{ \begin{array}{rl}
t_{j}, & \tn{if } j < i, \\
0, & \tn{if } j = i, \\
t_{j-1}, & \tn{if } j > i,
\end{array} \right.} \qaq \\
& d^{n} := \sum_{i=0}^{n} (-1)^{i} \delta^{n}_{i}.
\end{align*}
Then $(A^{\bl}, d^{\bl})$ is a complex of algebras. Following \cite[Preamble to Lemma 3.2.1]{TMF} we define a functor:
\begin{align*}
\Cx: \PSTkR \lra \Dm(\PSTkR)
\end{align*}
by setting, for any $n \in \Z$, for any $\shfF \in \PSTkR$ and for any smooth scheme $U$ over $k$, 
\begin{align*}
\Cnx(\shfF)(U) := \shfF(U \x \Delta^{n}),
\end{align*}
and observing that the maps:
\begin{align*}
\pd_{n} := \shfF(\id_{U} \x \Spec(d^{n}))
\end{align*}
turn $\Cnx(\shfF)$ into a chain complex of sheaves with transfers, since $(A^{\bl}, d^{\bl})$ is a complex of algebras. Clearly 
$\Cx$ restricts to a functor: 
\begin{align*}
\STtkR \lra \Dm(\STtkR),
\end{align*}
which we also denote by $\Cx$.
\end{definition}

\begin{notation}[\tbf{Voevodsky}] \label{Mhinot}
Let $X$ be a smooth scheme over $k$, let $\shfF \in \PSTkR$ a presheaf with transfers, and let $i$ be an integer. Following \cite[Preamble to Lemma 3.2.1]{TMF}, we use the notation:
\begin{itemize}
\item[\tn{(i)}] $M(X) := \Cx(L(X)) \in \Dm(\STtkR)$.
\item[\tn{(ii)}] $\hitau(\shfF) := \shfH^{-i}(\Cx(\shfF)) \in \STtkR$.  
\end{itemize}
If $\shfF := L(X)$, then we will write $\hitau(X)$ in place of $\hitau(L(X))$.
\end{notation}

\begin{proposition}[\tbf{Voevodsky}] \label{hinishi}
For any $\tau$ sheaf with transfers $\shfF \in \STtkR$ and any integer $i$, the sheaf $\hitau(\shfF)$ is homotopy-invariant. In particular, for any smooth scheme $X$ over $k$, the sheaf with transfers\tn{:} 
\begin{align*}
\hitau(X) = \shfH^{-i}(\Cx(L(X)))
\end{align*}
is homotopy-invariant, and so $M(X)$ is an object in $\DMettkR$.
\end{proposition}
\begin{prooff}
This is proved for the {\nsv} topology and $R = \Z$ in \cite[Lemma 3.2.1]{TMF}; the proof in the case of a general coefficient ring $R$ and for the {\etl} topology is similar.
\end{prooff}

\begin{proposition}[\tbf{Voevodsky}] \label{quotA}
Let $\catt$ be the minimal thick triangulated subcategory of $\Dm(\STtkR)$ containing all complexes of the form\tn{:} 
\begin{align*}
L(\pr_{X}): L(X \x \Ao) \lra L(X)
\end{align*}
for all smooth schemes $X$. The quotient category $\Dm(\STtkR)/\catt$ inherits a structure of tensor category from $(\Dm(\STtkR), \oxL)$.
\end{proposition}
\begin{prooff}
This is proved for the {\etl} topology and $R = \Z$ in \cite[Corollary 9.6]{MVW}. The proof in the case of a general coefficient ring $R$ and for the {\nsv} topology is similar; c.f. \cite[Definition 14.2]{MVW}.
\end{prooff}

\begin{proposition}[\tbf{Voevodsky}] \label{embquot}
Denote by\tn{:} 
\begin{align*}
i: \DMettkR \lra \Dm(\STtkR)
\end{align*}
the natural embedding functor and by\tn{:}
\begin{align*}
q: \Dm(\STtkR) \lra \Dm(\STtkR)/\catt
\end{align*}
 the natural quotient functor. Then the composition\tn{:}
\begin{align*}
q \circ i: \DMettkR \lra \Dm(\STtkR)/\catt
\end{align*}
is an equivalence of categories.
\end{proposition}
\begin{prooff}
That this composition is an equivalence is \cite[Proposition 3.2.3]{TMF} in the case of the {\nsv} topology and and $R = \Z$; the proof for a general coefficient ring $R$ and for the {\etl} topology is similar.
\end{prooff}

\begin{proposition}[\tbf{Voevodsky}] \label{DMeetktentri}
The category $\DMettkR$ is tensor-triangulated and idempotent complete. Consequently, if $R$ contains $\Q$ there exist symmetric powers on $\DMettkR$ in the sense of Definition \ref{symaltn}.
\end{proposition}
\begin{prooff}
By Proposition \ref{quotA}, $\Dm(\STtkR)/\catt$ inherits a structure of tensor category from $\Dm(\STtkR)$. Since $\DMettkR$ is equivalent to $\Dm(\STtkR)/\catt$ by Proposition \ref{embquot}, we see that $\DMettkR$ has a structure of tensor triangulated category. The tensor product in $\DMettkR$ works as follows: observing that $\Cx$ is an inverse equivalence to 
$f := q \circ i$, one sets for all objects $X^{\bl}$ and $Y^{\bl}$ in $\DMettkR$:
\begin{align*}
X^{\bl} \ox Y^{\bl} := \Cx\left( f \left( X^{\bl} \right) \ox f \left( Y^{\bl} \right) \right).
\end{align*}
The definition for the tensor product of morphisms in $\DMettkR$ is the same. That $\DMettkR$ is idempotent complete is \cite[Lemma 3.1.13]{TMF} in the {\nsv} case and for $R = \Z$; the proof is similar for general coefficients and in the {\etl} case.
\end{prooff}

\begin{remark}
By the proof of Proposition \ref{DMeetktentri}, the functor $q \circ i$ of Proposition \ref{embquot} is a tensor functor. Note however that the embedding $i: \DMettkR \lra \Dm(\STtkR)$ is \emph{not} a tensor functor; c.f.
\cite[Remark following Lemma 3.2.5]{TMF}. 
\end{remark}

\begin{notation} \label{lfuncnot}
We denote by:
\begin{align*}
\ell_{\tau}: \DMegk \lra \DMettkZ
\end{align*}
the functor which assigns to a complex of smooth schemes $(X^{i})_{i \in \Z}$ the complex consisting of the corresponding representing sheaves $(L(X^{i}))_{i \in \Z}$.
\end{notation}

\begin{proposition}[\tbf{Voevodsky}]
\begin{itemize}
\item[\tn{(i)}] The natural embedding\tn{:} 
\begin{align*}
i: \DMettkZ \lra \Dm(\STtkZ)
\end{align*}
possesses a left-adjoint\tn{:}
\begin{align*}
\RCx: \Dm(\STtkZ) \lra \DMettkZ
\end{align*}
which restricts to $\Cx$ on the subcategory $\STtkZ$.
\item[\tn{(ii)}] The functor $\ell_{\tau}: \DMegk \ra \DMettkZ$ is a tensor functor, and the following diagram commutes\tn{:}
\begin{align*}
\xymatrix{
\Hb(\Smck) \ar@{->}[rr]^(0.45){L} \ar@{->}[dd] & & \Dm(\STtkZ) \ar@{->}[dd]^(0.5){\RCx} \\
& & \\
\DMegk \ar@{->}[rr]_(0.45){\ell_{\tau}} & & \DMettkZ.}
\end{align*}
Moreover, $\RCx(L(X))$ is canonically isomorphic to $\Cx(X)$.
\end{itemize}
\end{proposition}
\begin{prooff}
This is proved in \cite[Theorem 3.2.6]{TMF} in the {\nsv} case. The proof of the {\etl} case is similar.
\end{prooff}

\section{Chow and Voevodsky motives with $\Q$-coefficients}

In this section, we recall the construction of a functor $\Phi$ from the category $\Chowek$ of effective Chow motives to the Voevodsky category $\DMeetkQ$ of motivic complexes. 

\begin{proposition}[\tbf{Voevodsky}] \label{etnissame}
The functor $\STEk \ra \STNk$ of Remark \ref{niset} induces an equivalence\tn{:}
\begin{align*}
\DMeetkQ \ilra \DMenkQ
\end{align*}
of tensor triangulated categories.
\end{proposition}
\begin{prooff}
This is \cite[Proposition 3.3.2]{TMF} or \cite[Theorem 14.30]{MVW}.
\end{prooff}

\begin{convention}
For the rest of this thesis, we work exclusively with the {\etl} topology. Since in subsequent we are mainly interested in the category $\DMeetkQ$, we could equivalently work over the Nisnevich topology, by Proposition \ref{etnissame}. However, since group schemes are naturally sheaves on the big {\etl} site over $\Speck$, this standpoint is more natural for us. 
\end{convention}

\begin{notation} \label{abcnot}
We denote by:
\begin{align*}
& a: \DMegk \lra \DMegk \ox \Q, \\
& b: \DMeetkZ \lra \DMeetkQ, \qaq \\
& c: \DMegk \ox \Q \lra \DMegk \sx \Q
\end{align*}
respectively, (a) the functor of Notation \ref{sxdef}, (b) the functor induced by the functor $\STtkZ \ra \STtkQ$ which is the identity on morphisms, and on objects sends a sheaf $\shfF$ of abelian groups to the sheaf $\shfF \ox \Q$ of vector spaces, and (c) the natural functor from a category to its idempotent completion.
\end{notation}

\begin{construction} \label{jconst}
We construct a functor: 
\begin{align*}
j: \DMegk \sx \Q \lra \DMeetkQ
\end{align*}
as follows. Denote by:
\begin{align*}
\wt{\ell}_{\et}: \DMegk \ox \Q \ra \DMeetkQ
\end{align*}
the $\Q$-linear extension of the composite functor $b \circ \ell_{\et}: \DMegk \ra \DMeetkQ$, where $\ell_{\et}$ is the functor of Notation \ref{lfuncnot}. Since $\DMeetkQ$ is idempotent complete by Proposition \ref{DMeetktentri}, the universal property of the idempotent completion $\DMegk \sx \Q$ implies there is a unique functor: 
\begin{align*}
j: \DMegk \sx \Q \lra \DMeetkQ
\end{align*}
such that the diagram:
\begin{align*}
\xymatrix{
\DMegk \ox \Q \ar@{->}[rr]^(0.5){\wt{\ell}_{\et}} \ar@{->}[dd]_(0.5){c} & & \DMeetkQ \\
& & \\
\DMegk \sx \Q \ar@{->}[uurr]_(0.5){j} & & }
\end{align*}
is commutative.
\end{construction}

\begin{remark} \label{jttens}
The functor\tn{:}
\begin{align*}
j: \DMegk \sx \Q \lra \DMeetkQ
\end{align*}
of Construction \ref{jconst} is tensor triangulated, because all functors involved in its construction are tensor triangulated.
\end{remark}

\begin{construction}[\tbf{Voevodsky}] \label{PsiXYconst} 
Let $X$ and $Y$ be two smooth projective varieties over a perfect field $k$. We recall here the construction given in \cite[Proposition 2.1.4]{TMF} of the functorial morphism:
\begin{align*}
\Psi_{X, Y}^{\Z}: \Mor_{\ChowZek}(X, Y) \lra \Mor_{\DMegk}(\Mgm(X), \Mgm(Y)),
\end{align*}
where $\ChowZek$ denotes the category of $\Z$-linear effective Chow motives favoured by Voevodsky. Firstly, we write:
\begin{align*}
r_{0}: X \x Y \lra X \x \Ao \x Y \qaq r_{1}: X \x Y \lra X \x \Ao \x Y
\end{align*}
for the morphisms obtained by restricting $X \x \Ao$ to $X \x \{ 0 \}$ and $X \x \{ 1 \}$, respectively. Now denote by:
\begin{align*}
g_{X, Y}: \VCork(X \x \Ao, Y) \lra \VCork(X, Y)
\end{align*}
the difference $(r_{0})^{*} - (r_{1})^{*}$ of the induced pullback morphisms on Voevodsky correspondences, and 
denote by:
\begin{align*}
f_{X, Y}: \VCork(X, Y) \lra \Mor_{\DMegk}(\Mgm(X), \Mgm(Y)) 
\end{align*}
the morphism determined by the functor $\Mgm: \Smck \ra \DMegk$ of Notation \ref{MgmNot}. Since the category $\DMegk$ is defined so that: 
\begin{align*}
\Mgm(\pr_{X}): X \x \Ao \lra X
\end{align*}
is an isomorphism, it follows that $f_{X, Y}$ factors through the cokernel of $g_{X, Y}$, yielding a morphism:
\begin{align*}
\psi_{X, Y}: \coker(g_{X, Y}) \lra \Mor_{\DMegk}(\Mgm(X), \Mgm(Y)). 
\end{align*}
But by the definition of rational equivalence, there is a canonical homomorphism
$c_{X, Y}: \coker(g_{X, Y}) \lra \CCork^{0}(X, Y)$, which is an isomorphism, by the deep result \cite[Theorem 7.1]{BCC}. Let us now define:
\begin{align*}
\Psi_{X, Y}^{\Z}: \CCork^{0}(X, Y) \lra \Mor_{\DMegk}(\Mgm(X), \Mgm(Y)) 
\end{align*}
to be the composition $\psi_{X, Y} \circ c_{X, Y}^{-1}$. 
\end{construction} 

\begin{notation} \label{PsiXYnot} 
Observing that $\Chowek = \ChowZek \sx \Q$, let us denote by:
\begin{align*}
\Psi_{X, Y}: \Mor_{\Chowek}(X, Y) \lra \Mor_{\DMegk \sx \Q}(\Mgm(X), \Mgm(Y))
\end{align*}
the morphism $\Psi_{X, Y}^{\Z} \sx \Q$, where $\Psi_{X, Y}^{\Z}$ was defined in Construction \ref{PsiXYconst}.
\end{notation} 

\begin{proposition}[\tbf{Voevodsky}] \label{psimon}
The morphism\tn{:}
\begin{align*}
\Psi_{X, Y}: \Mor_{\Chowek}(X, Y) \lra \Mor_{\DMegk \sx \Q}(X, Y)
\end{align*}
of Notation \ref{PsiXYnot} induces a functor\tn{:}
\begin{align*}
\mathsmap{\Psi}{\Chowek}{\DMegk \sx \Q}
{\begin{array}{rc}
\tn{\small{Objects:}} & (X, p) \\
\tn{\small{Morphisms:}} & f
\end{array}}
{\begin{array}{c}
(\scf{X}, \scf{p}) \\
\Psi_{X, Y}(f), 
\end{array}}
\end{align*}
where the symbol $\scf{X}$ here denotes the object of $\Smck$ considered as a complex concentrated in degree zero. Moreover, $\Psi$ is additive and monoidal, and the diagram of functors\tn{:} 
\begin{align*}
\xymatrix{
\SmProjk \ar@{^{(}->}[rr] \ar@{->}[dd]_(0.5){h} & & \Smk \ar@{->}[dd]^(0.5){\Mgm} \\
& & \\
\Chowek \ar@{->}[rr]_(0.5){\Psi} & & \DMegk \sx \Q}
\end{align*}
commutes, meaning there is a canonical isomorphism $F_{\tn{gm}}: \Psi \circ h \ra \Mgm$.
\end{proposition}
\begin{prooff}
In \cite[Proposition 2.1.4]{TMF} it is proved that the $\Z$-linear counterpart $\Psi^{\Z}$ of $\Psi$ of Construction \ref{PsiXYconst} is a functor, and also that the associated diagram of functors is commutative. That $\Psi^{\Z}$ is additive and monoidal is \cite[\S 2.2 ``Relations to Chow motives'']{TMF}. The corresponding statements for our $\Q$-linear functor $\Psi$ follow formally by applying $\sx \Q$.
\end{prooff}

\begin{proposition} \label{psisymcomm}
The functor\tn{:} 
\begin{align*}
\Psi: \Chowek \ra \DMegk \sx \Q
\end{align*}
of Proposition \ref{psimon} commutes with $\Symn$ and $\Sym$, for all non-negative integers $n$, in the sense that the diagrams of functors\tn{:}
\begin{align*}
\xymatrix{
\Chowek \ar@{->}[rr]^(0.45){\Psi} \ar@{->}[dd]_(0.5){\Symn} & & 
\DMegk \sx \Q \ar@{->}[dd]^(0.5){\Symn} \\
& & \\
\Chowek \ar@{->}[rr]_(0.45){\Psi} & & \DMegk \sx \Q} \qaq 
\xymatrix{
\Chowek \ar@{->}[rr]^(0.45){\Psi} \ar@{->}[dd]_(0.5){\Sym} & & 
\DMegk \sx \Q \ar@{->}[dd]^(0.5){\Sym} \\
& & \\
\Chowek \ar@{->}[rr]_(0.45){\Psi} & & \DMegk \sx \Q} 
\end{align*}
both commute up to natural isomorphism. 
\end{proposition}
\begin{prooff}
Since $\Psi$ is an additive tensor functor by Proposition \ref{psimon}, it follows by Propositions \ref{symaltcomm} that it commutes with symmetric powers. 
\end{prooff}

\begin{notation} \label{Phifuncnot}
Let us denote by:
\begin{align*}
\Phi: \Chowek \lra \DMeetkQ
\end{align*}
the composite of the functors $j: \DMegk \sx \Q \ra \DMeetkQ$ from Construction \ref{jconst} and 
$\Psi: \Chowek \ra \DMegk \sx \Q$ of Proposition \ref{psimon}.
\end{notation}

\begin{proposition}[\tbf{Voevodsky}] \label{phisymcomm}
The functor\tn{:}
\begin{align*}
\Phi: \Chowek \lra \DMeetkQ
\end{align*}
of Notation \ref{Phifuncnot} is an additive tensor functor, and for all non-negative integers $n$, the diagrams:
\begin{align*}
\xymatrix{
\Chowek \ar@{->}[rr]^(0.45){\Phi} \ar@{->}[dd]_(0.5){\Symn} & & 
\DMeetkQ \ar@{->}[dd]^(0.5){\Symn} \\
& & \\
\Chowek \ar@{->}[rr]_(0.45){\Phi} & & \DMeetkQ} \qaq 
\xymatrix{
\Chowek \ar@{->}[rr]^(0.45){\Phi} \ar@{->}[dd]_(0.5){\Sym} & & 
\DMeetkQ \ar@{->}[dd]^(0.5){\Sym} \\
& & \\
\Chowek \ar@{->}[rr]_(0.45){\Phi} & & \DMeetkQ} 
\end{align*}
both commute up to natural isomorphism. Moreover, there is a natural isomorphism $F: M \ra \Phi \circ h$ of functors. 
\end{proposition}
\begin{prooff}
Since: 
\begin{align*}
& j: \DMegk \sx \Q \ra \DMeetkQ \qaq \\
& \Psi: \Chowek \ra \DMegk \sx \Q
\end{align*}
are both additive tensor functors by Remark \ref{jttens} and Proposition \ref{psimon}, respectively, so too is their composition $\Phi$. It now follows from Proposition \ref{symaltcomm} that $\Phi$ commutes with symmetric powers. 
The natural isomorphism $F: \Phi \circ h \ra M$ comes from applying the functor $j$ to the natural isomorphism 
$F_{\tn{gm}}: \Psi \circ h \ra \Mgm$ of Proposition \ref{psimon}.
\end{prooff}
\chapter{Motives of commutative group schemes}

Let $G$ be a \emph{homotopy invariant commutative group scheme} over a perfect field $k$, by which we mean that the associated {\etl} scheme $\shtG$ with transfers is homotopy invariant (see Definition \ref{hominvGdef}). In this chapter we construct, using work of
Barbieri-Viale/Kahn \cite{BVK}, Spie{\ss}/Szamuely \cite{SS}, and Suslin/Voevodsky \cite{SV}, a morphism:
\begin{align*}
\vp_{G}: M(G) \lra \Sym(\Mof(G))
\end{align*}
in the category $\DMeetkQ$ of triangulated motives, which is functorial in $G$ and compatible with products of group schemes. The object $\Mof(G)$ is defined as follows: in the fourth section a construction is given of a certain {\etl} sheaf $\shtG$ with transfers, whose underlying ordinary sheaf is $\shtG = \Mor_{\Schk}(\pul{X}, G)$; we set 
$\Mof(G) := \shtG$, regarded as a complex in $\DMeetkQ$ concentrated in degree zero. Let us outline the structure of the chapter.

\tbf{(i)} In the first section of this chapter we review two equivalent definitions of a group scheme $G$ over a general base scheme $S$. Firstly, we can think of $G$ as an honest $S$-scheme equipped with extra data in the form of various morphisms satisfying commutative diagrams reflecting the axioms of a group. Alternatively, one can regard a group scheme as a representable functor from the category of $S$-schemes to the category of groups. We recall the result, due to Grothendieck, that such a functor is a sheaf on the big {\etl} site over $\SchS$. This observation is important because it is most natural when working in the category of {\etl} motives to regard group schemes as {\etl} sheaves. 

\tbf{(ii)} Although we are ultimately concerned only with group schemes over the spectrum of a perfect field, it will be necessary for us to discuss the $n^{th}$ symmetric power $\Sn(X/S)$ over a general base scheme $S$, which is the subject of the second section. 

\tbf{(iii)} In the third section, we recall the construction, due to Suslin and Voevodsky \cite[\S 6]{SV}, of the morphism:
\begin{align*}
\tXYW: X \lra \Sd(Y).
\end{align*}
This morphism, which is associated to a degree $d$ correspondence $W$ in the product $X \x Y$ of two smooth separated $k$-schemes, will be used in subsequent sections to construct a canonical morphism $\alpha_{G}: M(G) \ra \Mof(G)$ in the category $\DMeetkQ$. 

\tbf{(iv)} In the fourth section, we recall, following Spieß and Szamuely \cite{SS}, that the associated ordinary {\etl} sheaf $\wt{G}$ has a natural structure of sheaf with transfers. More precisely, there is a sheaf $\shtG$ with transfers extending $\wt{G}$. Moreover, Orgogozo \cite{Org} has proven that if $G$ is a semiabelian, then $\wt{G}$ is homotopy invariant, hence defines an object in $\DMeetkQ$.

\tbf{(v) \& (vi)} In the fifth section, we recall the construction of a canonical morphism $\gamma_{G}: L(G) \ra \shtG$. Under the assumption that $\wt{G}$ is a homotopy invariant sheaf with transfers, we use $\gamma_{G}$ in the sixth section to construct a canonical morphism $\alpha_{G}: M(G) \ra \Mof(G)$. Again both morphisms are taken from Spieß and Szamuely \cite{SS}.

\tbf{(vii) \& (viii)} We use the morphism $\alpha_{G}: M(G) \ra \Mof(G)$ in the seventh section to construct a canonical morphism $\vp_{G}: M(G) \lra \Sym(\Mof(G))$. In the eighth section, we show that $\vp_{G}$ also preserves products of group schemes.\\\\
We emphasize that, with the exception of the construction of the morphism $\vp_{G}: M(G) \lra \Sym(\Mof(G))$ and associated functoriality results in the eighth section, none of the material in this chapter is original. In particular: the presentations of group schemes and quotients by group schemes in the first and second sections follow closely those found in \cite[Chapter 0]{Mum2} of the book of Fogarty, Kirwan and Mumford; all material concerning symmetric powers in various categories is classical, and the third section on the morphism $\tXYW$ is simply an expanded account of a construction \cite[Preamble to Theorem 6.8]{SV} of Suslin and Voevodsky. As mentioned above, Sections four to six are largely revision of material from Spieß and Szamuely \cite{SS}.

\begin{remark}
Throughout this chapter, $k$ will be a perfect field and $G$ and $H$ will represent quasi-projective commutative group schemes over $k$. All motives are taken to be in the category $\DMeetkQ$. Since  $\DMeetkQ$ is idempotent complete and $\Q$-linear by Proposition \ref{DMeetktentri}, the symmetric power functors $\Symn$ are well-defined in this category. 
\end{remark}

\section{Group schemes}

\begin{notation} \label{Xshfnot}
Throughout this section $S$ will be an arbitrary base scheme. For any $S$-scheme $X$, the presheaf of sets:
\begin{align*}
\Mor_{\SchS}(\pul{T}, X): (\SchS)^{\op} \lra \Sets
\end{align*}
will be denoted by $\wt{X}$.
\end{notation}

\begin{definition}[\tbf{\cite[Definition 0.1]{Mum2}}] \label{gpschdef}
A \emph{group scheme} over $S$ is an $S$-scheme $G$, together with an \emph{identity section} $e_{G}: S \ra G$, \emph{inverse morphism} $i_{G}: G \ra G$, and \emph{multiplication law} $m_{G}: G\x_{S} G \ra G$ such that the diagrams: 
\begin{align*}
& \tn{\tbf{(a)}} \quad 
\xymatrix{
G \x_{S} G \x_{S} G \ar@{->}[dd]_(0.5){m_{G} \x 1_{G}} \ar@{->}[rr]^(0.55){1_{G} \x m_{G}} & & 
G \x_{S} G \ar@{->}[dd]^(0.5){m_{G}} \\
& & \\
G \x_{S} G \ar@{->}[rr]_(0.55){m_{G}} & & G} \qquad  
\tn{\tbf{(b)}} \quad 
\xymatrix{
G \ar@{->}[d]_(0.5){\pi_{G}} \ar@{->}[rrr]^(0.5){\Delta_{G/S}} & & & G \x_{S} G \ar@<0.5ex>@/^3pt/[dd]^{1_{G} \x_{S} i_{G}} 
\ar@<-0.5ex>@/_2pt/[dd]_{i_{G} \x_{S} 1_{G}} \\
S \ar@{->}[d]_(0.5){e_{G}} & & & \\
G \ar@{<-}[rrr]_(0.5){m_{G}} & & & G \x_{S} G} \\
& \tn{\tbf{(c)}} \quad \xymatrix{
G \ar@{->}[r]^(0.4){\simeq} \ar@{->}[d]_(0.4){\simeq} \ar@{->}[ddrrr]^(0.5){1_{G}} & S \x_{S} G \ar@{->}[rr]^{e_{G} \x_{S} 1_{G}} & &
G \x_{S} G \ar@{->}[dd]^(0.5){m_{G}} \\
G \x_{S} S \ar@{->}[d]_{1_{G} \x_{S} e_{G}} & & & \\
G \x_{S} G \ar@{->}[rrr]_(0.5){m_{G}} & & & G}
\end{align*}
commute. Here, $\pi_{G}: G \ra S$ and $\Delta_{G/S}: G \ra G \x_{S} G$ are the structure and diagonal morphisms, respectively. We say that $G$ is a \emph{commutative group scheme} if the diagram:
\begin{align*}
\xymatrix{
G \x_{S} G \ar@{->}[dd]_(0.5){\tau_{G}} \ar@{->}[rr]^(0.55){m_{G}} & & G \\
& & \\
G \x_{S} G \ar@{->}[uurr]_(0.55){m_{G}} & &} 
\end{align*}
is also commutative, where $\tau_{G}: G \x_{S} G \ra G \x_{S} G$ is the morphism exchanging the factors. 
\end{definition}

\begin{remark} 
The functor $\SchS \ra \Fun((\SchS)^{\op}, \Sets)$ sending $X$ to $\wt{X}$ preserves products and is fully faithful, by the Yoneda Lemma. Consequently, we can equivalently define a group scheme over $S$ as a representable functor from $(\SchS)^{\op}$ to the category $\Gp$ of abstract groups. If the functor factors through the subcategory of abelian groups, then we call it a \emph{commutative group scheme}. 
\end{remark}

\begin{proposition}[\tbf{Grothendieck}]  \label{Gsheaf}
Let $G$ be a group scheme over $S$. Then $\wt{G}$ is a sheaf both on the big {\etl} and on the big flat sites $(\SchS)_{\et}$ and $(\SchS)_{\fl}$ over $S$.
\end{proposition}
\begin{prooff}
\cite[Expos{\'e} VII, \S 2]{SGA4t2}.
\end{prooff}

\section{Symmetric powers of a quasiprojective scheme}

The purpose of this section is to define, for any non-negative integer $n$ and quasiprojective $k$-scheme $X$, the $n^{\tn{th}}$ symmetric power $\Sn(X)$.

\subsection{Quotients of schemes by group actions}

\begin{definition}[\tbf{\cite[Definition 0.3]{Mum2}}] \label{actdef}
Let $X$ be a scheme over $S$ and $G$ be a group scheme over $S$ with identity section 
$e_{G}: \Spec(S) \ra G$ and group law $m_{G}: G \x_{S} G \ra G$. An \emph{action} of $G$ on $X$ is a morphism of $S$-schemes $\rho: G \x_{S} X \ra X$ such that the diagrams: 
\begin{align*}
\xymatrix{
G \x_{S} G \x_{S} X \ar@{->}[rr]^(0.5){1_{G} \x_{S} \rho} \ar@{->}[dd]_(0.5){m_{G} \x_{S} 1_{X}} & & 
G \x_{S} X \ar@{->}[dd]^(0.5){\rho} \\
& & \\
G \x_{S} X \ar@{->}[rr]_(0.5){\rho} & & X} \qaq
\xymatrix{
X  \ar@{->}[dd]_(0.5){1_{X}} \ar@{->}[rr]^(0.5){\simeq} & & S \x_{S} X \ar@{->}[dd]^(0.55){e_{G} \x_{S} 1_{X}} \\
& & \\
X \ar@{<-}[rr]_(0.5){\rho} & & G \x_{S} X}
\end{align*}
commute.  
\end{definition}

\begin{remark} \label{gpact}
Let $\Gamma$ be an abstract group and denote by $\Gamma^{\tn{cs}}$ the associated constant group scheme over $S$. Then we define an action of $\Gamma$ on an $S$-scheme $X$ to be an action $\rho: \Gamma^{\tn{cs}} \x_{S} X \ra X$ in the sense of Definition \ref{actdef}. Equivalently, an action of $\Gamma$ on $X$ can be regarded as a morphism of monoids:
\begin{align*}
\rho: \Gamma \lra \End_{\SchS}(X).
\end{align*}
We will somtimes abuse notation by writing $g$ in place of $\rho(g)$.
\end{remark}

\begin{notation}[\tbf{\cite[Definition 0.4]{Mum2}}]  \label{Psinot}
Let $X$ be a scheme over $S$ and let $G$ be a group scheme over $S$. Denote by 
$\Psi = \Psi_{G, X}: G \x_{S} X \ra X \x_{S} X$ the unique morphism, induced by the universal property of the fibre product, such that the diagram:
\begin{align*}
\xymatrix{
& G \x_{S} X \ar@/^1.2pc/@{->}[ddr]^(0.5){p_{X}} \ar@{->}[d]^(0.5){\Psi}  
\ar@/_1.2pc/@{->}[ddl]_(0.5){\rho} \ar@{->}[d]^(0.5){\Psi}& \\
& X \x_{S} X \ar@{->}[dr]_(0.5){p_{2}} \ar@{->}[dl]^(0.5){p_{1}} & \\
X & & X}
\end{align*}
is commutative.
\end{notation}

\begin{definition}[\tbf{Geometric quotient -- see \cite[Definition 0.6]{Mum2}}]  \label{quotdef}
Let $X$ be a scheme over $S$, let $G$ be a group scheme over $S$ and let $\rho: G \x X \ra X$ be an action of $G$ on $X$. Given a pair $(Y, q)$ consisting of an $S$-scheme $Y$ and morphism $q: X \ra Y$, we say that $Y$ is a \emph{(geometric) quotient of $X$ by $G$ with canonical projection $q$} if the following conditions are satisfied: 
\begin{itemize}
\item[\tn{(i)}] The diagram: 
\begin{align*}
\xymatrix{
G \x_{S} X \ar@{->}[rr]^(0.5){\rho} \ar@{->}[dd]_(0.5){p_{X}} & & X \ar@{->}[dd]^(0.5){q} \\
& & \\
X \ar@{->}[rr]_(0.5){q} & & Y} 
\end{align*}
commutes. 
\item[\tn{(ii)}] The scheme-theoretic image of the morphism $\Psi: G \x_{S} X \ra X \x_{S} X$ of Notation \ref{Psinot} is $X \x_{Y} X$.
\item[\tn{(iii)}] The map of topological spaces underlying $q$ is open and surjective.
\item[\tn{(iv)}] The induced sheaf morphism $q^{\#}: \shfO_{Y} \ra q_{*} \sOX$ identifies $\shfO_{Y}$ with the subsheaf $\sOX^{G}$ of $\sOX$ of $G$-invariant sections. 
\end{itemize}
In this case we write $Y = X/G$. 
\end{definition}

\begin{proposition}[\tbf{Mumford}] \label{catquot}
Let $\rho: G \x_{S} X \ra X$ be an $S$-group scheme action, suppose that the quotient $X/G$ exists, and denote by $q = q_{X/G}: X \ra X/G$ the associated quotient map. Then for any morphism $r:X \ra Z$ of $S$-schemes such that the diagram\tn{:} 
\begin{align*}
\xymatrix{
G \x_{S} X \ar@{->}[rr]^(0.5){\rho} \ar@{->}[dd]_(0.5){p_{X}} & & X \ar@{->}[dd]^(0.5){r} \\
& & \\
X \ar@{->}[rr]_(0.5){r} & & Z} 
\end{align*}
commutes, there exists a unique morphism $s: X/G \ra Z$ such that $r = s \circ q$. Consequently, the quotient $X/G$ is unique up to unique isomorphism. 
\end{proposition}
\begin{prooff}
\cite[Proposition 0.1]{Mum2}.
\end{prooff}

\begin{remark} 
The universal condition in Proposition \ref{catquot} is the definition of the \emph{categorical quotient of $X$ by $G$}, due to Mumford \cite[Definition 0.5]{Mum2}. The proposition therefore says that the geometric quotient of Definition \ref{quotdef} is also a categorical quotient.
\end{remark}

\begin{remark} \label{unpropGact}
Let $\rho: \Gamma \x_{S} X \ra X$ denote an action of an abstract group on an $S$-scheme $X$, such that $X/\Gamma$ exists. Then Proposition \ref{catquot} implies that $X/\Gamma$ possesses the following universal property. For any morphism of $S$-schemes $f: X \ra Y$ such that $f \circ g = f$ for all $g \in \Gamma$, there exists a unique morphism $h: X/G \ra Y$ such that $h \circ q_{X/G} = f$. In other words, there exists an $h$ such that the upper triangle in the diagram:
\begin{align*}
\xymatrix{
X \ar@{->}[rr]^(0.5){q_{X/G}} \ar@{->}[dd]_(0.5){g} \ar@{->}[ddrr]^(0.5){f} & & 
X/G \ar@{->}[dd]^(0.5){h} \\
& & \\
X \ar@{->}[rr]_(0.5){f} & & Y} 
\end{align*}
commutes, provided the lower one does, for all $g \in \Gamma$.
\end{remark}

\begin{proposition} \label{qpquot}
Let $\Gamma$ be a finite group acting on a quasiprojective $S$-scheme $X$. Then the quotient $X/\Gamma$ exists.
\end{proposition}
\begin{prooff}
This follows from \cite[Expos{\'e} V, Proposition 1.8]{SGA1}.
\end{prooff}

\subsection{Symmetric powers}

For a positive integer $n$, we denote by $\Sym(n)$ the symmetric group on $n$ letters. Given an $S$-scheme $X$, we define in this subsection the $n^{\tn{th}}$ symmetric power $\Sn(X/S)$.

\begin{notation}
Let $X$ be a quasiprojective scheme over $S$ and let $n$ be a positive integer. Then $\Sym(n)$ acts on 
\begin{align*}
X^{n} := \ub{X \x_{S} \cdots \x_{S} X}_{n \tn{ copies}}
\end{align*}
by permutation of the factors. We define $\Sn(X/S) := X^{n}/\Sym(n)$, which we call the \emph{$n^{\tn{th}}$ symmetric power of $X$}. We denote the quotient $X^{n} \ra \Sn(X/S)$ by $q_{X}^{n}$ in place of $q_{X/\Sn(X/S)}$. We also set $\tn{S}^{0}(X/S) := S$ and define $q_{X}^{0}$ to be the identity morphism $1_{S}: S \ra S$.
\end{notation}

\begin{remark} \label{unpropSact}
Let $X$ be a quasiprojective $S$-scheme and let $n$ be a positive integer. Then the universal property of 
Remark \ref{unpropGact} reads as follows. Let $f: X \ra Y$ be a morphism of $S$-schemes such that the lower triangle in the diagram:
\begin{align*}
\xymatrix{
X^{n} \ar@{->}[rr]^(0.5){q_{X}^{n}} \ar@{->}[dd]_(0.5){\sigma} \ar@{->}[ddrr]^(0.5){f} & & 
\Sn(X/S) \ar@{->}[dd]^(0.5){h} \\
& & \\
X^{n} \ar@{->}[rr]_(0.5){f} & & Y} 
\end{align*}
commutes, for all $\sigma \in \Sym(n)$. Then there exists a unique morphism $h: \Sn(X/S) \ra Y$ such that the upper triangle commutes as well. Using this property, the symbol $\Sn$ can be extended to a functor from the category of quasiprojective $S$-schemes to the category of $S$-schemes. Indeed, if $f:X \ra Y$ is a morphism of $S$-schemes then 
$q^{n}_{Y} \circ f^{n} \circ \sigma = q^{n}_{Y} \circ f^{n}$, for all $\sigma \in \Sym(n)$. This induces a unique morphism:
\begin{align*}
\Sn(f): \Sn(X/S) \lra \Sn(Y/S)
\end{align*}
such that the diagram:
\begin{align*}
\xymatrix{
X^{n} \ar@{->}[rr]^(0.5){q_{X}^{n}} \ar@{->}[dd]_(0.5){f^{n}} & & \Sn(X/S) \ar@{->}[dd]^(0.5){\Sn(f)} \\
& & \\
Y^{n} \ar@{->}[rr]_(0.5){q_{Y}^{n}} & & \Sn(Y/S)} 
\end{align*}
commutes.
\end{remark}

\begin{remark}
Let $n$ be a positive integer. Recall the construction of the $n^{\tn{th}}$ symmetric power in the category of sets: if $X$ is a set, then $\Sn(X)$ is defined to be the set of all unordered $n$-tuples of elements of $X$. Formally, it is the quotient of $X^{n}$ by the equivalence relation identifying any two $n$-tuples whose elements are identical up to a permutation $\sigma \in \Sym(n)$. We denote the image of a typical $n$-tuple $(x_{1}, \ldots, x_{n}) \in X^{n}$ under the obvious projection $X^{n} \ra \Sn(X)$ by $x_{1} x_{2} \cdots x_{n}$ (juxtaposition of entries). We also set 
$\tn{S}^{0}(X) := S$.
\end{remark}

\begin{remark} \label{symcomp}
Let $X$ be a quasiprojective scheme over $S$, and let $n$ be a positive integer. Define for any $S$-scheme $T$ a morphism of sets: 
\begin{align*}
\mathsmap{\psi}{\wt{X}(T)^{n}}{\wt{\Sn(X/S)}(T)}{(f_{1}, \ldots, f_{n})}{q_{X}^{n} \circ f,}
\end{align*}
where $q_{X}^{n}: X^{n} \ra \Sn(X/S)$ is the canonical projection and $f$ is the unique map, induced by the universal property of the fibre product $X^{n}$ over $S$, such that the diagram:
\begin{align*}
\xymatrix{
T \ar@{->}[d]_(0.5){f} \ar@{->}[dr]^(0.5){f_{i}} & \\
X^{n} \ar@{->}[r]_(0.5){\pr_{i}} & X}
\end{align*}
commutes for all $1 \le i \le n$. Then by the universal property of the symmetric power of sets, there is a unique morphism: 
\begin{align*}
\vp: \Sn \left( \wt{X}(T) \right) \lra \wt{\Sn(X/S)}(T)
\end{align*}
such that the diagram:
\begin{align*}
\xymatrix{
\wt{X}(T)^{n} \ar@{->}[ddrr]_(0.5){\psi} \ar@{->}[rr]^(0.5){q^{n}_{\wt{X}(T)}} & & 
\Sn \left( \wt{X}(T) \right) \ar@{->}[dd]^(0.5){\vp} \\
& & \\
& & \wt{\Sn(X/S)}(T)}
\end{align*}
commutes.
\end{remark}

\begin{notation}
We use the notation:
\begin{align*}
\Sinf(Y/S) := \coprod_{d \ge 0} \Sd(Y/S)
\end{align*}
for any quasiprojective $S$-scheme $Y$.
\end{notation}

\begin{remark} \label{schememon}
For any quasiprojective $S$-scheme $Y$, the $S$-scheme $\Sinf(Y/S)$ is a \emph{monoid scheme}, in the sense that the associated presheaf:
\begin{align*}
\wt{\Sinf(Y/S)}: \SchS \lra \Sets
\end{align*}
factors through the category of monoids. Indeed, for any $S$-scheme $T$ and any pair of non-negative integers $d$ and $e$, there is a natural composition law:
\begin{align*}
\mathsmaptwo{\Sd\left( \wt{Y}(T) \right) \x \Se \left( \wt{Y}(T) \right)}{\tn{S}^{d+e} \left( \wt{Y}(T) \right)}
{\Big( (x_{1} \cdots x_{d}), (y_{1} \cdots y_{e}) \Big)}
{x_{1} \cdots x_{d} \, y_{1} \cdots y_{e},}
\end{align*}
which by Remark \ref{symcomp} induces a monoid structure: 
\begin{align*}
\Sd(Y/S) \x \Se(Y/S) \lra \tn{S}^{d+e}(Y/S)
\end{align*}
on the presheaf $\Sinf(Y/S)$. Consequently, for any $S$-scheme $X$, the set\tn{:}
\begin{align*}
\Mor_{\SchS}\left(X, \Sinf(Y/S) \right)
\end{align*}
is a monoid.
\end{remark} 

\begin{construction} \label{sigmade}
Let $d$ and $e$ be non-negative integers. We define a morphism of functors:
\begin{align*}
\Sigma^{d,e}: \Sd \circ \Se \lra \Sde
\end{align*}
as follows. If $X$ is a quasiprojective $S$-scheme, then the diagram:
\begin{align*}
\xymatrix{
\left( X^{e} \right)^{d} \ar@{->}[rr]^(0.5){(q_{X}^{e})^{d}} \ar@{=}[dd] & & 
\Se(X/S)^{d} \ar@{->}[rr]^(0.5){q_{\Se(X/S)}^{d}} \ar@{->}[dd]^(0.5){f} & & 
\Sd(\Se(X/S)) \ar@{->}[dd]^(0.5){\Sigma^{d,e}(X/S)} \\
& & & & \\
X^{de} \ar@{->}[rr]_(0.5){q_{X}^{de}} & & \Sde(X/S) \ar@{=}[rr] & & \Sde(X/S).} 
\end{align*}
defines a morphism $\Sigma^{d,e}(X/S): \Sde(X/S) \ra \Sd(\Se(X/S))$. Here, the vertical morphism $f$ in the middle is the 
$d$-fold monoid composition:
\begin{align*}
\Se(X/S)^{d} = \ub{\Se(X/S) \x \cdots \x \Se(X/S)}_{d \tn{ copies}} \lra \Sde(X/S)
\end{align*}
of Remark \ref{schememon}, and $\Sigma^{d,e}(X/S)$ is defined to be the unique morphism, induced by the universal property of 
$\Sd(\Se(X/S))$, such that the right-hand square commutes. Since the diagram is constructed using universal properties, it is easy to see that the morphism $\Sigma^{d,e}(X/S)$ is functorial in $X$. That is, for any morphism $f: X \ra Y$ of quasiprojective $k$-schemes, the diagram: 
\begin{align*}
\xymatrix{
\Sd(\Se(X/S)) \ar@{->}[rr]^(0.5){\Sigma^{d,e}(X/S)} \ar@{->}[dd]_(0.5){\Sd(\Se(f))} & & 
\Sde(X/S) \ar@{->}[dd]^(0.5){\Sde(f)} \\
& & \\
\Sd(\Se(Y/S)) \ar@{->}[rr]_(0.5){\Sigma^{d,e}(Y/S)} & & \Sde(Y/S)}
\end{align*}
commutes.
\end{construction}

\begin{notation} \label{symGmult}
Let $G$ be a commutative quasiprojective group scheme over $S$, let $n$ be a non-negative integer, and denote by 
$m_{G}^{n}: G^{n} \ra G$ the $n$-fold group law. We define $\cj{m}^{n}_{G}: \Sn(G/S) \ra G$ to be the induced multiplication law on $\Sn(G/S)$. That is, it is the unique morphism making the diagram:
\begin{align*}
\xymatrix{
G^{n} \ar@{->}[rr]^(0.4){q_{G}^{n}} \ar@{->}[ddrr]_(0.5){m_{G}^{n}} & & 
\Sn(G/S) \ar@{->}[dd]^(0.5){\cj{m}^{n}_{G}} \\
& & \\
& & G}
\end{align*}
commute.
\end{notation}

\begin{proposition} \label{sdemult}
Let $G$ be a commutative quasiprojective group scheme over $S$, and let $d$ and $e$ be non-negative integers. Then the following diagram commutes\tn{:}
\begin{align*}
\xymatrix{
\Sd(\Se(G/S)) \ar@{->}[rr]^(0.5){\Sd(\cj{m}^{e}_{G}/S)} \ar@{->}[dd]_(0.5){\Sigma^{d,e}(G/S)} & & 
\Sd(G/S) \ar@{->}[dd]^(0.5){\cj{m}^{d}_{G}} \\
& & \\
\Sde(G/S) \ar@{->}[rr]_(0.5){\cj{m}^{de}_{G}} & & G.}
\end{align*}
\end{proposition}
\begin{prooff}
This can be checked by a diagram chase, using the universal properties of the morphisms $\Sigma^{d,e}(X)$ and 
$\cj{m}^{d}_{G}$ described in Construction \ref{sigmade} and Notation \ref{symGmult}, respectively. Alternatively, for an arbitrary $S$-scheme $T$, it is enough by Remark \ref{symcomp} to show that the diagram:
\begin{align*}
\xymatrix{
\Sd \left( \Se(\wt{G}(T)) \right) \ar@{->}[rr]^(0.55){\Sd \left( \cj{m}^{e}_{G}/S \right)(T)} 
\ar@{->}[dd]_(0.5){\Sigma_{G}^{d,e}(T)} & & 
\Sd \left( \wt{G}(T) \right) \ar@{->}[dd]^(0.5){\cj{m}^{d}_{G}(T)} \\
& & \\
\Sde \left( \wt{G}(T) \right) \ar@{->}[rr]_(0.5){\cj{m}^{de}_{G}(T)} & & \wt{G}(T)}
\end{align*}
commutes, which is clear.
\end{prooff}

\section{The morphism $\tXYW$}

Let $k$ be a perfect field, let $X$ and $Y$ be two smooth schemes over $k$, and let $W$ be a smooth integral subscheme of $X \x Y$ which is finite and surjective over a connected component of $X$. Denote by $d$ be the degree of the resulting projection $W \ra X$. We recall the construction, due to Suslin and Voevodsky 
\cite[Preamble to Theorem 6.8]{SV}, of the map:
\begin{align*}
\tXYW: X \lra \Sd(Y)
\end{align*}
appearing in the paper \cite[Lemma 3.2]{SS} of Spie{\ss} and Szamuely.

\begin{proposition} \label{WSPec}
Let $W$ and $X$ be $k$-schemes, with $W$ integral and $X$ smooth, and suppose $f: W \ra X$ is a finite flat morphism of constant degree $d$. Then $\shfA := f_{*}(\shfO_{W})$ is a locally free sheaf on $X$ of $\shfO_{X}$-algebras of constant rank $d$, and 
$W = \Spec(\shfA)$.
\end{proposition}
\begin{prooff}
The result is implied by \cite[Corollary 1.3.2]{EGA2}, since the finite morphism $f$ is affine.
\end{prooff}

\begin{notation} \label{symaltshf}
Let $X$ be a scheme over a perfect field $k$, let $\shfF$ be a locally free sheaf of $\sOX$-algebras on $X$, and let $n$ be a non-negative integer. The following are standard definitions.
\begin{itemize}
\item[\tn{(i)}] The \emph{tensor algebra sheaf} $\tn{T}(\shfF)$ associated to $\shfF$ is the direct sum 
$\bos_{n \ge 0} \shfF^{\ox n}$, of sheaves of $\sOX$-modules, with the usual convention $\shfF^{\ox 0} := \sOX$. 
\item[\tn{(ii)}] The \emph{exterior algebra sheaf} $\bigwedge \shfF$ is the sheaf attached to the presheaf 
$U \mapsto \tn{T}(\shfF)(U)/I(U)$, where $I(U)$ is the two-sided ideal $I$ generated by all elements of the form $x \ox x$, with $x \in \shfF(U)$. The image of the composition:
\begin{align*}
\shfF^{\ox n} \lra \tn{T}(\shfF) \lra \bigwedge \shfF,
\end{align*}
where the morphism on the right is the natural projection, is denoted $\bigwedge^{n}(\shfF)$, and called the \emph{$n^{\tn{th}}$ exterior power} of $\shfF$.
\item[\tn{(iii)}] By $(\shfF^{\ox n})^{\Sym(n)}$ we denote the subpresheaf of $\shfF^{\ox n}$ invariant under the action of the symmetric group $\Sym(n)$ on $n$ letters. 
\end{itemize}
Following Suslin and Voevodsky \cite[Preamble to Theorem 6.8]{SV}, we observe that the algebra multiplication law on $\shfF^{\ox n}$ induces a bilinear pairing:
\begin{align*}
\bfm_{\shfF, n}: (\shfF^{\ox n})^{\Sym(n)} \ox_{\sOX} \bigwedge^{n} \shfF \lra \bigwedge^{n} \shfF
\end{align*}
of $\sOX$-modules.
\end{notation}

\begin{construction}[\tbf{Suslin and Voevodsky}] \label{tXWconst}
We recall a construction of Suslin and Voevodsky \cite[Preamble to Theorem 6.8]{SV}. Let $W$ and $X$ be $k$-schemes, with $W$ integral and $X$ smooth, let $f: W \ra X$ be a finite surjective morphism, write:
\begin{align*}
d:= [k(W): k(X)],
\end{align*}
and set $\shfA := f_{*}(\shfO_{W})$, so that $W = \Spec(\shfA)$, by Proposition \ref{WSPec}. We construct a morphism:
\begin{align*}
t_{X, W}: X \lra \Sd(W/X)
\end{align*}
of $k$-schemes as follows. Firstly, we note that it is shown in \cite[Preamble to Theorem 6.8]{SV} that one may assume the morphism $f$ is flat and of constant degree $d$. By Notation \ref{symaltshf}, we have a bilinear pairing of $\sOX$-modules:
\begin{align*}
\bfm_{\shfA, n}: (\shfA^{\ox d})^{\Sym(d)} \ox_{\sOX} \bigwedge^{d} \shfA \lra \bigwedge^{d} \shfA.
\end{align*}
This pairing induces, by adjunction of $\Hom$ and $\ox$ in the category of $\sOX$-modules, a representation:
\begin{align*}
\rho_{\shfA}: (\shfA^{\ox d})^{\Sym(d)} \lra \End_{\sOX}\left( \bigwedge^{d} \shfA \right) = \sOX,
\end{align*}
where the isomorphism holds because $\shfA$ is a locally free $\sOX$-algebra of degree $d$. Now, one can check that:
\begin{align*}
\Spec(\sOX) \simeq X \qaq \Spec(\Symd(\shfA)) \simeq \Sd(\Spec(\shfA)) \simeq \Sd(W/X),
\end{align*}
where $\Spec$ denotes the sheaf-spec functor. Applying $\Spec$ to $\rho_{\shfA}$ therefore gives:
\begin{align*}
t_{X, W} := \Spec(\rho_{\shfA}): X \lra \Sd(W/X),
\end{align*}
as required. 
\end{construction}

\begin{construction}[\tbf{Suslin and Voevodsky}] \label{thetacon}
We recall a construction of Suslin and Voevodsky \cite[Preamble to Theorem 6.8]{SV}. Let $X$ and $Y$ be smooth schemes over a perfect field $k$ with $X$ connected, and let $W$ be an integral closed subscheme of $X \x Y$, finite and surjective over a component of $X$, so that $W \in \VCork(X, Y)$. Denote by:
\begin{align*}
\iota_{W}: W \lra X \x Y
\end{align*}
the associated inclusion, and write $d$ for the degree of the natural projection $p_{X} \circ \iota_{W}: W \ra X$. One constructs a morphism\tn{:}
\begin{align*}
\tXYW: X \lra \Sd(Y)
\end{align*}
as follows. Intuitively, we proceed as follows: since $W \sseq X \x Y$ has degree $d$, it associates to every point 
point $x$ in $X$ a collection of $d$ points $y_{1}, \ldots, y_{d}$ in $Y$. The morphism $\tXYW$ sends $x$ to the unordered $d$-tuple in $\Sd(Y)$ determined by the $y_{i}$. More formally, consider the composite mapping:
\begin{align*}
g: X \sxra{t_{X, W}} & \Sd(W/X) \sxra{\Sd(\iota_{W}/X)} \Sd((X \x Y)/X) \\
\sxra{\lambda_{X, Y}} & X \x_{\Speck} \Sd(Y),
\end{align*}
where:
\begin{align*}
\lambda_{X, Y}: \Sd((X \x Y)/X) \lra X \x_{\Speck} \Sd(Y)
\end{align*}
is the canonical base-change isomorphism associated to $\Sd$, recalling that $\Sd(Y) =  \Sd(Y/\Speck)$. Then $g$ is a morphism of $X$-schemes, hence induces a morphism: 
\begin{align*}
\tXYW: X \lra \Sd(Y)
\end{align*}
of $k$-schemes. 
\end{construction}

\section{The induced presheaf $\shtG$ with transfers of a group scheme $G$}

\begin{construction}[\tbf{Spieß and Szamuely}] \label{betamorph}
Let $X$ and $Y$ be smooth schemes over a perfect field $k$, let $G$ be a commutative group scheme over $k$, and suppose that $c := \sum n_{i} Z_{i} \in \VCork(X, Y)$ is a Voevodsky correspondence. We recall from 
\cite[Lemma 3.2]{SS} the construction of a morphism: 
\begin{align*}
\beta(c) = \beta_{X, Y, G}(c): \Mor_{\Schk}(Y, G) \lra \Mor_{\Schk}(X, G).
\end{align*}
Firstly, suppose that $c = Z$ is the correspondence associated to a single closed integral subscheme 
$Z \sseq X \x Y$, and let $d$ denote the degree of the induced projection $Z \ra X$. Then for 
$g \in \Mor_{\Schk}(Y, G)$, we define $\beta(c)(g) \in \Mor_{\Schk}(X, G)$ to be the composition:
\begin{align*}
\beta(c)(g): X \sxra{\tXYZ} \Sym^{d}(Y) \sxra{\Sym^{d}(g)} \Sym^{d}(G) \slra{\cj{m}_{G}^{d}} G, 
\end{align*}
where $\cj{m}_{G}^{d}$ denotes canonical morphism of group schemes induced by the $d$-fold summation map 
$G^{d} \ra G$, as in Notation \ref{symGmult}. Since $G$ is a commutative group scheme, $\Mor_{\Schk}(X, G)$ is an abelian group, and so we can define for general $c = \sum n_{i} Z_{i} \in \VCork(X, Y)$ a general correspondence 
$\beta(c) := \sum_{i} n_{i} \cdot \beta(Z_{i})$. Thus:
\begin{align*}
\beta(c)(g) = \sum_{i} n_{i} \cdot \beta(Z_{i})(g) = 
\sum_{i} n_{i} \cdot \left( \cj{m}_{G}^{d_{i}} \circ \Sym^{d_{i}}(g) \circ \theta_{X, Y}^{Z_{i}} \right), 
\end{align*}
where $d_{i}$ denotes the degree of the projection $Z_{i} \ra X$.
\end{construction}

\begin{proposition}[\tbf{Spieß and Szamuely}] \label{Gtrans}
Let $G$ be a commutative group scheme over a perfect field $k$. Then $\wt{G} := \Mor_{\Schk}(\pul{X}, G)$ is a sheaf with transfers on the big {\etl} site. More precisely, 
\begin{align*}
\mathsmap{\sht{G}}{\Smck}{\Ab}
{\begin{array}{rc}
\tn{\small{Objects}:} & X \\
\tn{Morphisms:} & c \in \VCork(X, Y)
\end{array}}
{\begin{array}{c}
\sh{G}(X) \\
\beta(c): \sh{G}(Y) \ra \sh{G}(X)
\end{array}}
\end{align*}
defines a sheaf with transfers such that the diagram of functors\tn{:}
\begin{align*}
\xymatrix{
(\Smck)^{\op} \ar@{->}[rr]^(0.6){\sht{G}} & & \Ab \\
& & \\
(\Smk)^{\op} \ar@{->}[uurr]_(0.5){\sh{G}} \ar@{->}[uu]^(0.5){(\scd)^{\op}} & &}
\end{align*}
commutes. 
\end{proposition}
\begin{prooff}
See \cite[Lemma 3.2]{SS}.
\end{prooff}

\begin{convention}
Let $G$ be a commutative group scheme $G$ over $k$. From now on, we agree to denote the associated {\etl} sheaf with transfers $\shtG$ simply by $\wtG$, to ease the notation. 
\end{convention}

\begin{definition} \label{hominvGdef}
Let $G$ be a commutative group scheme $G$ over $k$. We say that $G$ is \emph{homotopy invariant}
if the associated {\etl} sheaf with transfers $\wtG$ is homotopy invariant.
\end{definition}

\begin{proposition}[\tbf{Orgogozo}] \label{savhinv}
A semiabelian variety over $k$ is a homotopy invariant commutative group scheme. 
\end{proposition}
\begin{prooff}
\cite[Lemma 3.3.1]{Org}.
\end{prooff}

\section{The sheaf morphism $\gamma_{G}: L(G) \ra \wtG$}

\begin{construction}[\tbf{Spieß and Szamuely}]  \label{gammaX}
Let $G$ be a commutative group scheme over a perfect field $k$, and let $X$ be a smooth $k$-scheme. We recall from 
\cite[\S 3]{SS} the morphism:
\begin{align*}
\gamma_{G}(X): L(G)(X) \lra \wtG(X)
\end{align*}
of abelian groups. By Construction \ref{betamorph}, we have for all $c \in \VCork(X, G)$ a morphism:
$\beta_{X, G, G}(c)(\id_{G}) \in \Mor_{\Schk}(X, G)$. We now define:
\begin{align*}
\mathsmap{\gamma_{G}(X)}{\VCork(X, G)}{\Mor_{\Schk}(X, G)}{c}{\beta_{X, G, G}(c)(\id_{G}).}
\end{align*}
We will show below in Proposition \ref{gammdef} that $\gamma_{G}(X)$ is natural in $X$.
\end{construction}

\begin{proposition} \label{gammdef}
The morphism $\gamma_{G}(X)$ of Construction \ref{gammaX} defines a morphism\tn{:} 
\begin{align*}
\gamma_{G}: L(G) \lra \wtG
\end{align*}
of {\etl} sheaves with transfers.
\end{proposition}
\begin{prooff}
Suppose $c \in \VCork(X, Y)$ and $d \in \VCork(Y, G)$. Since $\wtG$ is a sheaf, we have that 
$\wtG(d \circ c) = \wtG(c) \circ \wtG(d)$, which means that:
\begin{align}
\beta_{X, G, G}(d \circ c) = \beta_{X, Y, G}(c) \circ \beta_{Y, G, G}(d). \label{bcomprel}
\end{align}
Now, to show that $\gamma_{G}$ is a morphism of {\etl} sheaves with transfers, we must show that the diagram:
\begin{align*}
\xymatrix{
\VCork(Y, G) \ar@{->}[rr]^(0.5){\gamma_{G}(Y)} \ar@{->}[dd]_(0.5)(0.5){c^{*} = \VCork(c, G)} & & 
\Mor_{\Schk}(Y, G) \ar@{->}[dd]^(0.5){\wtG(c) = \beta_{X, Y, G}(c)}\\
& & \\
\VCork(X, G) \ar@{->}[rr]_(0.5){\gamma_{G}(X)} & & \Mor_{\Schk}(X, G)}
\end{align*}
commutes. Since $c^{*}(d) = d \circ c \in \VCork(X, G)$, we have:
\begin{align}
(\gamma_{G}(X) \circ c^{*})(d) = \gamma_{G}(X)(d \circ c) = \beta_{X, G, G}(d \circ c)(\id_{G}).
\label{bcd}
\end{align}
On the other hand, 
\begin{align}
\nonumber \left( \wtG(c) \circ \gamma_{G}(Y) \right)(d) & = 
\wtG(c) \big( \beta_{Y, G, G}(d)(\id_{G}) \big) \\
& = \beta_{X, Y, G}(c) \big( \beta_{Y, G, G}(d)(\id_{G}) \big).
\label{bccd}
\end{align}
By the relation (\ref{bcomprel}), we see that (\ref{bcd}) and (\ref{bccd}) agree.
\end{prooff}

\begin{proposition} \label{gammanatG}
The morphism $\gamma_{G}$ is natural in $G$ in the following sense. Let $f: G \ra H$ be a morphism of commutative group schemes over $k$. Then the diagram\tn{:}
\begin{align*}
\xymatrix{
L(G) \ar@{->}[rr]^(0.5){\gamma_{G}} \ar@{->}[dd]_(0.5){L(f)} & & \wtG \ar@{->}[dd]^(0.5){\wtf} \\
& & \\
L(H) \ar@{->}[rr]_(0.5){\gamma_{H}} & & \wtH}
\end{align*} 
commutes.
\end{proposition}
\begin{prooff}
Let $X$ be a smooth $k$-scheme. We wish to show that the diagram:
\begin{align*}
\xymatrix{
\VCork(X, G) \ar@{->}[rr]^(0.5){\gamma_{G}(X)} \ar@{->}[dd]_(0.5){L(f)(X)} & & 
\Mor_{\Schk}(X, G) \ar@{->}[dd]^(0.5){\wtf(X)} \\
& & \\
\VCork(X, H) \ar@{->}[rr]_(0.5){\gamma_{H}(X)} & & \Mor_{\Schk}(X, H)}
\end{align*} 
commutes. We test commutativity of the square on an arbitrary element $c$ of $\VCork(X, G)$. By linearity of $\gamma_{G}(X)$, we can assume that $c$ is given by a closed integral subscheme $U \sseq X \x G$. Denoting by $d$ the degree of the associated projection $U \ra X$, we have:
\begin{align}
\nonumber \big( \wtf(X) \circ \gamma_{G}(X) \big)(U) & = \wtf(X) \big( \beta_{X, G, G}(U)(\id_{G}) \big) \\ 
& = \wtf(X) \big( \cj{m}^{d}_{G} \circ \theta_{X, G}^{U} \big) = f \circ \cj{m}^{d}_{G} \circ \theta_{X, G}^{U}. \label{gnat1}
\end{align}
On the other hand:
\begin{align}
\nonumber \big( \gamma_{H}(X) \circ L(f) \big)(U) & = \gamma_{H}(X) \big( \Gamma(f) \circ U \big) \\
\nonumber & = \beta_{X, H, H}(\Gamma(f) \circ U)(\id_{H}) \\
& = \cj{m}^{d}_{H} \circ \theta_{X, H}^{\Gamma(f) \circ U}. \label{gnat2}
\end{align}
Observe that $d := d(\Gamma(f) \circ U) = d(\Gamma(f)) \cdot d(U) = d(U)$, since $d(\Gamma(f)) = 1$. Now, proving the right-hand side of (\ref{gnat1}) is equal the right-hand side of (\ref{gnat2}) means proving that the diagram:
\begin{align*}
\xymatrix{
X \ar@{->}[rr]^(0.4){\theta_{X, G}^{U}} \ar@{->}[ddrr]_(0.5){\theta_{X, H}^{\Gamma(f) \circ U}} & & 
\Sym^{d}(G) \ar@{->}[rr]^(0.6){\cj{m}^{d}_{G}} \ar@{->}[dd]^(0.5){\Sym^{d}(f)} & & G \ar@{->}[dd]^(0.5){f} \\
& & & & \\
& & \Sym^{d}(H) \ar@{->}[rr]_(0.6){\cj{m}^{d}_{H}} & & H}
\end{align*}
commutes. But the square on the right commutes, because $f$ is a morphism of group schemes and $\cj{m}^{d}_{G}$ is $d$-fold summation in the groups $G$ and $H$, and the triangle on the left commutes by functoriality of the morphism 
$\tXYW$. It follows that $\wtf(X) \circ \gamma_{G}(X) = \gamma_{H}(X) \circ L(f)$.
\end{prooff}

\section{The morphism $\alpha_{G}: M(G) \ra \Mof(G)$}

\begin{remark}
Let $G$ be a homotopy invariant commutative group scheme over a perfect field $k$. In this section 
we recall the construction of Spie{\ss}/Szamuely \cite[\S 3]{SS} of a certain functorial morphism $\alpha_{G}: M(G) \ra \Mof(G)$. We emphasise that all definitions and results are due to Spie{\ss} and Szamuely: we merely rephrase them in our own notation and expand upon the details.
\end{remark}

\begin{notation}
Let $\shfF \in \STEkQ$ be a homotopy-invariant {\etl} sheaf with transfers. We denote by $[\shfF]_{0}$ the complex in $\DMeetkQ$ consisting of $\shfF$ concentrated in degree zero. 
\end{notation}

\begin{notation} \label{mofnot}
In the special case $\shfF = \wtG$, with $G$ a homotopy invariant commutative group scheme over $k$, we denote by $\Mof(G)$ the object $[\wtG]_{0}$ in $\DMeetkQ$.
\end{notation}

\begin{notation} \label{etaGHdef}
Let $G$ and $H$ be commutative group schemes over a perfect field $k$, and let $p_{G}: G \x H \ra G$ and 
$p_{H}: G \x H \ra H$ denote the canonical projection morphisms. Let us denote by:
\begin{align*}
\eta_{G, H}: \Mof(G \x H) \lra \Mof(G)\os \Mof(H)
\end{align*}
the morphism in $\DMeetkQ$ defined by:
\begin{align*}
\mathsmap{\eta_{G, H}(U)}{\Mor_{\Schk}(U, G \x H)}{\Mor_{\Schk}(U, G) \os \Mor_{\Schk}(U, H)}{f}
{(p_{G} \circ f, p_{H} \circ f),}
\end{align*}
for all smooth $k$-schemes $U$. It is a canonical isomorphism, by the universal property of the fibre product. 
\end{notation}

\begin{remark} \label{projGHcomm}
From the definition of $\eta_{G, H}$ in Notation \ref{etaGHdef}, we note the functorial relations $p_{\Mof(G)} \circ \eta_{G, H} = \Mof(p_{G})$ and $p_{\Mof(H)} \circ \eta_{G, H} = \Mof(p_{H})$ hold. Indeed, for any
\begin{align*}
f \in \Mor_{\Schk}(U, G \x H) = \Mof(G \x H)(U), 
\end{align*}
it is immediate that:
\begin{align*}
\left( p_{\Mof(G)}(U) \circ \eta_{G, H}(U) \right)(f) = p_{G} \circ f = \Mof(p_{G})(U)(f).
\end{align*}
The same calculation with $H$ in place of $G$ proves the other relation. 
\end{remark}

\begin{notation} \label{augmapnot}
Let that $\shfF \in \STEkQ$ be a homotopy invariant sheaf with transfers. We denote by\tn{:}
\begin{align*}
a_{\shfF}: \Cx(\shfF) \lra [\shfF]_{0},
\end{align*}
the \emph{augmentation map}, which is defined to be the identity morphism in degree zero and the zero morphism elsewhere. (Note that homotopy invariance of $\shfF$ is necessary here to be able to conclude that the augmentation map is a bona fide map of chain complexes.) In the special case $\shfF = \wtG$, with $G$ a homotopy invariant commutative group scheme over $k$, we write:
\begin{align*}
a_{G}: \Cx(\wtG) \lra \Mof(G),
\end{align*}
in place of $a_{\wtG}: \Cx(\wtG) \lra [\wtG]_{0}$.
\end{notation}

\begin{lemma}[\tbf{Spie{\ss}/Szamuely}] \label{augqi}
Suppose that $\shfF \in \STEkQ$ is homotopy invariant sheaf with transfers. Then the augmentation map\tn{:}
\begin{align*}
a_{\shfF}: \Cx(\shfF) \lra [\shfF]_{0}
\end{align*}
of Notation \ref{augmapnot} is\tn{:} 
\begin{itemize}
\item[\tn{(i)}] A quasi-isomorphism in $\DMeetkQ$.
\item[\tn{(ii)}] Functorial in $\shfF$, in the following sense. If $\shfG \in \STEkQ$ is a homotopy invariant sheaf with transfers and $f: \shfF \ra \shfG$ is a morphism of sheaves with transfers, then the diagram\tn{:}
\begin{align*}
\xymatrix{
\Cx(\shfF) \ar@{->}[rr]^(0.5){a_{\shfF}} \ar@{->}[dd]_(0.5){\Cx(f)} & & 
[\shfF]_{0} \ar@{->}[dd]^(0.5){[f]_{0}} \\
& & \\
\Cx(\shfG) \ar@{->}[rr]_(0.5){a_{\shfG}} & & [\shfG]_{0}}
\end{align*}
commutes.
\end{itemize}
\end{lemma}
\begin{prooff}
Commutativity of the diagram in (ii) is immediate from the definitions. To prove (i), we note that, as observed in \cite[\S 2]{SS}, the homotopy invariance of $\shfF$, together with the fact that $\An$ is isomorphic to the affine $n$-simplex $\Delta^{n}$, imply that: 
\begin{align*}
\Cnx(\shfF)(X) = \shfF(X \x \Delta^{n}) = \shfF(X),
\end{align*}
for all smooth $k$-schemes $X$ and non-negative integers $n$. Thus $\Cnx(\shfF) = \shfF$ for all $n$, and is the constant simplicial sheaf defined by $\shfF$, the cohomology of which is equal to $\shfF$ in degree zero and zero otherwise. 
\end{prooff}

\begin{notation} \label{alphanot}
Let $G$ be a homotopy invariant commutative group scheme over a perfect field $k$. We denote by:
\begin{align*}
\alpha_{G}: M(G) \lra \Mof(G)
\end{align*}
the morphism in $\DMeetkQ$ obtained by composing: 
\begin{align*}
\Cx(\gamma_{G}): M(G) = \Cx(L(G)) \lra \Cx(\wtG)
\end{align*}
with the augmentation map $a_{G}: \Cx(\wtG) \ra \Mof(G)$.
\end{notation}

\begin{proposition}[\tbf{Spie{\ss}/Szamuely}] \label{alphafunc}
The morphism $\alpha_{G}: M(G) \ra \Mof(G)$ of Notation \ref{alphanot} is functorial, in the sense that for any morphism $f: G \ra H$ of
homotopy invariant commutative group schemes over a perfect field $k$, the diagram\tn{:}
\begin{align*}
\xymatrix{
M(G) \ar@{->}[rr]^(0.5){\alpha_{G}} \ar@{->}[dd]_(0.5){M(f)} & & 
\Mof(G) \ar@{->}[dd]^(0.5){\Mof(f)} \\
& & \\
M(H) \ar@{->}[rr]_(0.5){\alpha_{H}} & & \Mof(H)}
\end{align*}
commutes. 
\end{proposition}
\begin{prooff}
By definition of $\alpha_{G}$, the above diagram may be written in expanded form as:
\begin{align*}
\xymatrix{
\Cx(L(G)) \ar@{->}[rr]^(0.5){\Cx(\gamma_{G})} \ar@{->}[dd]_(0.5){\Cx(L(f))} & & 
\Cx(\wtG) \ar@{->}[rr]^(0.5){a_{\wtG}} \ar@{->}[dd]_(0.5){\Cx \left( \wtf \right)} & & 
\Mof(G) \ar@{->}[dd]^(0.5){\Mof(f)} \\
& & & & \\
\Cx(L(H)) \ar@{->}[rr]_(0.5){\Cx(\gamma_{H})} & & \Cx(\wtH) \ar@{->}[rr]_(0.5){a_{\wtH}} & & \Mof(H).}
\end{align*}
The left-hand square commutes, since it is obtained by applying the functor $\Cx$ to the commutative diagram 
of Proposition \ref{gammanatG} which describes the naturality of $\gamma_{G}$ in $G$, and the right-hand square is commutative by Lemma \ref{augqi} Part (ii). Thus the whole diagram commutes. Note that this is implicit in \cite[\S 3]{SS}. 
\end{prooff}

\section{The morphism $\vp_{G}: M(G) \ra \Sym(\Mof(G))$}

\begin{notation} \label{deltasymnot}
Let $G$ be a homotopy invariant commutative group scheme over $k$. We denote by $\Delta_{G}^{n}: G \ra G^{n}$ the $n$-fold diagonal morphism, and recall the symbols: 
\begin{align*}
& \Sym(\Mof(G)) := \bos_{n \ge 0} \Symn(\Mof(G)) \qaq \\
& \pi_{G}^{n}: \Mof(G)^{\ox n} \lra \Symn(\Mof(G)) 
\end{align*}
of Notation \ref{symnot} and Definition \ref{symaltn}, respectively. We recall further from Notation \ref{tensym} that the zeroth tensor and symmetric powers of all objects and morphisms are equal to $\one$ and $1_{\one}$, respectively, so that in particular:
\begin{align*}
M(G)^{\ox 0} = \Mof(G)^{\ox 0} = \Sym^{0}(\Mof(G)).
\end{align*}
Consistent with this, we adopt the new conventions that $G^{0} = \Speck$ and that $\Delta^{0}_{G}$ is equal to the structure morphism $\pi_{G}: G \ra \Speck$.
\end{notation}

\begin{notation} \label{vpGnnot}
Let $G$ be a homotopy invariant commutative group scheme over a perfect field $k$ and $n$ a non-negative integer. We define $\vp_{G}^{n}: M(G) \ra \Symn(\Mof(G))$ to be the composite morphism:
\begin{align*}
\vp_{G}^{n}: M(G) \sxra{M(\Delta_{G}^{n})} M(G)^{\ox n} \sxra{\alpha_{G}^{\ox n}} \Mof(G)^{\ox n}
\sxra{\pi_{G}^{n}} \Symn(\Mof(G))
\end{align*}
in the category $\DMeetkQ$. 
\end{notation}

\begin{remark} \label{zeroonevp}
From the conventions on zeroth powers in Notation \ref{deltasymnot}, it follows that $\vp_{G}^{0}$ is equal to $M(\pi_{G}): M(G) \ra \one$, where $\pi_{G}: G \ra \Speck$ denotes the structure morphism. Moreover, since $\Delta^{1}_{G}: G \ra G$ is equal to the identity $1_{G}$, and $\Sym(\Mof(G)) = \Mof(G)$ with 
$\pi_{G}^{1} = 1_{\Mof(G)}$, it follows that $\vp_{G}^{1}$ is equal to $\alpha_{G}: M(G) \ra \Mof(G)$.
\end{remark}

\begin{proposition} \label{vpGfactorsnicely}
Let $G$ be a homotopy invariant commutative group scheme over a perfect field $k$, and let $n$ be a non-negative integer. 
Then the diagram\tn{:}
\begin{align*} 
\xymatrix@C=20pt{   
M(G) \ar@{->}[dd]_(0.5){M(\Delta^{n}_{G})} \ar@{->}[rr]^(0.42){\vp^{n}_{G}} & & 
\Symn(\Mof(G)) \ar@{->}[dd]^(0.5){\iota_{G}^{n}} \\
& & \\
M(G)^{\ox n} \ar@{->}[rr]_(0.5){\alpha_{G}^{\ox n}} & & \Mof(G)^{\ox n}} 
\end{align*} 
commutes in the category $\DMeetkQ$.
\end{proposition}
\begin{prooff}
For any permutation $\sigma \in \Sym(n)$, the diagonal morphism $\Delta^{n}_{G}: G \ra G^{n}$ in the category of $k$-schemes satisfies 
$\sigma_{G^{n}} \circ \Delta^{n}_{G} = \Delta^{n}_{G}$. Since $\sigma_{M(G)} = M(\sigma_{G})$, it follows that
$\sigma_{M(G)^{\ox n}} \circ M(\Delta^{n}_{G}) = M(\Delta^{n}_{G})$ in the category $\DMeetkQ$. Since moreover clearly 
$\sigma_{\Mof(G)^{\ox n}} \circ \alpha_{G}^{\ox n} = \alpha_{G}^{\ox n} \circ \sigma_{M(G)^{\ox n}}$, we have now: 
\begin{align*} 
\sigma_{\Mof(G)^{\ox n}} \circ \alpha_{G}^{\ox n} \circ M(\Delta^{n}_{G}) & = 
\alpha_{G}^{\ox n} \circ \sigma_{M(G)^{\ox n}} \circ M(\Delta^{n}_{G}) \\
& = \alpha_{G}^{\ox n} \circ M(\Delta^{n}_{G}).
\end{align*} 
By the universal property of the morphism $\iota_{G}^{n}$ in Proposition \ref{symunprop}, this identity implies the existence of a unique morphism $\vsig^{n}_{G}: M(G) \ra \Symn(\Mof(G))$ such that the diagram:
\begin{align*} 
\xymatrix@C=20pt{   
M(G) \ar@{->}[dd]_(0.5){M(\Delta^{n}_{G})} \ar@{->}[rr]^(0.42){\vsig^{n}_{G}} & & 
\Symn(\Mof(G)) \ar@{->}[dd]^(0.5){\iota_{G}^{n}} \\
& & \\
M(G)^{\ox n} \ar@{->}[rr]_(0.5){\alpha_{G}^{\ox n}} & & \Mof(G)^{\ox n}} 
\end{align*} 
commutes, giving the relation $\iota_{G}^{n} \circ \vsig^{n}_{G} = \alpha_{G}^{\ox n} \circ M(\Delta^{n}_{G})$. In view of the identity $\pi_{G}^{n} \circ \iota_{G}^{n} = 1_{\Symn(\Mof(G))}$, applying $\pi_{G}^{n}$ to both sides of this relation yields:
\begin{align*}
\vsig^{n}_{G} = \pi_{G}^{n} \circ \alpha_{G}^{\ox n} \circ M(\Delta^{n}_{G}) = \vp^{n}_{G}.
\end{align*}
This completes the proof.
\end{prooff}

\begin{definition}[\tbf{Kimura}] \label{finitekimdim}
Suppose that $G$ is a commutative group scheme over a perfect field $k$. We say that $G$ is \emph{of finite Kimura dimension} if the $n^{\tn{th}}$ symmetric power $\Symn(\Mof(G))$ vanishes for almost all $n$.
\end{definition}

\begin{notation}  \label{vpGnot}
Suppose that $G$ is a homotopy invariant commutative group scheme over a perfect field $k$ of finite Kimura dimension. For all non-negative integers $n$ write $i_{G}^{n}: \Symn(\Mof(G)) \lra \Sym(\Mof(G))$ for the canonical inclusion. We denote by:
\begin{align*}
\vp_{G}: M(G) \lra \Sym(\Mof(G))
\end{align*}
the sum of morphisms $M(G) \sxra{\vp_{G}^{n}} \Symn(\Mof(G)) \sxra{i_{G}^{n}} \Sym(\Mof(G))$ over all $n \ge 0$ 
in the additive category $\DMeetkQ$. Note that the assumption of finite Kimura dimensionality is necessary for the sum to be well-defined.
\end{notation}

\begin{proposition} \label{vpnfunc}
The morphisms $\vp_{G}^{n}$ and $\vp_{G}$ are both functorial, in the sense that for any morphism 
$f: G \ra H$ of homotopy invariant commutative group schemes of finite Kimura dimension over $k$, the diagrams\tn{:}
\begin{align*}
\tn{(a)} \quad \xymatrix{
M(G) \ar@{->}[rr]^(0.4){\vp_{G}^{n}} \ar@{->}[dd]_(0.5){M(f)} & & 
\Symn(\Mof(G)) \ar@{->}[dd]_(0.5){\Symn(\Mof(f))} \\
& & \\
M(H) \ar@{->}[rr]_(0.4){\vp_{H}^{n}} & & \Symn(\Mof(H))} \quad \tn{and} \quad 
\tn{(b)} \quad \xymatrix{
M(G) \ar@{->}[rr]^(0.4){\vp_{G}} \ar@{->}[dd]_(0.5){M(f)} & & 
\Sym(\Mof(G)) \ar@{->}[dd]_(0.5){\Sym(\Mof(f))} \\
& & \\
M(H) \ar@{->}[rr]_(0.4){\vp_{H}} & & \Sym(\Mof(H))}
\end{align*} 
commute.
\end{proposition}
\begin{prooff}
To prove that Diagram (a) commutes, we note that by the definition of $\vp_{G}^{n}$, it may be written in expanded form as:
\begin{align*}
\xymatrix{
M(G) \ar@{->}[rr]^(0.45){M(\Delta_{G}^{n})} \ar@{->}[dd]_(0.5){M(f)} & & 
M(G)^{\ox n} \ar@{->}[rr]^(0.5){\alpha_{G}^{\ox n}} \ar@{->}[dd]_(0.5){M(f)^{\ox n}} & & 
\Mof(G)^{\ox n} \ar@{->}[rr]^(0.45){\pi_{G}^{n}} \ar@{->}[dd]_(0.5){\Mof(f)^{\ox n}} & & 
\Symn(\Mof(G)) \ar@{->}[dd]^(0.5){\Symn(\Mof(f))} \\
& & & & & & \\
M(H) \ar@{->}[rr]_(0.45){M(\Delta_{H}^{n})} & & 
M(H)^{\ox n} \ar@{->}[rr]_(0.5){\alpha_{H}^{\ox n}} & & 
\Mof(H)^{\ox n} \ar@{->}[rr]_(0.45){\pi_{H}^{n}}  & & 
\Symn(\Mof(H)).}
\end{align*}
The left-hand square commutes by the fact that the diagonal morphism satisfies a universal property given by the fibre product of schemes, and the right-hand square commutes by the universal property of the $n^{\tn{th}}$ symmetric power. The middle square commutes because it is the $n^{\tn{th}}$ tensor power of the diagram in Proposition \ref{alphafunc} describing the functoriality of $\alpha_{G}$. Therefore the whole diagram commutes. 

The commutativity of Diagram (b) now follows from the commutativity of Diagram (a), as follows. Firstly, observe that the diagram:
\begin{align*}
\xymatrix{
M(G) \ar@{->}[rr]^(0.4){\vp_{G}^{n}} \ar@{->}[dd]_(0.5){M(f)} & & 
\Symn(\Mof(G)) \ar@{->}[rr]^(0.5){i_{G}^{n}} \ar@{->}[dd]^(0.5){\Symn(\Mof(f))} & & 
\Sym(\Mof(G)) \ar@{->}[dd]^(0.5){\Sym(\Mof(f))} \\
& & \\
M(H) \ar@{->}[rr]_(0.4){\vp_{H}^{n}} & & \Symn(\Mof(H)) \ar@{->}[rr]_(0.5){i_{H}^{n}} & & \Sym(\Mof(H))} 
\end{align*}
commutes. Indeed, the square on the left is Diagram (a), which we have just shown to be commutative. Further, since $G$ has finite Kimura dimension, the direct sum $\Sym(\Mof(G))$ is finite, and is hence also a categorical product. It now follows directly from the universal property of the product that the square on the left commutes. We have therefore:
\begin{align*}
\Sym(\Mof(f)) \circ \left( i^{n}_{G} \circ \vp^{n}_{G} \right) = 
\left( i^{n}_{H} \circ \vp^{n}_{H} \right) \circ M(f).
\end{align*}
This relation implies that the equality (*) holds in the calculation:
\begin{align*}
\Sym(\Mof(f)) \circ \vp_{G} & =
\Sym(\Mof(f)) \circ \left( \sum_{n \ge 0} i^{n}_{G} \circ \vp^{n}_{G} \right) \\
& = \sum_{n \ge 0} \Bigg(  \Sym(\Mof(f)) \circ \left( i^{n}_{G} \circ \vp^{n}_{G} \right) \Bigg) \\
& \lgeqt{(*)} \sum_{n \ge 0} \Bigg( \left( i^{n}_{H} \circ \vp^{n}_{H} \right) \circ M(f) \Bigg) \\
& = \left( \sum_{n \ge 0} \left( i^{n}_{H} \circ \vp^{n}_{H} \right) \right) \circ M(f) = \vp_{H} \circ M(f)\tn{;}
\end{align*}
all other equalities hold either by definition of $\vp_{G}$ or by bilinearity of the composition of morphisms in the additive category $\DMeetkQ$. Note that the sums in the above computation are finite, since $G$ and $H$ have finite Kimura dimension. We conclude that the Diagram (b) is commutative. 
\end{prooff}

\section{Multiplicativity of $\vp_{G}$}

\begin{remark}
Let us recall some notation from Chapter $1$. Let $m$ and $n$ be non-negative integers with sum $N$, let $\Sym(N)$ denote the symmetric group on $N$ letters, and let $\sigma \in \Sym(N)$. Let further $X, Y$ and $Z$ be objects in a 
$\Q$-linear idempotent complete tensor category $\catc$. We have the following morphisms.  
\begin{itemize}
\item[\tn{(i)}] The map $\sigma_{Z}: Z^{\ox n} \ra Z^{\ox n} $ of Notation \ref{sigmanot} permuting factors by $\sigma$.
\item[\tn{(ii)}] The map
$s^{n}_{Z} := \frac{1}{n!} \left( \sum_{\sigma \in \Sym(n)} \sigma_{Z} \right): Z^{\ox n} \ra Z^{\ox n}$ of Notation \ref{ansn}.
\item[\tn{(iii)}] The morphisms: 
\begin{align*}
& \xi_{X, Y}: \Sym(X \os Y) \ra \Sym(X) \ox \Sym(Y), \\
& \xi_{X, Y}^{m,n}: \Sym^{N}(X \os Y) \ra \Symm(X) \ox \Symn(Y) \qaq \\
& \xi_{X, Y}^{N}: \Sym^{N}(X \os Y) \lra \bos_{m+n = N} \Sym^{m}(X) \ox \Sym^{n}(Y)
\end{align*}
of Notation \ref{xinot}.
\end{itemize}
For us, $\catc$ will be the category $\DMeetkQ$.
\end{remark}

\begin{remark} \label{Qfactcat}
Since $\DMeetkQ$ is by construction a $\Q$-linear category, the functor $M: \Smck \ra \DMeetkQ$ of
Notation \ref{Mhinot} factors through $\Smck \ox \Q$. 
\end{remark}

\begin{lemma} \label{SmnMGH}
Let $X$ be a smooth scheme of finite type over a perfect field $k$. Then for any non-negative integer $n$, the morphisms\tn{:}
\begin{align*}
& s^{n}_{M(X)}: M(X)^{\ox n} \lra M(X)^{\ox n} \qaq \\
& M(s^{n}_{X}): M(X)^{\ox n} \lra M(X)^{\ox n}
\end{align*}
in the category $\DMeetkQ$ agree up to the canonical isomorphism of $M(X^{n})$ and $M(X)^{\ox n}$. Here, the morphism $s^{n}_{X}: X^{n} \ra X^{n}$ is regarded as belonging to the $\Q$-linear category $\Smck \ox \Q$, following Remark \ref{Qfactcat}. 
\end{lemma}
\begin{prooff}
If $n = 0$ then both morphisms are equal to $1_{\one}: \one \ra \one$. Suppose then that $n$ is positive, and let
$\sigma$ be a permutation in $\Sym(n)$. Then by definition of the tensor product in the category of motives $\DMeetkQ$, the diagram:
\begin{align*}
\xymatrix{
M(X^{n}) \ar@{->}[rr]^(0.48){M(\sigma_{X})} \ar@{->}[dd]_{\simeq} & & 
M(X^{n}) \ar@{->}[dd]_{\simeq} \\
& & \\
M(X)^{\ox n} \ar@{->}[rr]_(0.48){\sigma_{M(X)}} & & M(X)^{\ox n}}
\end{align*}
commutes. Summing this diagram over all $\sigma \in \Sym(n)$ and multiplying by $(1/n!)$ thus gives a commutative diagram:  
\begin{align*}
\xymatrix{
M(X^{n}) \ar@{->}[rr]^(0.48){\sum_{\sigma} M(\sigma_{X})} 
\ar@{->}[dd]_{\simeq} & & M(X^{n}) \ar@{->}[dd]_{\simeq} \\
& & \\
M(X)^{\ox n} \ar@{->}[rr]_(0.48){\sum_{\sigma} \sigma_{M(X)}} & & M(X)^{\ox n}.}
\end{align*}
Now by definition, $s^{n}_{M(X)} = (1/n!)\sum_{\sigma} \sigma_{M(X)}$, and so the above diagram implies that 
$s^{n}_{M(X)}$ and $(1/n!) \sum_{\sigma} M(\sigma_{X})$ agree up to the canonical isomorphism of $M(X^{n})$ and $M(X)^{\ox n}$. But since: 
\begin{align*}
M: \Smck \ox \Q \lra \DMeetkQ
\end{align*}  
is a $\Q$-linear functor, it follows that:
\begin{align*}
\frac{1}{n!} \sum_{\sigma} M(\sigma_{X}) = M \left( \frac{1}{n!} \sum_{\sigma} \sigma_{X} \right) = M(s^{n}_{X}),
\end{align*}  
completing the proof.
\end{prooff}

\begin{proposition} \label{vpmultprop}
Let $G$ and $H$ be homotopy invariant commutative group schemes of finite Kimura dimension over a perfect field $k$. Then\tn{:} 
\begin{itemize}
\item[\tn{(i)}] The diagram\tn{:}
\begin{align*}
\xymatrix{
M(G \x H) \ar@{->}[dd]^(0.5){\simeq} \ar@{->}[rr]^(0.45){\vp_{G \x H}} & & 
\Sym(\Mof(G \x H)) \ar@{->}[rr]^(0.47){\Sym(\eta_{G, H})} & & 
\Sym(\Mof(G) \os \Mof(H)) \ar@{->}[dd]^(0.5){\xi_{G, H}} \\
& & & & \\
M(G) \ox M(H) \ar@{->}[rrrr]_(0.45){\vp_{G} \ox \vp_{H}} & & & & \Sym(\Mof(G)) \ox \Sym(\Mof(H))} 
\end{align*}
commutes. 
\item[\tn{(ii)}] For any positive integer $r$, the map $\vp_{G^{r}}$ is an isomorphism if and only if $\vp_{G}$ is an isomorphism.
\end{itemize}
\end{proposition}
\begin{prooff}
We begin by observing that since both $\xi_{G, H}$ and $\eta_{G, H}$ are isomorphisms, the commutativity of the diagram in statement (i) immediately implies statement (ii).\\\\
To prove (i), it suffices to prove that the diagram: 
\begin{align*}
\xymatrix{
M(G \x H) \ar@{->}[dd]^(0.5){\simeq} \ar@{->}[rr]^(0.4){\vp_{G \x H}^{N}} & & 
\Sym^{N}(\Mof(G \x H)) \ar@{->}[rr]^(0.5){\Sym^{N}(\eta_{G, H})} & & 
\Sym^{N}(\Mof(G) \os \Mof(H)) \ar@{->}[dd]^(0.5){\xi_{G, H}^{m,n}} \\
& & & & \\
M(G) \ox M(H) \ar@{->}[rrrr]_(0.4){\vp_{G}^{m} \ox \vp_{H}^{n}} & & & & \Symm(\Mof(G)) \ox \Symn(\Mof(H))} 
\end{align*}
commutes, for all non-negative integers $m$ and $n$ with sum $N$. This diagram may be written in expanded form by glueing, from left to right, the four subdiagrams:
\begin{align*}
\tn{(a)} \qquad \quad \xymatrix{
M(G \x H) \ar@{->}[rrr]^(0.45){M(\Delta_{G \x H}^{N})} 
\ar@{->}[dd]^(0.5){\simeq} & & & 
M(G \x H)^{\ox N} \ar@{->}[d]^{s^{N}_{M(G \x H)}} \\
& & & M(G \x H)^{\ox N} \ar@{->}[d]^(0.5){M(p_{G})^{\ox m} \ox M(p_{H})^{\ox n}} \\
M(G) \ox M(H) \ar@{->}[rrr]_(0.45){M(\Delta_{G}^{m}) \ox M(\Delta_{H}^{n})} & & & 
M(G)^{\ox m} \ox M(H)^{\ox n},} 
\end{align*}

\begin{samepage}

\begin{align*}
\tn{(b)} \quad \xymatrix{
M(G \x H)^{\ox N} \ar@{->}[rr]^(0.5){\alpha_{G \x H}^{\ox N}} \ar@{->}[dd]_{s^{N}_{M(G \x H)}} & & 
\Mof(G \x H)^{\ox N} \ar@{->}[dd]^{s^{N}_{\Mof(G \x H)}} \\
& & \\
M(G \x H)^{\ox N} \ar@{->}[rr]^(0.5){\alpha_{G \x H}^{\ox N}} 
\ar@{->}[dd]_(0.5){M(p_{G})^{\ox m} \ox M(p_{H})^{\ox n}} & & \Mof(G \x H)^{\ox N} 
\ar@{->}[dd]^(0.5){\Mof(p_{G})^{\ox m} \ox \Mof(p_{H})^{\ox n}} \\
& & \\
M(G)^{\ox m} \ox M(H)^{\ox n} \ar@{->}[rr]_(0.47){\alpha_{G}^{\ox m} \ox \alpha_{H}^{\ox n}} & & 
\Mof(G)^{\ox m} \ox \Mof(H)^{\ox n},} \\
\end{align*}

\begin{align*}
\tn{(c)} \quad \xymatrix{
\Mof(G \x H)^{\ox N} \ar@{->}[dd]_{s^{N}_{\Mof(G \x H)}} \ar@{->}[rrr]^(0.45){\pi_{G \x H}^{N}} & & & 
\Sym^{N} \big( \Mof(G \x H) \big) \ar@{->}[dd]^{\xi^{m,n}_{\Mof(G \x H)}} \\
& & & \\
\txt{$\Mof(G \x H)^{\ox m} \ox$ \\ $\Mof(G \x H)^{\ox n}$}  
\ar@{->}[dd]_(0.5){p_{\Mof(G)}^{\ox m} \ox p_{\Mof(H)}^{\ox n}} 
\ar@{->}[rrr]_(0.45){\pi_{\Mof(G \x H)}^{m} \ox \pi_{\Mof(G \x H)}^{n}} & & & 
\txt{$\Sym^{m} \big( \Mof(G \x H) \big) \ox$ \\ $\Sym^{n} \big( \Mof(G \x H) \big)$}
\ar@{->}[dd]^(0.5){\txt{\scs{$\Symm(\Mof(p_{G})) \ox$} \\ \scs{$\Symn(\Mof(p_{H}))$}}} \\
& & & \\
\Mof(G)^{\ox m} \ox \Mof(H)^{\ox n} \ar@{->}[rrr]_(0.45){\pi_{G}^{m} \ox \pi_{H}^{n}} & & & \Symm(\Mof(G)) \ox \Symn(\Mof(H)),} 
\end{align*}
and
\begin{align*}
\tn{(d)} \: \xymatrix @C=5pt {
\Sym^{N} \big( \Mof(G \x H) \big) \ar@{->}[dd]_{\xi^{m,n}_{\Mof(G \x H)}} \ar@{->}[rr]^{\Sym^{N}(\eta_{G, H})} & & \Sym^{N} \big( \Mof(G) \os \Mof(H) \big) \ar@{->}[dd]^{\xi^{m,n}_{\Mof(G) \os \Mof(H)}} 
\ar@/_6.5pc/@{->}[dddd]_(0.3){\xi_{G, H}^{m,n}} \\
& & \\
\txt{$\Sym^{m} \big( \Mof(G \x H) \big) \ox$ \\ $\Sym^{n} \big( \Mof(G \x H) \big)$}
\ar@{->}[dd]_(0.5){\txt{\scs{$\Symm(\Mof(p_{G})) \ox$} \\ \scs{$\Symn(\Mof(p_{H}))$}}} & & 
\txt{$\Sym^{m} \big( \Mof(G) \os \Mof(H) \big) \ox$ \\ $\Sym^{n} \big( \Mof(G) \os \Mof(H) \big)$}
\ar@{->}[dd]^(0.5){\txt{\scs{$\Symm(p_{\Mof(G)}) \ox$} \\ \scs{$\Symn(p_{\Mof(H)})$}}} \\
& & \\
\Symm(\Mof(G)) \ox \Symn(\Mof(H)) \ar@{=}[rr] & & \Symm(\Mof(G)) \ox \Symn(\Mof(H)).}
\end{align*}

\end{samepage}
We check in turn that each diagram commutes.\\\\
\tbf{(a)} By Lemma \ref{SmnMGH}, the morphisms $s^{N}_{M(G)}$ and $M(s^{N}_{G})$ agree up to canonical isomorphism, and so the first subdiagram is just the application of the functor $M$ to the diagram of $k$-schemes:
\begin{align*}
\xymatrix{
G \x H \ar@{->}[rr]^(0.45){\Delta_{G \x H}^{N}} \ar@{->}[ddrr]_(0.45){\Delta_{G}^{m} \x \Delta_{H}^{n}} & &
(G \x H)^{N} \ar@{->}[d]^{s^{N}_{G \x H}} \\
& & (G \x H)^{N} \ar@{->}[d]^(0.5){p_{G}^{m} \x p_{H}^{n}} \\
& & G^{m} \x H^{n}}
\end{align*}
considered in the category $\Smck \ox \Q$, following Remark \ref{Qfactcat}. We check this commutes by considering a $T$-point $(g, h) \in (G \x H)(T)$, where $T$ is a $k$-scheme. Observe first that:
\begin{align*}
s^{N}_{G \x H} \big( \ob{(g,h), \ldots, (g,h)}^{N \tn{ copies}} \big) & = 
\frac{1}{N!} \left( \sum_{\sigma \in \Sym(N)} \big( (g,h), \ldots, (g,h) \big) \right) \\
& = \frac{1}{N!} \left( N! \cdot \big( (g,h), \ldots, (g,h) \big) \right) \\
& = \big( (g,h), \ldots, (g,h) \big).
\end{align*}
Let us emphasize here that the summation and multiplication with $N!$ and $1/N!$ take place in $\Smck \ox \Q$, and have nothing to do with the commutative group scheme laws on $G$ and $H$. From the above calculation, it follows that:
\begin{align*}
\Big( (p_{G}^{m} \x p_{H}^{n}) \circ s^{N}_{G \x H} \circ \Delta_{G \x H}^{N} \Big) \big( (g,h) \big) & = 
\Big( (p_{G}^{m} \x p_{H}^{n}) \circ s^{N}_{G \x H} \Big) \big( \ob{(g,h), \ldots, (g,h)}^{N \tn{ copies}} \big) \\
& = \left( p_{G}^{m} \x p_{H}^{n} \right) \big( (g,h), \ldots, (g,h) \big) \\
& = \big( \ub{g, \ldots, g}_{m \tn{ copies}}, \ub{h, \ldots, h}_{n \tn{ copies}} \big) \\
& = \big( \Delta_{G}^{m} \x \Delta_{H}^{n} \big)\big( (g,h) \big).
\end{align*}
Thus the first square commutes.\\\\
\tbf{(b)} The top square: 
\begin{align*}
\xymatrix{
M(G \x H)^{\ox N} \ar@{->}[rr]^(0.5){\alpha_{G \x H}^{\ox N}} \ar@{->}[dd]_{s^{N}_{M(G \x H)}} & & 
\Mof(G \x H)^{\ox N} \ar@{->}[dd]^{s^{N}_{\Mof(G \x H)}} \\
& & \\
M(G \x H)^{\ox N} \ar@{->}[rr]^(0.5){\alpha_{G \x H}^{\ox N}} & & \Mof(G \x H)^{\ox N}}
\end{align*}
of the second diagram commutes because the morphisms $s^{N}_{X}$ are each a sum of permutations of the factors of the respective products, and one has more generally that:
\begin{align*}
\xymatrix{
X^{\ox N} \ar@{->}[rr]^(0.5){f^{\ox N}} \ar@{->}[dd]_{\sigma_{X}} & & 
Y^{\ox N} \ar@{->}[dd]^{\sigma_{Y}} \\
& & \\
X^{\ox N} \ar@{->}[rr]_(0.5){f^{\ox N}} & & Y^{\ox N}}
\end{align*}
commutes for any objects $X$ and $Y$, any morphism $f: X \ra Y$ in an additive $\Q$-linear idempotent complete tensor category $\catc$, and any permutation $\sigma \in \Sym(n)$. Since 
$M(G \x H)^{\ox N} \simeq M(G \x H)^{\ox m} \ox M(G \x H)^{\ox n}$, the bottom square may be rewritten as:
\begin{align*}
\xymatrix{
M(G \x H)^{\ox m} \ox M(G \x H)^{\ox n} \ar@{->}[rr]^(0.5){\alpha_{G \x H}^{\ox N}} 
\ar@{->}[dd]_(0.5){M(p_{G})^{\ox m} \ox M(p_{H})^{\ox n}} & & \Mof(G \x H)^{\ox m} \ox \Mof(G \x H)^{\ox n}
\ar@{->}[dd]^(0.5){\Mof(p_{G})^{\ox m} \ox \Mof(p_{H})^{\ox n}} \\
& & \\
M(G)^{\ox m} \ox M(H)^{\ox n} \ar@{->}[rr]_(0.47){\alpha_{G}^{\ox m} \ox \alpha_{H}^{\ox n}} & & 
\Mof(G)^{\ox m} \ox \Mof(H)^{\ox n}}
\end{align*}
But this is the tensor product of the two squares:
\begin{align*}
\xymatrix{
M(G \x H)^{\ox m} \ar@{->}[rr]^(0.5){\alpha_{G \x H}^{\ox m}} \ar@{->}[dd]_(0.5){M(p_{X})^{\ox m}} & & 
\Mof(G \x H)^{\ox m}\ar@{->}[dd]^(0.5){\Mof(p_{X})^{\ox m}} \\
& & \\
M(X)^{\ox m} \ar@{->}[rr]_(0.47){\alpha_{G}^{\ox m}} & & \Mof(X)^{\ox m}} 
\end{align*}
where $X$ is equal to $G$ or $H$, both of which commute by the functoriality of the morphism $\alpha_{?}$ described in Proposition \ref{alphafunc}. \\\\
\tbf{(c)} The top square of the third diagram commutes by definition of the morphism $\xi^{m,n}_{\Mof(G \x H)}$ (see Notation \ref{xinot}). The bottom square is the product of the two diagrams: 
\begin{align*}
\xymatrix{
\Mof(G \x H)^{\ox m} \ar@{->}[dd]_(0.5){p_{\Mof(G)}^{\ox m}} \ar@{->}[rr]^(0.45){\pi_{\Mof(G \x H)}^{m} } & &  
\Sym^{m} \big( \Mof(G \x H) \big) \ar@{->}[dd]^(0.5){\Symm(\Mof(p_{G}))} \\
& & \\
\Mof(G)^{\ox m} \ar@{->}[rr]_(0.45){\pi_{G}^{m}} & & \Symm(\Mof(G))} 
\end{align*}
and
\begin{align*}
\xymatrix{
\Mof(G \x H)^{\ox n} \ar@{->}[dd]_(0.5){p_{\Mof(H)}^{\ox n}} \ar@{->}[rr]^(0.45){\pi_{\Mof(G \x H)}^{n} } & &  
\Symn \big( \Mof(G \x H) \big) \ar@{->}[dd]^(0.5){\Symn(\Mof(p_{H}))} \\
& & \\
\Mof(H)^{\ox n} \ar@{->}[rr]_(0.45){\pi_{H}^{n}} & & \Symn(\Mof(H)),}
\end{align*}
both of which commute by Proposition \ref{symnfunc}.\\\\
\tbf{(d)} Firstly, the equality of morphisms:
\begin{align*}
\xi_{G, H}^{m,n} = \big( \Symm(p_{\Mof(G)}) \ox \Symn(p_{\Mof(H)}) \big) \circ \xi^{m,n}_{\Mof(G) \os \Mof(H)}
\end{align*}
holds by definition of the map $\xi^{m,n}_{G, H}$ in Notation \ref{xinot}. The top square is:
\begin{align*}
\xymatrix @C=5.5pc {
\Sym^{N} \big( \Mof(G \x H) \big) \ar@{->}[dd]_{\xi^{m,n}_{\Mof(G \x H)}} \ar@{->}[rr]^{\Sym^{N}(\eta_{G, H})} & & 
\Sym^{N} \big( \Mof(G) \os \Mof(H) \big) \ar@{->}[dd]_{\xi^{m,n}_{\Mof(G) \os \Mof(H)}} \\
& & \\
\txt{$\Sym^{m} \big( \Mof(G \x H) \big) \ox$ \\ $\Sym^{n} \big( \Mof(G \x H) \big)$} 
\ar@{->}[rr]_{\Sym^{m}(\eta_{G, H}) \ox \Sym^{n}(\eta_{G, H})} & & 
\txt{$\Sym^{m} \big( \Mof(G) \os \Mof(H) \big) \ox$ \\ $\Sym^{n} \big( \Mof(G) \os \Mof(H) \big)$,}}
\end{align*}
which is commutative because $\eta_{G, H}$ is a canonical isomorphism, and the bottom square is a product of the two diagrams:
\begin{align*}
\xymatrix{
\Sym^{m} \big( \Mof(G \x H) \big) \ar@{->}[rr]^(0.45){\Symm(\eta_{G, H})} \ar@{->}[dd]_(0.5){\Symm(\Mof(p_{G}))} & & 
\Sym^{m} \big( \Mof(G) \os \Mof(H) \big) \ar@{->}[dd]^(0.5){\Symm(p_{\Mof(G)})} \\
& & \\
\Symm(\Mof(G))\ar@{=}[rr] & & \Symm(\Mof(G))} 
\end{align*}
and
\begin{align*}
\xymatrix{
\Symn \big( \Mof(G \x H) \big) \ar@{->}[rr]^(0.45){\Symn(\eta_{G, H})} \ar@{->}[dd]_(0.5){\Symn(\Mof(p_{H}))} & & 
\Symn \big( \Mof(G) \os \Mof(H) \big) \ar@{->}[dd]^(0.5){\Symn(p_{\Mof(H)})} \\
& & \\
\Symn(\Mof(H))\ar@{=}[rr] & & \Symn(\Mof(H)),}
\end{align*}
which commute by Remark \ref{projGHcomm}. This completes the proof.
\end{prooff}
\chapter{The Albanese scheme}

Let $k$ be a perfect field and $X$ a scheme over $k$. In this chapter we will review the definitions and main properties of the {\SeAV} $\AXz$ and Albanese scheme $\AlbXk$ of $X$. We emphasize that none of the material in this chapter is due to the author. These concepts will be used in the next chapter to show that the morphism:
\begin{align*}
\vp_{A}: M(A) \lra \Sym(\Mof(A))
\end{align*}
of Notation \ref{vpGnot} an isomorphism, when $A$ is an abelian variety over $k$, by first proving the special case 
$A = \AlbzCk$, where $C$ is a smooth proper curve. We now outline the contents of each section. 

\tbf{(i)} In the first section of this chapter we give, in a general categorical context, a definition of a universal morphism. The universal properties defining both the {\SeAV} and the Albanese scheme are examples of universal morphisms, so the terminology is a convenient generalisation. 

\tbf{(ii) \& (iii)} In the second and third sections, we recall the definitions of a semiabelian variety and of a scheme-theoretic torsor, respectively.  

\tbf{(iv)} In the fourth section, we review the concept of the \emph{\SeAV} $\AXz$ of $X$. The {\SeAV} was first defined by Serre in \cite{Se1} to be a semiabelian variety satisfying a certain universal property. It was shown to exist when $X$ satisfies the following condition: 
\begin{condition} \label{existcond}
The scheme $X$ is reduced, separated and of finite type over the perfect field $k$. 
\end{condition}

\tbf{(v)} In the fifth section we review the definition of the \emph{Albanese scheme $\AlbXk$} of $X$ due to Ramachandran \cite[\S 1]{Ra}, based on a suggestion of Serre. This definition is given in terms of a second universal property, and hence is logically independent from that of the {\SeAV}. 

\tbf{(vi)} In the sixth section we construct, for a connected $k$-scheme $X$ such that $\AlbXk$ exists, a natural morphism:
\begin{align*}
p_{\Zcs}: \wAlbXk \lra \Zcs
\end{align*}
of commutative group schemes over $k$, where $\Zcs$ is the constant group scheme associated to the abelian group $\Z$. We then show that the neutral component $\AlbzXk$ of $\AlbXk$ is equal to the kernel of $p_{\Zcs}$, which allows us, to write down a canonical isomorphism between the neutral component $\AlbzXk$ of $\AlbXk$ and the Serre-Albanese variety $\AXz$. This isomorphism exists provided $\AlbXk$ exists and $X$ is connected. 

\tbf{(vii)} The seventh section lists those properties of the Albanese scheme which we will use in subsequent chapters. In particular, we will recall in Proposition \ref{hzalbC} the classical isomorphism:
\begin{align*}
\bC: \hzet(M(C)) \lra \wAlbCk
\end{align*}
proven in different forms by Lichtenbaum \cite{Li1} and Suslin and Voevodsky \cite[Theorem 3.1]{SV}. The isomorphism $\bC$ will be used in the next chapter to help prove that $\vp_{A}: M(A) \ra \Sym(\Mof(A))$ is an isomorphism. 

\tbf{(viii)} As a by-product of the construction of $\bC$, as well as the other properties of the Albanese scheme listed in the seventh section, we prove in the eight section that 
\begin{align*}
\vp_{T}: M(T) \ra \Sym(\Mof(T))
\end{align*}
is an isomorphism for any torus $T$ over a perfect field $k$ (Proposition \ref{vptorusisom}).

\section{Universal morphisms}

The universal properties defining both $\AXz$ and $\AlbXk$ are both examples of the following abstract definition, which will allow us to streamline the exposition. 

\begin{definition}
Suppose that $\catc$ is a category, $\catd$ is a subcategory of $\catc$ and $X$ is an object of $\catc$. Then a morphism $u_{X}: X \ra Y$, with $Y$ an object of $\catd$, is said to be a \emph{universal morphism from $X$ into the category $\catd$} if, given any morphism 
$\vp: X \ra Z$ in $\catc$ with $Z$ an object of $\catd$, there exists a unique morphism $\psi: Y \ra Z$ in $\catd$ such that the diagram: 
\begin{align*}
\xymatrix{
X \ar@{->}[dr]_(0.45){\vp} \ar@{->}[rr]^(0.45){u_{X}} & & 
Y \ar@{->}[dl]^(0.45){\psi} \\
& Z &}
\end{align*}
commutes.
\end{definition}

\section{Semiabelian varieties}

\begin{definition} \label{alggpsch}
Let $G$ be a group scheme over $S$. We say that $G$ is an \emph{algebraic group scheme over $S$} if it is of finite type over $S$, and $G$ is a \emph{locally algebraic group scheme over $S$} if it is locally of finite type over $S$.
\end{definition}

\begin{definition} \label{savdef}
We say that an algebraic group scheme $G$ over $S$ is \emph{semiabelian} if it is an extension of an abelian group scheme $A$ by a torus $T$. Thus we have a short exact sequence: 
\begin{align*}
1 \lra \wt{T} \lra \wt{G} \lra \wt{A} \lra 0 
\end{align*}
in the category of sheaves of abelian groups on the big {\etl} site on $\SchS$. Observe that $G$ is an algebraic group scheme. 
\end{definition}

\begin{definition} \label{largegroupdef}
A locally algebraic group scheme $G$ over $S$ is \emph{a semiabelian group scheme} if its neutral component $G^{0}$ is a semiabelian algebraic group scheme. Following Ramachandran \cite[Definition 1.2]{Ra}, we refer to a semiabelian locally algebraic group scheme over $S$ as \emph{large group scheme over $S$}.
\end{definition}

\section{Torsors and pointed schemes}

\begin{definition} \label{torsdef}
Let $X$ be a scheme over $S$ which is faithfully flat and locally of finite type. Let $G$ be a group scheme over $S$ and let $\rho: G \x_{S} X \ra X$ be an action of $G$ on $X$. Then we call $X$ a \emph{$G$-torsor} if the morphism of $S$-schemes $\Psi_{G, X}: G \x_{S} X \ra X \x_{S} X$ of Notation \ref{Psinot} is an isomorphism. 
\end{definition} 

\begin{definition}
A \emph{pointed $S$-scheme} is a pair $(X, \xo)$ consisting of an $S$-scheme $X$ and an $S$-valued point 
$\xo: S \lra X$. A \emph{morphism $(X, \xo) \ra (Y, \yo)$ of pointed $S$-schemes} is a morphism of $S$ schemes compatible with $\xo$ and $\yo$. The pointed $S$-schemes form a category. Given a pointed $S$-scheme $(X, \xo)$ we define $\xoT \in X(T)$ to be the composite morphism: 
\begin{equation}
T \ilra S \x_{S} T \sxra{\xo \x_{S} T} X \x_{S} T \slra{\pr_{X}} X,
\end{equation}
for any $S$-scheme $T$. 
\end{definition}

\begin{remark}
We will always consider a group scheme $G$ over $S$ to be pointed by the neutral element $e_{G}: S \ra G$ of $G(S)$.
\end{remark}

\begin{remark} \label{pointdefmorph}
Let $(X, \xo)$ be a pointed $S$-scheme, and let $G$ be a group scheme over $S$ acting on $X$. It follows from Definition \ref{torsdef} that the composite morphism: 
\begin{align*}
\vp_{\xo}: G \ilra G \x_{S} S \sxra{1_{G} \x x_{0}} G \x_{S} X \slra{\Psi} X \x_{S} X \slra{\pr_{1}} X
\end{align*}
is an isomorphism of $S$-schemes. 
\end{remark}

\begin{definition} \label{torsmorphdef}
Let $G$ and $H$ be $S$-group schemes and let $X$ and $Y$ are $G$- and $H$-torsors, respectively. A \emph{morphism of torsors $(X, G) \ra (Y, H)$} consists of a morphism $\vp: G \ra H$ of $S$-group schemes, together with a morphism $f: X \ra Y$ of $S$-schemes such that the diagram: 
\begin{align*}
\xymatrix{
G \x_{S} X \ar@{->}[dd]_(0.5){\Psi_{G, X}} \ar@{->}[rr]^(0.5){\vp \x_{S} f} & & 
H \x_{S} Y \ar@{->}[dd]^(0.5){\Psi_{H, Y}}\\
& & \\
X \x_{S} X \ar@{->}[rr]_(0.5){f \x_{S} f} & & Y \x_{S} Y} 
\end{align*}
commutes. 
\end{definition} 

\section{The {\SeAV}} \label{SeAVsec}

\begin{remark}
From now on all schemes will be based over the spectrum of the perfect field $k$.
\end{remark}

\begin{definition} \label{SAVT}
Let $X$ be a scheme over the perfect field $k$. Suppose that there exists an algebraic semiabelian group scheme $A_{X}^{0}$ over $k$ and an $A_{X}^{0}$-torsor $A_{X}^{1}$, together with a universal morphism $\muXo: X \lra A_{X}^{1}$ from $X$ into the category of torsors on semiabelian schemes. Then we call $A_{X}^{0}$ the \emph{{\SeAV}} of $X$, and we call $A_{X}^{1}$ the \emph{Serre-Albanese torsor}.
\end{definition}

\begin{remark}
We make the following convention to ease the notation. Let $X$ be any $k$-scheme. Whenever we say that 
``the {\SeAV} $\AXz$ exists'', we mean ``there exists a triplet $(\AXz, \AXo,\muXo)$ which satisfies Definition \ref{SAVT}''.
\end{remark}

\begin{remark}
Unpacked, the definition means the following. Given any triplet $(G, T, \vp)$ with $G$ an algebraic semiabelian group scheme, $T$ a $G$-torsor and $X \slra{\vp} T$ a morphism of $k$-schemes, there exists a unique morphism
\begin{align*}
\psi = (\psi_{0}, \psi_{1}): (A_{X}^{0}, A_{X}^{1}) \lra (G, T)
\end{align*}
of $G$-torsors such that the diagram:
\begin{align*}
\xymatrix{
X \ar@{->}[rr]^(0.45){\muXo} \ar@{->}[dr]_(0.45){\vp} & & 
A_{X}^{1} \ar@{->}[dl]^(0.45){\psi_{1}} \\
& T &} 
\end{align*}
commutes. 
\end{remark}

\begin{proposition} \label{seavexists}
Let $k$ be a perfect field and $X$ a $k$-scheme satisfying Condition \ref{existcond}. Then the {\SeAV} exists and is unique up to canonical isomorphism.
\end{proposition}
\begin{prooff}
The universal property defining the triplet $(\AXz, \AXo,\muXo)$ guarantees its uniqueness. Its existence is proven in \cite[\S 4, Theorem 5]{Se1}.
\end{prooff}

\begin{proposition} \label{varnotorsor}
Let $(X, \xo)$ be a pointed $k$-scheme such that the {\SeAV} $\AXz$ of $X$ exists, set
$\yo := \muXo \circ \xo: \Speck \ra \AXo$, and let $\vp_{\yo}: \AXz \ra \AXo$ be the morphism of Remark \ref{pointdefmorph}. Then the composite mapping\tn{:}
\begin{align*}
\muXz := \vp_{\yo}^{-1} \circ \muXo: X \lra \AXz
\end{align*}
is a universal morphism from the object $(X, \xo)$ in the category of pointed $k$-schemes into the subcategory of semiabelian varieties.
\end{proposition}
\begin{prooff}
The universality of $\muXz$ is inherited from the universality of $\muXo$.
\end{prooff}

\section{The Albanese scheme} \label{ASsec} 

\begin{remark}
In what follows, all sheaves will be on the big {\etl} site $\ket$ over $\Speck$. We recall 
(Notation \ref{Xshfnot}) that $\wt{X}$ denotes the {\etl} sheaf represented by a $k$-scheme $X$, and that $\cs{\Gamma}$ denotes the constant group scheme over $k$ associated to an abstract group $\Gamma$ (Remark \ref{gpact}).
\end{remark}

\begin{notation}
We write $W \mapsto \Z[W]$ for the functor $\Sets \lra \Ab$ taking a set to the free abelian group generated by its elements.
\end{notation}

\begin{notation} \label{ZXnot} 
For any scheme $X$ over $\Speck$, let us denote by $\ZX$ the sheafification of the presheaf:
\begin{align*}
\mathsmap{\Z[X]}{\ket}{\Ab}{T}{\Z[X(T)].}
\end{align*}
We write $\sigma_{X}: \Z[X] \lra \ZX$ for the associated morphism of presheaves.
\end{notation}

\begin{definition} \label{albschdef}
Let $X$ be a $k$-scheme. Suppose that there exists a large group scheme $\AlbXk$, together with 
a universal morphism:
\begin{align*}
\muX: \ZX \lra \wAlbXk
\end{align*}
in the category $\Sh(\ket)$ of sheaves on $\ket$ from $\ZX$ into the subcategory of $\Sh(\ket)$ consisting of sheaves represented by large group schemes (Definition \ref{largegroupdef}). Then we call $\AlbXk$ the \emph{Albanese scheme of $X$} and $\muX$ the \emph{universal morphism associated with $X$}.
\end{definition}

\begin{remark}
Let $X$ be any $k$-scheme. To ease notation, we adopt the following convention. Whenever we say 
``the Albanese scheme $\AlbXk$ exists'', we mean ``the pair $(\AlbXk, \muX)$ of Definition \ref{albschdef} exists''.
\end{remark}

\begin{proposition} \label{albschemeexists}
Let $X$ be a scheme over $k$ satisfying Condition \ref{existcond}. Then the Albanese scheme exists and is unique up to canonical isomorphism. 
\end{proposition}
\begin{prooff}
The universal property defining the pair $(\AlbXk, \muX)$ guarantees its uniqueness, by the standard argument. See \cite[Theorem 1.11]{Ra} for a proof of the existence of a representable sheaf $\wAlbXk$ satisfying the universal property described in Definition \ref{albschdef}, but for the flat instead of the {\etl} topology; indeed this is Ramachandran's definition of the Albanese scheme. However, by Proposition \ref{Gsheaf}, due to Grothendieck, the presheaf on $\Schk$ associated to the scheme $\AlbXk$ is also a sheaf for the {\etl} topology.
\end{prooff}

\begin{remark} \label{Albfuncrem} 
The universal property defining the Albanese scheme implies that the assignment $X \mapsto \AlbXk$ of objects is in fact a functor from the category of all $k$-schemes satisfying Condition \ref{existcond} to the category of large group schemes. Indeed, a map $f: X \ra Y$ of $k$-schemes induces a map on representable sheaves of sets $\wt{f}: \wt{X} \ra \wt{Y}$ extending $\Z$-linearly to a map of abelian presheaves $\Z[f]: \Z[X] \ra \Z[Y]$. Sheafifying now gives a morphism $\Z^{f}: \ZX \ra \Z^{Y}$, and the universal property of $\muX: \ZX \ra \AlbXk$ then induces a map $\Alb_{f/k}: \AlbXk \ra \AlbXk$ such that the diagram: 
\begin{align*}
\xymatrix{
\ZX \ar@{->}[rr]^(0.45){\muX} \ar@{->}[dd]_{\Z^{f}} & & 
\AlbXk \ar@{->}[dd]^{\Alb_{f/k}} \\
& & \\
\ZY \ar@{->}[rr]_(0.45){\muY} & & \AlbYk}
\end{align*}
commutes. Preservation of composition and identity maps follows from the uniqueness of $\Alb_{f/k}$.
\end{remark}

\section{The connected components of the Albanese scheme}

We recall that for a connected scheme $X$ over $k$ such that $\AlbXk$ exists, there exists a natural morphism of commutative $k$-group schemes:
\begin{align*}
p_{\Zcs}: \wAlbXk \lra \wtZcs,
\end{align*}
the fibres of which are the connected components of $\AlbXk$. In particular, the neutral component $\AlbzXk$ of $\AlbXk$ is equal to the kernel of $p_{\Zcs}$. Following Ramachandran \cite{Ra}, we recall in Section \ref{SeAVsec} that this observation implies that $\AlbzXk$ is canonically isomorphic to the {\SeAV} $\AXz$.

\begin{remark} \label{Zspeckrem}
Observe that the constant group scheme $\wtZcs$ is the sheafification $\Z^{\Speck}$ of the presheaf $\Z[\Speck]$. It is a large group scheme, since its neutral component is the trivial semiabelian variety $\Speck$.
\end{remark}

\begin{notation} \label{degXnot} 
For any $k$-scheme $X$ with structure morphism denoted by $\pi_{X}: X \ra \Speck$, let $\degX: \ZX \ra \wtZcs$ be the morphism of sheaves induced by the ``summing coefficients'' morphism: 
\begin{align*}
\mathsmap{d_{X}(T)}{\Z[X](T)}{\Z[\Speck](T)}{\sum n_{i} x_{i}}{\sum n_{i},}
\end{align*}
of presheaves, where the summation on the right-hand side represents $\sum n_{i}$ copies of the structure morphism 
$T \ra \Speck$. Equivalently, $\degX$ may be defined as the morphism: 
\begin{align*}
\Z^{\pi_{X}}: \ZX \lra \Z^{\Speck} = \wtZcs
\end{align*}
obtained by applying the functor $\Z^{?}$ to the structure map $\pi_{X}$. In case $X = \Speck$, we have 
$\Z^{\Speck} = \wtZcs$ by Remark \ref{Zspeckrem}, and the degree map is equal to $1_{\wtZcs}: \wtZcs \ra \wtZcs$.
\end{notation}

\begin{proposition} \label{AlbSpeck}
Let $k$ be a perfect field. The Albanese scheme $\Alb_{\Speck/k}$ of $\Speck$ is canonically isomorphic to $\Zcs$. More precisely, the unique morphism $f$ such that the diagram\tn{:}
\begin{align*}
\xymatrix{
\wtZcs \ar@{->}[rr]^(0.42){\mu_{\Speck}} \ar@{->}[ddr]_(0.5){1_{\wtZcs}} 
& & \wAlb_{\Speck/k} \ar@{->}[ddl]^(0.5){f} \\
& & \\
& \wtZcs &}
\end{align*}
commutes is an isomorphism with inverse $\mu_{\Speck}$.
\end{proposition}
\begin{prooff}
The commutativity of the above diagram means that $f$ is a left inverse for $\mu_{\Speck}$. To see that it is also a right-inverse, observe that the outer two triangles of the diagram: 
\begin{align*}
\xymatrix{
\wtZcs  \ar@{->}[dd]_(0.5){\mu_{\Speck}} \ar@{->}[ddrr]_(0.5){1_{\wtZcs}}  
\ar@/^0.9pc/@{->}[ddrrrr]_(0.5){\mu_{\Speck}} \ar@/^1.1pc/@{->}[ddrrrrrr]^(0.5){1_{\wtZcs}} & & & & & & \\
& & & & & & \\
\wAlb_{\Speck/k} \ar@{->}[rr]_(0.58){f} & & \wtZcs  \ar@{->}[rr]_(0.42){\mu_{\Speck}} & & 
\wAlb_{\Speck/k} \ar@{->}[rr]_(0.6){f} & & \wtZcs}
\end{align*}
are identical and commute by definition of $f$, and the middle triangle is obviously commutative. Thus the whole diagram commutes. Applying the universal property of $\Alb_{\Speck/k}$ to the subdiagram consisting of the left-hand and middle triangles, we see that both $\mu_{\Speck} \circ f$ and $1_{\Alb_{\Speck/k}}$ are morphisms such that the subdiagram commutes, hence they are equal, so $f$ is a right inverse to $\mu_{\Speck}$. 
\end{prooff}

\begin{notation} \label{jxonnot} 
Let $(X, \xo)$ be a pointed $k$-scheme. For all $n \in \Z$ we denote by $\jxon: \wt{X} \ra \Z[X]$ the morphism given by
\begin{align*}
\mathsmap{\jxon(T)}{\wt{X}(T)}{\Z[X](T)}{x}{x - n \cdot \xoT,}
\end{align*}
for all $k$-schemes $T$, with the convention $\jxo := \jxoo$. Further, we define $\ixon$ to be the composite morphism:
\begin{align*}
\ixon: \wt{X} \slra{\jxon} \Z[X] \slra{\sigma_{X}} \ZX \slra{\muX} \wAlbXk
\end{align*}
of sheaves of sets on $\ket$, with the convention $\ixo := \ixo^{1}$.
\end{notation}

\begin{notation}
For a locally algebraic $k$-scheme $X$, we denote by $\pioX$ the \emph{scheme of connected components of $X$} and by 
$q_{X}: X \ra \pioX$ the associated canonical morphism. By \cite[Chapter I, \S 4, Proposition 6.5]{DG}, this has the following universal property. For every $k$-scheme $E$ which is {\etl} over $k$ and morphism $f: X \ra E$, there exists a unique morphism $g: \pioX \ra E$ such that $g \circ q_{X} = f$. Moreover, the morphism $q_{X}$ is faithfully flat, and its fibres are the connected components of $X$.
\end{notation}

\begin{notation} \label{pZcsnot} 
Let us denote by $p_{\Zcs}$ the morphism induced by the universal property of $\wAlbXk$ in the diagram\tn{:}
\begin{align*}
\xymatrix{
\ZX \ar@{->}[rr]^(0.5){\muX} \ar@{->}[ddr]_(0.5){\degX} 
& & \wAlbXk \ar@{->}[ddl]^(0.5){p_{\Zcs}} \\
& & \\
& \wtZcs &}
\end{align*}
In view of the functoriality of the morphism $\muX$ (Remark \ref{Albfuncrem}) and the fact that 
$\degX = \Z^{\pi_{X}}$ (Notation \ref{degXnot}), we see that the morphism $p_{\Zcs}$ is equal to $\Alb_{\pi_{X}/k}$.
\end{notation} 

\begin{proposition} \label{pZcsdef}
Let $X$ be a connected $k$-scheme such that $\AlbXk$ exists. Then the scheme $\pio(\AlbXk)$ of connected components of $\AlbXk$ is canonically isomorphic to $\Zcs$, and the morphism $p_{\Zcs}$ of Notation \ref{pZcsnot} coincides with the canonical morphism\tn{:}
\begin{align*}
\wt{q} := \wt{q_{\AlbXk}}: \wAlbXk \lra \wt{\pio(\AlbXk)} = \wtZcs. 
\end{align*} 
Therefore, the connected components of $\AlbXk$ are indexed by the integers. 
\end{proposition}
\begin{prooff}
The existence of a canonical isomorphism between $\pio(\AlbXk)$ and $\Zcs$ is a consequence of 
\cite[Definition 1.5]{Ra}, as stated in \emph{loc. cit.}, Equation (4). The equality $p_{\Zcs} = \wt{q}$ follows from the universal property of $\AlbXk$ and from the universal property of $\Z^{\pioX}$ discussed in \cite[Lemma 1.4]{Ra}.
\end{prooff}

\begin{notation} \label{albncon}
Given an integer $n$ and a connected $k$-scheme $X$ such that $\AlbXk$ exists, we denote the connected component of 
$\AlbXk$ corresponding to $n$ by $\AlbnXk$. 
\end{notation}

\begin{corollary} \label{betagamma}
Let $X$ be a connected $k$-scheme such that the Albanese scheme $\AlbXk$ exists. Then\tn{:}
\begin{itemize}
\item[\tn{(i)}] The neutral component $\AlbzXk$ and kernel $\AlbXk^{0}$ of $p_{\Zcs}$ are identical. 
\item[\tn{(ii)}] For any integer $n$, the morphism $\ixon$ of Notation \ref{jxonnot} factors through $\AlbXk^{1-n}$. In particular, $\ixo = \ixo^{1}$ factors through the natural inclusion \\ $i_{X}: \wAlbzXk \ra \wAlbXk$.
\end{itemize}
\end{corollary}
\begin{prooff}
\tbf{(i)} Indeed, the map $p_{\Zcs}: \wAlbXk \ra \wtZcs$ from Proposition \ref{pZcsdef} is a morphism of (large) group schemes, hence preserves the neutral element.\\\\
\tbf{(ii)} Consider the diagram:
\begin{equation*}
\xymatrix{
\wt{X} \ar@{->}[rr]^(0.5){\jxon} \ar@/^2pc/@{->}[rrrrrr]^(0.5){\ixon} \ar@{->}[ddrrrr]_(0.5){\gammaX^{n}} & & 
\Z[X] \ar@{->}[rr]^(0.5){\sigma_{X}} \ar@{->}[ddrr]^(0.5){d_{X}} & & 
\ZX \ar@{->}[rr]^(0.5){\muX} \ar@{->}[dd]^(0.5){\degX} & & 
\wAlbXk \ar@{->}[ddll]^(0.5){p_{\Zcs}} \\
& & & & & & \\
& & & & \wtZcs, & &}
\end{equation*}
where $\gammaX^{n} \defeq d_{X} \circ \jxon$, so that the left-hand triangle commutes. Now, $p_{\Zcs}$ is the morphism induced by universality of $\wAlbXk$, so the right-hand triangle also commutes. Since the middle triangle commutes by construction of $d_{X}$ and $\degX$, so too does the large outer one. Let $[n]$ denote the connected component of the group scheme $\Zcs$ corresponding to the integer $n$. By definition of the morphism $\jxon$, 
\begin{align*}
\gammaX^{n} = d_{X} \circ \jxon \quad \tn{ lands in } \quad [1-n], 
\end{align*}
for all $n \in \Z$, which implies that:
\begin{align*}
\im(p_{\Zcs} \circ \ixon) = \im(\gammaX^{n}) \sseq [1-n],
\end{align*}
and so it follows that $\ixon$ factors through $\AlbXk^{1-n}$.
\end{prooff}

\begin{remark}
Let $X$ be a connected $k$-scheme such that $\AlbXk$ exists. It is shown in \cite[Corollary 1.10]{Ra} that the neutral component $\AlbzXk$ of $\AlbXk$ is canonically isomorphic to the {\SeAV} of $X$ over $k$, and $\AlbXk^{1}$ is canonically isomorphic to $A^{1}_{X}$.
\end{remark}

\section{Main properties of $\AlbXk$}

\begin{remark} 
We recall from Chapter 1 some notation concerning motives. Let $X$ be a smooth scheme over $k$, let $\shfF \in \STEkQ$ be an {\etl} sheaf with transfers, and let $i$ be an integer. Then: 
\begin{itemize}
\item[\tn{(i)}] $L(X)$ denotes (Notation \ref{Lfuncnot}) the sheaf $\VCork(\pul{X}, X)$. 
\item[\tn{(ii)}] $\Cx(\shfF)$ denotes (Definition \ref{Cxdef}) the complex in $\Dm(\STtkQ)$ with $n^{\tn{th}}$ entry:
\begin{align*}
\Cnx(\shfF)(U) := \shfF(U \x \Delta^{n}).
\end{align*}
Here, $\Delta^{n}$ denotes the spectrum of the $k$-algebra:
\begin{align*}
A^{n} := k[T_{0}, \ldots, T_{n}]\Big/ \left\langle \left( \sum_{i=0}^{n} t_{i} \right) - 1 \right\rangle,
\end{align*}
and the $n^{\tn{th}}$ differential of $\Cx(\shfF)$ is $\pd_{n} := \shfF(\id_{U} \x \Spec(d^{n}))$, where $d^{n}$ is the alternating sum of the maps:
\begin{align*}
& \mathsmap{\delta^{n}_{i}}{A^{n}}{A^{n-1}}{t_{j}}
{\left\{ \begin{array}{rl}
t_{j}, & \tn{if } j < i, \\
0, & \tn{if } j = i, \\
t_{j-1}, & \tn{if } j > i,
\end{array} \right.} 
\end{align*}
\item[\tn{(iii)}] $M(X)$ is defined (Notation \ref{Mhinot}) to be the motive $\Cx(L(X)) \in \Dm(\STtkQ)$.
\item[\tn{(iv)}] $\hiet(\shfF)$ is defined to be (Notation \ref{Mhinot}) the cohomology sheaf 
$\shfH^{-i}(\Cx(\shfF)) \in \STtkQ$.  
\end{itemize}
Recall further that we write $\hzet(X)$ in place of $\hzet(L(X))$.
\end{remark}

\begin{remark} \label{hznisspeckrem}
By homotopy invariance of $\Speck$, it follows that $\hiet(\one)$ is equal to $\wtZcs$ if $i = 0$, and vanishes otherwise.
\end{remark}

\begin{construction} \label{consthX}
Let $X$ be a smooth scheme over $k$. We will construct a morphism $h_{X}: \ZX \ra \hzet(M(X))$ of {\etl} sheaves with transfers. Firstly, observe that $L(X)(\pul{X} \x_{k} \Delta^{0}) = L(X)(\pul{X})$, since $\Delta^{0} = \Speck$, and so $\hzet(M(X))$ is the sheafification of the presheaf with transfers:
\begin{align*}
\mathsmaptwo{\Smck}{\Ab}{T}{L(X)(T) \big/ \im(d^{X}_{1}(T)).}
\end{align*}
We define now the morphism $h^{\prm \prm}_{X}: \wtX \ra \hzet(M(X))$ to be the composition:
\begin{align*}
\wtX \sxra{\Gamma} L(X) \slra{p} L(X) \big/ \im(d^{X}_{1}) \slra{s} \hzet(M(X)).
\end{align*}
Here, $\Gamma$ is the functor assigning to a morphism $f: T \ra X$ of $k$-schemes its graph, $p$ denotes the natural surjection and $s$ denotes the canonical map from a presheaf to its sheafification. The map $h^{\prm \prm}$ then extends $\Z$-linearly from $\wtX$ to a morphism:
$h_{X}^{\prm}: \Z[X] \lra \hzet(M(X))$. By the universal property of the sheafification functor, there is now a unique morphism:
\begin{align*}
h_{X}: \ZX \lra \hzet(M(X))
\end{align*}
such that the diagram:
\begin{align*}
\xymatrix{
\Z[X] \ar@{->}[ddr]_(0.5){h_{X}^{\prm}}  \ar@{->}[rr]^(0.5){\sigma_{X}} & & 
\ZX \ar@{->}[ddl]^(0.5){h_{X}} \\
& & \\
& \hzet(M(X)) & }
\end{align*}
is commutative. 
\end{construction}

\begin{notation} \label{picdivnot}
Let $X$ be a $k$-scheme. Following \cite[Definition 2.2]{Kl}, we recall the definition of the Picard functor 
$\PicXk: \Schk \ra \Ab$: 
\begin{align*}
\PicXk(T) := \Pic(X_{T})/\Pic(T),
\end{align*}
where $\Pic(X_{T})$ denotes the usual Picard group of isomorphism classes of line bundles on $X_{T}$.
Following \cite[Definition 3.6]{Kl}, we denote by $\DivXk$ the abelian presheaf $\Schk \ra \Ab$ defined by:
\begin{align*}
\DivXk(T) := \{ \tn{effective divisors $D$ on $X_{T}/T$} \},
\end{align*}
for all $k$-schemes $T$. Further, following \cite[Definition 4.6]{Kl}, we denote by 
$a_{X/k}: \DivXk \ra \PicXk$ the \emph{Abel map} of functors defined by: 
\begin{align*}
\mathsmap{a_{X/k}(T)}{ \DivXk(T)}{\PicXk(T)}{D}{\shfO_{X_{T}}(D)}
\end{align*}
for all $k$-schemes $T$.
\end{notation}

\begin{proposition} \label{picdivrep}
Let $X$ be a projective $k$-scheme. Then the presheaves $\wPicXk: \Schk \ra \Ab$ and 
$\wDivXk: \Schk \ra \Ab$ of Notation \ref{picdivnot} are both representable, and are hence both sheaves on the big {\etl} site. Moreover, the scheme $\PicXk$ representing $\wPicXk$ is a large group scheme. 
\end{proposition}
\begin{prooff}
Representability of $\wDivXk$ is \cite[Theorem 3.7]{Kl} and representability of $\wPicXk$ is \cite[Theorem 4.8]{Kl}, where it is also shown that the representing scheme $\PicXk$ is locally of finite type. The neutral component $\PicoXk$ is a semiabelian variety over $k$, by \cite[\S 2.1]{Ra}. Hence $\PicXk$ is a large group scheme. Sheafhood of both $\wDivXk$ and $\wPicXk$ now follows immediately from their representability by Proposition \ref{Gsheaf}.
\end{prooff}

\begin{proposition} \label{albpic}
Let $C$ be a smooth projective curve over a perfect field $k$ and denote by $\Gamma^ {\p}: \ZC \ra  \DivCk$ the morphism of {\etl} sheaves induced by the morphism $\Gamma: \Z[C] \ra  \DivCk$ of presheaves sending a formal 
$\sum n_{i} \cdot f_{i}$ of morphisms $f_{i}: T \ra C$ to the associated formal sum of the graphs $\Gamma(f_{i}) \sseq T \x C$. Then the composite morphism\tn{:}
\begin{align*}
\ZC \slra{\Gamma^ {\p}} \wDivCk \slra{\wt{a_{X/k}}} \wPicCk
\end{align*}
is a morphism of abelian sheaves, and the morphism $\psi_{C}: \AlbCk \ra \PicCk$ induced by the universal property of the Albanese scheme is an isomorphism of large group schemes. 
\end{proposition}
\begin{prooff}
\cite[Theorem 6.1]{Ra}.
\end{prooff}

\begin{remark} \label{hznismof}
Let $X$ be any $k$-scheme such that $\AlbXk$ exists. Then $\hzet(\Mof(\AlbXk)) = \wAlbXk$, since $\Mof(\AlbXk)$ is by definition a complex concentrated in degree zero, and $\wAlbXk$ is an {\etl} sheaf with transfers, by Proposition \ref{Gtrans}, because $\AlbXk$ is a commutative group scheme, by Proposition \ref{picdivrep}.
\end{remark}

\begin{construction} \label{bXmap}
Let $(X, \xo)$ be a pointed $k$-scheme. We recall the construction of Szamuely and Spie{\ss} \cite[postscript to Lemma 3.2]{SS} of a canonical morphism from $\hzet(M(X))$ to $\wAlbXk$. Let: 
\begin{align*}
\alpha_{\AlbXk}: M(\AlbXk) \ra \Mof(\AlbXk)
\end{align*}
be the morphism defined in \ref{alphanot}. Observe that the composite morphism: 
\begin{align*}
M(X) \sxra{M(\ixo)} M(\AlbXk) \sxra{\alpha_{\AlbXk}} \Mof(\AlbXk),
\end{align*}
factors through $\hzet(M(X))$, considered as a complex concentrated in degree zero, since both $M(X)$ and $\Mof(\AlbXk)$ are complexes concentrated in non-positive degrees, and the latter is concentrated in degree zero. This yields a morphism:
\begin{align*}
b_{X}: \hzet(M(X)) \lra \wAlbXk,
\end{align*}
in $\STEkQ$ such that the diagram: 
\begin{align*}
\xymatrix{
\hzet(M(X)) \ar@{->}[ddrrr]^(0.5){b_{X}} \ar@{->}[dd]_(0.5){\hzet(M(\ixo))} & & & \\
& & & \\
\hzet(M(\AlbXk)) \ar@{->}[rrr]_(0.55){\hzet(\alpha_{\AlbXk})} & & & \wAlbXk.}
\end{align*}
is commutative. Here we have used the fact that $\hzet(\Mof(\wAlbXk)) = \wAlbXk$, again due to the fact that $\Mof(\wAlbXk)$ is a complex concentrated in degree zero. 
\end{construction}

\begin{proposition} \label{bXcan}
Let $(X, \xo)$ be a pointed scheme over $k$. The morphism\tn{:}
\begin{align*}
b_{X}: \hzet(M(X)) \lra \wAlbXk
\end{align*}
of Construction \ref{bXmap} is canonical in the sense that for any $k$-scheme morphism $f: X \ra Y$, the diagram\tn{:}
\begin{align*}
\xymatrix{
\hzet(M(X)) \ar@{->}[dd]_(0.5){\hzet(M(f))} \ar@{->}[rr]^(0.55){b_{X}} & & 
\wAlbXk \ar@{->}[dd]^(0.5){\wAlb_{f/k}} \\
& & \\
\hzet(M(Y)) \ar@{->}[rr]_(0.55){b_{Y}} & & \wAlb_{Y/k}}
\end{align*}
commutes. 
\end{proposition}
\begin{prooff}
By Construction \ref{bXmap}, the morphism $b_{X}$ is equal to the composition of $\hzet(M(\ixo))$, which is canonical by the universal property of the Albanese scheme, and $\hzet(\alpha_{\AlbXk})$, which is canonical by Proposition \ref{alphafunc}.
\end{prooff}

\begin{corollary} \label{bXZcscomm}
Let $(X, \xo)$ be a pointed scheme over $k$ with structure morphism $\pi_{X}: X \ra \Speck$. Then the diagram\tn{:}
\begin{align*}
\xymatrix{
\hzet(M(X)) \ar@{->}[ddr]_(0.5){\hzet(M(\pi_{X}))} \ar@{->}[rr]^(0.55){b_{X}} & & 
\wAlbXk \ar@{->}[ddl]^(0.5){p_{\Zcs}} \\
& & \\
& \wtZcs &}
\end{align*}
is commutative. 
\end{corollary} 
\begin{prooff}
By Remark \ref{hznisspeckrem} and Proposition \ref{AlbSpeck}, respectively, we have that $\hzet(M(\Speck))$ and $\Alb_{\Speck/k}$ are both equal to $\wtZcs$, and it is easy to see that:
\begin{align*}
b_{\Speck} = 1_{\wtZcs}: \wtZcs \lra \wtZcs.
\end{align*}
Recalling that $p_{\Zcs}$ is equal to $\wAlb_{\pi_{X}/k}$ by Notation \ref{pZcsnot}, the result now follows by applying Proposition \ref{bXcan} to the structure morphism $\pi_{X}$.
\end{prooff}

\begin{proposition} \label{hzalbC}
Let $(C, \xo)$ be a smooth connected pointed curve over $k$. Then the morphism\tn{:}
\begin{align*}
\bC: \hzet(M(C)) \lra \wAlbCk
\end{align*}
of Construction \ref{bXmap} is an isomorphism of abelian {\etl} sheaves with transfers. 
\end{proposition}
\begin{prooff}
See \cite[\S 3]{Li1} for the construction of a canonical isomorphism:
\begin{align*}
\bC^{\p}: \hzet(M(C)) \ra \wPicCk.
\end{align*}
(Note that a similar result is proven in \cite[Theorem 3.1]{SV} for a smooth irreducible affine scheme of relative dimension one over a normal affine base scheme $S$). One can check that our $\bC$ is the composition of $\bC^{\p}$ with the map $\psi_{C}: \AlbCk \ra \PicCk$ defined in Proposition \ref{albpic}.
\end{prooff}

\begin{remark}
If $X$ is not a curve, then $\hzet(M(X))$ and $\wAlbXk$ are not isomorphic in general. Indeed, if $X=A$ is an abelian variety, one sees that: 
\begin{align*}
\hzet(M(A))(k) \neq \wAlb_{A/k}(k),
\end{align*}
as follows. Firstly, we will see in Proposition \ref{AlbGknice} below that $\Alb_{A/k} = A \x_{k} \Zcs$, hence 
$A(k) = \Albz_{A/k}(k)$ is the kernel of the degree map: 
\begin{align*}
\deg: \wAlb_{A/k}(k) \lra \Z
\end{align*}
on $k$-points. On the other hand, in \cite[\S 0]{Bl1} it is shown that if $I$ is the kernel of the degree map:
\begin{align*}
\deg: \hzet(M(A))(k) = \CH_{0}(A) \lra \Z,
\end{align*}
where $\CH_{0}(A)$ is the Chow group of zero cycles on $A$, then (i) $I$ is generated by zero cycles of the form $(a) - (0_{A})$, where 
$a$ is a $k$-point of $A$ and $0_{A}$ is the identity, and (ii) there is a short exact sequence:
\begin{align*}
0 \lra I^{*2} \lra I \slra{f} A(k) \ra 0,
\end{align*}
where $f$ is the map sending $(a) - (0_{A}) \in I$ to $a \in A(k)$ and $I^{*2}$ denotes the second Pontrjagin power of $I$ (see Definition \ref{pontproddef} for the Pontrjagin product). Moreover, if $k$ is algebraically closed and $A$ has dimension $n$, then $I^{* n+1} = 0$, but $I^{*m}$ does not vanish in general for $1 \le m \le n$. and so:
\begin{align*}
A(k) \simeq I/I^{*2} \not\simeq I,
\end{align*}
implying that $\hzet(M(A))(k) \neq \wAlb_{A/k}(k)$.
\end{remark}

\begin{construction} \label{AlbGkmorphi}
Let $G$ be a semiabelian variety over $k$. We construct a canonical morphism of commutative groups schemes:
\begin{align*}
\chi_{G}: \Alb_{G/k}\lra \Zcs \x G.
\end{align*}
Denote by $\tn{sum}_{G}: \Z^{G} \ra \wt{G}$ the morphism of abelian {\etl} sheaves induced by the morphism of presheaves $\tn{sum}_{G}^{\p}: \Z[G] \ra \wtG$ which is defined, for all connected $k$-scheme $T$, by:
\begin{align*}
\mathsmap{\tn{sum}_{G}^{\p}(T)}{\Z[G](T)}{\wtG(T)}{\ub{\sum_{i} n_{i} \cdot f_{i}}_{\tn{\tiny{formal sum}}}}
{\ub{\sum_{i} n_{i} \cdot f_{i}}_{\tn{\tiny{group law sum}}}.}
\end{align*}
Let $p_{G}$ denote the unique morphism, induced by the universal property of $\Alb_{G/k}$, such that the left-hand diagram: 
\begin{align*}
\xymatrix{
\Z^{G} \ar@{->}[ddr]_(0.5){\tn{sum}_{G}}  \ar@{->}[rr]^(0.45){\mu_{G}} & & 
\wAlb_{G/k} \ar@{->}[ddl]^(0.5){p_{G}} \\
& & \\
& \wtG & } \qaq
\xymatrix{
\Z^{G} \ar@{->}[ddr]_(0.5){\deg_{G}} \ar@{->}[rr]^(0.45){\mu_{G}} & & 
\wAlb_{G/k} \ar@{->}[ddl]^(0.5){p_{\Zcs}} \\
& & \\
& \wtZcs & }
\end{align*}
commutes, and recall from Proposition \ref{pZcsdef} that $p_{\Zcs}$ is the unique morphism such that the right-hand diagram commutes, where $\deg_{G}: \Z^{G} \ra \wtZcs$ is the degree morphism from Notation \ref{degXnot}. We define $\chi_{G}: \Alb_{G/k} \ra \wtZcs \x G$ to be the unique morphism, induced by the universal property of the fibre product, such that the diagram: 
\begin{align*}
\xymatrix{
& \Alb_{G/k} \ar@/^1.2pc/@{->}[ddr]^(0.5){p_{\Zcs}} \ar@{->}[d]^(0.5){\chi_{G}}  
\ar@/_1.2pc/@{->}[ddl]_(0.5){p_{G}} & \\
& G \x \wtZcs \ar@{->}[dr]_(0.5){\pi_{2}} \ar@{->}[dl]^(0.5){\pi_{1}} & \\
G & & \wtZcs}
\end{align*}
is commutative. 
\end{construction}

\begin{proposition} \label{AlbGknice}
Let $G$ be a semiabelian variety over a perfect field $k$, and let\tn{:}
\begin{align*}
& p_{\Zcs}: \Alb_{G/k} \lra \Zcs, \qquad p_{G}: \Alb_{G/k} \lra G \qaq \\
& \chi_{G}: \Alb_{G/k} \lra \Zcs \x G
\end{align*}
be the morphisms of Construction \ref{AlbGkmorphi}. Then\tn{:}
\begin{itemize}
\item[\tn{(i)}] The morphism $\chi_{G}$ is an isomorphism of {\etl} sheaves with transfers. 
\item[\tn{(ii)}] The morphism $p_{G}$ splits the morphism $\iota_{e}: G \ra \AlbGk$ of Notation \ref{jxonnot}.
\end{itemize} 
\end{proposition}
\begin{prooff}
\tbf{(i)} \cite[Corollary 1.12, Part (v)]{Ra}.\\\\
\tbf{(ii)} The middle and right-hand triangles in the diagram:
\begin{equation*}
\xymatrix{
\wt{G} \ar@{->}[rr]_(0.5){j_{e}} \ar@{->}[ddrrrr]_(0.5){1_{\wtG}} \ar@/^1.2pc/@{->}[rrrrrr]^(0.5){\iota_{e}} & & 
\Z[G] \ar@{->}[rr]_(0.55){\sigma_{G}} \ar@{->}[ddrr]^(0.5){\tn{sum}_{G}^{\p}} & & 
\Z^{G} \ar@{->}[rr]_(0.45){\mu_{G}} \ar@{->}[dd]^(0.5){\tn{sum}_{G}} & & 
\wAlbGk \ar@{->}[ddll]^(0.5){p_{G}} \\
& & & & & & \\
& & & & \wtG, & &}
\end{equation*}
both commute by Construction \ref{AlbGkmorphi}, and it is clear from the definitions of the morphisms $j_{e}$ and $\tn{sum}_{G}^{\p}$ that the left-hand triangle commutes also. Since the morphism $\iota_{e}: G \ra \AlbGk$ is equal to the composite of the three maps along the top row by its definition in Notation \ref{jxonnot}, the whole diagram commutes. 
\end{prooff}

\section{The Voevodsky motive of a torus}

\begin{proposition} \label{symosym1} 
Let $k$ be a perfect field and let $\Gm$ denote the multiplicative group scheme over $k$. Then $\Symn(\Mof(\Gm))$  vanishes for $n \ge 2$, and so we have the identities\tn{:}
\begin{itemize}
\item[\tn{(i)}] $\Sym(\Mof(\Gm)) = \one \os \Mof(\Gm)$ and 
\item[\tn{(ii)}] $\hzet(\Sym(\Mof(\Gm))) = \wt{\Zcs} \os \wtGm$
\end{itemize}
in the categories $\DMeetkQ$ and $\STEkQ$, respectively. 

In particular, $\Gm$ has finite Kimura dimension, implying that any torus $T$ has finite Kimura dimension. 
\end{proposition}
\begin{prooff}
By \cite[Theorem 4.1]{MVW}, $\Mof(\Gm)$ is isomorphic to $\tato$. It now follows from \cite[Lemma 2.5.10]{Big}, together with the fact that $\Alt^{n}(\one)$ vanishes for $n \ge 2$, that: 
\begin{align*}
\Symn(\Mof(\Gm)) = \Symn(\tato) = \Altn(\one(1))[n] = \Altn(\one)(n)[n] = 0,
\end{align*}
which gives the first identity. The second now follows by additivity of the functor $\hzet$, the fact that 
$\hzet(\one) = \wtZcs$, and the definition of $\Mof(G)$ as the complex concentrated in degree zero coming from the sheaf $\wtG$.
\end{prooff}

\begin{proposition} \label{bGmhophiGm}
Let $k$ be a perfect field and let $\Gm$ denote the multiplicative group scheme over $k$. Then the diagram\tn{:}
\begin{align*}
\xymatrix{
\hzet(M(\Gm)) \ar@{->}[dd]_(0.5){\hzet(\vp_{\Gm})} \ar@{->}[rr]^(0.55){b_{\Gm}} & & 
\wAlb_{\Gm/k} \ar@{->}[dd]^(0.5){\chi_{\Gm}} \\
& & \\
\wt{\Zcs} \os \wt{\Gm} \ar@{<-}[rr]_(0.55){\eta_{\Gm, \Zcs}} & & \wt{\Zcs \x \Gm}}
\end{align*}
is commutative, hence $\hzet(\vp_{\Gm})$ is an isomorphism since the other three maps are isomorphisms. Here we recall that $\vp_{\Gm}: M(\Gm) \ra \Sym(\Mof(\Gm))$ is the morphism of Notation \ref{vpGnot}, and $\eta_{\Gm, \Zcs}: \wt{\Zcs \x \Gm} \lra \wt{\Zcs} \os \wt{\Gm}$ is the canonical isomorphism of Notation \ref{etaGHdef}.
\end{proposition}
\begin{prooff}
For $i = 1, 2$, let us denote by $\pr_{i}$ the canonical projection of $\wt{\Zcs} \os \wt{\Gm}$ onto the $i^{\tn{th}}$ summand. Since the finite direct sum $\wt{\Zcs} \os \wt{\Gm}$ is a product in the category $\DMeetkQ$, to show that the diagram commutes it suffices to prove that:
\begin{align}
\pr_{i} \circ \eta_{\Gm, \Zcs} \circ \chi_{\Gm} \circ b_{\Gm} = \pr_{i} \circ \hzet(\vp_{\Gm}), \label{projeq}
\end{align}
for $i = 1, 2$. Recall the maps $p_{\Zcs}: \Alb_{G/k} \ra \Zcs$ and $p_{\Gm}: \Alb_{G/k} \ra \Gm$ introduced in Construction \ref{AlbGkmorphi}. It follows directly from the definitions of the morphisms $\chi_{\Gm}$ and 
$\eta_{\Gm, \Zcs}$ that the diagram: 
\begin{align*}
\xymatrix{
& & & & \wt{\Zcs} \ar@{<-}[d]^(0.5){\pr_{1}} \\
\wAlb_{\Gm/k} \ar@{->}[rr]_(0.5){\chi_{\Gm}} \ar@/^1.1pc/@{->}[rrrru]^(0.55){p_{\Zcs}} 
\ar@/_1.1pc/@{->}[rrrrd]_(0.55){p_{\Gm}} & &
\wt{\Zcs \x \Gm} \ar@{->}[rr]_(0.5){\eta_{\Gm, \Zcs}} & & \wt{\Zcs} \os \wt{\Gm} \\
& & & & \wt{\Gm} \ar@{<-}[u]_(0.5){\pr_{2}}}
\end{align*}
commutes. Moreover, since the morphism $\vp_{\Gm}$ is by definition the sum of $\vp_{\Gm}^{0}$ and $\vp_{\Gm}^{1}$, the diagram:
\begin{align*}
\xymatrix{
& & \hzet(M(\Gm)) \ar@{->}[dd]^(0.5){\hzet(\vp_{\Gm})} \ar@/^1.1pc/@{->}[ddrr]^(0.5){\hzet(\vp_{\Gm}^{0})} 
\ar@/_1.1pc/@{->}[ddll]_(0.5){\hzet(\vp_{\Gm}^{1})} & & \\
& & & & \\
\wt{\Zcs} \ar@{<-}[rr]_(0.45){\pr_{1}} & & \wt{\Zcs} \os \wt{\Gm} \ar@{->}[rr]_(0.57){\pr_{2}} & & \wt{\Gm}}
\end{align*}
also commutes. By the commutativity of the above two diagrams, Relation (\ref{projeq}) reduces to:
\begin{align}
p_{\Zcs} \circ b_{\Gm} & = \hzet(\vp^{0}_{\Gm}), \qquad \tn{for $i = 0$,} \quad \tn{and} \label{prelone} \\
p_{\Gm} \circ b_{\Gm} & = \hzet(\vp^{1}_{\Gm}), \qquad \tn{for $i = 1$.} \label{preltwo}
\end{align}
To see that the first of these identities holds, observe that by Remark \ref{zeroonevp}, the morphism $\vp^{0}_{\Gm}$ is equal to  $M(\pi_{\Gm}): M(\Gm) \ra \one$, where $\pi_{\Gm}: \Gm \ra \Speck$ is the structure morphism of $\Gm$, so $\hzet(\vp^{0}_{\Gm}) = \hzet(M(\pi_{\Gm}))$; hence (\ref{prelone}) follows from Corollary \ref{bXZcscomm}. 

To verify the second identity, we note again by Remark \ref{zeroonevp} that $\vp^{1}_{\Gm}$ is equal to the morphism
$\alpha_{\Gm}: M(\Gm) \ra \Mof(\Gm) $ of Notation \ref{alphanot}, and so:
\begin{align}
\hzet(\vp^{1}_{\Gm}) = \hzet(\alpha_{\Gm}). \label{vpalpha}
\end{align} 
Consider now the diagram: 
\begin{align*}
\xymatrix{
\hzet(M(\Gm)) \ar@{->}[rrr]^(0.45){\hzet(M(\iota_{e}))} \ar@{->}[ddrrr]_(0.5){b_{\Gm}}
\ar@{->}[dd]_(0.5){\hzet(\alpha_{\Gm})} & & &
\hzet(M(\Alb_{\Gm/k})) \ar@{->}[dd]^(0.5){\hzet(\alpha_{\Alb_{\Gm/k}})} \\
& & & \\
\hzet(\Mof(\Gm)) \ar@{->}[rrr]_(0.45){\hzet(\Mof(\iota_{e}))} \ar@{=}[dd] & & & 
\hzet(\Mof(\Alb_{\Gm/k})) \ar@{=}[dd] \\
& & & \\
\wt{\Gm} \ar@{->}[rrr]_(0.45){\iota_{e}} & & & \wAlb_{\Gm/k}.}
\end{align*}
The perimeter of the top square is commutative by functoriality of the morphism $\alpha_{\Gm}$ shown in Proposition \ref{alphafunc}, and the upper-right triangle therein by definition of the morphism $b_{\Gm}$. The bottom square is commutative because $\Mof(\Gm)$ and $\Mof(\Alb_{\Gm/k})$ are complexes concentrated in degree zero. Combining the commutativity of the diagram with (\ref{vpalpha}), we see that: 
\begin{align*}
b_{\Gm} = \iota_{e} \circ \hzet(\alpha_{\Gm}) = \iota_{e} \circ \hzet(\vp^{1}_{\Gm}),
\end{align*}
which implies the identity (\ref{preltwo}), since $p_{\Gm}$ splits $\iota_{e}$ by Proposition \ref{AlbGknice}. This completes the proof.
\end{prooff}

\begin{corollary} \label{vpGmcor}
Let $k$ be a perfect field, and let $\Gm$ be the multiplicative group scheme over $k$. Then the morphism\tn{:}
\begin{align*}
\vp_{\Gm}: M(\Gm) \lra \Sym(\Mof(\Gm))
\end{align*}
of Notation \ref{vpGnot} is an isomorphism in the category $\DMeetkQ$.
\end{corollary}
\begin{prooff}
We note first that the morphism $\vp_{\Gm}$ exists, by Proposition \ref{symosym1}, since $\Gm$ has finite Kimura dimension. Now, by combining \cite[Lemma 2.13]{MVW} with \cite[Theorem 4.1]{MVW}, one sees that $M(G)$ is isomorphic to $\one \os \Mof(\Gm)$, and is hence cohomologically concentrated in degree zero, since both summands are. By Proposition \ref{symosym1}, $\Sym(\Mof(\Gm))$ is also isomorphic to $\one \os \Mof(\Gm)$ hence also cohomologically concentrated in degree zero. To see that $\vp_{\Gm}$ is an isomorphism, it therefore suffices to prove that $\hzet(\vp_{\Gm})$ is an isomorphism. But this was shown in Proposition \ref{bGmhophiGm}.
\end{prooff}

\begin{proposition} \label{vptorusisom}
Let $k$ be a perfect field, and let $T$ be a torus over $k$. Then the morphism\tn{:}
\begin{align*}
\vp_{T}: M(T) \lra \Sym(\Mof(T))
\end{align*}
of Notation \ref{vpGnot} is an isomorphism in the category $\DMeetkQ$.
\end{proposition}
\begin{prooff}
Observe that by Proposition \ref{symosym1} the morphism $\vp_{T}$ exists, since $T$ has finite Kimura dimension. By Proposition \ref{etldes}, we may assume that the field $k$ is algebraically closed, and hence that the torus $T$ is split, say isomorphic to the $r^{\tn{th}}$ power $\Gmr$ of the multiplicative group scheme, for some positive integer $r$. By
Proposition \ref{vpmultprop}, the map $\vp_{\Gmr}: M(\Gmr) \ra \Sym(\Mof(\Gmr))$ is an isomorphism provided
$\vp_{\Gm}: M(\Gm) \ra \Sym(\Mof(\Gm))$ is an isomorphism, but this we know from Corollary \ref{vpGmcor}
\end{prooff}
\chapter{The Voevodsky motive of an abelian variety}

Let $A$ be an abelian variety over a perfect field $k$. In this chapter we show that the morphism:
\begin{align*}
\vp_{A}: M(A) \lra \Sym(\Mof(A))
\end{align*}
of Notation \ref{vpGnot} is an isomorphism in the category $\DMeetkQ$, by showing that it is equal to $\Phi(\psi_{A})$, where $\Phi: \Chowek \ra \DMeetkQ$ is the functor of Notation \ref{Phifuncnot}, and:
\begin{align*}
\psi_{A}: h(A) \ilra \Sym(\ho(A))
\end{align*}
is the classical decomposition of the motive $h(A)$ of an abelian variety in the category $\Chowek$. Let us describe the structure of this proof in more detail. 

\tbf{(i)} In the first section, we review the construction of an isomorphism: 
\begin{align*}
\psi_{A}: h(A) \ilra \Sym(\ho(A))
\end{align*}
in the category $\Chowek$. Versions of this classical result have been published by many authors: Grothendieck, Kleiman and Lieberman \cite[Appendix to \S 2]{GKL}, Deninger and Murre \cite{DeMu}, K{\"u}nnemann \cite{Ku1}, and Shermenev \cite{Sh1}. In particular, we stress that none of the mathematical content of this section is due to the author. We simply restate the work of the above-cited authors in our own notation. 

\tbf{(ii)} In the second section, we use the decomposition $h(C) = \hz(C) \os \ho(C) \os h^{2}(C)$ of Remark \ref{hodef} of the Chow motive $h(C)$ of a smooth projective curve $C$ to understand the object $\Mof(\AlbzCk)$ in $\DMeetkQ$, where we recall that $\Mof(\AlbzCk)$ denotes the complex consisting of the {\etl} sheaf with transfers $\wAlbzCk$ concentrated in degree zero. More precisely, we construct an isomorphism: 
\begin{align}
\hzet((\Phi \circ \ho)(C)) \lra \wAlbzCk. \label{CAlbC}
\end{align}

\tbf{(iii)} In the third section, we construct a functorial morphism:
\begin{align*}
D(A): (\Phi \circ \ho)(A) \lra \Mof(A)
\end{align*}
in $\DMeetkQ$, which we show is an isomorphism. The proof proceeds, by first reducing to the case where $A$ is the Serre-Albanese variety $\AlbzCk$ of a smooth projective curve $C$. To show that $D(\AlbzCk)$ is an isomorphism, we 
use the isomorphism (\ref{CAlbC}). 

\tbf{(iv)} In the fourth section, we use the fact $D(A): (\Phi \circ \ho)(A) \ra M(A)$ is an isomorphism to show that 
the morphisms:
\begin{align*}
& \Phi(\psi_{A}): (\Phi \circ h)(A) \lra \Sym((\Phi \circ \ho)(A)) \qaq \\
& \vp_{A}: M(A) \lra \Sym(\Mof(A))
\end{align*}
agree up to isomorphism. Since $\psi_{A}: h(A) \lra \Sym(\ho(A))$ is an isomorphism in $\Chowek$, this implies that $\vp_{A}: M(A) \ra \Sym(\Mof(A))$ is an isomorphism in $\DMeetkQ$.

\section{The Chow motive of an abelian variety}

Let $A$ be an abelian variety of dimension $g$ over a perfect field $k$. The aim of this section is to recall the constructions of K{\"u}nnemann \cite{Ku1} of a direct sum decomposition $h(A) = \bos_{0 \le n \le 2g} h^{n}(A)$ of the Chow motive $h(A)$ of $A$, together with an isomorphism:
\begin{align*}
\psi_{A}^{n}: \hn(A) \lra \Symn(\ho(A)).
\end{align*}
Recall that we use the covariant definition of the category $\Chowk$ (see Convention \ref{covcon}), whereas K{\"u}nnemann's Chow motives are contravariant. The two categories are equivalent via the functor which is the identity on objects and transposes morphisms, so we denote contravariant Chow motives by $\tChowk$, and note that the definition of all morphisms will be the transpose of the corresponding morphisms of K{\"u}nnemann.

\begin{definition}[\tbf{Pontrjagin product}] \label{pontproddef}
Let $A$ be an abelian variety over a perfect field $k$, and let $\CH(A)$ denote the direct sum $\bos_{r \ge 0} \CH_{r}(A)$. We recall that the \emph{Pontrjagin product} on $\CH(A) \ox \Q$ is defined by:
\begin{align*}
\mathsmap{- \: * \: -}{ \big( \CH(A) \ox \Q \big) \x \big( \CH(A) \ox \Q \big)}{\CH(A) \ox \Q}{(\alpha, \beta)}
{(m_{A})_{*}(\alpha \x \beta),}
\end{align*}
where $(m_{A})_{*}: \CH(A \x A) \ox \Q \ra \CH(A) \ox \Q$ denotes the map induced by the multiplication morphism $m_{A}: A \x A \ra A$.
\end{definition}

\begin{remark}
Let $A$ be an abelian variety over a perfect field $k$ and let $n$ be a positive integer. Recall that $\Sym(n)$ denotes the symmetric group on $n$ letters. Let $\sigma \in \Sym(n)$. We recall from Notation \ref{sigmanot} that:
\begin{align*}
\sigma_{A}: A^{n} \lra A^{n}
\end{align*}
denotes the morphism permuting the factors by $\sigma$.
\end{remark}

\begin{remark} 
Recall the category $\CSmProjk$ of \emph{Chow correspondences}, whose objects are smooth projective varieties over $k$, and whose morphisms are given by:
\begin{align*}
\Mor_{\CSmProjk}(X, Y) := \bos_{i} \CH_{d_{X_{i}}}(X_{i} \x Y) \ox \Q,
\end{align*}
where the $X_{i}$ are the connected components of $X$ and $d_{X_{i}} = \dim(X_{i})$. 
\end{remark}

\begin{notation}[\tbf{K{\"u}nnemann}] \label{pislambda}
Let $A$ be an abelian variety over a perfect field $k$ of dimension $g$. We recall the following morphisms from 
\cite[\S 3]{Ku1}. Since K{\"u}nnemann works in $\tChowk$ but we work in $\Chowk$, our morphisms are the transpose of his. 
\begin{itemize}
\item[\tn{(i)}] From \cite[\S 3.1, Equation (8)]{Ku1}, we have:
\begin{align*}
p_{n} = p_{n, A} & := \frac{1}{(2g - n)!} \log \left( \pt{[\Gamma(1_{A}) ]} \right)^{*(2g - n)} \\
& = \frac{1}{(2g - n)!} 
\left\{ \sum_{n=1}^{\infty} \left( \pt{[\Gamma(1_{A})]} - \pt{[\Gamma(\ve_{A})]} \right)^{*n} \right\},
\end{align*}
where $\ve_{A}: A \ra A$ denotes the constant morphism sending all points to the neutral element $e_{A}$. Note that \cite[Theorem 1.4.1]{Ku1} implies the terms in this expansion vanish for large enough $n$, and so this logarithm gives a well-defined element of $\CH_{g}(A \x A) \ox \Q$. We use the notation $p_{n}$ instead of K{\"u}nnemann's $\pi_{n}$, to avoid confusion with the canonical projector $\pi^{n}_{X}: X^{\ox n} \ra \Symn(X)$ associated to the $n^{\tn{th}}$ symmetric power.
\item[\tn{(ii)}] From \cite[\S 2.6]{Ku1}, we have:
\begin{align*}
s_{n} := \frac{1}{n!} \sum_{\sigma_{A} \in \Sym(n)} [\Gamma(\sigma_{A})],
\end{align*}
which is an element of $\CH_{ng}(A^{n} \x A^{n}) \ox \Q$.
\item[\tn{(iii)}] For elements  $a_{1}, \ldots, a_{n}$ of the multiplicative abelian group $\mu_{2}$ with two elements, let 
$\chi(a_{1}, \ldots, a_{n}) \in \{ 1, -1 \}$ be the ``sign'' of the product $a_{1} a_{2} \cdots a_{n}$, considered as a rational number, and denote by: 
\begin{align*}
\mult_{a_{1}, \ldots, a_{n}}: A^{n}\lra A^{n}
\end{align*}
the morphism restricting to multiplication by $a_{i}$ in the $i^{\tn{th}}$ factor. We recall \cite[\S 2.6]{Ku1}: 
\begin{align*}
\lambda_{n} := \frac{1}{2^{n}} \sum_{(a_{1}, \ldots, a_{n}) \in \mu_{2}^{n}} 
\chi(a_{1}, \ldots, a_{n})[ \Gamma(\mult_{a_{1}, \ldots, a_{n}}) ],
\end{align*}
which is again an element of $\CH_{ng}(A^{n} \x A^{n}) \ox \Q$.
\end{itemize}
The above morphisms are defined for all $0 \le n \le 2g$. Note that for any non-negative integer $n$:
\begin{align*}
\CH_{ng}(A^{n} \x A^{n}) \ox \Q & = \Mor_{\CSmProjk}(A^{n}, A^{n}) \\ 
& = \Mor_{\Chowek}(h(A^{n}), h(A^{n})),
\end{align*}
and so the above morphisms may be regarded either as elements of $\CSmProjk$ or $\Chowek$.
\end{notation}

\begin{proposition}[\tbf{Deninger and Murre}, \tbf{K{\"u}nnemann}] \label{pidecomp}
Let $A$ be a $g$-dimensional abelian variety over a perfect field $k$. Then the morphisms $p_{n}$ of Notation \ref{pislambda} give a 
decomposition\tn{:}
\begin{align*}
\sum_{n = 0}^{2g} p_{n} = [\Gamma(1_{A})] \in \CH_{g}(A \x A) \ox \Q,
\end{align*} 
and satisfy the following relations.
\begin{itemize}
\item[\tn{(i)}] For all $0 \le m \neq n \le 2g$, we have $p_{n}^{2} = p_{n}$ and $p_{n} \circ p_{m} = 0$.
\item[\tn{(ii)}] For any  homomorphism $f: A \ra B$ of abelian varieties over $k$,
\begin{align*}
\Gamma(f) \circ p_{n, A} = p_{n, B} \circ \Gamma(f).
\end{align*}
Therefore $f$ induces a morphism $\hn(f): \hn(A) \ra \hn(B)$.
\item[\tn{(iii)}] For all $0 \le n \le 2g$, we have $\pt{p_{n}} = p_{2g - n}$.
\end{itemize}
\end{proposition}
\begin{prooff}
The decomposition $\sum_{n = 0}^{2g} p_{n} = [\Gamma(1_{A})]$ as well as Relation (i) were first proved in \cite[Theorem 3.1]{DeMu}. Relation (ii) is \cite[Proposition 3.3]{DeMu}, and Relation (iii) is shown in \cite[Theorem 3.1.1]{Ku1}. Although both these sources work in the contravariant category $\tChowek$, one sees immediately that the results still hold in $\Chowek$. Indeed, since transposition is additive and $[\pt{\Gamma(1_{A})}] = [\Gamma(1_{A})]$, the decomposition is also valid in $\Chowek$, and the relations (i), (ii) and (iii) are visibly the transposes of the ones written down in \cite[Theorem 3.1.1]{Ku1}.
\end{prooff}

\begin{notation} \label{hnAnot}
Let $A$ be a $g$-dimensional abelian variety over a perfect field $k$. We denote by $\hn(A)$ the effective Chow motive
$(A, p_{n}, 0)$, for all $0 \le n \le 2g$. Note that since $p_{n}$ is a projector by Proposition \ref{pidecomp} Part (i), the diagram:
\begin{align*}
\xymatrix{
A \ar@{->}[dd]_(0.5){[\Gamma(1_{A})]} \ar@{->}[rr]^(0.5){p_{n}} & & A \ar@{->}[dd]^(0.5){p_{n}} \\
& & \\
A \ar@{->}[rr]_(0.5){p_{n}} & & A}
\end{align*}
commutes, and so $p_{n}$ defines a morphism from $h(A)$ to $\hn(A)$, which we denote by $p_{n}(A): h(A) \ra \hn(A)$.
A similar diagram shows that $p_{n}$ defines a morphism from $\hn(A)$ to $h(A)$, which we write as
$j_{n}(A): \hn(A) \ra h(A)$. Further, we use the symbol: 
\begin{align*}
\iota_{\hn(A)}: \hn(A) \lra \bos_{n}^{2g} \hn(A)
\end{align*}
for the canonical inclusion.
\end{notation}

\begin{proposition} \label{jipifuncrem}
The assignment $A \mapsto \hn(A)$ on objects naturally extends to an additive functor\tn{:}
\begin{align*}
\hn: \AbVk \lra \Chowek.
\end{align*}
The canonical embedding $j_{n}(A)$ and projection $p_{n}(A)$ of Notation \ref{hnAnot} may be regarded as morphisms of functors $j_{n}: \hn \ra h$ and $p_{n}: h \ra \hn$. Moreover, one has the additivity property\tn{:}
\begin{align*}
\ho(A \x B) \simeq \ho(A) \os \ho(B)
\end{align*}
for $n=1$.
\end{proposition} 
\begin{prooff}
Let $f: A \ra B$ be a morphism of abelian varieties. Then by Proposition \ref{pidecomp} Part (ii), the diagram:
\begin{align*}
\xymatrix{
A \ar@{->}[dd]_{p_{n, A}} \ar@{->}[rr]^{f} & & B \ar@{->}[dd]^{p_{n, B}} \\
& & \\
A \ar@{->}[rr]_{f} & & B}
\end{align*}
commutes. This implies, firstly, that $f$ can naturally be considered as a morphism from $\hn(A)$ to $\hn(B)$, which we denote by $\hn(f): \hn(A) \ra \hn(B)$. It is easy to check that $\hn$ respects compositions and preserves the identity, and hence defines a functor: 
\begin{align*}
\hn: \AbVk \lra \Chowek.
\end{align*}
Secondly, recalling that $j_{n}(A)$ is by definition $p_{n, A}$, the commutativity of this diagram also implies that the morphism $j_{n}(A): \hn(A) \ra h(A)$ of Notation \ref{hnAnot} is functorial, i.e. that $j_{n}: \hn \ra h$ is a morphism of functors. For the same reason, $p_{n}(A): h(A) \ra \hn(A)$ is also a functor. The additivity property 
$\ho(A \x B) \simeq \ho(A) \os \ho(B)$ follows from the motivic K{\"u}nneth formula discussed in the proof of 
\cite[Lemma 3.2.1]{Ku1}. 
\end{prooff}

\begin{proposition}[\tbf{Murre}] \label{hoCAlbisom} 
Let $(C, \xo)$ be a smooth, proper and connected curve over a field $k$. Then\tn{:}
\begin{align*}
\ho(\ixo): \ho(C) \lra \ho(\AlbzCk)
\end{align*}
is a functorial isomorphism in $\Chowek$. 
\end{proposition}
\begin{prooff}
\cite[Theorem 3.1]{Mur}.
\end{prooff}

\begin{proposition}[\tbf{Classical decomposition of $h(A)$}] \label{dsisomchow}
Let $A$ be a $g$-dimensional abelian variety over a perfect field $k$. Then the morphism\tn{:}
\begin{align*}
f(A) := \left( \sum_{n=0}^{2g} \iota_{\hn(A)} \circ p_{n} \right): h(A) \lra \bos_{n=0}^{2g} \hn(A)
\end{align*}
is an isomorphism of effective Chow motives. 
\end{proposition}
\begin{prooff}
As noted in the introduction, this result is classical and has been published in various forms by many authors: Grothendieck, Kleiman and Lieberman \cite[Appendix to \S 2]{GKL}, Deninger and Murre \cite{DeMu}, K{\"u}nnemann \cite{Ku1}, and Shermenev \cite{Sh1}. Let us simply note that the result follows immediately from Proposition \ref{pidecomp}, which gives a decomposition $[\Gamma(1_{A})] = \sum_{n=0}^{2g} p_{n}$ in:
\begin{align*}
\CH_{g}(A \x A) \ox \Q = \Mor_{\Chowek}(A, A)
\end{align*}
into mutually orthogonal idempotents.
\end{prooff}

\begin{proposition}[\tbf{K{\"u}nnemann}] \label{slambdarel}
Let $A$ be a $g$-dimensional abelian variety over a perfect field $k$. Then the morphisms $p_{n}$, $s_{n}$ and $\lambda_{n}$ of Notation \ref{pislambda} satisfy the following relations. 
\begin{itemize}
\item[\tn{(i)}] $s_{n}^{2} = s_{n}$, $\lambda_{n}^{2} = \lambda_{n}$ and 
$\lambda_{n} \circ s_{n} = s_{n} \circ \lambda_{n}$.
\item[\tn{(ii)}] $p_{n, A^{n}} \circ s_{n} = s_{n} \circ p_{n, A^{n}}$ and 
$p_{n, A^{n}} \circ \lambda_{n} = \lambda_{n} \circ p_{n, A^{n}}$.
\item[\tn{(iii)}] $\lambda_{n} \circ p_{n, A^{n}} = p_{1, A}^{\ox n}$.
\item[\tn{(iv)}] $s_{n} \circ h(\Delta^{n}_{A}) = h(\Delta^{n}_{A})$.
\item[\tn{(v)}] $p_{n, A^{n}} \circ h(\Delta^{n}_{A}) = h(\Delta^{n}_{A}) \circ p_{n, A}$.
\end{itemize}
\end{proposition}
\begin{prooff}
Relations (i), (ii) and (iii) are proved in \cite[Lemma 3.2.1]{Ku1} for contravariant Chow motives, but one sees immediately that our relations are the transposes of the ones given there. Relation (iv) is obvious, and Relation (v) follows by applying Proposition \ref{pidecomp}, Part (ii) to the morphism $\Delta^{n}_{A}: A \ra A^{n}$.
\end{prooff}

\begin{remark} \label{symnchowA}
Let $A$ be a $g$-dimensional abelian variety over a perfect field $k$. Observe that the morphism 
$s_{n}: A^{n} \ra A^{n}$ of Notation \ref{pislambda} agrees with the morphism $s_{X}^{n}$ of Notation 
\ref{ansn}. Now by Remark \ref{symnchow}, we have: 
\begin{align*}
\Symn(\hof(A)) = \Symn((A, p_{1, A}, 0)) = (A^{n}, s_{n} \circ p_{1, A}^{\ox n}, 0),
\end{align*}
and further that:
\begin{align*}
s_{n} = \pi^{n}_{\hof(A)}: \hof(A)^{\ox n} \lra \Symn(\hof(A)) 
\end{align*}
is the canonical projection. Moreover, as remarked in \cite[Lemma 3.2.1]{Ku1} Part (iv), one can use the first three relations of 
Proposition \ref{slambdarel} to rewrite $s_{n} \circ p_{1, A}^{\ox n}$ as 
$\lambda_{n} \circ s_{n} \circ p_{n, A^{n}}$, and hence:
\begin{align*}
\Symn(\hof(A)) = \Symn \big( (A, p_{1, A}, 0) \big) = (A^{n}, \lambda_{n} \circ s_{n} \circ p_{n, A^{n}}, 0),
\end{align*}
for all $0 \le n \le 2g$.
\end{remark} 

\begin{notation} \label{Pnnot}
Let $A$ be a $g$-dimensional abelian variety over a perfect field $k$. For an integer $0 \le n \le 2g$, we set: 
\begin{align*}
P_{n}(A) := p_{n, A^{n}} \circ \lambda_{n} \circ [\Gamma(\Delta^{n}_{A})] \in \CH_{g}(A \x A^{n}) \ox \Q.
\end{align*}
Note that this is nothing but the transpose of the Chow cycle denoted by $\Phi_{n}$ in \cite[Theorem 3.3.1]{Ku1}.
\end{notation} 

\begin{proposition}[\tbf{K{\"u}nnemann}] \label{Pnconst}
Let $A$ be a $g$-dimensional abelian variety over a perfect field $k$ and let $0 \le n \le 2g$ be an integer. The Chow cycle\tn{:}
\begin{align*}
P_{n}(A) = p_{n, A^{n}} \circ \lambda_{n} \circ [\Gamma(\Delta^{n}_{A})] \in \CH_{g}(A \x A^{n}) \ox \Q.
\end{align*}
of Notation \ref{Pnnot} defines an isomorphism $P_{n}(A): \hn(A) \ra \Symn(\ho(A))$ of effective Chow motives.
\end{proposition}
\begin{prooff}
It is shown in \cite[Theorem 3.3.1]{Ku1} that the morphism denoted there by $\Phi_{n}$, which in our notation is $\pt{P_{n}}$, is an isomorphism in the category $\tChowek$. Hence $P_{n}(A)$ is an isomorphism in $\Chowek$.
\end{prooff}

\begin{proposition} \label{kunnvssteve}
Let $A$ be a $g$-dimensional abelian variety over a perfect field $k$. Then the diagram\tn{:}
\begin{align*}
\xymatrix{
h(A) \ar@{->}[rr]^(0.45){h(\Delta^{n}_{A})} \ar@{->}[dd]_(0.5){p_{n}} & & 
h(A)^{\ox n} \ar@{->}[rr]^(0.47){p_{1}(A)^{\ox n}} & & \ho(A)^{\ox n} \ar@{->}[dd]^(0.5){s_{n}} \\
& & & & \\
\hn(A) \ar@{->}[rrrr]_(0.45){P_{n}(A)} & & & & \Symn(\ho(A))}
\end{align*}
of morphisms in $\Chowek$ is commutative for all integers $0 \le n \le 2g$. 
\end{proposition}
\begin{prooff}
All steps in the calculation below follow either from relations (i) -- (v) of Proposition \ref{slambdarel}, or the idempotency relation:
\begin{align}
p_{n, A^{n}}^{2} = p_{n, A^{n}} \label{piidemptwo}
\end{align}
shown in Proposition \ref{pidecomp}. Recalling that $p_{1}(A) = p_{1, A}$, we have along the upper-right path:
\begin{align*}
s_{n} \circ p_{1, A}^{\ox n} \circ h(\Delta^{n}_{A}) 
& \lgeqt{(iii)} s_{n} \circ \lambda_{n} \circ p_{n, A^{n}} \circ h(\Delta^{n}_{A}) \\
& \lgeqt{(i)} \lambda_{n} \circ s_{n} \circ p_{n, A^{n}} \circ h(\Delta^{n}_{A}) \\
& \lgeqt{(ii)} \lambda_{n} \circ p_{n, A^{n}} \circ s_{n} \circ h(\Delta^{n}_{A}) \\
& \lgeqt{(iv)} \lambda_{n} \circ p_{n, A^{n}} \circ h(\Delta^{n}_{A}).
\end{align*}
Along the lower-left path:
\begin{align*}
P_{n}(A) \circ p_{n, A} & = p_{n, A^{n}} \circ \lambda_{n} \circ h(\Delta^{n}_{A}) \circ p_{n, A} \\
& \lgeqt{(v)} p_{n, A^{n}} \circ \lambda_{n} \circ p_{n, A^{n}} \circ h(\Delta^{n}_{A}) \\
& \lgeqt{(ii)} p_{n, A^{n}} \circ p_{n, A^{n}} \circ \lambda_{n} \circ h(\Delta^{n}_{A}) \\
& \lgeqt{(\ref{piidemptwo})} p_{n, A^{n}} \circ \lambda_{n} \circ h(\Delta^{n}_{A}) \\
& \lgeqt{(ii)} \lambda_{n} \circ p_{n, A^{n}} \circ h(\Delta^{n}_{A}).
\end{align*}
This completes the proof.
\end{prooff}

\begin{notation} \label{chowpsinot}
Let $A$ be a $g$-dimensional abelian variety over a perfect field $k$, and write: 
\begin{align*}
i_{\ho(A)}^{n}: \Symn(\ho(A)) \ra \bos_{n = 0}^{2g} \Symn(\ho(A))
\end{align*}
for the canonical inclusion. We define $\psi_{A}^{n}: h(A) \ra \Symn(\ho(A))$ to be the composite morphism:
\begin{align*}
\psi_{A}^{n}: h(A) \sxra{h(\Delta^{n}_{A})} h(A)^{\ox n} & \sxra{p_{1}(A)^{\ox n}} \ho(A)^{\ox n} \slra{s_{n}} \Symn(\ho(A)), 
\end{align*}
and denote the sum over $0 \le n \le 2g$ of all composite morphisms:
\begin{align*}
h(A) \sxra{\psi_{A}^{n}} \Symn(\ho(A)) \sxra{i_{\ho(A)}^{n}} \bos_{n = 0}^{2g} \Symn(\ho(A))
\end{align*}
by $\psi_{A}: h(A) \ra \bos_{n = 0}^{2g} \Symn(\ho(A))$.
\end{notation}

\begin{proposition}[\tbf{Classical decomposition of $h(A)$, constructive}] \label{kunnchowisom} 
Let $A$ be a $g$-dimensional abelian variety over a perfect field $k$. Then the morphism\tn{:}
\begin{align*}
\psi_{A}: h(A) \lra \bos_{n = 0}^{2g} \Symn(\ho(A)) 
\end{align*}
of Notation \ref{chowpsinot} is an isomorphism of effective Chow motives.
\end{proposition}
\begin{prooff}
As noted in the proof of Proposition \ref{dsisomchow}, this result has been published in various forms by many authors. The proof below is simply a paraphrasing of this work using our notation. By Proposition \ref{kunnvssteve}, we have: 
\begin{align}
\nonumber \psi^{n}_{A} & = i_{\ho(A)}^{n} \circ s_{n} \circ p_{1}(A)^{\ox n} \circ h(\Delta^{n}_{A}) \\
& = i_{\ho(A)}^{n} \circ P_{n}(A) \circ p_{n}. \label{psirelone}
\end{align}
Let $0 \le n \le 2g$, set $P(A) := \bos_{n=0}^{2g} P_{n}(A)$, and consider the diagram:
\begin{align*}
\xymatrix{
\hn(A) \ar@{->}[rrr]^(0.45){P_{n}(A)} \ar@{->}[dd]_(0.5){\iota_{\hn(A)}} & & &
\Symn(\ho(A)) \ar@{->}[dd]^(0.5){\iota_{\Symn(\ho(A))}} \\
& & & \\
\bos_{n=0}^{2g} \hn(A) \ar@{->}[rrr]_(0.45){P(A)} & & & \bos_{n=0}^{2g} \Symn(\ho(A)),}
\end{align*}
in which the maps $\iota_{\hn(A)}$ and $\iota_{\Symn(\ho(A))}$ are the canonical inclusions. By definition, the direct sum $P(A)$ is the unique morphism such that the diagram commutes. Moreover $P(A)$ is an isomorphism, since $P_{n}(A)$ is an isomorphism, by Proposition \ref{Pnconst}. By the relation (\ref{psirelone}) and the commutativity of the above diagram, we see:
\begin{align*}
\psi_{A} & = \sum_{n = 0}^{2g} \Big( \iota_{\Symn(\ho(A))} \circ P_{n}(A) \circ p_{n} \Big) 
= \sum_{n = 0}^{2g} \Big( P(A) \circ \iota_{\hn(A)} \circ p_{n} \Big) \\
& = P(A) \circ \sum_{n = 0}^{2g} \Big( \iota_{\hn(A)} \circ p_{n} \Big) 
= P(A) \circ f(A),
\end{align*}
where the last equality follows by definition of the isomorphism: 
\begin{align*}
f(A): h(A) \ra \os \hn(A)
\end{align*}
of Proposition \ref{dsisomchow}. Both $P(A)$ and $f(A)$ are isomorphisms, and so it follows that $\psi_{A}$ is an isomorphism.
\end{prooff}

\section{The Voevodsky motive of a curve}

\begin{convention} \label{etlQconv}
In keeping with Chapter 1, all motivic complexes will be considered in the category $\DMeetkQ$ of {\etl} motives with $\Q$-coefficients. 
\end{convention}

\begin{remark} 
We review some notation, definitions and results from previous chapters which will be used below.\\\\
\tbf{(i)} From Chapter 1, we recall the functors: 
\begin{itemize}
\item[\tn{(a)}] $h: \SmProjk \ra \Chowek$ of Definition \ref{hfunc}, 
\item[\tn{(b)}] $\Phi: \Chowek \ra \DMeetkQ$ of Notation \ref{Phifuncnot}, and 
\item[\tn{(c)}] $M: \Smk \ra \DMeetkQ$ of Notation \ref{Mhinot} and Proposition \ref{hinishi}. 
\end{itemize}
By Proposition \ref{phisymcomm} $\Phi$ is an additive tensor functor, and there is a natural isomorphism of functors 
$F: \Phi \circ h \ra M$.\\\\ 
\tbf{(ii)} From Chapter 2, we have in Notation \ref{mofnot} the object $\Mof(G)$ in $\DMeetkQ$ attached to any homotopy invariant commutative $k$-group scheme $G$. In particular, $\Mof(A)$ is defined for any abelian variety $A$ over $k$, since abelian varieties (even semiabelian varieties) are homotopy invariant, by Proposition \ref{savhinv}. Thus we may regard $\Mof$ as a functor:
\begin{align*}
\Mof: \AbVk \lra \DMeetkQ.
\end{align*} 
\noindent \tbf{(iii)} From Proposition \ref{hzalbC} in Chapter 3, we have attached to any pointed curve $C$ over a field $k$ an isomorphism: 
\begin{align*}
b_{C}: \hzet(M(C)) \ra \wAlbCk
\end{align*} 
of {\etl} sheaves with transfers.   
\end{remark} 

\begin{proposition} \label{calcHik}
Let $C$ be a smooth, proper and connected curve over a field $k$. Then we have\tn{:}
\begin{align*}
\hiet((\Phi \circ h^{0})(C)) & = \left\{
\begin{array}{ll}
\Zcs, & \tn{if } i = 0, \\
0, & \tn{if } i \neq 0,
\end{array} \right.   \qaq \\
\hiet((\Phi \circ h^{2})(C)) & = \left\{
\begin{array}{ll}
\wtGm, & \tn{if } i = 1, \\
0, & \tn{if } i \neq 1.
\end{array} \right. 
\end{align*}
\end{proposition}
\begin{prooff}
By Proposition \ref{hdecomp} part (iii) we have a decomposition $h(\Pjo) = \onec \os \Lf$ of effective Chow motives. Since $\Phi \circ h$ and $M$ are isomorphic functors and $\Phi$ is additive, this yields:
\begin{align*}
\onev \os \onev(1)[2] = M(\Pjo) = (\Phi \circ h)(\Pjo) = \Phi(\onec) \os \Phi(\Lf),
\end{align*}
where we use temporarily subscripts to distinguish between $\onec = h(\Speck)$ in the category $\Chowek$ and $\onev = M(\Speck)$ in the category $\DMeetkQ$. Since $\Phi$ is also monoidal, the identity object for the tensor product is preserved: $\Phi(\onec) = \onev$, which implies that $\Phi(\Lf) = \onev(1)[2]$. On the other hand, by Proposition \ref{hdecomp} we have $h^{0}(C) = \onec$ and $h^{2}(C) = \Lf$, and so:
\begin{itemize}
\item[\tn{(i)}] $(\Phi \circ h^{0})(C) = \onev$ and 
\item[\tn{(ii)}] $(\Phi \circ h^{2})(C) = \onev(1)[2]$.
\end{itemize}
From (i), it now follows that:
\begin{align*}
\hiet((\Phi \circ h^{0})(C)) = \hiet(\onev) = 
\left\{ \begin{array}{ll}
\Zcs, & \tn{if } i = 0, \\
0, & \tn{if } i \neq 0,
\end{array} \right.
\end{align*}
as required. To verify the computation for $\hiet((\Phi \circ h^{2})(C))$, observe that since $\onev(1)[1]$ is isomorphic to $\Gm$, by 
\cite[Theorem 4.1]{MVW}, we have:
\begin{align*}
\hiet((\Phi \circ h^{2})(C)) = \hiet(\onev(1)[2]) = \hiet(\Gm[1]) = \left\{
\begin{array}{ll}
\wtGm, & \tn{if } i = 1, \\
0, & \tn{if } i \neq 1,
\end{array} \right. 
\end{align*}
where the last equality follows since $\Gm[1]$ is a complex concentrated in degree minus one. 
\end{prooff}

\begin{corollary} \label{hojonekern}
Let $(C, \xo)$ be a smooth connected pointed curve over an algebraically closed field $k$, and denote by 
$\pi_{C}: C \ra \Speck$ the structure morphism. Then\tn{:}
\begin{align*}
0 \lra \hzet((\Phi \circ \ho)(C)) & \sxra{\hzet((\Phi \circ j^{1})(C))} 
\hzet((\Phi \circ h)(C)) \\
& \sxra{\hzet((\Phi \circ h)(\pi_{C}))} \Zcs \lra 0
\end{align*}
is a short exact sequence of abelian {\etl} sheaves with transfers.  
\end{corollary}
\begin{prooff}
By Proposition \ref{imkernp}, the kernel and image of the idempotent morphism: 
\begin{align*}
p := \xo \circ \pi_{C}: h(C) \lra h(C)
\end{align*}
of effective Chow motives are given by the morphisms in the sequence:
\begin{align}
\ho(C) \os h^{2}(C) \sxra{j^{1}_{C} \os j^{2}_{C}} h(C) \sxra{h(\pi_{C})} \onec. \label{jonejtwo}
\end{align}
We now apply the (additive) composition of functors $\hzet \circ \Phi$ to the above sequence. Since $\hzet((\Phi \circ h^{2}))(C) = 0$ by Proposition \ref{calcHik}, this kills the direct summand $h^{2}(C)$ in (\ref{jonejtwo}), and we obtain the required short exact sequence:
\begin{align*}
0 \lra \hzet((\Phi \circ \ho)(C)) & \sxra{\hzet((\Phi \circ j^{1})(C))} 
\hzet((\Phi \circ h)(C)) \\
& \sxra{\hzet((\Phi \circ h)(\pi_{C}))} \Zcs \lra 0
\end{align*}
in the abelian category $\STEkQ$.
\end{prooff}

\begin{construction} \label{boconst}
Let $(C, \xo)$ be a smooth connected pointed curve over an algebraically closed field $k$. By Corollary \ref{betagamma}, the sheaf $\wAlbzCk$ associated to $\AlbzCk$ is the kernel of the degree map $p_{\Zcs}: \wAlbCk \ra \Zcs$ of Notation \ref{pZcsnot}. We therefore have a short exact sequence:
\begin{align*}
0 \lra \wAlbzCk \slra{i_{C}} \wAlbCk \slra{\delta_{C}} \Zcs \lra 0
\end{align*}
of abelian {\etl} sheaves with transfers. We construct a morphism:
\begin{align*}
p_{C}: \AlbCk \lra \AlbzCk
\end{align*}
splitting $i_{C}$, as follows. Denote by $\pi_{C}: C \ra \Speck$ the structure morphism. Then $\delta_{C}$ is the same as:
\begin{align*}
\wAlb_{\pi_{C}/k}: \wAlbCk \lra \wAlb_{\Speck/k} = \Zcs,
\end{align*}
which is split by the morphism $\wAlb_{\xo/k}: \Zcs \lra \wAlbCk$. This splitting induces a direct sum decomposition 
$\wAlbCk = \wAlbzCk \os \Zcs$ and hence also the desired morphism $p_{C}: \wAlbCk \ra \wAlbzCk$ splitting $i_{C}$.
\end{construction} 

\begin{notation} \label{bcoac}
Let $(C, \xo)$ be a smooth connected pointed curve over an algebraically closed field $k$. We define the composite morphism:
\begin{align*}
\bCo: \hzet(M(C)) \slra{\bC} \wAlbCk \slra{p_{C}} \wAlbzCk,
\end{align*}
where $\bC$ is the morphism from Construction \ref{bXmap} and $p_{C}$ is the morphism from Construction \ref{boconst}. We define also:
\begin{align*}
a_{C}: \wAlbzCk \slra{i_{C}} \wAlbCk \slra{\bC^{-1}} \hzet(M(C)),  
\end{align*}
and observe that $\bCo \circ a_{C}$ is the identity on $\AlbzCk$ by definition. 
\end{notation} 

\begin{proposition} \label{calbisom}
Let $(C, \xo)$ be a smooth projective pointed curve over an algebraically closed field $k$. Then the composite morphism\tn{:}
\begin{align*}
\hzet((\Phi \circ \ho)(C)) & \sxra{\hzet(\Phi \circ j_{1})(C)} \hzet((\Phi \circ h)(C)) \\
& \sxra{\hzet(F(C))} \hzet(M(C)) \sxra{\bCo} \wAlbzCk
\end{align*}
is an isomorphism of {\etl} sheaves with transfers. 
\end{proposition}
\begin{prooff}
By Corollary \ref{betagamma}, $\wAlbzCk$ is the kernel of the morphism $p_{\Zcs}: \AlbCk \ra \Zcs$ of Notation \ref{pZcsnot}. Thus the bottom row of the diagram:
\begin{align*}
\xymatrix{
0 \ar@{->}[rr] & & \wAlbzCk \ar@{=}[dd] \ar@{->}[rr]^(0.45){a_{C}} & & 
\hzet(M(C)) \ar@{->}[dd]^{\bC} \ar@{->}[rr]^(0.6){\hzet(M(\pi_{C}))} & &
\Zcs \ar@{=}[dd] \ar@{->}[rr] & & 0 \\
& & & & & & & & \\
0 \ar@{->}[rr] & & \wAlbzCk \ar@{->}[rr]_(0.5){i_{C}} & & 
\wAlbCk \ar@{->}[rr]_(0.55){p_{\Zcs}} & & \Zcs \ar@{->}[rr] & & 0}
\end{align*}
is a short exact sequence of {\etl} sheaves with transfers, and by definition of the morphism $a_{C}$ in Notation \ref{bcoac}, the left-hand square commutes. By Corollary \ref{bXZcscomm}, the right-hand square is also commutative. Hence the whole diagram commutes, and so the top row is also a short exact sequence. Combining this with Corollary \ref{hojonekern}, it follows that both rows in the diagram:
\begin{align*}
\xymatrix{
0 \ar@{->}[r] & 
\hzet((\Phi \circ \ho)(C)) \ar@{->}[rrr]^(0.5){\hzet(\Phi(j^{1}(C)))} \ar@{-->}[dd]^{\simeq}_{f} & & &
\hzet((\Phi \circ h)(C)) \ar@{->}[rrr]^(0.6){\hzet((\Phi \circ \ho)(\pi_{C})) } 
\ar@{->}[dd]^{\hzet(F(C))}_{\simeq} & & & 
\Zcs \ar@{->}[r] \ar@{=}[dd] & 0 \\
& & & & & & & & \\
0 \ar@{->}[r] & \wAlbzCk \ar@{->}[rrr]_(0.45){a_{C}} \ar@<0.1ex>@/^1.0pc/@{<-}[rrr]^(0.45){\bCo} & & & 
\hzet(M(C)) \ar@{->}[rrr]_(0.6){\hzet(M(\pi_{C}))} & & &
\Zcs \ar@{->}[r] & 0 }
\end{align*}
are short exact sequences. Since $F: \Phi \circ h \ra M$ is an isomorphism of functors, the vertical map in the centre is an isomorphism and the right-hand square of the diagram is commutative, inducing an isomorphism: 
\begin{align*}
f: \hzet((\Phi \circ \ho)(C)) \lra \wAlbzCk 
\end{align*}
on the kernels such that the left-hand square also commutes, as indicated. But as remarked in Notation \ref{bcoac}, the morphism $\bCo$ splits $a_{C}$, and so it follows that the isomorphism $f$ is equal to the composition of morphisms
$\bCo \circ \hzet(F(C)) \circ \hzet((\Phi \circ j^{1})(C))$, completing the proof. 
\end{prooff}

\section{The morphism $D(A): (\Phi \circ \ho)(A) \ra \Mof(A)$}

\begin{notation} \label{GAdef}
Let $A$ be an abelian variety over $k$. We define the morphism $D(A): (\Phi \circ \ho)(A) \ra \Mof(A)$ in $\DMeetkQ$ to be the composition:
\begin{align*}
D(A): (\Phi \circ \ho)(A) \sxra{\Phi(j_{1}(A))} (\Phi \circ h)(A) \sxra{F(A)} M(A) \slra{\alpha_{A}} \Mof(A).
\end{align*}
We recall that $j_{1}(A): \ho(A) \ra h(A)$ is the direct summand embedding of Notation \ref{hnAnot}, the symbol $F$ denotes the natural isomorphism $F: \Phi \circ h \ra M$ of functors $\SmProjk \ra \DMeetkQ$ of Proposition \ref{phisymcomm}, and  the morphism $\alpha_{A}: M(A) \ra \Mof(A)$ is defined in Notation \ref{alphanot}.
\end{notation}

\begin{proposition} \label{GAfunc}
The morphism $D(A): (\Phi \circ \ho)(A) \ra \Mof(A)$ of Notation \ref{GAdef} is functorial in the following sense. If 
$f: A \ra B$ is a morphism of abelian varieties, then the diagram\tn{:} 
\begin{align*}
\xymatrix{
(\Phi \circ \ho)(A) \ar@{->}[dd]_{(\Phi \circ \ho)(f)} \ar@{->}[rr]^{D(A)} & & 
\Mof(A) \ar@{->}[dd]^{\Mof(f)} \\
& & \\
(\Phi \circ \ho)(B) \ar@{->}[rr]_{D(B)} & & \Mof(B)}
\end{align*}
commutes. 
\end{proposition}
\begin{prooff}
The morphism $D(A)$ is the composition of the three morphisms (i) $\Phi(j_{1}(A))$, (ii) $F(A)$ and (iii) $\alpha_{A}$. The first of these is functorial since $\Phi$ is a functor and $j_{i}(A)$ is functorial, by Proposition \ref{jipifuncrem}, the second because $F$ is a natural transformation, and the third by Proposition \ref{alphafunc}.
\end{prooff}

\begin{definition}
Recall that a homomorphism $f: A \ra B$ of abelian varieties is an \emph{isogeny} if it is surjective, and its kernel is a finite group scheme. 
\end{definition}

\begin{remark} \label{isogrem}
Since tensoration with $\Q$ kills torsion, all morphisms in the $\Q$-linear category $\left( \AbVk \right) \ox \Q$ (see Definition \ref{sxdef}) coming from isogenies in $\AbVk$ are isomorphisms.
\end{remark}

\begin{lemma}[\tbf{Classical result}] \label{albofac}
Let $A$ be an abelian variety over an infinite field $k$. Then there exists a smooth, proper curve $C$ over $k$ and a second abelian variety $B$ such that $\AlbzCk$ is isogenous to the product $A \x B$.
\end{lemma}
\begin{prooff}
By \cite[Theorem 10.1]{Mi2}, there exists a smooth proper curve $C$ over $k$ such that $A$ is a quotient of $\Jac(C)$, where $\Jac(C)$ is the Jacobian variety of $C$. By \cite[\S 6]{Mi2}, the Jacobian is canonically isomorphic to $\AlbzCk$, and so we have a surjective homomorphism: 
\begin{align*}
q: \AlbzCk \lra A
\end{align*}
of abelian varieties. Since $\left( \AbVk \right) \ox \Q$ is an abelian category, the morphism $q$ has a kernel $B := \ker(q) \sseq \AlbzCk$. Since $\left( \AbVk \right) \ox \Q$ is moreover a semisimple category by the Poincar{\'e}-Weil reducibility theorem (see \cite[Theorem 12.1]{Mi4}), there is a product decomposition $A^{\p} \x B = \AlbzCk$, where $A$ is isogenous to $A^{\p}$. Hence $A \x B$ is isogenous to $\AlbzCk$.
\end{prooff}

\begin{proposition} \label{albored}
Suppose that for every smooth, proper curve $C$ over an algebraically closed field $k$, the morphism\tn{:}
\begin{align*}
D(\AlbzCk): (\Phi \circ \hof)(\AlbzCk) \lra \Mof(\AlbzCk) 
\end{align*}
of motives is an isomorphism. Then $D(A)$ is an isomorphism for all abelian varieties $A$ over $k$. 
\end{proposition}
\begin{prooff}
Let $A$ be an abelian variety over $k$. Since $k$ is algebraically closed it is infinite, so we may invoke Lemma \ref{albofac} and choose a smooth proper curve $C$, a second abelian variety $B$ over $k$ and an isogeny:
\begin{align*}
f: A \x B \lra \AlbzCk.
\end{align*}
By Remark \ref{isogrem}, $f$ is an isomorphism in the category $\left( \AbVk \right) \ox \Q$. Now, the functor $\Mof: \AbVk \ra \DMeetkQ$ plainly factors via $\left( \AbVk \right) \ox \Q$. We have then two functors: 
\begin{align*}
& \Mof: \left( \AbVk \right) \ox \Q \lra \DMegk \sx \Q \\
& \Phi: \Chowek \lra \DMeetkQ
\end{align*}
which are both additive, by Notation \ref{etaGHdef} and Proposition \ref{phisymcomm}, respectively. Since moreover the functor $\ho: \left( \AbVk \right) \ox \Q \ra \Chowek \ox \Q$ is additive by Proposition \ref{jipifuncrem}, the composite functor 
$\Phi \circ \ho$ is also additive.  This implies that the left-most morphisms in both rows of the commutative diagram:
\begin{align*}
\xymatrix{
\Mof(A)  \ar@{<-}[dd]_{D(A) \os D(B)} \os \Mof(B) \ar@{->}[r]^(0.52){\simeq} & 
\Mof(A \x B) \ar@{->}[rr]_(0.48){\simeq}^(0.48){\Mof(f)} & &
\Mof(\AlbzCk) \ar@{<-}[dd]_{\simeq}^{D(\AlbzCk)} \\
& & & & \\
(\Phi \circ \ho)(A) \os (\Phi \circ \ho)(B) \ar@{->}[r]_(0.57){\simeq} & 
(\Phi \circ \ho)(A \x_{k} B) \ar@{->}[rr]_(0.5){(\Phi \circ \ho)(f)} & &
(\Phi \circ \ho)(\AlbzCk).}
\end{align*}
are isomorphisms, as indicated. Moreover, the isogeny $f: A \x B \ra \AlbzCk$ becomes an isomorphism in $\left( \AbVk \right) \ox \Q$, and so $\Mof(f)$ and $(\Phi \circ \ho)(f)$ are also isomorphisms. Since $D(\AlbzCk) = D(A \x B)$ is an isomorphism, it follows that $D(A) \os D(B)$ is an isomorphism, and hence that $D(A)$ and $D(B)$ are isomorphisms. 
\end{prooff}

\begin{proposition} \label{GAisom}
Let $A$ be an abelian variety over a perfect field $k$. Then the morphism\tn{:}
\begin{align*}
D(A): (\Phi \circ \ho)(A) \lra \Mof(A)
\end{align*}
of Notation \ref{GAdef} is an isomorphism in $\DMeetkQ$. 
\end{proposition}
\begin{prooff}
By the {\etl} reduction result Proposition \ref{etldes}, we may assume that $k$ is an algebraically closed field, and by Proposition \ref{albored}, it suffices to prove the result for the {\SeAV} $\AlbzCk$ of a smooth proper curve $C$ over $k$. Since $k$ is algebraically closed, we may fix a $k$-rational point $\xo$ of $C$. The objects $\Mof(A)$ and $(\Phi \circ \ho)(A)$ of $\DMeetkQ$ are both complexes concentrated in degree zero, the former by definition and the latter because $\ho(A) = (A, p_{1}, 0)$ is an effective Chow motive, which is sent to a complex concentrated in degree zero by definition of the functor $\Phi$ in Notation \ref{Phifuncnot}. Recalling from Remark \ref{hznismof} that 
$\hzet(\Mof(\AlbzCk)) = \wAlbzCk$, it therefore suffices to prove that: 
\begin{align*}
\hzet(D(\AlbzCk)): \hzet((\Phi \circ \ho)(\AlbzCk)) \lra \hzet(\Mof(\AlbzCk)) = \wAlbzCk
\end{align*}
is an isomorphism in the category $\STEkQ$. Consider now the diagrams:
\begin{align*}
\tn{(a)} \quad \xymatrix{
\hzet((\Phi \circ \ho)(C)) \ar@{->}[rrrr]^(0.43){\hzet((\Phi \circ j^{1})(C))} 
\ar@{->}[dd]_(0.5){\hzet((\Phi \circ \ho)(\ixo))} & & & &
\hzet((\Phi \circ h)(C)) \ar@{->}[dd]_(0.5){\hzet((\Phi \circ h)(\ixo))} \\
& & & & \\
\hzet((\Phi \circ \ho)(\AlbzCk)) \ar@{->}[rrrr]_(0.5){\hzet((\Phi \circ j^{1})(\AlbzCk))} & & & & 
\hzet((\Phi \circ h)(\AlbzCk))} 
\end{align*}
and
\begin{align*}
\tn{(b)} \quad \xymatrix{
\hzet((\Phi \circ h)(C)) \ar@{->}[rrr]^(0.55){\hzet(F(C))} \ar@{->}[dd]_(0.5){\hzet((\Phi \circ h)(\ixo))} & & & 
\hzet(M(C)) \ar@{->}[dd]^(0.5){\hzet(M(\ixo))} \\
& & & \\
\hzet((\Phi \circ h)(\AlbzCk)) \ar@{->}[rrr]_(0.55){\hzet(F(\AlbzCk))} & & & \hzet(M(\AlbzCk))}
\end{align*}
and
\begin{align*}
\tn{(c)} \qquad \qquad \qquad \qquad \xymatrix{
\hzet(M(C)) \ar@{->}[ddrrr]^(0.5){\bC^{\circ}} \ar@{->}[dd]_(0.5){\hzet(M(\ixo))} & & & \\
& & & \\
\hzet(M(\AlbzCk)) \ar@{->}[rrr]_(0.55){\hzet(\alpha_{\AlbzCk})} & & & \wAlbzCk.}
\end{align*}
Diagram (a) is commutative, by the functoriality of $\ixo$ from Proposition \ref{hoCAlbisom}. Diagram (b) is commutative because $F: \Phi \circ h \ra M$ is a morphism of functors, and Diagram (c) is commutative by definition of the natural transformation $b_{C}$ of Construction \ref{bXmap}. Thus the large diagram obtained by glueing Diagrams (a), (b) and (c) together, from left to right, is commutative. Moreover, the composition:
\begin{align*}
\bCo \circ \hzet(F(C)) \circ \hzet(\Phi \circ j_{1})(C): \hzet((\Phi \circ \ho)(C)) \lra \AlbzCk
\end{align*}
along the top of this large diagram is an isomorphism, by Proposition \ref{calbisom}, as is the left-most edge:
\begin{align*}
\hzet((\Phi \circ h)(\ixo)): \hzet((\Phi \circ \ho)(C)) \lra \hzet((\Phi \circ \ho)(\AlbzCk))
\end{align*}
by Proposition \ref{hoCAlbisom}. The commutativity of the diagram now implies that the composition: 
\begin{align*}
\hzet(D(A)) = \hzet(\alpha_{\AlbzCk}) \circ \hzet(F(\AlbzCk)) \circ \hzet((\Phi \circ j^{1})(\AlbzCk))
\end{align*}
along the bottom is also an isomorphism. This completes the proof. 
\end{prooff}

\section{The isomorphism $\vp_{A}: M(A) \ra \Sym(\Mof(A))$}

\begin{lemma} \label{GAsquarecomm}
Let $A$ be an abelian variety over a perfect field $k$. Then the diagram\tn{:}
\begin{align*}
\xymatrix{
(\Phi \circ h)(A) \ar@{->}[dd]_(0.5){F(A)} \ar@{->}[rr]^(0.5){\Phi(p_{1}(A))} & & 
(\Phi \circ \ho)(A) \ar@{->}[dd]^(0.5){D(A)} \\
& & \\
M(A) \ar@{->}[rr]_(0.5){\alpha_{A}} & & \Mof(A)} 
\end{align*}
is commutative. Here, $p_{1}(A): h(A) \ra \ho(A)$ is the morphism defined in Notation \ref{hnAnot}.
\end{lemma}
\begin{prooff}
This follows directly from the definition of $D(A)$ as the composite morphism: 
\begin{align*}
D(A): (\Phi \circ \ho)(A) \sxra{\Phi(j_{1}(A))} (\Phi \circ h)(A) \sxra{F(A)} M(A) \slra{\alpha_{A}} \Mof(A),
\end{align*}
of Notation \ref{GAdef}, together with the fact that $p_{1}(A)$ splits $j_{1}(A)$. 
\end{prooff}

\begin{remark} \label{Afinkimdim}
Let $A$ be a $g$-dimensional abelian variety over a perfect field $k$, and let $n$ be an integer larger than $2g$. By \cite[Theorem 3.1.1]{Ku1}, we have $\Symn(h^{1}(A)) = 0$ and by Proposition \ref{GAisom}, $(\Phi \circ \ho)(A)$ is isomorphic to $\Mof(A)$. Since $\Phi$ commutes with $\Symn$, we see:
\begin{align*}
\Symn(\Mof(A)) \simeq \Symn((\Phi \circ \ho)(A)) \simeq \Phi \left( \Symn(\ho(A)) \right) = 0.
\end{align*}
Thus $\Mof(A)$ has finite Kimura dimension (see Definition \ref{finitekimdim}).
\end{remark}

\begin{proposition} \label{vpAisom}
Let $A$ be an abelian variety over a perfect field $k$. Then the morphism\tn{:}
\begin{align*}
\vp_{A}: M(A) \lra \Symn(\Mof(A))
\end{align*}
of Notation \ref{vpGnot} is an isomorphism in the category $\DMeetkQ$.
\end{proposition}
\begin{prooff}
Recall that the morphism $\vp_{A}: M(A) \ra \Sym(\Mof(A))$ of Notation \ref{vpGnot} is defined to be the sum 
$\vp_{A} := \sum_{n} i^{n}_{\Mof(A)} \circ \vp_{A}^{n}$, where $\vp_{A}^{n}$ is the composite morphism
$\pi_{A}^{n} \circ \alpha_{A}^{\ox n} \circ M(\Delta_{A}^{n})$, and:
\begin{align*}
i^{n}_{\Mof(A)}: \Symn(\Mof(A)) \lra \bos_{n} \Symn(\Mof(A))
\end{align*}
is the canonical inclusion. Note that $\vp_{A}$ is well-defined, since $\Mof(A)$ has finite Kimura dimension, by Remark \ref{Afinkimdim}. Moreover, the direct sum:
\begin{align*}
\Sym(\Mof(A)) = \bos_{n \ge 0} \Symn(\Mof(A))
\end{align*}
is finite, and is hence a categorical product. Now, observe that he composition of the morphisms along the bottom row in the diagram:
\begin{align*}
\xymatrix{
(\Phi \circ h)(A) \ar@{->}[dd]_(0.5){F(A)} \ar@{->}[rr]_(0.46){(\Phi \circ h)(\Delta_{A}^{n})} 
\ar@/^1.5pc/@{->}[rrrrr]^(0.5){\Phi(\psi_{A}^{n})} & & 
(\Phi \circ h)(A)^{\ox n} \ar@{->}[dd]_(0.5){F(A)^{\ox n}} \ar@{->}[rr]_(0.5){\Phi(p_{1}(A))^{\ox n}} 
& & (\Phi \circ \ho)(A)^{\ox n} \ar@{->}[dd]_(0.5){D(A)^{\ox n}} \ar@{->}[r]_(0.45){\Phi(s_{n})} & 
\Symn((\Phi \circ \ho)(A)) \ar@{->}[dd]_(0.5){\Symn(D(A))} \\
& & & & & \\
M(A) \ar@{->}[rr]^(0.46){M(\Delta_{A}^{n})} \ar@/_1.5pc/@{->}[rrrrr]_(0.5){\vp_{A}^{n}}
& & M(A)^{\ox n} \ar@{->}[rr]^(0.5){\alpha_{A}^{\ox n}} & & 
\Mof(A)^{\ox n} \ar@{->}[r]^(0.45){\pi^{n}_{A}} & \Symn(\Mof(A))} 
\end{align*}
is equal to $\vp_{A}^{n}$ by definition, and the composition along the top row is equal to $\Phi(\psi_{A}^{n})$, by the definition of the morphism $\psi_{A}^{n}$ in Notation \ref{chowpsinot}. The left-hand square commutes by the fact that $F: \Phi \circ h \ra M$ is a morphism of tensor functors. By Remark \ref{symnchowA}, the canonical projection: 
\begin{align*}
\pi^{n}_{\hof(A)}: \hof(A)^{\ox n} \lra \Symn(\hof(A))
\end{align*}
in $\Chowek$ is equal to $s_{n}$. Applying the monoidal functor:
\begin{align*}
\Phi: \Chowek \ra \DMeetkQ
\end{align*}
to $s_{n}$, and recalling that $\Symn$ commutes with $\Phi$ by Corollary \ref{phisymcomm}, it follows that $\Phi(s_{n})$ is equal to the canonical projection:
\begin{align*}
\pi^{n}_{(\Phi \circ \hof)(A)}: (\Phi \circ \hof)(A)^{\ox n} \lra \Symn((\Phi \circ \hof)(A)).
\end{align*}
Therefore the right-hand square commutes, by functoriality of $\pi^{n}_{X}: X^{\ox n} \ra \Symn(X)$. The middle square is the $n^{\tn{th}}$ tensor power of the commutative diagram of Lemma \ref{GAsquarecomm}, and hence is itself commutative. Hence the whole diagram commutes. Since $F: \Phi \circ h \ra M$ is an isomorphism of functors, the two left-most vertical maps are isomorphisms. Since $D(A)$ is an isomorphism by Proposition \ref{GAisom}, the two right-most vertical maps are also isomorphisms. Furthermore, setting $\Sigma(A) := \os_{n} \Symn(D(A))$ to ease the notation, the diagram:
\begin{align*}
\xymatrix{
\Symn((\Phi \circ \ho)(A)) \ar@{->}[dd]_(0.5){\Symn(D(A))} \ar@{->}[rr]^(0.45){i_{(\Phi \circ \ho(A))}^{n}} & & 
\bos_{n} \Symn((\Phi \circ \ho)(A)) \ar@{->}[dd]^(0.5){\Sigma(A)} \\
& & \\
\Symn(\Mof(A)) \ar@{->}[rr]_(0.45){i_{\Mof(A)}^{n}} & & \bos_{n} \Symn(\Mof(A))}
\end{align*}
commutes by the universal property of the direct sum, and the vertical maps are isomorphisms, because $D(A)$ is an isomorphism. Glueing together the above diagrams produces a larger diagram consisting of four squares, whose commutativity is expressed by the relation:
\begin{align*}
i^{n}_{\Mof(A)} \circ \vp_{A}^{n} \circ F(A) = \Sigma(A) \circ i_{(\Phi \circ \ho(A))}^{n} \circ \Phi(\psi_{A}^{n}).
\end{align*}
Using bilinearity of composition of morphisms and the fact that neither $F(A)$ nor $\Sigma(A)$ depends on $n$, it follows from this relation that: 
\begin{align*}
\vp_{A} \circ F(A) & = \left( \sum_{n} i^{n}_{\Mof(A)} \circ \vp_{A}^{n} \right) \circ F(A) 
= \sum_{n} \left( i^{n}_{\Mof(A)} \circ \vp_{A}^{n} \circ F(A) \right) \\
& = \sum_{n} \left( \Sigma(A) \circ i_{(\Phi \circ \ho(A))}^{n} \circ \Phi(\psi_{A}^{n}) \right) 
= \Sigma(A) \circ \sum_{n} \left( i_{(\Phi \circ \ho(A))}^{n} \circ \Phi(\psi_{A}^{n}) \right) \\
& = \Sigma(A) \circ \Phi(\psi_{A}).
\end{align*}
We have already seen that $\Sigma(A)$ and $F(A)$ are isomorphisms, and it was shown in Proposition \ref{kunnchowisom} that $\psi_{A}$ is also an isomorphism. It follows that $\vp_{A}$ is an isomorphism. 
\end{prooff}
\chapter{The motive of a semiabelian variety}

Let $k$ be a perfect field, and let $G$ be a semiabelian variety over $k$ of rank $r$, by which we mean that $G$ is an extension of an abelian variety $A$ by a rank $r$ torus $T$. The main result of this chapter, and of the thesis as a whole, is Theorem \ref{mainresult}, in which it is proved that if $G$ has rank one, then there exists an isomorphism:
\begin{align}
M(G) \simeq \Sym(\Mof(G)) \label{isorankone}
\end{align}
in the category $\DMeetkQ$. The author was unable to prove, even in the rank one case, that the morphism 
$\vp_{G}: M(G) \ra \Sym(\Mof(G))$ of Notation \ref{vpGnot} is an isomorphism. However, for a general semiabelian variety $G$ of arbitrary rank $r$ we have the following two partial results. 
\begin{itemize}
\item[\tn{(a)}] It is shown in Proposition \ref{trivextcase} in the first section that $\vp_{G}$ is an isomorphism if $G$ is the trivial extension $A \x T$ of $A$ by $T$. Indeed, this a direct consequence of the isomorphy of $\vp_{T}$ and $\vp_{A}$ and multiplicativity of $\vp_{G}$, all of which were shown in previous chapters. 
\item[\tn{(b)}] Secondly, we have the following conditional result. Under the (unproven!) assumption that a certain square:
\begin{align} 
\xymatrix@C=14pt{   
M(H)(1)[2] \ar@{->}[dd]_(0.5){\vp_{H}(1)[2]} \ar@{->}[rr]^(0.52){\pd_{\Aff(G), H} \circ M(s_{0})} & & 
M(G)[1] \ar@{->}[dd]^(0.5){\vp_{G}[1]} \\
& & \\
\Sym(\Mof(H))(1)[2] \ar@{->}[rr]_(0.53){-f^{\bl}[1]} & & \Sym(\Mof(G))[1]} \label{techdiag}
\end{align} 
in the category $\DMeetkQ$ commutes, it is shown that for $G$ of arbitrary rank $r$, the morphism $\vp_{G}$ is an isomorphism. Let us emphasize however that since the author is unable to prove that (\ref{techdiag}) commutes, the only unconditional result we have is the much weaker statement (\ref{isorankone}), that there exists an isomorphism in the rank one case. 
\end{itemize}
The square (\ref{techdiag}) arises as the rightmost portion of a larger diagram: 
\begin{align} 
\xymatrix@C=14pt{   
M(G) \ar@{->}[rr]^(0.5){M(g)} \ar@{->}[dd]_(0.5){\vp_{G}} & & 
M(H) \ar@{->}[rr]^(0.45){1_{M(H)} \cup c_{1}(\Aff(G))} \ar@{->}[dd]^(0.5){\vp_{H}} & & 
M(H)(1)[2] \ar@{->}[dd]^(0.5){\vp_{H}(1)[2]} \ar@{->}[rr]^(0.52){\pd_{\Aff(G), H} \circ M(s_{0})} & & 
M(G)[1] \ar@{->}[dd]^(0.5){\vp_{G}[1]} \\
& & & & & & \\
\Sym(\Mof(G)) \ar@{->}[rr]_(0.49){g^{\bl}} & & 
\Sym(\Mof(H)) \ar@{->}[rr]_(0.43){h^{\bl}} & & 
\Sym(\Mof(H))(1)[2] \ar@{->}[rr]_(0.53){-f^{\bl}[1]} & & \Sym(\Mof(G))[1]} \label{pseudomorphtri}
\end{align} 
in the category $\DMeetkQ$, the rows of which are both exact triangles. A brief explanation of the construction of the diagram (\ref{pseudomorphtri}), will allow us both to appreciate why the unknown commutativity of the square (\ref{techdiag}) is the obstruction to $\vp_{G}$ being an isomorphism, and will also provide a convenient framework to outline the structure of this chapter. 

\tbf{(i)} In the first section, we show that $\vp_{G}$ is an isomorphism in the case where $G$ is the trivial extension 
$A \x T$ of $A$ by $T$.

\tbf{(ii)} The object $H$ in Diagram (\ref{pseudomorphtri}) is a semiabelian variety of rank $r-1$, sitting in a short exact sequence: 
\begin{align} 
1 \lra \Gm \slra{f} G \slra{g} H \lra 0. \label{GmGHses}
\end{align} 
That such a short exact sequence always exists is shown in the second section. 

\tbf{(iii)} In the third section, we construct, given a short exact sequence: 
\begin{align*} 
0 \lra \shfF \slra{f} \shfG \slra{g} \shfH \lra 0
\end{align*} 
in a general $\Q$-linear tensor abelian category $\cata$, an exact triangle:
\begin{align*}
\Sym(\shfH) \ox \shfF & \sxra{f^{\bl}} \Sym(\shfG) \sxra{g^{\bl}} \Sym(\shfH) \sxra{h^{\bl}} \shfF \ox \Sym(\shfH)[1]
\end{align*}
in the derived category $\Dm(\cata)$.

\tbf{(iv)} In the fourth section, we specialise this to the case where $\cata = \STEkQ$ and where the short exact sequence is (\ref{GmGHses}), to obtain an exact triangle:
\begin{align}
\nonumber \Sym(\Mof(H))(1)[1] \sxra{f^{\bl}} \Sym(\Mof(G)) & \sxra{g^{\bl}} \Sym(\Mof(H)) \\
& \sxra{h^{\bl}} \Sym(\Mof(H))(1)[2] \label{bottomtri}
\end{align}
in the derived category $\Dm(\STEkQ)$, and hence also in its localisation $\DMeetkQ$. Note that the bottom row of (\ref{pseudomorphtri}) is the exact triangle obtained by left-shifting (\ref{bottomtri}).

\tbf{(v)} The fifth section recalls the construction of the \emph{Euler triangle}:
\begin{align*} 
M(G) \sxra{M(g)} M(H) \sxra{1_{M(H)} \cup  c_{1}(\Aff(G))} M(H)(1)[2] \sxra{\pd_{\Aff(G), H} \circ M(s_{0})} M(G)[1]
\end{align*} 
in the category $\DMeetkQ$, which forms the top row of (\ref{pseudomorphtri}). This construction was given initially by to Huber/Kahn 
\cite[\S 7]{HK} in the line bundle case, and was subsequently written down for arbitrary vector bundles by D{\'e}glise \cite[\S 1.2]{Deg1}. For a general vector bundle $p: E \ra X$ in the category of $k$-schemes, the Euler triangle is the special case of the Gysin triangle construction associated to the closed immersion given by the zero section. In our case, the vector bundle in question is the associated line bundle $\cj{g}: \Aff(G) \ra H$ obtained by viewing the short exact sequence (\ref{GmGHses}) as a $\Gm$-bundle over $H$. 

\tbf{(vi)} We begin the sixth section by recalling, given a smooth $k$-scheme $X$, the definition of the canonical isomorphism:
\begin{align*} 
c_{1} := c_{1}^{X}: \Pic(X) \lra \Mor_{\DMeetkQ}(M(X), \tatt),
\end{align*} 
via which one obtains the first Chern class $c_{1}(\shfL) \in \Mor_{\DMeetkQ}(M(X), \tatt)$ of a line bundle $\shfL$ on $X$. The case of interest to us is the first Chern class: 
\begin{align*} 
c_{1}(\Aff(G)) \in \Mor_{\DMeetkQ}(M(H), \tatt),
\end{align*} 
attached to the associated line bundle $\cj{g}: \Aff(G) \ra H$ of the $\Gm$-bundle $G$ over $H$ defined by the short exact sequence (\ref{GmGHses}). The remainder of the section is devoted to proving that $c_{1}(\Aff(G))$ satisfies certain commutative diagrams, working towards a proof that the middle square in (\ref{pseudomorphtri}) is commutative.

\tbf{(vii)} In the seventh section, we use what we know about the diagram (\ref{pseudomorphtri}) to prove the conditional result labelled (b) above, as well as the existence of the isomorphism (\ref{isorankone}). Using the results of the first six sections, we show that the left-hand and middle squares commute, and since the rank of $H$ is one less than the rank of $G$, we can use induction to assume that $\vp_{H}$ and $\vp_{H}(1)[2]$ are isomorphisms (in the base  case $r=1$, then $H$ has rank zero and is thus an abelian variety, so $\vp_{H}$ is an isomorphism by the main result of Chapter 4). If the right-hand square of Diagram (\ref{pseudomorphtri}) also commutes, then $\vp_{G}$ is also an isomorphism, by an elementary property of triangulated categories. In the base case of the induction, where $H$ is an abelian variety and $G$ a rank one extension of $H$, isomorphy of $\vp_{H}$ and $\vp_{H}(1)[2]$ allows us to conclude that there exists an isomorphism (\ref{isorankone}), not necessarily $\vp_{G}$, even if we do not know that the third square commutes.

\tbf{(viii)} The final section is an appendix, in which it is shown that the symmetric power sheaf $\Symn(\shfF)$ is equal to the sheafification of the presheaf $\Symnp(\shfF)$ defined by:
\begin{align*}
\Symnp(\shfF)(U) := \im \left( \frac{1}{n!} \sum_{\sigma \in \Sym(n)} \sigma_{\shfF(U)}: 
\shfF(U)^{\ox n} \lra \shfF(U)^{\ox n} \right),
\end{align*} 
which by definition is the $n^{\tn{th}}$ symmetric power $\Symn(\shfF(U))$ of the $\Q$-vector space $\shfF(U)$. This observation is useful because it allows us to use element notation to define morphisms of sheaves of $\Q$-vector spaces with transfers of the form 
$\vp: \Symn(\shfF) \ra \shfG$ or $\vp: \shfF \ra \Symn(\shfG)$ in a meaningful way.

\section{The product of a torus and an abelian variety}

\begin{proposition} \label{trivextcase}
Let $k$ be a perfect field, let $A$ and $T$ be an abelian variety and a torus over $k$, respectively, and let $G$ be the product $A \x T$. Then the morphism\tn{:}
\begin{align*}
\vp_{G}: M(G) \lra \Sym(\Mof(G))
\end{align*}
of Notation \ref{vpGnot} is an isomorphism in the category $\DMeetkQ$.
\end{proposition}
\begin{prooff}
Since $G = A \x T$, it follows from Proposition \ref{vpmultprop} that $\vp_{G}$ is an isomorphism in $\DMeetkQ$ if and only if $\vp_{A}$ and $\vp_{T}$ are isomorphisms. But this was shown in Theorem \ref{vpAisom} and Theorem \ref{vptorusisom}, respectively.
\end{prooff}

\section{An auxiliary short exact sequence of \\ semiabelian varieties}

\begin{definition} \label{savrankdef}
Let $G$ be a semiabelian variety over a perfect field $k$. Then the \emph{rank} of $G$ is defined to be the rank of the torus $T$ in the short exact sequence: 
\begin{align*}
1 \lra T \lra G \lra A \lra 0
\end{align*}
of {\etl} sheaves in $\SE(k, \Z)$ defining $G$. We adopt the convention that the rank of the trivial torus $T = \Speck$ is zero. Thus $G$ is an abelian variety if and only if it has rank zero.
\end{definition}

\begin{construction} \label{secondsesconst}
Let $G$ be a semiabelian variety of rank $r$ over an algebraically closed field $k = \cj{k}$, sitting in a short exact sequence 
$1 \lra \Gmr \slra{j} G \slra{p} A \lra 0$ in the category $\SE(k, \Z)$ of sheaves on the big {\etl} site over $\Speck$, where $A$ is an abelian variety. We construct a second short exact sequence: 
\begin{align}
1 \lra \Gm \slra{f} G \slra{g} H \lra 0, \label{GmGHsesdef}
\end{align}
in $\SE(k, \Z)$, as follows, by defining $f$ to be the composite of the morphism $j$ with the embedding $\Gmkb \ra \Gmkb^{r}$ into the first factor, and letting $g: G \ra H$ be the cokernel of $f$. Observe that since $\Q$-vector spaces are flat modules, the natural functor 
$\SE(k , \Z) \ra \SE(k, \Q)$ is exact, and so we can also view (\ref{GmGHsesdef}) as a short exact sequence in the category $\STEkQ$. 
\end{construction}

\begin{lemma} \label{secondseslemm}
Let $G$ be a semiabelian variety of rank $r$ over an algebraically closed field $k = \cj{k}$, sitting in a short exact sequence 
$1 \ra \Gmr \sra{j} G \sra{p} A \ra 0$. Then the commutative group scheme $H$ in the associated short exact sequence\tn{:}
\begin{align*}
1 \lra \Gm \slra{f} G \slra{g} H \lra 0
\end{align*}
of Construction \ref{secondsesconst} is a semiabelian variety of rank $r-1$.
\end{lemma}
\begin{prooff}
As in Construction \ref{secondsesconst}, let $i_{1}: \Gm \ra \Gmr$ denote the first factor embedding morphism. We define also 
$i^{\p}: \Gm^{r-1} \ra \Gmr$ to be the embedding into the last $r-1$ factors, and denote by $f^{\p}: \Gm^{r-1} \ra H$ the composite morphism $g \circ j \circ i^{\p}$. By the Yoneda Lemma, since we are working with sheaves of $\Q$-vector spaces, we may use element notation to discuss morphisms. In particular, we have:
\begin{align*}
& \mathsmap{i_{1}}{\Gm}{\Gmr}{x}{(x, 0, \ldots, 0)} \qaq \\
& \mathsmap{i^{\p}}{ \Gm^{r-1}}{\Gmr}{(x_{1}, \ldots, x_{r-1})}{(0, x_{1}, \ldots, x_{r-1}),}
\end{align*}
(Note that although we are working with the multiplicative group $\Gm$, we use the additive notation ``$0$'' instead of the multiplicative notation ``$1$'' to keep the notation readable in the coset calculations below). Since the morphism 
$g: G \ra H$ in Construction \ref{secondsesconst} is defined to be the cokernel of $f$, we may write $G/f(\Gm)$ in place of $H$, and think of $g$ as the canonical projection $G \ra G/f(\Gm)$. Denote by $g^{\p}$ the morphism:
\begin{align*}
\mathsmap{g^{\p}}{G/f(\Gm)}{A}{y + f(\Gm)}{p(y).}
\end{align*}
This is well-defined, for if $y = f(x) \in f(\Gm)$, then: 
\begin{align*}
p(y) = (p \circ f)(x) = (p \circ j \circ i_{1})(x) = 0, 
\end{align*}
recalling that $f = j \circ i_{1}$. Consider the diagram:
\begin{align*} 
\xymatrix@C=20pt@R=20pt{
1 \ar@{->}[rr] & & \Gmr \ar@{->}[rr]^(0.5){j} \ar@{<-}[dd]_(0.5){i^{\p}} & & 
G \ar@{->}[rr]^(0.5){p} \ar@{->}[dd]^(0.5){g} & & A \ar@{=}[dd] \ar@{->}[rr] & & 0 \\
& & & & & & & & \\
1 \ar@{->}[rr] & & \Gm^{r-1} \ar@{->}[rr]_(0.42){f^{\p}} & & G/f(\Gm) \ar@{->}[rr]_(0.6){g^{\p}} & & 
A \ar@{->}[rr] & & 0.}
\end{align*}
It suffices to show the bottom row is a short exact sequence. By the definitions of $f^{\p}$ and $g^{\p}$, we have:
\begin{align*} 
\im(f^{\p}) & = \{ j(0, x_{1}, \ldots, x_{r}) + f(\Gm) \mid x_{i} \in \Gmr \} \qaq \\
\ker(g^{\p}) & = \{ y + f(\Gm) \mid y \in G \tn{ and } p(y) = 0 \}.
\end{align*}
Since $p \circ j = 0$, we have that $p(j(0, x_{1}, \ldots, x_{r}))$ vanishes for all $x_{i} \in \Gm$, and so $\im(f^{\p})$ is contained in $\ker(g^{\p})$. Suppose then that $y + f(\Gm) \in \ker(g^{\p})$, so that $p(y) = 0$. Then exactness of the top row means that $y = j(x_{1}, \ldots, x_{r})$, for $x_{i} \in \Gm$. Since:
\begin{align*} 
y = j(x_{1}, \ldots, x_{r}) = \ub{j(x_{1}, 0, \ldots, 0)}_{\in f(\Gm)} + j(0, x_{2}, \ldots, x_{r}),
\end{align*}
it follows that: 
\begin{align*} 
y + f(\Gm) = j(0, x_{2}, \ldots, x_{r}) + f(\Gm).
\end{align*}
Thus $\im(f^{\p}) = \ker(g^{\p})$, and the proof is complete.
\end{prooff}

\section{An exact triangle associated to $\Sym(\shfG)$}

\begin{notation} \label{Aquotnot}
To ease the notation we denote throughout this section, we make the following assignments. The $\Q$-linear abelian tensor category 
$\STEkQ$ will be written $\cata$, and $\catt$ will be the thick triangulated subcategory of the derived category $\Dm(\cata)$ described in Proposition \ref{quotA}. Further, $\catq$ will denote the quotient category $\Dm(\cata)/\catt$. Further, we will write:
\begin{align*}
Q: \Dm(\cata) \lra \catq
\end{align*}
for the natural quotient functor. Abusing notation, we will often use the same symbol to denote both an object in $\cata$ and the associated complex, concentrated in degree zero, in the derived category $\Dm(\cata)$. 
\end{notation}

\begin{remark} \label{symQcomp}
Recall that $\Dm(\cata)$ and $\catq$ are tensor categories, and that both the functor $Q$, as well as the ``concentrated in degree zero'' functor $\cata \ra \Dm(\cata)$ are tensor-preserving functors. As a result, we have the following observations, valid for any objects $\shfF$ and $\shfG$ in $\cata$ and any non-negative integer $n$:
\begin{itemize}
\item[\tn{(i)}] $[\shfF]_{0} \ox [\shfG]_{0}$ and $[\shfF \ox \shfG]_{0}$ are isomorphic objects of $\Dm(\cata)$ and
\item[\tn{(ii)}] $Q(\shfF \ox \shfG)$ and $Q(\shfF) \ox Q(\shfG)$ are isomorphic objects of $\catq$.
\end{itemize}
As a result, 
\begin{itemize}
\item[\tn{(iii)}] $[\Symn(\shfF)]_{0}$ and $\Symn([\shfF]_{0})$ are isomorphic objects of $\Dm(\cata)$ and
\item[\tn{(iv)}] $\Symn(Q(\shfF))$ and $Q(\Symn(\shfF))$ are isomorphic objects of $\catq$, 
\end{itemize}
by Proposition \ref{symaltcomm}.
\end{remark}

\begin{notation} \label{symninot}
Let $i$ and $n$ be non-negative integers, and let $f: \shfF \ra \shfG$ be a monomorphism in $\cata$. For $0 \le i \le n$, let us denote by $t^{n, i}_{\shfG}$ the morphism: 
\begin{align*}
t^{n, i}_{\shfG} := f^{\ox i} \ox 1_{\shfG}^{\ox (n-i)}: \shfF^{\ox i} \ox \shfG^{\ox (n-i)} \lra \shfG^{\ox n}.
\end{align*}
We now define:
\begin{align*}
\Symn_{i}(\shfG) := \left\{ 
\begin{array}{ll}
\im \left( s^{n}_{\shfG} \circ t^{n, i}_{\shfG} \right),  & \tn{for } 0 \le i \le n, \\
0, & \tn{for } i > n.
\end{array} \right.
\end{align*}
Here, $s^{n}_{\shfG}: \shfG^{\ox n} \ra \shfG^{\ox n}$ is the morphism of Notation \ref{ansn}.
\end{notation}

\begin{construction} \label{symfilt}
Let the situation be as in Notation \ref{symninot}, and suppose $i$ is an integer $0 \le i \le n$. Observe that the triangle: 
\begin{align*}
\xymatrix{
\shfF^{\ox (i+1)} \ox \shfG^{\ox (n-i-1)} \ar@{->}[rrr]^(0.53){1_{\shfF}^{\ox i} \ox f \ox 1_{\shfG}^{\ox (n-i-1)}} 
\ar@{->}[ddr]_(0.5){t^{n, i+1}_{\shfG}} & & &
\shfF^{\ox i} \ox \shfG^{\ox (n-i)} \ar@{->}[ddll]^(0.5){t^{n, i}_{\shfG}} \\
& & & \\
& \shfG^{\ox n} & &}
\end{align*}
commutes. Thus $t^{n, i+1}_{\shfG}$ factors through $t^{n, i}_{\shfG}$, implying that:
$\Symn_{i+1}(\shfG) = \im \left( s^{n}_{\shfG} \circ t^{n, i+1}_{\shfG} \right)$ is a subobject of 
$\Symn_{i}(\shfG) = \im \left( s^{n}_{\shfG} \circ t^{n, i}_{\shfG} \right)$. Let us denote the resulting morphism by:
\begin{align*}
f^{n, i}_{\shfG}: \Symn_{i+1}(\shfG) \lra \Symn_{i}(\shfG).
\end{align*}
In view of the identities $t^{n, 0}_{\shfG} = 1_{\shfG}$ and $t^{n, n}_{\shfG} = f^{\ox n}$, it follows that 
$\Symn_{0}(\shfG) = \Symn(\shfG)$ and $\Symn_{n}(\shfG) = \Symn(\shfF)$, and so we have a descending filtration:
\begin{align*}
\Symn(\shfG) = \Symn_{0}(\shfG) & \suse \Symn_{1}(\shfG) \suse \Symn_{2}(\shfG) \suse \cdots \\
& \suse \Symn_{n}(\shfG) = \Symn(\shfF) \suse \Symn_{n+1}(\shfG) = 0
\end{align*}
of the $n^{\tn{th}}$ symmetric power $\Symn(\shfG)$ of $\shfG$.
\end{construction}

\begin{notation} \label{Gnifilt}
Let the situation be as in Notation \ref{symninot}. For $0 \le i \le n$, we denote by $\shfG^{n}_{i}$ the preimage of $\Symn_{i}(\shfG)$ under the canonical projection $\pi^{n}_{\shfG}: \shfG^{\ox n} \ra \Symn(\shfG)$. If $i > n$, then we define $\shfG^{n}_{i}$ to be equal to zero. The filtration on $\Symn(\shfG)$ of Construction \ref{symfilt} therefore induces a filtration:
\begin{align*}
\shfG^{\ox n} = \shfG^{\ox n}_{0} & \suse \shfG^{\ox n}_{1} \suse \shfG^{\ox n}_{2} \suse \cdots \\
& \suse \shfG^{\ox n}_{n} = \shfF^{\ox n} \suse \shfG^{\ox n}_{n+1} = 0
\end{align*}
of the $n^{\tn{th}}$ tensor power $\shfG^{\ox n}$ of $\shfG$. We denote by 
$\wt{f^{n, i}}_{\shfG}: \shfG^{n}_{i+1} \ra \shfG^{n}_{i}$ the successive inclusion morphisms, and by 
$\pi^{n, i}_{\shfG}: \shfG^{n}_{i} \ra \Symn_{i}(\shfG)$ the restriction of $\pi^{n}_{\shfG}$ to the subobject $\shfG^{n}_{i}$ of $\shfG^{\ox n}$. It is immediate that the diagram:
\begin{align*}
\xymatrix{
\shfG^{n}_{i+1} \ar@{->}[dd]_(0.5){\pi^{n, i+1}_{\shfG}} \ar@{->}[rr]^(0.5){\wt{f^{n, i}}} & & 
\shfG^{n}_{i} \ar@{->}[dd]^(0.5){\pi^{n, i}_{\shfG}} \\
& & \\
\Symn_{i+1}(\shfG) \ar@{->}[rr]_(0.5){f^{n, i}} & & \Symn_{i}(\shfG)}
\end{align*}
is commutative.
\end{notation}

\begin{notation} \label{symndots}
Let $X$ be a $\Q$-vector space. We associate to any $n$-tuple of elements $x_{1}, \ldots, x_{n} \in X$ the element: 
\begin{align*}
x_{1} \ssbl \cdots \ssbl x_{n} := \sum_{\sigma \in \Sym(n)} x_{\sigma(1)} \ox \cdots \ox x_{\sigma(n)} \in X^{\ox n}.
\end{align*}
Thus the canonical projection morphism $\pi^{n}_{X}: X^{\ox n} \ra \Symn(X)$ of Definition 1.1.7 is given by:
\begin{align*}
\mathsmap{\pi^{n}_{X}}{X^{\ox n}}{\Symn(X)}{x_{1} \ox \cdots \ox x_{n}}
{\frac{1}{n!} (x_{1} \sbl \cdots \sbl x_{n}),}
\end{align*}
and the $x_{1} \sbl \cdots \sbl x_{n}$ generate $\Symn(X)$.
\end{notation}

\begin{remark} \label{conccatconv} 
We often define a morphism $\vp: \shfG \ra \shfH$ of presheaves of $\Q$-vector spaces with transfers using element notation: 
\begin{align}
\mathsmap{\vp}{\shfG}{\shfH}{x}{f(x),} \label{setnotvp}
\end{align}
as if it were a morphism of vector spaces. By convention, this notation means that $\vp(U)(x)$ is defined to be $f(x)$ for all sections $U$, and we tacitly assert that the family of maps $\vp(U)$ satisfies the condition to form a map of presheaves. Denote by $\shfG^{\dg}$ the sheafification of $\shfG$. We will also take the notation (\ref{setnotvp}) as a means of defining the unique morphism 
$\vp^{\dg}: \shfG^{\dg} \ra \shfH$ induced by the universal property of sheafification. We use this convention in two settings -- in both instances, $\shfF$ and $\shfG$ are sheaves in $\cata$, and $n$ is a non-negative integer: \\\\
\tbf{(i)} Since $\shfF^{\ox n}$ is the sheafification of the presheaf $\shfF^{\oxpn}: U \mapsto \shfF(U)^{\ox n}$, the above convention tells us to take the notation:
\begin{align*}
\mathsmap{\vp}{\shfF^{\ox n}}{\shfG}{x_{1} \ox \cdots \ox x_{n}}{f(x_{1}, \ldots, x_{n})}
\end{align*}
to mean the unique sheaf morphism $\vp: \shfF^{\ox n} \ra \shfG$ induced by the morphism of presheaves defined by: 
\begin{align*}
\mathsmap{\vp(U)}{\shfF(U)^{\ox n}}{\shfG(U)}{x_{1} \ox \cdots \ox x_{n}}{f(x_{1}, \ldots, x_{n}),}
\end{align*}
for all sections $U$. In a similar manner, one can use element notation to define sheaf morphisms of the form $\vp: \shfF \ra \shfG^{\ox n}$  
\\\\
\tbf{(ii)} By Proposition \ref{symnaltdescrip}, the $n^{\tn{th}}$ symmetric power $\Symn(\shfF)$ is the sheafification of the presheaf 
$\Symnp(\shfF)$ defined by:
\begin{align*}
\Symnp(\shfF)(U) := \im \left( \frac{1}{n!} \sum_{\sigma \in \Sym(n)} \sigma_{\shfF(U)}: 
\shfF(U)^{\ox n} \lra \shfF(U)^{\ox n} \right),
\end{align*} 
which by definition is the $n^{\tn{th}}$ symmetric power $\Symn(\shfF(U))$ of the $\Q$-vector space $\shfF(U)$. The above convention therefore tells us to take the notation:
\begin{align*}
\mathsmap{\vp}{\Symn(\shfF)}{\shfG}{x_{1} \sbl \cdots \sbl x_{n}}{f(x_{1}, \ldots, x_{n})}
\end{align*}
to mean the unique sheaf morphism $\vp: \Symn(\shfF) \ra \shfG$ induced by the morphism of presheaves defined by: 
\begin{align*}
\mathsmap{\vp(U)}{\Symn(\shfF(U))}{\shfG(U)}{x_{1} \sbl \cdots \sbl x_{n}}{f(x_{1}, \ldots, x_{n}),}
\end{align*}
for all sections $U$. In a similar manner, one can use element notation to define sheaf morphisms of the form $\vp: \shfF \ra \Symn(\shfG)$. 
\end{remark}

\begin{construction} \label{quotfiltmorph}
Let $i$ and $n$ be non-negative integers with $0 \le i \le n$, and let:
\begin{align*}
0 \lra \shfF \slra{f} \shfG \slra{g} \shfH \lra 0
\end{align*}
be a short exact sequence in $\cata$. We define a morphism:
\begin{align*}
q^{n, i}_{\shfG}: \Sym^{i}(\shfF) \ox \Sym^{n-i}(\shfH) \lra \Symn_{i}(\shfG)/\Symn_{i+1}(\shfG) 
\end{align*}
using element notation, as follows:
\begin{align*}
\mathsmap{q^{n, i}_{\shfG}}{\Sym^{i}(\shfF) \ox \Sym^{n-i}(\shfH)}{\Symn_{i}(\shfG)/\Symn_{i+1}(\shfG) }
{(x_{1} \sbl \cdots \sbl x_{i}) \ox (z_{i+1} \sbl \cdots \sbl z_{n})}
{\left[ f(x_{1}) \sbl \cdots \sbl f(x_{i}) \sbl y_{i+1} \sbl \cdots \sbl y_{n} \right].}
\end{align*} 
We re-emphasize, following Remark \ref{conccatconv}, that the element notation in fact describes $q^{n, i}_{\shfG}$ on an arbitrary section $U$. Here, the $x_{k}$ are elements of $\shfF$, the $z_{k}$ are elements of $\shfH$, and the $y_{k}$ are arbitrarily chosen preimages of the $z_{k}$ in $\shfG$. The definition of $q^{n, i}_{\shfG}$ does not depend on the choice of preimage $y_{k}$ of $z_{k}$, for if $g(y_{k}) = z_{k} = g(y_{k}^{\p})$, then $y_{k} - y_{k}^{\p}$ is an element $f(x_{k})$ of $\ker(g) = \im(f)$, and it follows that: 
\begin{align*}
& f(x_{1}) \ssbl \cdots \ssbl f(x_{i}) \ssbl y_{i+1} \ssbl \cdots \ssbl y_{k} \ssbl \cdots \ssbl y_{n} 
- f(x_{1}) \ssbl \cdots \ssbl f(x_{i}) \ssbl y_{i+1} \ssbl \cdots \ssbl y_{k}^{\p} \ssbl \cdots \ssbl y_{n} \\ 
& = f(x_{1}) \ssbl \cdots \ssbl f(x_{i}) \ssbl y_{i+1} \ssbl \cdots \ssbl (y_{k} - y_{k}^{\p}) \ssbl \cdots \ssbl y_{n} \\ 
& = f(x_{1}) \ssbl \cdots \ssbl f(x_{i}) \ssbl y_{i+1} \ssbl \cdots \ssbl f(x_{k}) \ssbl \cdots \ssbl y_{n} \in \Symn_{i+1}(\shfG).
\end{align*}
Thus $q^{n, i}_{\shfG}$ is well-defined on sections. From this it follows that it is functorial, hence constitutes a morphism of sheaves.
\end{construction}

\begin{proposition} \label{qniisom}
Let $0 \ra \shfF \sra{f} \shfG \sra{g} \shfH \ra 0$ be a short exact sequence in $\cata$, and let $i$ and $n$ be non-negative integers with 
$0 \le i \le n$. The morphism\tn{:}
\begin{align*}
q^{n, i}_{\shfG}: \Sym^{i}(\shfF) \ox \Sym^{n-i}(\shfH) \lra \Symn_{i}(\shfG)/\Symn_{i+1}(\shfG) 
\end{align*}
of Construction \ref{quotfiltmorph} is an isomorphism. 
\end{proposition}
\begin{prooff}
Recalling Remark \ref{conccatconv}, it suffices to prove the statement in the category of $\Q$-vector spaces. This is proved for alternating powers in \cite[Chapter XIX, Proposition 1.3]{Lang}. The proof for symmetric powers is similar.
\end{prooff}

\begin{notation} \label{sesFGH}
Let $0 \lra \shfF \slra{f} \shfG \slra{g} \shfH \lra 0$ be a short exact sequence in $\cata$, and let $i$ and $n$ be non-negative integers with $0 \le i \le n$. For all non-negative integers $n$ and $i$ with $0 \le i \le n$, let us denote by: 
\begin{align*}
g^{n, i}: \Symn_{i}(\shfG) \lra \Sym^{i}(\shfF) \ox \Sym^{n-i}(\shfH)
\end{align*}
the composition of the natural projection: 
\begin{align*}
\Symn_{i}(\shfG) \lra \Symn_{i}(\shfG)/\Symn_{i+1}(\shfG)
\end{align*}
with the inverse $(q^{n, i}_{\shfG})^{-1}$ of the isomorphism $q^{n, i}_{\shfG}$ of Construction \ref{quotfiltmorph}. We therefore have a short exact sequence:
\begin{align*}
\nonumber 0 \lra \Symn_{i+1}(\shfG) & \slra{f^{n, i}} \Symn_{i}(\shfG) 
\slra{g^{n, i}} \Sym^{i}(\shfF) \ox \Sym^{n-i}(\shfH) \lra 0 
\end{align*}
in the abelian category $\cata$.
\end{notation}

\begin{remark} \label{gniexplicit}
From the definition of the morphism $q^{n, i}_{\shfG}$ in Construction \ref{quotfiltmorph}, we see that the morphism $g^{n, i}$ of Notation \ref{sesFGH} is given by: 
\begin{align*}
\mathsmap{g^{n, i}}{\Symn_{i}(\shfG)}{\Sym^{i}(\shfF) \ox \Sym^{n-i}(\shfH)}
{f(x_{1}) \sbl \cdots \sbl f(x_{i}) \sbl y_{i+1} \sbl \cdots \sbl y_{n}}
{x_{1} \sbl \cdots \sbl x_{i} \ox g(y_{i+1}) \sbl \cdots \sbl g(y_{n}),}
\end{align*}
where the $x_{k}$ are in $\shfF$ and the $y_{k}$ are in $\shfG$. In particular, in the case $i = 0$ we have: 
\begin{align*}
\mathsmap{g^{n, 0}}{\Symn(\shfG)}{\Sym^{n}(\shfH)}{y_{1} \sbl \cdots \sbl y_{n}}
{g(y_{1}) \sbl \cdots \sbl g(y_{n}),}
\end{align*}
which is nothing but the morphism: 
\begin{align*}
\Symn(g): \Symn(\shfG) \lra \Symn(\shfH)
\end{align*}
induced by the projection $g: \shfG \ra \shfH$.
\end{remark}

\begin{remark} \label{iotadefrem}
Let $\shfG$ be an object in $\cata$, and let $n$ be a non-negative integer. Let us recall the canonical embedding 
$\iota^{n}_{\shfG}: \Symn(\shfG) \ra \shfG^{\ox n}$ of Definition \ref{symaltn}. This is nothing but the inclusion of the image 
$\Symn(\shfG)$ of the morphism $s^{n}_{\shfG}: \shfG^{\ox n} \ra \shfG^{\ox n}$ into $\shfG^{\ox n}$, and so concretely:  
\begin{align*}
\mathsmap{\iota^{n}_{\shfG}}{\Symn(\shfG)}{\shfG^{\ox n}}{y_{1} \sbl \cdots \sbl y_{n}}
{y_{1} \sbl \cdots \sbl y_{n}.}
\end{align*}
Moreover, by the universal property of the categorical image, any morphism $g: \shfG \ra \shfH$ in $\cata$ induces a diagram:
\begin{align*}
\xymatrix{
\Symn(\shfG) \ar@{->}[rr]^(0.5){\Symn(g)} \ar@{->}[dd]_(0.5){\iota^{n}_{\shfG}} & & 
\Symn(\shfH) \ar@{->}[dd]^(0.5){\iota^{n}_{\shfH}} \\
& & \\
\shfG^{\ox n} \ar@{->}[rr]_(0.5){g^{\ox n}} & & \shfH^{\ox n},}
\end{align*}
which is commutative.
\end{remark}

\begin{lemma} \label{sescommFGH}
Let $0 \lra \shfF \slra{f} \shfG \slra{g} \shfH \lra 0$ be a short exact sequence in $\cata$, and let $n$ be a non-negative integer. Then the diagram\tn{:}
\begin{align*} 
\xymatrix@C=15pt{
0 \ar@{->}[rr] & & \Symn_{1}(\shfG) \ar@{->}[rr]^(0.5){f^{n,0}} \ar@{->}[d]_(0.5){g^{n,1}} & & 
\Symn(\shfG) \ar@{->}[rr]^(0.5){g^{n}} \ar@{->}[d]^(0.5){\iota^{n}_{\shfG}} & & 
\Symn(\shfH) \ar@{->}[dd]^(0.5){\iota^{n}_{\shfH}} \ar@{->}[rr] & & 0 \\
& & \shfF \ox \Sym^{n-1}(\shfH) \ar@{->}[d]_(0.5){1_{\shfF} \ox \iota^{n-1}_{\shfH}} & & 
\shfG^{\ox n} \ar@{->}[d]^(0.5){1_{\shfG} \ox g^{\ox n-1}} & & & & \\
0 \ar@{->}[rr] & & 
\shfF \ox \shfH^{\ox n-1} \ar@{->}[rr]_(0.5){f \ox 1_{\shfH^{\ox n-1}}} & & 
\shfG \ox \shfH^{\ox n-1} \ar@{->}[rr]_(0.58){g \ox 1_{\shfH^{\ox n-1}}} & & \shfH^{\ox n} \ar@{->}[rr] & & 0}
\end{align*}
describes a morphism of short exact sequences in the category $\cata$. 
\end{lemma}
\begin{prooff}
\tbf{(i)} Noting that $g^{n} = \Symn(g)$, the commutativity of the square on the right is equivalent to the identity 
$g^{\ox n} \circ \iota^{n}_{\shfG} = \iota^{n}_{\shfH} \circ \Symn(g)$, which is nothing but the functoriality of the morphism 
$\iota^{n}_{X}: \Symn(X) \ra X^{\ox n}$ expressed by the commutative diagram in Remark \ref{iotadefrem}.

\tbf{(ii)} To see that the left-hand square commutes, consider a typical generator $f(x_{1}) \sbl y_{n} \sbl \cdots \sbl y_{n}$ of 
$\Symn_{1}(\shfG) \sseq \Symn(\shfG)$. We compute first along the the bottom left-hand path of the left-hand square. Recall first from Remark \ref{gniexplicit} that $g^{n,1}$ is given by:
\begin{align*}
\mathsmap{g^{n,1}}{\Symn_{1}(\shfG)}{\shfF \ox \Sym^{n-1}(\shfH)}{f(x_{1}) \sbl y_{2} \sbl \cdots \sbl y_{n}}
{x_{1} \ox g(y_{2}) \sbl \cdots \sbl g(y_{n}).}
\end{align*}
Along the bottom left-hand path, we therefore have:
\begin{align}
\nonumber \left( \big( f \ox 1_{\shfH^{\ox n-1}} \big) \circ \big( 1_{\shfF} \ox \iota^{n-1}_{\shfH} \big) \circ g^{n,1} \right) &
\Big( f(x_{1}) \ssbl y_{2} \ssbl \cdots \ssbl y_{n} \Big) \\
\nonumber  & = \left( \big( f \ox 1_{\shfH^{\ox n-1}} \big) \circ \big( 1_{\shfF} \ox \iota^{n-1}_{\shfH} \big) \right) 
\Big( x_{1} \ox \left[ g(y_{2}) \ssbl \cdots \ssbl g(y_{n}) \right] \Big) \\
\nonumber  & = \big( f \ox 1_{\shfH^{\ox n-1}} \big) \Big( x_{1} \ox \left[ g(y_{2}) \ssbl \cdots \ssbl g(y_{n}) \right] \Big) \\
& = f(x_{1}) \ox \left[ g(y_{2}) \ssbl \cdots \ssbl g(y_{n}) \right]. \label{bottomleftleft}
\end{align}
To compute along the the top right-hand path of the left-hand square, we begin with a small calculation. Writing $y_{1} := f(x_{1})$ for convenience, observe that if $\sigma$ is a permutation in $\Sym(n)$ with $\sigma(1) \neq 1$, then $y_{1} = y_{\sigma(k)}$ for $2 \le k \le n$, and so $g(y_{\sigma(k)}) = g(y_{1}) = g(f(x_{1})) = 0$. This observation justifies the equality (*) in the computation: 
\begin{align}
\nonumber \sum_{\sigma \in \Sym(n)} y_{\sigma(1)} \ox & g(y_{\sigma(2)}) \ox \cdots \ox g(y_{\sigma(n)}) \\
\nonumber & \lgeqt{(*)} \sum_{\substack{\sigma \in \Sym(n) \\ \sigma(1) = 1}} y_{1} \ox g(y_{\sigma(2)}) \ox \cdots \ox g(y_{\sigma(n)}) \\
\nonumber & = \sum_{\sigma \in \Sym(n-1)} f(x_{1}) \ox g(y_{\sigma(2)}) \ox \cdots \ox g(y_{\sigma(n)}) \\
\nonumber & = f(x_{1}) \ox \left( \sum_{\sigma \in \Sym(n-1)} g(y_{\sigma(2)}) \ox \cdots \ox g(y_{\sigma(n)}) \right) \\
& = f(x_{1}) \ox \left[ g(y_{2}) \ssbl \cdots \ssbl g(y_{n}) \right] \tn{;} \label{symnreleq}
\end{align}
the other equalities being plain. Along the top-right of the left-hand square we therefore have: 
\begin{align}
\nonumber
\nonumber \left( \big(1_{\shfG} \ox g^{\ox n-1} \big) \circ \iota^{n}_{\shfG} \circ f^{n, 0} \right) 
& \Big( f(x_{1}) \ssbl y_{2} \ssbl \cdots \ssbl y_{n} \Big) \\
\nonumber & = \left( \big(1_{\shfG} \ox g^{\ox n-1} \big) \circ \iota^{n}_{\shfG} \right) 
\Big( f(x_{1}) \ssbl y_{2} \ssbl \cdots \ssbl y_{n} \Big) \\
\nonumber & = \left( \big(1_{\shfG} \ox g^{\ox n-1} \big) \right) \Big( f(x_{1}) \ssbl y_{2} \ssbl \cdots \ssbl y_{n} \Big) \\
\nonumber & = \left( \big(1_{\shfG} \ox g^{\ox n-1} \big) \right) 
\left( \sum_{\sigma \in \Sym(n)} y_{\sigma(1)} \ox y_{\sigma(2)} \ox \cdots \ox y_{\sigma(n)} \right) \\
\nonumber & = \sum_{\sigma \in \Sym(n)} y_{\sigma(1)} \ox g(y_{\sigma(2)}) \ox \cdots \ox g(y_{\sigma(n)}) \\
& = f(x_{1}) \ox \left[ g(y_{2}) \ssbl \cdots \ssbl g(y_{n}) \right], \label{toprightleft}
\end{align}
where the last equality follows from Relation (\ref{symnreleq}). Comparing (\ref{bottomleftleft}) with (\ref{toprightleft}), we see that the left-hand square commutes. This completes the proof.
\end{prooff}

\begin{construction} \label{conmorphrem}
Let $E := 0 \ra \shfF \slra{f} \shfG \slra{g} \shfH \ra 0$ be a short exact sequence in the category $\cata$. We recall that the \emph{connecting morphism} $h: \shfH \ra \shfF[1]$ in the derived category $\Dm(\cata)$ associated to $E$ is the morphism given by the roof: 
\begin{align*}
\xymatrix@C=5pt{
& & \left[ \shfF \slra{f} \shfG \right] \ar@{->}[ddrr]^(0.5){[1_{\shfF}]_{-1}} \ar@{->}[ddll]_(0.5){[g]_{0}} & & \\
& & & & \\
\shfH & & & & \shfF[1],}
\end{align*}
where $[g]_{0}$ denotes the morphism which is equal to $g$ in degree zero and zero in all other degrees, and $[1_{\shfF}]_{-1}$ is defined analogously. Then the sequence of morphisms:
\begin{align*}
\shfF \slra{f} \shfG \slra{g} \shfH \slra{h} \shfF[1]
\end{align*}
in $\Dm(\cata)$ forms an exact triangle, and by \cite[Corollary 10.7.5]{Wei1} the connecting morphism is functorial in the following sense. If the diagram: 
\begin{align*} 
\xymatrix@C=15pt{
0 \ar@{->}[rr] & & \shfF \ar@{->}[rr]^(0.52){f} \ar@{->}[dd]_(0.5){\vp} & & 
\shfG \ar@{->}[rr]^(0.5){g} \ar@{->}[dd]^(0.5){\psi} & & 
\shfH \ar@{->}[dd]^(0.5){\xi} \ar@{->}[rr] & & 0 \\
& & & & & & & & \\
0 \ar@{->}[rr] & & \shfF^{\p} \ar@{->}[rr]_(0.52){f^{\p}} & & 
\shfG^{\p} \ar@{->}[rr]_(0.5){g^{\p}} & & \shfH^{\p} \ar@{->}[rr] & & 0}
\end{align*} 
is a morphism of short exact sequences, then the diagram:
\begin{align*} 
\xymatrix@C=20pt{
\shfF \ar@{->}[rr]^(0.52){f} \ar@{->}[dd]_(0.5){\vp} & & 
\shfG \ar@{->}[rr]^(0.5){g} \ar@{->}[dd]^(0.5){\psi} & & 
\shfH \ar@{->}[dd]^(0.5){\xi} \ar@{->}[rr]^(0.45){h} & & \shfF[1] \ar@{->}[dd]^(0.5){\vp[1]} \\
& & & & & & \\
\shfF^{\p} \ar@{->}[rr]_(0.52){f^{\p}} & & \shfG^{\p} \ar@{->}[rr]_(0.5){g^{\p}} & & 
\shfH^{\p} \ar@{->}[rr]_(0.45){h^{\p}} & & \shfF^{\p}[1]}
\end{align*} 
is a morphism of exact triangles. 
\end{construction}

\begin{notation} \label{Xinot}
Let $0 \lra \shfF \slra{f} \shfG \slra{g} \shfH \lra 0$ be a short exact sequence in $\cata$, and let $i$ and $n$ be integers with $0 \le i \le n$. We denote by $X^{n}_{i}$ the two-term complex:
\begin{align*}
X^{n}_{i} := \left[ \Symn_{i+1}(\shfG) \sxra{f^{n, i}} \Symn_{i}(\shfG) \right] 
\end{align*}
in $\Dm(\cata)$, concentrated in degrees minus one and zero.
\end{notation}

\begin{remark} \label{symnitriangle}
Let $0 \lra \shfF \slra{f} \shfG \slra{g} \shfH \lra 0$ be a short exact sequence in $\cata$, and let $i$ and $n$ be non-negative integers with $0 \le i \le n$. Then by Construction \ref{conmorphrem}, the short exact sequence of Notation \ref{sesFGH} in the abelian category $\cata$ induces a connection morphism:
\begin{align*}
h^{n, i}: \Sym^{i}(\shfF) \ox \Sym^{n-i}(\shfH) \lra \Symn_{i+1}(\shfG)[1]
\end{align*}
in the derived category $\Dm(\cata)$, which is represented by the roof:
\begin{align*}
\xymatrix{
& X^{n}_{i} \ar@{->}[ddrr]^(0.5){[1_{\Symn_{i+1}(\shfG)}]_{-1}} \ar@{->}[ddl]_(0.5){[g^{n, i}]_{0}} & & \\
& & & \\
\Sym^{i}(\shfF) \ox \Sym^{n-i}(\shfH) & & & \Symn_{i+1}(\shfG)[1],}
\end{align*}
where $X^{n}_{i}$ is defined in Notation \ref{Xinot}, $[1_{\Symn_{i+1}(\shfG)}]_{-1}$ denotes the morphism $X^{n}_{i} \ra \Symn_{i+1}(\shfG)[1]$ which is the identity in degree minus one and zero elsewhere, and the morphism 
$[g^{n, i}]_{0}$ is similarly defined. We therefore have an exact triangle:
\begin{align*}
\Symn_{i+1}(\shfG) \slra{f^{n, i}} \Symn_{i}(\shfG) \slra{g^{n, i}} \Sym^{i}(\shfF) \ox 
\Sym^{n-i}(\shfH) \slra{h^{n, i}} \Symn_{i+1}(\shfG)[1]
\end{align*}
in the derived category $\Dm(\cata)$, and hence also an exact triangle: 
\begin{align*}
Q(\Symn_{i+1}(\shfG)) \sxra{Q(f^{n, i})} Q(\Symn_{i}(\shfG)) & \sxra{Q(g^{n, i})} \Sym^{i}(Q(\shfF)) \ox \Sym^{n-i}(Q(\shfH)) \\
& \sxra{Q(h^{n, i})} Q(\Symn_{i+1}(\shfG))[1]
\end{align*}
in the quotient category $\catq$.
\end{remark} 

\begin{proposition} \label{gn1isom}
Let $0 \lra \shfF \slra{f} \shfG \slra{g} \shfH \lra 0$ be a short exact sequence in $\cata$, let $n$ be a positive integer, assume that the functor $Q: \cata \ra \catq$ of Notation \ref{Aquotnot} is a tensor functor, and suppose $\Sym^{n}(Q(\shfF))$ vanishes for all $n \ge 2$. 
Then $Q(\Symn_{2}(\shfG))$ also vanishes, and the morphism\tn{:}
\begin{align*}
Q(g^{n, 1}): Q(\Symn_{1}(\shfG)) \lra Q(\shfF) \ox \Sym^{n-1}(Q(\shfH)), 
\end{align*}
coming from the morphism $g^{n, 1}$ of Notation \ref{sesFGH}, is an isomorphism. 
\end{proposition}
\begin{prooff}
Consider in the quotient category $\catq$ the exact triangle:
\begin{align}
\nonumber
Q(\Symn_{i+1}(\shfG)) \sxra{Q(f^{n, i})} Q(\Symn_{i}(\shfG)) & \sxra{Q(g^{n, i})} \Sym^{i}(Q(\shfF)) \ox \Sym^{n-i}(Q(\shfH)) \\
& \sxra{Q(h^{n, i})} Q(\Symn_{i+1}(\shfG))[1] \label{extriFGHtwo}
\end{align}
from Remark \ref{symnitriangle}, defined for all integers $i$ and $n$ satisfying $0 \le i \le n$. If $n = i = 1$, then 
$Q(\Symn_{i+1}(\shfG)) =  Q(\Sym^{1}_{2}(\shfG)) = 0$, and the second and third terms are both equal to $Q(\shfF)$. The exact triangle therefore reads:
\begin{align*}
0 \sxra{Q(f^{1, 1})} Q(\shfF) \sxra{Q(g^{1, 1})} Q(\shfF) \sxra{Q(h^{1, 1})} 0,
\end{align*}
implying that $Q(g^{1, 1})$ is an isomorphism (indeed, it is the identity on $\shfF$, by Remark \ref{gniexplicit}). Suppose now that 
$2 \le i \le n$. Then $\Sym^{i}(Q(\shfF)) = 0$, by assumption, and so the third term in the exact triangle (\ref{extriFGHtwo}) vanishes. The triangle therefore reads:
\begin{align*}
Q(\Symn_{i+1}(\shfG)) \sxra{Q(f^{n, i})} Q(\Symn_{i}(\shfG)) \sxra{Q(g^{n, i})} 0 \sxra{Q(h^{n, i})} Q(\Symn_{i+1}(\shfG))[1],
\end{align*}
implying that $Q(f^{n, i})$ is an isomorphism for all $2 \le i \le n$. This gives a chain of isomorphisms:
\begin{align*}
\Symn(Q(\shfF)) = Q(\Symn_{n}(\shfG)) & \sxra{Q(f^{n, n-1})} Q(\Symn_{n-1}(\shfG)) \sxra{Q(f^{n, n-2})} \cdots \\
\cdots & \sxra{Q(f^{n, 2})} Q(\Symn_{2}(\shfG)).
\end{align*}
But since $n$ is at least two, $\Symn(Q(\shfF))$ vanishes by assumption, and so $Q(\Symn_{2}(\shfF))$ also vanishes. It follows that for 
$n \ge 2$ and $i = 1$, the exact triangle (\ref{extriFGHtwo}) looks like:
\begin{align*}
0 \sxra{Q(f^{n, 1})} Q(\Symn_{1}(\shfG)) & \sxra{Q(g^{n, 1})} Q(\shfF) \ox \Sym^{n-1}(Q(\shfH)) \sxra{Q(h^{n, 1})} 0.
\end{align*}
Therefore $Q(g^{n, 1})$ is an isomorphism.
\end{prooff}

\begin{notation} \label{fngnhn} 
Let $0 \lra \shfF \slra{f} \shfG \slra{g} \shfH \lra 0$ be a short exact sequence in $\cata$, let $n$ be a non-negative integer, and suppose that $\Sym^{n}(Q(\shfF)) = 0$, for all $n \ge 2$. To aid in defining the following notation, observe that in the case $i = 0$, the exact triangle in $\catq$ of Remark \ref{symnitriangle} looks like:
\begin{align}
Q(\Symn_{1}(\shfG)) \sxra{Q(f^{n, 0})} \Symn(Q(\shfG)) & \sxra{Q(g^{n, 0})} \Sym^{n}(Q(\shfH)) 
\sxra{Q(h^{n, 0})} Q(\Symn_{1}(\shfG))[1]. \label{extrizero}
\end{align}
\tbf{(i)} We denote by $f^{n}$ the composite morphism:
\begin{align*}
f^{n}: Q(\shfF) \ox \Sym^{n-1}(Q(\shfH)) \sxra{(Q(g^{n, 1}))^{-1}} Q(\Symn_{1}(\shfG)) \sxra{f^{n, 0}} \Symn(Q(\shfG)),
\end{align*}
where $Q(g^{n, 1})$ is the isomorphism of Proposition \ref{gn1isom}. 
\tbf{(ii)} To simplify the notation, let us denote the morphism $Q(g^{n, 0})$ by:
\begin{align*}
g^{n}: \Symn(Q(\shfG)) \lra \Sym^{n}(Q(\shfH)).
\end{align*}
\tbf{(iii)} We denote by $h^{n}$ the composite morphism:
\begin{align*}
h^{n}: \Sym^{n}(Q(\shfH)) \sxra{Q(h^{n, 0})} Q(\Symn_{1}(\shfG))[1] \sxra{Q(g^{n, 1})[1]} Q(\shfF) \ox \Sym^{n-1}(Q(\shfH))[1].
\end{align*}
Note that $f^{n}$, $g^{n}$ and $h^{n}$ are all morphisms in the category $\catq$, and further that $h^{1} = Q(h^{1, 0}) = Q(h)$, since $g^{1,1}$ is the identity $\shfF \ra \shfF$. 
\end{notation}

\begin{notation} \label{finaltriangle}
Let $0 \lra \shfF \slra{f} \shfG \slra{g} \shfH \lra 0$ be a short exact sequence in $\cata$, let $n$ be a non-negative integer, and suppose that $\Sym^{n}(Q(\shfF)) = 0$, for all $n \ge 2$. Since the triangle (\ref{extrizero}) in Notation \ref{fngnhn} is exact, and $Q(g^{n, 1})$ is an isomorphism by Proposition \ref{gn1isom}, it follows that:
\begin{align*}
\Sym^{n-1}(Q(\shfH)) \ox Q(\shfF) & \sxra{f^{n}} \Symn(Q(\shfG)) \sxra{g^{n}} \Sym^{n}(Q(\shfH)) 
\sxra{h^{n}} \shfF \ox \Sym^{n-1}(Q(\shfH))[1]
\end{align*}
is also an exact triangle in $\catq$. Since exactness is preserved by direct summation, the triangle: 
\begin{align*}
\Sym(Q(\shfH)) \ox Q(\shfF) & \sxra{f^{\bl}} \Sym(Q(\shfG)) \sxra{g^{\bl}} \Symn(Q(\shfH)) \sxra{h^{\bl}} 
Q(\shfF) \ox \Sym(Q(\shfH))[1]
\end{align*}
is also exact. Here, $f^{\bl}$ denotes the direct sum of $f^{n}$ over all non-negative integers $n$, and the morphisms $g^{\bl}$ and $h^{\bl}$ are defined similarly.
\end{notation}

\begin{proposition} \label{midsquarebottom}
Let $0 \lra \shfF \slra{f} \shfG \slra{g} \shfH \lra 0$ be a short exact sequence in $\cata$, let $n$ be a positive integer, and suppose that 
$\Sym^{n}(Q(\shfF)) = 0$, for all $n \ge 2$. Then the diagram\tn{: }
\begin{align*}
\xymatrix{
\Symn(Q(\shfH)) \ar@{->}[rr]^(0.42){h^{n}} \ar@{->}[dd]_(0.5){\iota^{n}_{Q(\shfH)}} & & Q(\shfF) \ox \Sym^{n-1}(Q(\shfH))[1]
\ar@{->}[dd]^(0.5){\iota^{n-1}_{Q(\shfH)} \ox 1_{Q(\shfF)[1]}} \\
& & \\
Q(\shfH)^{\ox n} \ar@{->}[rr]_(0.4){h^{1} \ox 1_{Q(\shfH)^{\ox n-1}}} & & Q(\shfF) \ox Q(\shfH)^{\ox n-1}[1]}
\end{align*}
commutes in the category $\catq$. Here, $h^{1}$ and $h^{n}$ are the morphisms of Notation \ref{fngnhn}.
\end{proposition}
\begin{prooff}
By Lemma \ref{sescommFGH}, the diagram\tn{:}
\begin{align*} 
\xymatrix@C=15pt{
0 \ar@{->}[rr] & & \Symn_{1}(\shfG) \ar@{->}[rr]^(0.5){f^{n,0}} \ar@{->}[d]_(0.5){g^{n,1}} & & 
\Symn(\shfG) \ar@{->}[rr]^(0.5){g^{n}} \ar@{->}[d]^(0.5){\iota^{n}_{\shfG}} & & 
\Symn(\shfH) \ar@{->}[dd]^(0.5){\iota^{n}_{\shfH}} \ar@{->}[rr] & & 0 \\
& & \shfF \ox \Sym^{n-1}(\shfH) \ar@{->}[d]_(0.5){1_{\shfF} \ox \iota^{n-1}_{\shfH}} & & 
\shfG^{\ox n} \ar@{->}[d]^(0.5){1_{\shfG} \ox g^{\ox n-1}} & & & & \\
0 \ar@{->}[rr] & & 
\shfF \ox \shfH^{\ox n-1} \ar@{->}[rr]_(0.5){f \ox 1_{\shfH^{\ox n-1}}} & & 
\shfG \ox \shfH^{\ox n-1} \ar@{->}[rr]_(0.58){g \ox 1_{\shfH^{\ox n-1}}} & & \shfH^{\ox n} \ar@{->}[rr] & & 0}
\end{align*}
commutes. The connecting homomorphisms (c.f. Construction \ref{conmorphrem}) of the top and bottom rows, respectively, are given by 
$h^{n, 0}$ and $h \ox 1_{\shfH^{\ox n-1}}$, respectively, and so it follows by the functoriality of connecting homomorphisms that the diagram: 
\begin{align*}
\xymatrix{
\Symn(\shfH) \ar@{->}[rr]^(0.5){h^{n,0}} \ar@{->}[dd]_(0.5){\iota^{n}_{\shfH}} & & \Symn_{1}(\shfG)[1]
\ar@{->}[d]^(0.5){g^{n,1}[1]} \\
& & \shfF \ox \Sym^{n-1}(\shfH)[1] \ar@{->}[d]^(0.5){1_{\shfF} \ox \iota^{n-1}_{\shfH}[1]} \\
\shfH^{\ox n} \ar@{->}[rr]_(0.4){h \ox 1_{\shfH^{\ox n-1}}} & & \shfF \ox \shfH^{\ox n-1}[1]}
\end{align*}
commutes in the category $\Dm(\cata)$. Applying the functor $Q$ to this diagram gives a diagram in the category $\catq$:
\begin{align*}
\xymatrix{
\Symn(Q(\shfH)) \ar@{->}[rr]^(0.5){Q(h^{n,0})} \ar@{->}[dd]_(0.5){\iota^{n}_{Q(\shfH)}} & & Q(\Symn_{1}(\shfG))[1]
\ar@{->}[d]^(0.5){Q(g^{n,1})[1]} \\
& & Q(\shfF) \ox \Sym^{n-1}(Q(\shfH))[1] \ar@{->}[d]^(0.5){1_{Q(\shfF)} \ox \iota^{n-1}_{Q(\shfH)}[1]} \\
Q(\shfH)^{\ox n} \ar@{->}[rr]_(0.4){Q(h) \ox 1_{Q(\shfH)^{\ox n-1}}} & & Q(\shfF) \ox Q(\shfH)^{\ox n-1}[1].}
\end{align*}
Recalling from Notation \ref{fngnhn} that $h^{n}$ is by definition the composition $Q(g^{n, 1})[1] \circ Q(h^{n,0})$, and that 
$h^{1} = Q(h)$, we see that this is the desired commutative square.
\end{prooff}

\section{An exact triangle associated to $\Sym(\Mof(G))$}

\begin{remark} 
Let $k$ be a perfect field, and let $G$ be a homotopy invariant commutative group scheme (Definition \ref{hominvGdef}). Recall from Notation \ref{mofnot} that $\Mof(G)$ denotes the complex in $\DMeetkQ$, concentrated in degree zero, consisting of the {\etl} sheaf with transfers associated to $G$. 
\end{remark}

\begin{convention} \label{Qvecconv}
Since we work exclusively we sheaves of $\Q$-vectors spaces, we employ from now on the following conventions. 
\begin{itemize}
\item[\tn{(i)}] All commutative group schemes $G$ will be regarded as representable ordinary {\etl} sheaf of $\Q$-vector spaces. 
\item[\tn{(ii)}] We will write simply $\Pic(X)$ in place of $\Pic(X) \ox \Q$ for the $\Q$-vector space associated to the Picard group. 
\end{itemize}
\end{convention} 

\begin{notation} \label{exacttrisav}
Let $k = \cj{k}$ be an algebraically closed field, let $n$ be a positive integer, and recall the following assignments from Notation \ref{Aquotnot}: 
\begin{itemize}
\item[\tn{(i)}] $\cata := \STEkQ$, so
\item[\tn{(ii)}] $\Dm(\cata) = \Dm(\STEkQ)$ and 
\item[\tn{(iii)}] $\catq := \Dm(\STEkQ)/\catt = \DMeetkQ$, where $\catt$ is the minimal thick subcategory of $\Dm(\STEkQ)$ described in 
Proposition \ref{quotA}.
\end{itemize}
Recall also the convention that an object of $\cata$ and the associated complex in $\Dm(\cata)$ concentrated in degree zero are denoted by the same symbol. Then we have $Q(\shfF) = \Mof(\shfF)$ for any {\etl} sheaf $\shfF$ with transfers. Now, by Proposition \ref{symosym1},
\begin{align*}
\Sym^{n}(\Mof(\Gm)) = 0 \: \: \tn{in the category} \: \: \DMeetkQ,
\end{align*}
for all $n \ge 2$. Recalling also that $\Mof(\Gm) = \tato$, by \cite[Theorem 4.1]{MVW}, we see that the short exact sequence $1 \lra \Gm \slra{f} G \slra{g} H \lra 0$ from Construction \ref{secondsesconst} induces, in the manner described in Notation \ref{finaltriangle}, exact triangles:
\begin{align*}
\Sym^{n-1}(\Mof(H))(1)[1] & \sxra{f^{n}} \Symn(\Mof(G)) \sxra{g^{n}} \Sym^{n}(\Mof(H)) 
\sxra{h^{n}} \Sym^{n-1}(\Mof(H))(1)[2] 
\end{align*}
and
\begin{align*}
\Sym(\Mof(H))(1)[1] \sxra{f^{\bl}} \Sym(\Mof(G)) \sxra{g^{\bl}} \Sym(\Mof(H)) \sxra{h^{\bl}} \Sym(\Mof(H))(1)[2]
\end{align*}
in the category $\DMeetkQ$.
\end{notation}

\begin{proposition} \label{hsquarecommprop}
Let $1 \lra \Gm \slra{f} G \slra{g} H \lra 0$ be the short exact sequence of semiabelian varieties of Construction \ref{secondsesconst}  
over an algebraically closed field $k = \cj{k}$, and let $n$ be a positive integer. Then the square\tn{:}
\begin{align*}
\xymatrix{
\Symn(\Mof(H)) \ar@{->}[rr]^(0.42){h^{n}} \ar@{->}[dd]_(0.5){\iota^{n}_{\Mof(H)}} & & 
\Sym^{n-1}(\Mof(H))(1)[2] \ar@{->}[dd]^(0.5){\iota^{n-1}_{\Mof(H)}(1)[2]} \\
& & \\
\Mof(H)^{\ox n} \ar@{->}[rr]_(0.4){1_{\Mof(H)^{\ox n-1}} \ox h^{1}} & & \Mof(H)^{\ox n-1}(1)[2]}
\end{align*} 
commutes. Here, $h^{1}$ and $h^{n}$ are the morphisms of Notation \ref{exacttrisav}.
\end{proposition}
\begin{prooff}
Noting that $\Sym^{n}(\Mof(\Gm)) = 0$ for all $n \ge 2$ in the category $\DMeetkQ$, by Proposition \ref{symosym1},
the result follows directly from Proposition \ref{midsquarebottom}.
\end{prooff}

\section{The Euler exact triangle}

\begin{notation}[\tbf{Ext groups, classical}] \label{extonot}
Let $i$ be a non-negative integer and let $Z$ be an object in an abelian category $\cata$ with enough injectives. Then we denote as usual by 
$\Exti_{\cata}(Z, \pul{X})$ the $i^{\tn{th}}$ right derived functor $R^{i} \Hom_{\cata}(Z, \pul{X})$ of the covariant functor 
$\Hom_{\cata}(Z, \pul{X})$. Now suppose 
\begin{align*}
\ve := 0 \ra X \ra Y \ra Z \ra 0
\end{align*}
is a short exact sequence in $\cata$. Then for any object $W$ in $\cata$, we have by \cite[\S 6.4]{Ha} a long exact sequence:
\begin{align*}
0 & \lra \Hom_{\cata}(Z, W) \lra \Hom_{\cata}(Y, W) \lra \Hom_{\cata}(X, W) \\ 
& \slra{\delta^{\ve}_{W}} \Exto_{\cata}(Z, W) \lra \Exto_{\cata}(Y, W) \lra \Exto_{\cata}(X, W) \lra \cdots.
\end{align*}
In particular, if $W := X$, we have a group homomorphism:
\begin{align*}
\delta^{\ve}_{X}: \Hom_{\cata}(X, X) \lra \Exto_{\cata}(Z, X).
\end{align*}
Denote further by $\YExt_{\cata}(Z, X)$ the set of all isomorphisms classes of short exact sequences of the form $0 \ra X \ra Y \ra Z \ra 0$. Such a short exact sequence is called an \emph{extension of $Z$ by $X$}. By \cite[Theorem 3.4.3]{Wei1}, the map:
\begin{align*}
\mathsmaptwo{\YExt_{\cata}(Z, X)}{\Exto_{\cata}(Z, X)}{\ve := \left( 0 \ra X \ra Y \ra Z \ra 0 \right)}{\delta^{\ve}_{X}(1_{X})}
\end{align*}
is an isomorphism of abelian groups, functorial in both variables. In light of this canonical identification, we will use exclusively the notation 
$\Exto_{\cata}(Z, X)$ instead of $\YExt_{\cata}(Z, X)$.
\end{notation}

\begin{notation} \label{kappanot}
As in Notation \ref{Aquotnot}, let $\cata$ denote the category $\STEkQ$, and let $X$ and $Z$ be objects in $\cata$. Given an extension: 
\begin{align*}
\ve := 0 \lra X \slra{f} Y \slra{g} Z \lra 0
\end{align*}
of $Z$ by $X$, let us define $\lambda_{Z, X}(\ve)$ to be the connection morphism $h: Z \ra X[1]$ of Construction \ref{conmorphrem} in the derived category $\Dm(\cata)$. By \cite[Corollary 10.7.5]{Wei1}, the map: 
\begin{align*}
\lambda_{Z, X}: \Exto_{\cata}(Z, X) \ra \Mor_{\Dm(\cata)}(Z, X[1])
\end{align*}
is a bijection, functorial in both arguments. In particular, $\Exto_{\cata}(Z, X)$ has a canonical abelian group structure and $\lambda_{Z, X}$ is an isomorphism of abelian groups. We denote by:
\begin{align*}
\kappa_{Z, X}: \Mor_{\Dm(\cata)}(Z, X[1]) \lra \Exto_{\cata}(Z, X)
\end{align*}
the inverse of the isomorphism $\lambda_{Z, X}$.
\end{notation}

\begin{remark}[\tbf{Functoriality of Ext groups, classical}] \label{unqextorem}
Let $\cata$ be an abelian category, let $0 \ra X \ra A \ra Z \ra 0$ be an element of $\Exto_{\cata}(Z, X)$, and let 
$f: W \ra Z$ be a morphism in $\cata$. By \cite[Corollary 20.3]{Mitch}, this induces in a canonical way an extension 
$0 \ra X \ra B \ra W \ra 0$ in $\Exto_{\cata}(W, X)$ and morphism $g: A \ra B$ such that the diagram: 
\begin{align*}
\xymatrix@C=20pt{   
0 \ar@{->}[rr] & & X \ar@{->}[rr] \ar@{=}[dd] & & B \ar@{->}[dd]^(0.5){g} \ar@{->}[rr] & & 
W \ar@{->}[rr] \ar@{->}[dd]^(0.5){f} & & 0 \\
& & & & & & & & \\
0 \ar@{->}[rr] & & X \ar@{->}[rr] & & A \ar@{->}[rr] & & Z \ar@{->}[rr] & & 0}
\end{align*}
commutes, and the right-hand square is a pullback. Thus the morphism $f: W \ra Z$ induces a map:
\begin{align*}
f^{*}: \Exto_{\cata}(Z, X) \lra \Exto_{\cata}(W, X).
\end{align*}
It is a homomorphism of abelian groups.
\end{remark} 

\begin{notation} \label{zshfFnot}
Let $\shfF \in \PSTkQ$ be a presheaf with transfers. We denote by:
\begin{align*}
z_{\shfF}: [\shfF]_{0} \lra \Cx(\shfF)
\end{align*}
The morphism in $\DMeetkQ$ which is the identity map $1_{\shfF}$ in $\PSTkQ$ in degree zero and zero in all other degrees. 
For a scheme $X$ over $k$, we write $z_{X}: L(X) \ra M(X)$ in place of $z_{L(X)}: [L(X)]_{0} \ra \Cx(L(X))$. Observe that the induced homomorphism of abelian groups:
\begin{align*}
z_{X}^{*}: \Mor_{\Dm(\STEkQ)}(M(X), \one(1)[2]) \lra \Mor_{\Dm(\STEkQ)}(L(X), \one(1)[2]) 
\end{align*}
is an isomorphism, since the right derived functor $R \Cx: \Dm(\STEkQ) \ra \DMeetkQ$ is left adjoint to the inclusion 
$\DMeetkQ \ra \Dm(\STEkQ)$. 
\end{notation}

\begin{notation} \label{aincnot}
For any object $\cly{M}$ in $\DMeetkQ$, we define $a_{\cly{M}}$ to be the natural identity morphism:
\begin{align*}
a_{\cly{M}}: \Mor_{\DMeetkQ}(\cly{M}, \tatt) \lra \Mor_{\Dm(\STEkQ)}(\cly{M}, \tatt)
\end{align*}
arising from the definition of $\DMeetkQ$ as the full subcategory of the derived category $\Dm(\STEkQ)$ with consisting of complexes with homotopy invariant cohomology sheaves. Plainly, $a_{\cly{M}}$ is functorial in $\cly{M}$. Note that for any smooth $k$-scheme $X$, the map $a_{M(X)}$ is an isomorphism, by \cite[Lemma 9.19]{MVW}.
\end{notation}

\begin{notation} \label{psivthnot}
Let $X$ be a scheme over a perfect field $k$. We denote by:
\begin{align*}
\psi_{X}: & \; \Exto_{\STEkQ}(L(X), \Gm) \lra \Hoet(X, \Gm)
\end{align*}
the canonical isomorphism of \cite[Lemma 6.23]{MVW} induced by the Yoneda isomorphism 
$\Mor_{\STEkQ}(L(X), \Gm) \ra \Gamma(X, \Gm)$. Let us write $\eta_{X}$ for the composite map:
\begin{align*}
\eta_{X}: H^{1}_{\et}(X, \Gm) & \sxra{\left( \psi_{X} \right)^{-1}} \Exto_{\STEkQ}(L(X), \Gm) \\
& \sxra{\lambda_{L(X), \Gm}} \Mor_{\Dm(\STEkQ)}(L(X), \one(1)[2]),
\end{align*} 
where $\lambda_{L(X), \Gm}$ is defined in Notation \ref{kappanot}. Note that $\eta_{X}$ is precisely the canonical identification of \cite[Proposition 3.3.1]{TMF}.
\end{notation} 

\begin{notation} \label{b1mapnot}
Let $X$ be a smooth scheme over a perfect field $k$. We denote by:
\begin{align*}
b_{1}: H^{1}_{\et}(X, \Gm) \lra \Mor_{\DMeetkQ}(M(X), \one(1)[2])
\end{align*}
the composition:
\begin{align*}
H^{1}_{\et}(X, \Gm) & \slra{\eta_{X}} \Mor_{\Dm(\STEkQ)}(L(X), \one(1)[2]) \sxra{(z_{X}^{*})^{-1}} 
\Mor_{\Dm(\STEkQ)}(M(X), \one(1)[2]) \\
& \sxra{(a_{M(X)})^{-1}} \Mor_{\DMeetkQ}(M(X), \one(1)[2]),
\end{align*} 
where $z_{X}^{*}$ is the morphism of Notation \ref{zshfFnot}, $a_{M(X)}$ that of Notation \ref{aincnot}, and $\eta_{X}$ that of 
Notation \ref{psivthnot}. Note that $b_{1}$ is an isomorphism, by \cite[Theorem 14.27]{MVW}.
\end{notation} 

\begin{notation}[\tbf{D{\'e}glise}] \label{ostnot}
Let $X$ be a smooth scheme over a perfect field $k$, let $\cly{M}$ and $\cly{N}$ be objects in $\DMeetkQ$, and let 
$f: M(X) \ra \cly{M}$ and $g: M(X) \ra \cly{N}$ be morphisms in $\DMeetkQ$. Following \cite[\S 1.2]{Deg1}, we denote by $f \cup g: M(X) \ra \cly{M} \ox \cly{N}$ the composite morphism:
\begin{align*}
f \cup g: M(X) \sxra{M(\Delta_{X})} M(X) \ox M(X) \sxra{f \ox g} \cly{M} \ox \cly{N},
\end{align*}
where $\Delta_{X}: X \ra X \x X$ denotes the diagonal morphism. Note that we use the notation ``$\cup$'' instead of D{\'e}glise's ``$\sx$'', to avoid confusion with the concept of the idempotent completion $\catc \sx \Q$ of the category $\catc \ox \Q$, defined in Notation \ref{sxdef}.
\end{notation}

\begin{construction}[\tbf{Classical, after Huber/Kahn and D{\'e}glise}] \label{eulertriconst}
Let $k$ be a perfect field, and let $p: E \ra X$ be a vector bundle of rank $n$ in the category of $k$-schemes, with zero section $s_{0}$. Denote by $E^{\x}$ the complement of the zero section, and by $j_{0}: E^{\x} \ra E$ the associated natural embedding. Via the zero section $s_{0}$, we regard $X$ as a closed subscheme of $E$. This induces a Gysin triangle, that may, by \cite[Example 1.25]{Deg1}, be written as: 
\begin{align} 
M(E^{\x}) \sxra{M(p \circ j_{0})} M(X) & \sxra{1_{M(X)} \cup c_{n}(E)} M(X)(n)[2n] 
\sxra{\pd_{E, X} \circ M(s_{0})} M(E^{\x})[1]. \label{MGHEuler}
\end{align} 
We have also use the fact, shown in \cite[Corollary 13.16]{MVW}, that the induced morphism: 
\begin{align*} 
M(p): M(E) \lra M(X)
\end{align*} 
is an isomorphism in the category $\DMeetkQ$, to replace the second term $M(E)$ with $M(X)$. Note that whereas D{\'e}glise's construction takes place in the category $\DMegk$, we shall consider this triangle as being in the category $\DMeetkQ$, by applying the composition of functors: 
\begin{align*}
\DMegk \slra{a} \DMegk \ox \Q \slra{c} \DMegk \sx \Q \slra{j} \DMeetkQ
\end{align*}
of Notation \ref{abcnot} and Construction \ref{jconst}.
\end{construction}

\begin{notation} \label{affsav}
Let $G$ be a semiabelian variety over an algebraically closed field $k = \cj{k}$, sitting in a short exact sequence 
$1 \lra \Gm \slra{f} G \slra{g} H \lra 0$ of the kind described in Construction \ref{secondsesconst}. Then we may view $G$ as a 
$\Gm$-torsor over $H$, and so we have an associated line bundle over $H$, which we write as:
\begin{align*}
\cj{g}: \Aff(G) \lra H.
\end{align*}
Let $j_{0}:G \ra \Aff(G)$ denote the natural embedding. Then the diagram:
\begin{align} 
\xymatrix{
G \ar@{->}[rr]^(0.5){j_{0}} \ar@{->}[ddrr]_(0.5){g} & & \Aff(G) \ar@{->}[dd]^(0.5){\cj{g}} \\
& & \\
& & H } \label{affGHdiag}
\end{align} 
commutes. Since the $\Gm$-torsor structure of $G$ over $H$ is determined by the extension $\cly{E}$, we will also write $\Aff(\cly{E})$ to denote the line bundle associated.
\end{notation}

\begin{notation} [\tbf{Chern classes after Huber/Kahn and D{\'e}glise}] \label{chernclassnot} 
Let $X$ be a smooth scheme over a perfect field $k$. Write: 
\begin{align*}
\vth_{X}: \Hoet(X, \Gm) \lra \Pic(X)
\end{align*}
for the canonical isomorphism of \cite[Proposition 4.9]{Mi1}, which attaches to the cohomology class of a $\Gm$-torsor $T$ over $X$ the associated line bundle $\Aff(T)$ of Notation \ref{affsav}, and recall from Convention \ref{Qvecconv} that both the source and target of $\vth_{X}$ are to be regarded as $\Q$-vector spaces. We now denote by: 
\begin{align*}
c_{1}^{X}: \Pic(X) \lra \Mor_{\DMeetkQ}(M(X), \one(1)[2])
\end{align*}
the concatenation $b_{1} \circ (\vth_{X})^{-1}$, where $b_{1}$ is the morphism of Notation \ref{b1mapnot}.
\end{notation} 

\begin{notation} \label{eulerseqsav}
Consider the line bundle $\cj{g}: \Aff(G) \ra H$ of Notation \ref{affsav}. Since the fibre is $\Ao$, we have 
$\Aff(G)^{\x} = G$. The Euler triangle (\ref{MGHEuler}) of Construction \ref{eulertriconst} for $\cj{g}: \Aff(G) \ra H$ therefore looks like: 
\begin{align*} 
M(G) \sxra{M(g)} M(H) \sxra{1_{M(H)} \cup  c_{1}(\Aff(G))} M(H)(1)[2] \sxra{\pd_{\Aff(G), H} \circ M(s_{0})} M(G)[1],
\end{align*} 
where we have used the identity $M(\cj{g} \circ j_{0}) = M(g)$ of the commutative diagram (\ref{affGHdiag}) of Notation \ref{affsav} to rewrite the first morphism.
\end{notation}

\section{Calculations with the first Chern class}

\begin{construction} \label{zetaMGmorphs} 
Let $G$ be a commutative group scheme over a perfect field $k$. We construct isomorphisms: 
\begin{align*}
\zeta_{\Mof(G)}: & \Mor_{\DMeetkQ}(\Mof(G), \tatt) \lra \Exto_{\STEkQ}(G, \Gm) \\
\zeta_{M(G)}: & \Mor_{\DMeetkQ}(M(G), \tatt) \lra \Exto_{\STEkQ}(L(G), \Gm)
\end{align*}
as follows. Recall that $\Gm$ is canonically isomorphic to $\tato$ in the category $\DMeetkQ$, by 
\cite[Theorem 4.1]{MVW}. We define $\zeta_{\Mof(G)}$ to be the composition:
\begin{align*}
\zeta_{\Mof(G)}: \Mor_{\DMeetkQ}(\Mof(G), \tatt) & \sxra{a_{\Mof(G)}} \Mor_{\Dm(\STEkQ)}(\Mof(G), \tatt) \\
& \sxra{\kappa_{G, \Gm}} \Exto_{\STEkQ}(G, \Gm) ,
\end{align*}
where $a_{\Mof(G)}$ is the isomorphism of Notation \ref{aincnot}, and $\kappa_{G, \Gm}$ is the isomorphism of Notation \ref{kappanot}. We define $\zeta_{M(G)}$ to be the composition:
\begin{align*}
\zeta_{M(G)}: \Mor_{\DMeetkQ}(M(G), \tatt) & \sxra{a_{M(G)}} \Mor_{\Dm(\STEkQ)}(M(G), \tatt) \\
& \slra{z_{G}^{*}} \Mor_{\Dm(\STEkQ)}(L(G), \tatt) \\
& \sxra{\kappa_{L(G), \Gm}} \Exto_{\STEkQ}(L(G), \Gm),
\end{align*}
where the map: 
\begin{align*}
z_{G}^{*}: \Mor_{\Dm(\STEkQ)}(M(G), \tatt) \ra \Mor_{\Dm(\STEkQ)}(L(G), \tatt)
\end{align*}
is the group homomorphism induced by the map $z_{G}: L(G) \ra M(G)$ in the category $\DMeetkQ$ defined in Notation \ref{zshfFnot}. It is an isomorphism, since the right derived functor $R \Cx: \Dm(\STEkQ) \ra \DMeetkQ$ is left adjoint to the inclusion 
$\DMeetkQ \ra \Dm(\STEkQ)$. 
\end{construction}

\begin{notation} \label{alphastarnot}
Let $G$ be a homotopy invariant commutative group scheme over a perfect field $k$. The morphism 
$\alpha_{G}: M(G) \ra \Mof(G)$ from Notation \ref{alphanot} induces a morphism: 
\begin{align*}
\mathsmap{\alpha_{G}^{*}}{\Mor_{\DMeetkQ}(\Mof(G), \tatt)}{\Mor_{\DMeetkQ}(M(G), \tatt)}
{f}{f \circ \alpha_{G}}
\end{align*}
of abelian groups. 
\end{notation}

\begin{notation} \label{betanot} 
Let $G$ be a homotopy invariant commutative group scheme over a perfect field $k$. We define $\beta_{G}$ to be the composite morphism:
\begin{align*}
\beta_{G}: L(G) \slra{z_{G}} M(G) \slra{\alpha_{G}} \Mof(G)
\end{align*}
in the category $\DMeetkQ$. Recalling that $\alpha_{G}$ is defined in Notation \ref{alphanot} to be the composition 
$a_{G} \circ \Cx(\gamma_{G})$, it follows that $\beta_{G} = a_{G} \circ \Cx(\gamma_{G}) \circ z_{G}$.
\end{notation}

\begin{remark} \label{betaequalsgamma}
Recall the morphism $\gamma_{G}: L(G) \lra \wtG$ in the category $\STEkQ$ of Construction \ref{gammaX}. Treating the {\etl} sheaves with transfers $L(G)$ and $\wtG$ as complexes concentrated in degree zero, we may regard $\gamma_{G}$ as a morphism $L(G) \ra \Mof(G)$ in the category $\DMeetkQ$. We claim that $\gamma_{G}$ is equal to the morphism $\beta_{G}$ of Notation \ref{betanot}, i.e. that the diagram:
\begin{align*}
\xymatrix@C=20pt{ 
L(G) \ar@{->}[rr]^(0.5){z_{G}} \ar@{->}[dd]_(0.5){\gamma_{G}} & & 
M(G) \ar@{->}[dd]^(0.5){\Cx(\gamma_{G})} \\
& & \\
\Mof(G) \ar@{<-}[rr]_(0.5){a_{G}} & & \Cx(\wtG)}
\end{align*}
of complexes in $\DMeetkQ$ commutes. Since $L(G)$ and $\Mof(G)$ are complexes concentrated in degree zero,
this diagram commutes trivially in non-zero degrees, while in degree zero, it looks like:
\begin{align*}
\xymatrix@C=20pt{ 
L(G) \ar@{->}[rr]^(0.5){1_{L(G)}} \ar@{->}[dd]_(0.5){\gamma_{G}} & & 
L(G) \ar@{->}[dd]_(0.5){\gamma_{G}} \\
& & \\
\Mof(G) \ar@{<-}[rr]_(0.5){1_{\Mof(G)}} & & \Mof(G),}
\end{align*}
which plainly commutes. We will therefore regard $\beta_{G}: L(G) \ra \Mof(G)$ interchangeably as a morphism in $\DMeetkQ$ 
or in $\STEkQ$.
\end{remark} 

\begin{proposition} \label{alphbetacomm}
Let $G$ be a homotopy invariant commutative group scheme over a perfect field $k$. Then the diagram\tn{:} 
\begin{align*}
\xymatrix@C=20pt{ 
\Mor_{\DMeetkQ}(\Mof(G), \tatt)  \ar@{->}[rr]^(0.5){\alpha_{G}^{*}}  \ar@{->}[dd]_(0.5){\zeta_{\Mof(G)}} & & 
\Mor_{\DMeetkQ}(M(G), \tatt)  \ar@{->}[dd]^(0.5){\zeta_{M(G)}} \\
& & \\
\Exto_{\STEkQ}(G, \Gm) \ar@{->}[rr]_(0.5){\beta_{G}^{*}} & & \Exto_{\STEkQ}(L(G), \Gm) }
\end{align*}
is commutative. Here, $\zeta_{\Mof(G)}$ and $\zeta_{M(G)}$ are the morphisms of Construction \ref{zetaMGmorphs}, $\beta_{H}$ is the morphism of Notation \ref{betanot}, and $\alpha_{G}$ is the morphism of Notation \ref{alphanot}.
\end{proposition}
\begin{prooff}
By definition of the morphisms $\zeta_{\Mof(G)}$ and $\zeta_{M(G)}$ in Construction \ref{zetaMGmorphs}, the given diagram may be expanded as: 
\begin{align*}
\xymatrix@C=20pt{ 
\Mor_{\DMeetkQ}(\Mof(G), \tatt)  \ar@{->}[rr]^(0.5){\alpha_{G}^{*}} \ar@{->}[dd]_(0.5){a_{\Mof(G)}} 
\ar@/_9pc/@{->}[dddddd]^(0.5){\zeta_{\Mof(G)}} & & 
\Mor_{\DMeetkQ}(M(G), \tatt)  \ar@{->}[dd]^(0.5){a_{M(G)}} \ar@/^9pc/@{->}[dddddd]_(0.5){\zeta_{M(G)}} \\
& & \\
\Mor_{\Dm(\STEkQ)}(\Mof(G), \tatt)  \ar@{->}[rr]^(0.5){\alpha_{G}^{*}} \ar@{=}[dd] & & 
\Mor_{\Dm(\STEkQ)}(M(G), \tatt)  \ar@{->}[dd]^(0.5){z_{G}^{*}} \\
& & \\
\Mor_{\Dm(\STEkQ)}(\Mof(G), \tatt) \ar@{->}[rr]^{\beta_{G}^{*}} 
\ar@{->}[dd]_(0.5){\kappa_{\Mof(G), \Gm}} & & 
\Mor_{\Dm(\STEkQ)}(L(G), \tatt) \ar@{->}[dd]^(0.5){\kappa_{L(G), \Gm}} \\
& & \\
\Exto_{\STEkQ}(G, \Gm) \ar@{->}[rr]_{\beta_{G}^{*}} & & \Exto_{\STEkQ}(L(G), \Gm).}
\end{align*}
The top square of this diagram commutes by the functoriality of the morphism $a_{\cly{M}}$ of Notation \ref{aincnot}, the middle square commutes because the morphism $\beta_{G}$ is equal to $\alpha_{G} \circ z_{G}$, by its definition in Notation \ref{betanot}, and the bottom square commutes by the functoriality of the morphism $\kappa_{Z, X}$ of Notation \ref{kappanot}. 
\end{prooff}

\begin{notation} \label{omegavsig} 
Let $H$ be a commutative group scheme over a perfect field $k$. We denote by:
\begin{align*}
& \vsig_{H}: \Exto_{\STEkQ}(H, \Gm) \lra \Exto_{\SE(k, \Q)}(H, \Gm) \qaq \\
& \omega_{H}: \Exto_{\SE(k, \Q)}(H, \Gm) \lra \Hoet(H, \Gm)
\end{align*}
respectively, the morphism induced by the ``forget transfers'' functor: 
\begin{align*}
\fgt: \STEkQ \ra \SEkQ, 
\end{align*}
and the canonical isomorphism which considers an extension $\cly{E}$ of $H$ by $\Gm$ as a $\Gm$-torsor over $H$. Using the canonical isomorphy of the first derived functor and Cech cohomology groups, we describe $\omega_{H}$ as a map 
$\Exto_{\SE(k, \Q)}(H, \Gm) \lra \cHoet(H, \Gm)$, explicitly in terms of cocycles, as follows (c.f. \cite[\S 4]{Mi3} or \cite[\S 11]{Mi1}). Let:
\begin{align*}
\cly{E}: 1 \lra \Gm \slra{\alpha} G \slra{\beta} H \lra 0
\end{align*}
be an element of $\Exto_{\SE(k, \Q)}(H, \Gm)$, giving $G$ the structure of a $\Gm$-torsor over $H$. Choose a cover 
$\cly{U} = \{ U_{i} \ra H \}_{i \in I}$ trivialising this torsor, meaning that for all $s_{i} \in G(U_{i})$, the map
$\vp_{s_{i}}: \Gm|_{U_{i}} \ra G|_{U_{i}}$ given on all sections $V$ by:
\begin{align*}
\mathsmap{\vp_{s_{i}}(V)}{\Gm|_{U_{i}}(V)}{G|_{U_{i}}(V)}{g}{s_{i}|_{V} \cdot g}
\end{align*}
is an isomorphism of sheaves. Fixing any such set $\{ s_{i} \in G(U_{i}) \mid i \in I \}$, it follows from the isomorphy of the 
$\vp_{s_{i}}$ that for all $i, j \in I$ there are unique $g_{ij} \in \Gm(U_{ij})$ satisfying the equation:
\begin{align*}
 s_{j}|_{U_{ij}} = s_{i}|_{U_{ij}} \cdot g_{ij} \in G(U_{ij}).
\end{align*}
In particular, for all $i, j$ and $k \in I$, we have (omitting restriction notation for easy reading): 
\begin{align*}
 s_{ij} \cdot g_{ik} = s_{ik} = s_{j} g_{jk} = s_{i} \cdot g_{ij} \cdot g_{jk},
\end{align*}
implying that $g_{ik} = g_{ij} \cdot g_{jk}$, so that $g = (g_{ij})_{i,j \in I}$ represents a class $[g]$ in $\cHoet(\cly{U}, \Gm)$, hence in 
$\cHoet(H, \Gm)$. A different choice of the elements $s_{i}^{\prm} \in G(U_{i})$ yields different $g_{ij}^{\prm} \in \Gm(U_{ij})$, but the $1$-cocycles $g = (g_{ij})_{i,j \in I}$ and $g^{\prm} = (g_{ij}^{\prm})_{i,j \in I}$ are cohomologous (\cite[\S 4]{Mi3}). Similarly, an extension $\cly{E}^{\prm}$ of $H$ by $\Gm$ isomorphic to $\cly{E}$ yields via this procedure a 1-cocycle cohomologous to $g$. Thus 
$\cly{E} \mapsto [g]$ is a well-defined morphism $\omega_{H}: \Exto_{\SE(k, \Q)}(H, \Gm) \ra \cHoet(H, \Gm)$.
\end{notation} 

\begin{construction} \label{constR} 
Let $H$ be a commutative group scheme over a perfect field $k$ and denote by 
$Y_{H}: \Mor_{\SE(k, \Q)}(\ul{H}, \pul{X}) \ra \Gamma(H, \pul{X})$ the Yoneda transformation. We give here an explicit description of the morphism:
\begin{align*}
(R^{1} Y_{H})(\Gm): \Exto_{\SE(k, \Q)}(H, \Gm) \lra \cHoet(H, \Gm).
\end{align*}
In the following, we denote $(R^{1} Y_{H})(\Gm)$ by $R$ for brevity. To begin, let:
\begin{align*}
\cly{E}: 1 \lra \Gm \slra{\alpha} G \slra{\beta} H \lra 0
\end{align*}
be an element of $\Exto_{\SE(k, \Q)}(H, \Gm)$. By definition of the first right-derived functor, the last square in the diagram: 
\begin{align}
\xymatrix@C=10pt@R=15pt{ 
\Mor(\ul{H}, \ul{\Gm}) \ar@{->}[rr] \ar@{->}[dd]_(0.5){Y_{H}(\Gm)} & & \Mor(\ul{H}, \ul{G}) \ar@{->}[rr] \ar@{->}[dd]^(0.5){Y_{H}(G)} 
& & \Mor(\ul{H}, \ul{H}) \ar@{->}[rr]^(0.4){\delta_{1}} \ar@{->}[dd]^(0.5){Y_{H}(H)} & & 
\Exto_{\SE(k, \Q)}(\ul{H}, \ul{\Gm}) \ar@{->}[dd]^(0.5){(R^{1} Y_{H})(\Gm)} \\
& & & & & & \\
\Gamma(H, \ul{\Gm}) \ar@{->}[rr] & & \Gamma(H, \ul{G})  \ar@{->}[rr] & & \Gamma(H, \ul{H})  \ar@{->}[rr]_(0.47){\delta_{2}} & & \cHoet(H, \Gm)} \label{exth1diag}
\end{align} 
commutes. Here, the top and bottom rows are the come from the long exact sequences associated to Ext groups and Cech cohomology groups, respectively, all morphism sets are in the category $\SE(k, \Q)$, and the underline notation for the group schemes emphasizes their roles as sheaves. Now, the boundary map $\delta_{1}$ takes any morphism $\vp: H \ra H$ to the extension of $H$ by $\Gm$ defined by the top row of the diagram:
\begin{align}
\xymatrix@C=20pt@R=15pt{ 
0 \ar@{->}[rr] & & \Gm \ar@{->}[rr]^(0.5){\alpha^{\prm}} \ar@{=}[dd] & & 
G^{\prm} \ar@{->}[rr]^(0.5){\beta^{\prm}} \ar@{->}[dd]^(0.5){\gamma} & & H \ar@{->}[rr] \ar@{->}[dd]^(0.5){\vp} & & 0 \\
& & & & & & & & \\
0 \ar@{->}[rr] & & \Gm \ar@{->}[rr]_(0.5){\alpha} & & G \ar@{->}[rr]_(0.47){\beta} & & H \ar@{->}[rr] & & 0.}
\end{align} 
In this diagram, the right-hand square is the pushout of the morphisms $\beta$ and $\vp$, so in particular $G^{\prm}$ is given by the fibre product $G \x_{H, \vp} H$, where $H$ is considered as an $H$-scheme via $\vp$, and $\beta^{\prm}$ and $\gamma$ are the associated projections. In particular: 
\begin{align}
\delta_{1}(1_{H}) = \cly{E}. \label{delta1E}
\end{align}
The boundary map $\delta_{2}$ is defined as follows (c.f. \cite[Proof of Proposition 4.5]{Mi3}). Given a morphism $h \in H(H)$, let 
$\cly{U} = \{ U_{i} \ra H \}_{i \in I}$ be a cover on which $\beta(H)$ is locally surjective, and choose for each $i \in I$ elements 
$g_{i} \in G(U_{i})$ such that $\beta(U_{i})(g_{i}) = h|_{U_{i}}$. Then:
\begin{align*}
\beta \left( (g_{i}|_{U_{ij}})^{-1} (g_{j}|_{U_{ij}}) \right) = (h|_{U_{ij}})^{-1} h|_{U_{ij}} = 1,
\end{align*} 
so $(g_{i}|_{U_{ij}})^{-1} (g_{j}|_{U_{ij}}) \in \ker(\beta) = \im(\alpha)$, and we may define:
\begin{align*}
c_{ij} := \alpha^{-1} \big( (g_{i}|_{U_{ij}})^{-1} (g_{j}|_{U_{ij}}) \big).
\end{align*}
Then $c(h) = (c_{ij})_{i,j \in I}$ is a 1-cocycle, hence defines a cohomology class $[c(h)]$ in $\cHoet(\cly{U}, \Gm)$ and consequently also in $\cHoet(H, \Gm)$. The commutativity of the last square of (\ref{exth1diag}) and the (\ref{delta1E}) imply that:
\begin{align*}
R(\cly{E}) = [c(1_{H})]
\end{align*}
defining the morphism $R$ explicitly.
\end{construction} 

\begin{proposition} \label{omegahRagree}
Let $H$ be a commutative group scheme over a perfect field $k$. Then the morphisms\tn{:}
\begin{align*}
\omega_{H}: & \Exto_{\SE(k, \Q)}(H, \Gm) \lra \Hoet(H, \Gm) \qaq \\
R: & \Exto_{\SE(k, \Q)}(H, \Gm) \lra \Hoet(H, \Gm)
\end{align*}
from Notation \ref{omegavsig} and Construction \ref{constR} agree.
\end{proposition} 
\begin{prooff} 
Let $\cly{E}: 1 \ra \Gm \sra{\alpha} G \sra{\beta} H \ra 0$ be an element of $\Exto_{\SE(k, \Q)}(H, \Gm)$. By Construction \ref{constR}, we have $R(\cly{E}) = [c(1_{H})]$, where:
\begin{itemize}
\item[\tn{(i)}] $c(1_{H}) = (c_{ij})_{i,j \in I}$, for $c_{ij}$ defined by:
\item[\tn{(ii)}] $c_{ij} := \alpha^{-1} \big( (g_{i}|_{U_{ij}})^{-1} (g_{j}|_{U_{ij}}) \big)$, where:
\item[\tn{(iii)}] $g_{i} \in G(U_{i})$ are chosen such that $\beta(U_{i})(g_{i}) = 1_{H}|_{U_{i}}$, where: 
\item[\tn{(iv)}] $\cly{U} = \{ U_{i} \ra H \}_{i \in I}$ is a cover on which $\beta(H)$ is locally surjective.
\end{itemize}
Since the cover $\cly{U}$ permits local sections $g_{i}$ to the identity $1_{H}$, it also trivialises the $\Gm$-torsor $\cly{E}$. We may therefore use $\cly{U}$ as outlined in Notation \ref{omegavsig} to construct a 1-cocycle $g = (g_{ij})_{i,j \in I}$ yielding a cohomology class 
$\omega_{H}(\cly{E}) = [g] \in \cHoet(H, \Gm)$. The only choice this involves is the set of $s_{i} \in G(U_{i})$. We set the $s_{i}$ equal to the $g_{i}$ from (iii) above, so that the $g_{ij} \in \Gm(U_{ij})$ are uniquely determined by the relation 
$g_{j}|_{U_{ij}} = g_{i}|_{U_{ij}} \cdot g_{ij}$. Recalling that $\Gm$ acts on $G$ by the morphism $\alpha$, this implies: 
\begin{align*}
g_{ij} = \alpha^{-1}((g_{i}|_{U_{ij}})^{-1} g_{j}|_{U_{ij}}) = c_{ij},
\end{align*}
and so $g = (g_{ij}) = (c_{ij}) = c(1_{H})$, meaning that $\omega_{H}$ and $R$ agree.
\end{prooff}

\begin{proposition} \label{c1invprop}
Let $G$ be a commutative group scheme over a perfect field $k$. Then\tn{:} 
\begin{align*}
\begin{array}{rcl}
\vth_{G} \circ \psi_{G} \circ \zeta_{M(G)} & : & \Mor_{\DMeetkQ}(M(G), \tatt) \lra \Pic(G) \qaq \\
c_{1} := c_{1}^{G} & : & \Pic(G) \lra \Mor_{\DMeetkQ}(M(G), \tatt)
\end{array}
\end{align*}
are mutually inverse isomorphisms. Here, $c_{1}^{G}$ is the Chern class morphism defined in Notation \ref{chernclassnot}.
\end{proposition}
\begin{prooff}
We will show that composing the maps in either order gives the appropriate identity. We only need to unpack the definitions. 
Indeed: 
\begin{itemize}
\item[\tn{(i)}] By Notation \ref {chernclassnot}, $c_{1}^{G} = b_{1} \circ (\vth_{G})^{-1}$. 
\item[\tn{(ii)}] Denote by $d_{1}$ the composition $\vth_{G} \circ \psi_{G} \circ \zeta_{M(G)}$, for clarity. 
\item[\tn{(iii)}] By Construction \ref{zetaMGmorphs}, $\zeta_{M(G)} = \kappa_{L(G), \Gm} \circ z_{G}^{*} \circ a_{M(G)}$.
\item[\tn{(iv)}] By Notation \ref{b1mapnot}, $b_{1} = (a_{M(G)})^{-1} \circ (z_{G}^{*})^{-1} \circ \eta_{G}$.
\item[\tn{(v)}] By Notation \ref{psivthnot}, $\eta_{G} = \lambda_{L(G), \Gm} \circ \left( \psi_{G} \right)^{-1}$.
\item[\tn{(vi)}] Recall from Notation \ref{kappanot} that $\kappa_{L(G), \Gm}$ and $\lambda_{L(G), \Gm}$ are mutually inverse isomorphisms by definition.
\end{itemize}
We therefore have:
\begin{align*}
c_{1} \circ d_{1} & \lgeqt{(i) and (ii)} b_{1} \circ (\vth_{G})^{-1} \circ \vth_{G} \circ \psi_{G} \circ \zeta_{M(G)} \\ 
& = b_{1} \circ \psi_{G} \circ \zeta_{M(G)} \lgeqt{(iii)} b_{1} \circ \psi_{G} \circ \kappa_{L(G), \Gm} \circ z_{G}^{*} \circ a_{M(G)} \\
& \lgeqt{(iv)} 
(a_{M(G)})^{-1} \circ (z_{G}^{*})^{-1} \circ \eta_{G} \circ \psi_{G} \circ \kappa_{L(G), \Gm} \circ z_{G}^{*} \circ a_{M(G)} \\
& \lgeqt{(v)} 
(a_{M(G)})^{-1} \circ (z_{G}^{*})^{-1} \circ \lambda_{L(G), \Gm} \circ \left( \psi_{G} \right)^{-1} \circ \psi_{G} \circ \kappa_{L(G), \Gm} \circ z_{G}^{*} \circ a_{M(G)} \\
& = (a_{M(G)})^{-1} \circ (z_{G}^{*})^{-1} \circ \lambda_{L(G), \Gm} \circ \kappa_{L(G), \Gm} \circ z_{G}^{*} \circ a_{M(G)} \\
& \lgeqt{(vi)} (a_{M(G)})^{-1} \circ (z_{G}^{*})^{-1} \circ z_{G}^{*} \circ a_{M(G)} \\
& = 1 \in \Mor_{\DMeetkQ}(M(G), \tatt).
\end{align*}
Composition the other way round yields:
\begin{align*}
d_{1} \circ c_{1} & \lgeqt{(i) and (ii)} \vth_{G} \circ \psi_{G} \circ \zeta_{M(G)} \circ b_{1} \circ (\vth_{G})^{-1} \\ 
& \lgeqt{(iii)} \vth_{G} \circ \psi_{G} \circ \kappa_{L(G), \Gm} \circ z_{G}^{*} \circ a_{M(G)} \circ b_{1} \circ (\vth_{G})^{-1} \\
& \lgeqt{(iv)} 
\vth_{G} \circ \psi_{G} \circ \kappa_{L(G), \Gm} \circ z_{G}^{*} \circ a_{M(G)} \circ (a_{M(G)})^{-1} \circ (z_{G}^{*})^{-1} 
\circ \eta_{G} \circ (\vth_{G})^{-1} \\
& = \vth_{G} \circ \psi_{G} \circ \kappa_{L(G), \Gm} \circ \eta_{G} \circ (\vth_{G})^{-1} \\
& \lgeqt{(v)} \vth_{G} \circ \psi_{G} \circ \kappa_{L(G), \Gm} \circ \lambda_{L(G), \Gm} \circ \left( \psi_{G} \right)^{-1} \circ (\vth_{G})^{-1} \\
& \lgeqt{(vi)} \vth_{G} \circ \psi_{G} \circ \left( \psi_{G} \right)^{-1} \circ (\vth_{G})^{-1} \\
& = 1 \in \Pic(G).
\end{align*}
This completes the proof. 
\end{prooff}

\begin{proposition} \label{psiomegacomm}
Let $H$ be a commutative group scheme over a perfect field $k$. Then the diagram\tn{:} 
\begin{align}
\xymatrix@C=20pt@R=15pt{ 
\Exto_{\STEkQ}(H, \Gm) \ar@{->}[rr]^(0.47){\beta_{H}^{*}(\Gm)} \ar@{->}[dd]_(0.5){\vsig_{H}} & & 
\Exto_{\STEkQ}(L(H), \Gm) \ar@{->}[dd]^(0.5){\psi_{H}} \\
& & \\
\Exto_{\SE(k, \Q)}(H, \Gm) \ar@{->}[rr]_(0.55){\omega_{H}} & & \Hoet(H, \Gm)} \label{extlevdiag}
\end{align}
is commutative. Here, $\beta_{H}^{*}(\Gm)$ is the map on $\Exto$ groups induced by the morphism $\beta_{H}$ of Notation \ref{betanot}, $\psi_{H}$ is the morphism of Notation \ref{psivthnot}, and $\vsig_{H}$ and $\omega_{H}$ are the morphisms of Notation \ref{omegavsig}.
\end{proposition}
\begin{prooff}
Let us denote by $\beta_{H}^{*}: \Mor_{\STEkQ}(H, \pul{X}) \ra \Mor_{\STEkQ}(L(H), \pul{X})$ the natural transformation of functors defined for sheaf $\shfF$ in $\STEkQ$ by:
\begin{align*}
\mathsmap{\beta_{H}^{*}(\shfF)}{\Mor_{\STEkQ}(H, \shfF)}{\Mor_{\STEkQ}(L(H), \shfF)}{\vp}{\vp \circ \beta_{H},} 
\end{align*} 
which is induced by the morphism $\beta_{H}: L(H) \ra \Mof(H)$ of {\etl} sheaves with transfers. Observe here the abuse of Notation, in that we also use in Diagram (\ref{extlevdiag}) the symbol $\beta_{H}^{*}$ in the induced homomorphism:
\begin{align*}
\beta_{H}^{*}(\Gm): \Exto_{\STEkQ}(H, \Gm) \lra \Exto_{\STEkQ}(L(H), \Gm).
\end{align*} 
Further, we write 
\begin{align*}
& F_{1}: \Mor_{\SE(k, \Q)}(H, \pul{X}) \lra \Gamma(H, \pul{X}) \qaq \\
& F_{2}: \Mor_{\STEkQ}(L(H), \pul{X}) \lra \Gamma(H, \pul{X})
\end{align*} 
for the natural transformations defined, for any sheaf $\shfF$ in $\STEkQ$, by:
\begin{align*}
& \mathsmap{F_{1}(\shfF)}{\Mor_{\SE(k, \Q)}(H, \shfF)}{\Gamma(H, \shfF)}{\left( \vp: H \ra \shfF \right)}{\vp(H)(1_{H})} \qaq \\
& \mathsmap{F_{2}(\shfF)}{\Mor_{\STEkQ}(L(H), \shfF)}{\Gamma(H, \shfF)}{\left( \vp: L(H) \ra \shfF \right)}{\vp(H)(1_{H}^{\prm}).}
\end{align*} 
Here, $1_{H}$ and $1_{H}^{\prm}$ denote the identities $H \ra H$ in the categories $\Smk$ and $\Smck$, respectively. Consider now the diagram: 
\begin{align}
\xymatrix@C=20pt@R=15pt{ 
\Mor_{\STEkQ}(H, \pul{X}) \ar@{->}[rr]^(0.47){\beta_{H}^{*}} \ar@{->}[dd]_(0.5){\tn{fgt}} & & 
\Mor_{\STEkQ}(L(H), \pul{X}) \ar@{->}[dd]^(0.5){F_{2}} \\
& & \\
\Mor_{\SE(k, \Q)}(H, \pul{X}) \ar@{->}[rr]_(0.55){F_{1}} & & \Gamma(H, \pul{X}).} \label{homlevdiag}
\end{align}
We check that (\ref{homlevdiag}) is a commutative diagram of natural transformations by letting $\shfF$ be any sheaf in $\STEkQ$ and examining:
\begin{align}
\xymatrix@C=20pt@R=15pt{ 
\Mor_{\STEkQ}(H, \shfF) \ar@{->}[rr]^(0.47){\beta_{h}^{*}(\shfF)} \ar@{->}[dd]_(0.5){\tn{fgt}(\shfF)} & & 
\Mor_{\STEkQ}(L(H), \shfF) \ar@{->}[dd]^(0.5){F_{2}(\shfF)} \\
& & \\
\Mor_{\SE(k, \Q)}(H, \shfF) \ar@{->}[rr]_(0.55){F_{1}(\shfF)} & & \Gamma(H, \shfF).} \label{homlevdiagtwo}
\end{align}
Let $\vp: H \ra \shfF$ be an element of $\Mor_{\STEkQ}(H, \shfF)$ and observe that $\beta_{H}(H)(1_{H}^{\prm}) = 1_{H}$, simply because $\beta_{H}(H): \Mor_{\Smck}(H, H) \ra \Mor_{\Smk}(H, H)$ is a group homomorphism. We compute the image of $\vp$ along the top right path of Diagram (\ref{homlevdiagtwo}). For ease of reading, we omit the notation ``$F(\shfF)$'' from the morphisms $F$:
\begin{align*}
(F_{2} \circ \beta_{H}^{*})(\vp) & = F_{2}(\beta_{H}^{*}(\vp)) = F_{2}(\vp \circ \beta_{H}) \\
& = (\vp \circ \beta_{H})(H)(1_{H}^{\prm}) = \big( \vp(H) \circ \beta_{H}(H) \big) \Big(1_{H}^{\prm} \Big) \\
& = \vp(H) \Big( \beta_{H}(H) \big( 1_{H}^{\prm} \big) \Big) \\
& = \vp(H)(1_{H}).
\end{align*}
But along the bottom left of the diagram, we have: 
\begin{align*}
(F_{1} \circ \tn{fgt})(\vp) = F_{1}(\vp) = \vp(H)(1_{H}),
\end{align*}
proving that (\ref{homlevdiagtwo}) commutes, hence also the Diagram (\ref{homlevdiag}). Recalling that $\omega_{H} = (R^{1} Y_{H})(\Gm)$ by Proposition \ref{omegahRagree}, we see that Diagram (\ref{extlevdiag}) is nothing but the application of the first right-derived functor $R^{1}$ to diagram (\ref{homlevdiag}), evaluated over the section $\Gm$. This completes the proof.
\end{prooff}

\begin{proposition} \label{conealphah}
Let $1 \lra \Gm \slra{f} G \slra{g} H \lra 0$ be the short exact sequence of semiabelian varieties of Construction \ref{secondsesconst}
over an algebraically closed field $k = \cj{k}$, and let $\cj{g}: \Aff(G) \ra H$ be the associated line bundle of Notation \ref{affsav}. Then the diagram\tn{:}
\begin{align} 
\xymatrix@C=30pt@R=15pt{   
M(H) \ar@{->}[ddrr]^(0.48){c_{1}(\Aff(G))} \ar@{->}[dd]_(0.5){\alpha_{H}} & & \\
& & \\
\Mof(H) \ar@{->}[rr]_(0.48){h^{1}} & & \tatt} \label{c1h1diag}
\end{align} 
is commutative. Here, $h^{1}$ is the morphism defined in Notation \ref{fngnhn}.
\end{proposition}
\begin{prooff}
The commutativity of the diagram (\ref{c1h1diag}) is precisely the statement that $\alpha^{*}_{H}(h^{1}) = c_{1}(\Aff(G))$. To check this equality, consider the commutative diagram:
\begin{align*}
\xymatrix@C=20pt@R=15pt{ 
\Mor_{\DMeetkQ}(\Mof(H), \tatt)  \ar@{->}[rr]^(0.5){\alpha_{H}^{*}}  \ar@{->}[dd]_(0.5){\zeta_{\Mof(H)}} & & 
\Mor_{\DMeetkQ}(M(H), \tatt)  \ar@{->}[dd]^(0.5){\zeta_{M(H)}} \\
& & \\
\Exto_{\STEkQ}(H, \Gm) \ar@{->}[rr]^(0.5){\beta_{H}^{*}} \ar@{->}[dd]_(0.5){\vsig_{H}} & & 
\Exto_{\STEkQ}(L(H), \Gm) \ar@{->}[dd]^(0.5){\psi_{H}} \\
& & \\
\Exto_{\SE(k, \Q)} (H, \Gm) \ar@{->}[rr]_(0.58){\omega_{H}} & & \Hoet(H, \Gm)} 
\end{align*}
obtained by glueing the lower edge of the square from Proposition \ref{alphbetacomm} to the upper edge of the square from Proposition \ref{psiomegacomm}. Let $\cly{E}$ denote the extension: 
\begin{align*}
1 \lra \Gm \slra{f} G \slra{g} H \lra 0
\end{align*}
in the category $\STEkQ$. By definition of the morphism $\zeta_{\Mof(H)}$ in Construction \ref{zetaMGmorphs}, we have 
$\zeta_{\Mof(H)}(h^{1}) = \cly{E}$. Since the morphism $\vsig$ of Notation \ref{omegavsig} is just ``forget transfers'', one has 
$\vsig(\cly{E}) = \cly{E}$, and so it follows from the commutativity of the diagram that:
\begin{align}
\nonumber
(\psi_{H} \circ \zeta_{M(H)} \circ \alpha_{H}^{*})(h^{1}) & = (\omega_{H} \circ \vsig \circ \zeta_{\Mof(H)})(h^{1}) \\
& = (\omega_{H} \circ \vsig)(\cly{E}) = \omega_{H}(\cly{E}). \label{doublecomm}
\end{align}
Now, since $\omega_{H}(\cly{E})$ is the cohomology class of the $\Gm$-torsor $G$ over $H$ defined by $\cly{E}$, we have 
$\vth_{H}(\omega_{H}(\cly{E})) = \Aff(\cly{E})$, by definition of the morphism $\vth_{H}$ in Notation \ref{chernclassnot}. It therefore follows from relation (\ref{doublecomm}) that:
\begin{align}
(\vth_{H} \circ \psi_{H} \circ \zeta_{M(H)} \circ \alpha_{H}^{*})(h^{1}) = \vth_{H}(\omega_{H}(\cly{E})) = \Aff(\cly{E}). \label{psiaffGrel}
\end{align}
By Proposition \ref{c1invprop}, the first Chern class morphism:
\begin{align*}
c_{1}: \Pic(H) \lra \Mor_{\DMeetkQ}(\Mof(H), \tatt)
\end{align*}
is inverse to the composition $\vth_{H} \circ \psi_{H} \circ \zeta_{M(H)}$, and so it follows from relation (\ref{psiaffGrel}) that 
$\alpha_{H}^{*}(h^{1}) = c_{1}(\Aff(\cly{E}))$. This completes the proof.
\end{prooff}

\begin{proposition} \label{topmiddlesquare}
Let $1 \lra \Gm \slra{f} G \slra{g} H \lra 0$ be the short exact sequence of semiabelian varieties of Construction \ref{secondsesconst}
over an algebraically closed field $k = \cj{k}$, and let $\cj{g}: \Aff(G) \ra H$ be the associated line bundle of Notation \ref{affsav}. Then the diagram\tn {:}
\begin{align} 
\xymatrix@C=50pt{   
M(H) \ar@{->}[rr]^(0.45){1_{M(H)} \cup c_{1}(\Aff(G)) } \ar@{->}[dd]_(0.5){M(\Delta^{n}_{H})} & & 
M(H)(1)[2] \ar@{->}[dd]^(0.5){M(\Delta^{n-1}_{H})(1)[2]} \\
& & \\
M(H)^{\ox n} \ar@{->}[rr]_(0.45){1_{M(H)^{\ox n -1}} \ox c_{1}(\Aff(G))} & & M(H)^{\ox n-1}(1)[2]} \label{MHtopsquare}
\end{align} 
is commutative.
\end{proposition}
\begin{prooff}
By the definition of the symbol ``$\cup$'' in Notation \ref{ostnot}, 
\begin{align*} 
1_{M(H)} \cup c_{1}(\Aff(G) = \Big( 1_{M(H)} \ox c_{1}(\Aff(G)) \Big) \circ M(\Delta^{2}_{H}).
\end{align*} 
Since moreover $M(\Delta^{n-1}_{H})(1)[2] = M(\Delta^{n-1}_{H}) \ox 1_{\one(1)[2]}$, the composition along the top-right path of (\ref{MHtopsquare}) is given by:
\begin{align*} 
& M(\Delta^{n-1}_{H})(1)[2] \circ \Big( 1_{M(H)} \cup c_{1}(\Aff(G)) \Big) \\ 
& = \Big( M(\Delta^{n-1}_{H}) \ox 1_{\one(1)[2]} \Big) \circ \Big( 1_{M(H)} \ox c_{1}(\Aff(G)) \Big) 
\circ M(\Delta^{2}_{H}) \\
& = \Big( M(\Delta^{n-1}_{H}) \circ 1_{M(H)} \Big) \ox \Big( 1_{\one(1)[2]} \circ c_{1}(\Aff(G)) \Big) \circ M(\Delta^{2}_{H}) \\
& = \Big( M(\Delta^{n-1}_{H}) \ox c_{1}(\Aff(G)) \Big) \circ M(\Delta^{2}_{H}),
\end{align*} 
where the second equality follows from the fact that the composition of two tensor products of morphisms can be computed factorwise. The commutativity of (\ref{MHtopsquare}) is therefore equivalent to the commutativity of the perimeter of the diagram:
\begin{align*} 
\xymatrix@C=70pt{   
M(H) \ar@{->}[rr]^(0.45){M(\Delta^{2}_{H})} \ar@{->}[dd]_(0.5){M(\Delta^{n}_{H})} & & 
M(H)^{\ox 2} \ar@{->}[dd]^(0.5){M(\Delta^{n-1}_{H}) \ox c_{1}(\Aff(G))} \\
& & \\
M(H)^{\ox n} \ar@{->}[rr]_(0.45){\left( 1_{M(H)^{\ox n-1}} \right) \ox c_{1}(\Aff(G))} 
\ar@{<-}[uurr]^(0.55){\tn{$M(\Delta^{n-1}_{H}) \ox 1_{M(H)} \phn{XXXX}$}} & & M(H)^{\ox n-1}(1)[2]}
\end{align*} 
The upper triangle in this diagram is the application of the functor $M$ to the evidently commutative diagram of $k$-group schemes:
\begin{align*} 
\xymatrix@C=40pt{   
H \ar@{->}[rr]^(0.45){\Delta^{2}_{H}} \ar@{->}[dd]_(0.5){\Delta^{n}_{H}} & & H^{2} \\
& & \\
H^{n}, \ar@{<-}[uurr]_(0.5){\tn{$\Delta^{n-1}_{H} \x 1_{H}$}} & & }
\end{align*}
and the commutativity of the lower triangle is clear.
\end{prooff}

\section{The motive of a semiabelian variety}

\begin{proposition} \label{candiagramcomm}
Let $1 \lra \Gm \slra{f} G \slra{g} H \lra 0$ be the short exact sequence of Construction \ref{secondsesconst} over an algebraically closed 
field $k = \cj{k}$. Then the diagrams\tn{:}
\begin{align} 
\xymatrix@C=14pt{   
M(G) \ar@{->}[rr]^(0.5){M(g)} \ar@{->}[dd]_(0.5){\vp_{G}} & & M(H) \ar@{->}[dd]^(0.5){\vp_{H}} \\
& & \\
\Sym(\Mof(G)) \ar@{->}[rr]_(0.49){g^{\bl}} & & \Sym(\Mof(H))} \tn{ and }
\xymatrix@C=14pt{   
M(H) \ar@{->}[rr]^(0.45){1_{M(H)} \cup c_{1}(\Aff(G))} \ar@{->}[dd]_(0.5){\vp_{H}} & & 
M(H)(1)[2] \ar@{->}[dd]^(0.5){\vp_{H}(1)[2]} \\
& & \\
\Sym(\Mof(H)) \ar@{->}[rr]_(0.43){h^{\bl}} & & \Sym(\Mof(H))(1)[2]} \label{canondiag}
\end{align} 
are both commutative in the category $\DMeetkQ$.
\end{proposition}
\begin{prooff}
\tbf{(i)} Let $n$ be a positive integer. Then $g^{n}: \Symn(\Mof(G)) \ra \Symn(\Mof(H))$ is equal to $\Symn(g)$, as observed in Remark \ref{gniexplicit}, and so the direct sum $g^{\bl}$ is equal to to $\Sym(g)$. It follows from the functoriality of the morphism $\vp_{G}$ shown in Proposition \ref{vpnfunc} that the left-hand square commutes.\\\\
\tbf{(ii)} To see that the right-hand square commutes, we consider the following diagram:
\begin{align} 
\xymatrix@C=30pt{   
M(H) \ar@{->}[rr]^(0.45){1_{M(H)} \cup c_{1}(\Aff(G))} \ar@{->}[dd]^(0.5){M(\Delta^{n}_{H})} 
\ar@{}[ddrr]|(0.48){\tn{\large{\tbf{1.}}}} 
\ar@/_4pc/@{->}[dddddd]_(0.5){\vp^{n}_{H}} & & 
M(H)(1)[2] \ar@{->}[dd]_(0.5){M(\Delta^{n-1}_{H})(1)[2]} \ar@/^4.5pc/@{->}[dddddd]^(0.5){\vp^{n-1}_{H}(1)[2]} \\
& & \\
M(H)^{\ox n} \ar@{->}[rr]_(0.45){1_{M(H)^{\ox n -1}} \cup c_{1}(\Aff(G))} 
\ar@{}[ddrr]|(0.48){\tn{\large{\tbf{2.}}}} 
\ar@{->}[dd]^(0.5){\alpha^{\ox n}_{H}}_{\tn{\large{\tbf{4.} \phn{W}}}} & & 
M(H)^{\ox n-1}(1)[2] \ar@{->}[dd]_(0.5){\alpha^{\ox n-1}_{H}(1)[2]}^{\tn{\large{\phn{W} \tbf{5.}}}} \\
& & \\
\Mof(H)^{\ox n} \ar@{->}[rr]_(0.43){1_{\Mof(H)^{\ox n -1}} \ox h^{1}} \ar@{<-}[dd]^(0.5){\iota^{n}_{H}} 
\ar@{}[ddrr]|(0.48){\tn{\large{\tbf{3.}}}} & & 
\Mof(H)^{\ox n-1}(1)[2] \ar@{<-}[dd]_(0.5){\iota^{n-1}_{H}(1)[2]} \\
& & \\
\Sym^{n}(\Mof(H)) \ar@{->}[rr]_(0.43){h^{n}} & & \Sym^{n-1}(\Mof(H))(1)[2].} \label{largemidsquare}
\end{align} 
Squares 1 and 3 commute by Propositions \ref{topmiddlesquare} and \ref{hsquarecommprop}, respectively. Square 2 is the tensor product of the two diagrams: 
\begin{align*} 
\xymatrix@C=30pt{   
M(H)^{\ox n-1} \ar@{->}[rr]^(0.5){1_{M(H)^{\ox n -1}}} \ar@{->}[dd]_(0.5){\alpha^{\ox n-1}_{H}} & & 
M(H)^{\ox n-1}\ar@{->}[dd]^(0.5){\alpha^{\ox n-1}_{H}} \\
& & \\
\Mof(H)^{\ox n-1} \ar@{->}[rr]_(0.5){1_{\Mof(H)^{\ox n-1}}} & & \Mof(H)^{\ox n-1}} \qaq
\xymatrix@C=30pt{   
M(H) \ar@{->}[rr]^(0.48){c_{1}(\Aff(G))} \ar@{->}[dd]_(0.5){\alpha_{H}} & & \tatt \ar@{->}[dd]^(0.5){\id} \\
& & \\
\Mof(H) \ar@{->}[rr]_(0.48){h^{1}} & & \tatt.}
\end{align*} 
The square on the left is trivially commutative, and the one on the right commutes by Proposition \ref{conealphah}. Since moreover Squares 4 and 5 commute by Proposition \ref{vpGfactorsnicely}, it follows that the whole of Diagram (\ref{largemidsquare}) is commutative. By considering the perimeter of the subdiagram of (\ref{largemidsquare}) consisting of all except the third square, we obtain the relation:
\begin{align} 
\iota^{n-1}_{H}(1)[2] \circ \vp^{n-1}_{H}(1)[2] \circ \left( 1_{M(H)} \cup c_{1}(\Aff(G)) \right) = 
\left( 1_{\Mof(H)^{\ox n -1}} \ox h^{1} \right) \circ \iota^{n}_{H} \circ \vp^{n}_{H}. \label{perim}
\end{align} 
In view of the identity: 
\begin{align} 
\pi^{n-1}_{H} \circ \iota^{n-1}_{H} = 1_{\Sym^{n-1}(\Mof(H))}, \label{piiotan}
\end{align} 
composing Relation (\ref{perim}) on the left with $\pi^{n-1}_{H}(1)[2]$ yields the first and last equalities in calculation:
\begin{align*} 
\vp^{n-1}_{H}(1)[2] \circ \left( 1_{M(H)} \cup c_{1}(\Aff(G)) \right) & = 
\pi^{n-1}_{H}(1)[2] \circ \left( 1_{\Mof(H)^{\ox n -1}} \ox h^{1} \right) \circ \iota^{n}_{H} \circ \vp^{n}_{H} \\
& = \pi^{n-1}_{H}(1)[2] \circ \iota^{n-1}_{H}(1)[2] \circ h^{n} \circ \vp^{n}_{H} \\
& = h^{n} \circ \vp^{n}_{H}.
\end{align*} 
The middle equality follows from the commutativity of Square 3 in the diagram (\ref{largemidsquare}). Thus we have shown that the square: 
\begin{align} 
\xymatrix@C=30pt{   
M(H) \ar@{->}[rr]^(0.45){1_{M(H)} \cup c_{1}(\Aff(G))} \ar@{->}[dd]_(0.5){\vp^{n}_{H}} & & 
M(H)(1)[2] \ar@{->}[dd]^(0.5){\vp^{n-1}_{H}(1)[2]}\\
& & \\
\Sym^{n}(\Mof(H)) \ar@{->}[rr]_(0.43){h^{n}} & & \Sym^{n-1}(\Mof(H))(1)[2]} \label{compactsquare}
\end{align} 
is commutative. Recalling that $i^{n}_{\Mof(H)}: \Symn(\Mof(H)) \ra \Sym(\Mof(H)) = \os_{n} \Symn(\Mof(H))$ denotes the canonical inclusion, we observe that the diagram: 
\begin{align} 
\xymatrix@C=30pt{   
\Sym^{n}(\Mof(H)) \ar@{->}[rr]^(0.43){h^{n}} \ar@{->}[dd]_(0.5){i^{n}_{\Mof(H)}} & & 
\Sym^{n-1}(\Mof(H))(1)[2] \ar@{->}[dd]^(0.5){i^{n-1}_{\Mof(H)}(1)[2]} \\
& & \\
\Sym(\Mof(H)) \ar@{->}[rr]_(0.43){h^{\bl}} & & \Sym(\Mof(H))(1)[2]} \label{symnhdiag}
\end{align} 
is commutative, by the universal property of the direct sum. Glueing (\ref{symnhdiag}) along the bottom of 
(\ref{compactsquare}) produces a large diagram whose commutativity is expressed by the relation: 
\begin{align} 
i^{n-1}_{\Mof(H)}(1)[2] \circ \vp^{n-1}_{H}(1)[2] \circ \left( 1_{M(H)} \cup c_{1}(\Aff(G)) \right) = 
h^{\bl} \circ i^{n}_{\Mof(H)} \circ \vp^{n}_{H}. \label{ihrelone}
\end{align} 
Now, the morphism $\vp_{H}: M(H) \ra \Sym(\Mof(H))$ is defined in Notation \ref{vpGnot} to be the sum 
$\vp_{H} := \sum_{n} i^{n}_{\Mof(H)} \circ \vp_{H}^{n}$. Consequently:
\begin{align} 
\nonumber \vp_{H}(1)[2] & = \sum_{n} \left( i^{n-1}_{\Mof(H)} \circ \vp^{n-1}_{H} \right)(1)[2] \\
& = \sum_{n} \left( i^{n-1}_{\Mof(H)}(1)[2] \circ \vp^{n-1}_{H}(1)[2] \right). \label{ihreltwo}
\end{align} 
Using bilinearity of composition of morphisms and the fact that neither $1_{M(H)} \cup c_{1}(\Aff(G))$ nor 
$h^{\bl} = \os_{n} h^{n}$ depends on $n$, and it follows from Relations (\ref{ihrelone}) and (\ref{ihreltwo}) that:
\begin{align*} 
\vp_{H}(1)[2] \circ \left( 1_{M(H)} \cup c_{1}(\Aff(G)) \right) & = 
\left\{ \sum_{n} \left( i^{n-1}_{\Mof(H)}(1)[2] \circ \vp^{n-1}_{H}(1)[2] \right) \right\} \circ 
\left( 1_{M(H)} \cup c_{1}(\Aff(G)) \right) \\
& = \sum_{n} \left\{ \left( i^{n-1}_{\Mof(H)}(1)[2] \circ \vp^{n-1}_{H}(1)[2] \right) \circ 
\left( 1_{M(H)} \cup c_{1}(\Aff(G)) \right) \right\} \\
& = \sum_{n} \left( h^{\bl} \circ i^{n}_{\Mof(H)} \circ \vp^{n}_{H} \right) \\
& = h^{\bl} \circ \left\{ \sum_{n} \left( i^{n}_{\Mof(H)} \circ \vp^{n}_{H} \right) \right\} \\
& = h^{\bl} \circ \vp_{H},
\end{align*} 
and so the right-hand square in (\ref{canondiag}) commutes.
\end{prooff}

\begin{proposition}[\tbf{Classical}] \label{tricatisompropfund}
Let $\catt$ be a triangulated category, and suppose that\tn{:} 
\begin{align*} 
\xymatrix@C=14pt@R=14pt{   
A \ar@{->}[rr] \ar@{->}[dd]_(0.5){\vp} & & B \ar@{->}[rr] \ar@{->}[dd]_(0.5){\psi} & & 
C \ar@{->}[rr] \ar@{->}[dd]^(0.5){\xi} & & A[1] \ar@{->}[dd]^(0.5){\vp[1]}\\
& & & & & & \\
X \ar@{->}[rr] & & Y \ar@{->}[rr] & & Z \ar@{->}[rr] & & X[1]}
\end{align*}
is a commutative diagram, the rows of which are both exact triangles in $\catt$. Then if two of the three maps $\vp$, $\psi$ and $\xi$ are isomorphisms, the remaining map is also. 
\end{proposition}
\begin{prooff}
If $\vp$ and $\psi$ are isomorphisms, it follows by \cite[Proposition 1.1.20]{Nee} that $\xi$ is also an isomorphism. The other two cases follow by shifting the triangles.
\end{prooff}

\begin{lemma} \label{savfinkimdim}
A semiabelian variety $G$ over a perfect field $k$ has finite Kimura dimension (see Definition \ref{finitekimdim}).
\end{lemma} 
\begin{prooff}
To say that $G$ has finite Kimura dimension is to say that $\Symn(G)$ is isomorphic to the zero object for large enough $n$. By Proposition\ref{etldes}, this can be checked over an algebraic closure field $\cj{k}$ of $k$. We now proceed by induction on the rank $r$ of $G$. In the base case where $r$ is zero, $G$ is an abelian variety, which has finite Kimura dimension by Remark \ref{Afinkimdim}. Assume now that $r$ is positive, and consider the short exact sequence: 
\begin{align*} 
1 \lra \Gm \slra{f} G \slra{g} H \lra 0
\end{align*}
from Construction \ref{secondsesconst}, in which $H$ is a semiabelian variety of rank $r-1$. As observed in Notation \ref{exacttrisav}, this short exact sequence induces an exact triangle: 
\begin{align*}
\nonumber
\Sym^{n-1}(\Mof(H))(1)[1] & \sxra{f^{n}} \Symn(\Mof(G)) \sxra{g^{n}} \Sym^{n}(\Mof(H)) \\
& \sxra{h^{n}} \Sym^{n-1}(\Mof(H))(1)[2] 
\end{align*}
in the category $\DMeetkQ$, for any positive integer $n$. Since $H$ has rank $r-1$, it has finite Kimura dimension, by the inductive hypothesis. Thus for sufficiently large $n$, the exact triangle above reads:
\begin{align*}
0 \sxra{f^{n}} \Symn(\Mof(G)) \sxra{g^{n}} 0 \sxra{h^{n}} 0,
\end{align*}
implying that $\Sym^{n}(\Mof(G))$ vanishes.
\end{prooff}

\begin{theorem} \label{mainresult}
Let $G$ be a semiabelian variety over a perfect field $k$. Let $\cj{k}$ be a fixed algebraic closure of $k$, write $G_{\cj{k}}$ for the base extension $G \x_{k} \cj{k}$, and denote by\tn{:} 
\begin{align*}
1 \lra \Gm \slra{f} G_{\cj{k}} \slra{g} H \lra 0
\end{align*}
the short exact sequence in $\STE(\cj{k}, \Q)$ of Construction \ref{secondsesconst}. Then\tn{:}
\begin{itemize}
\item[\tn{(i)}] Suppose that the diagram\tn{:}
\begin{align*} 
\xymatrix@C=14pt{   
M(H)(1)[2] \ar@{->}[dd]_(0.5){\vp_{H}(1)[2]} \ar@{->}[rr]^(0.52){\pd_{\Aff(G), H} \circ M(s_{0})} & & 
M(G)[1] \ar@{->}[dd]^(0.5){\vp_{G}[1]} \\
& & \\
\Sym(\Mof(H))(1)[2] \ar@{->}[rr]_(0.53){-f^{\bl}[1]} & & \Sym(\Mof(G))[1]} 
\end{align*} 
is commutative. Then the morphism\tn{:}
\begin{align*}
\vp_{G}: M(G) \lra \Sym(\Mof(G))
\end{align*}
of Notation \ref{vpGnot} is an isomorphism in the category $\DMeetkQ$.
\item[\tn{(ii)}] If $G$ has rank one, then there exists an isomorphism $M(G) \simeq \Sym(\Mof(G))$ in the category $\DMeetkQ$, irrespective of whether the square in Part \tn{(i)} commutes.
\end{itemize}
\end{theorem}
\begin{prooff}
\tbf{(i)} Letting $r$ denote the rank of the semiabelian variety $G$, we note that by Lemma \ref{secondseslemm}, the object $H$ in the short exact sequence: 
\begin{align*}
1 \lra \Gm \slra{f} G_{\cj{k}} \slra{g} H \lra 0,
\end{align*}
is a semiabelian variety of rank $r-1$. By Proposition \ref{etldes}, it suffices to prove that the morphism $\vp_{G}$ is an isomorphism after extension to the designated algebraic closure $\cj{k}$ of $k$. For the remainder of the proof, we will therefore work exclusively over $\cj{k}$. Abusing notation for ease of reading, we will write $G$ in place of $G_{\cj{k}}$. Consider now the diagram:
\begin{align*} 
\xymatrix@C=14pt{   
M(G_{\cj{k}}) \ar@{->}[rr]^(0.5){M(g)} \ar@{->}[dd]_(0.5){\vp_{G}} & & 
M(H) \ar@{->}[rr]^(0.45){1_{M(H)} \cup c_{1}(\Aff(G))} \ar@{->}[dd]^(0.5){\vp_{H}} & & 
M(H)(1)[2] \ar@{->}[dd]^(0.5){\vp_{H}(1)[2]} \ar@{->}[rr]^(0.52){\pd_{\Aff(G), H} \circ M(s_{0})} & & 
M(G)[1] \ar@{->}[dd]^(0.5){\vp_{G}[1]} \\
& & & & & & \\
\Sym(\Mof(G)) \ar@{->}[rr]_(0.49){g^{\bl}} & & 
\Sym(\Mof(H)) \ar@{->}[rr]_(0.43){h^{\bl}} & & 
\Sym(\Mof(H))(1)[2] \ar@{->}[rr]_(0.53){-f^{\bl}[1]} & & \Sym(\Mof(G))[1],} 
\end{align*} 
in which the top row is Euler triangle of Notation \ref{eulerseqsav}, and the bottom row is obtained by left-shifting the exact triangle of Notation \ref{exacttrisav} once. The left-hand square and middle squares in this diagram commute, by Proposition \ref{candiagramcomm}, and the right-hand square commutes by hypothesis. We prove by induction on the rank $r$ of $G$ that $\vp_{G}$ is an isomorphism. If $r = 0$, then $G$ has rank zero, and is therefore an abelian variety. It follows that $\vp_{G}$ is an isomorphism, by Proposition \ref{vpAisom}. Assume then that $r$ is strictly positive. Then $\vp_{H}$, and hence also $\vp_{H}(1)[2]$, is an isomorphism by the inductive hypothesis, since $H$ is a semiabelian variety of rank $r-1$. By Proposition \ref{tricatisompropfund}, this implies that $\vp_{G}$ is an isomorphism.\\\\
\tbf{(ii)} If $G$ has rank one, then $H$ has rank zero, and is therefore an abelian variety. It follows that $\vp_{H}$, and hence also $\vp_{H}(1)[2]$, is an isomorphism, by Proposition \ref{vpAisom}. Consider the diagram:
\begin{align*} 
\xymatrix@C=14pt{   
M(G) \ar@{->}[rr]^(0.5){M(g)} \ar@{-->}[dd]_(0.5){\xi_{G}} & & 
M(H) \ar@{->}[rr]^(0.45){1_{M(H)} \cup c_{1}(\Aff(G))} \ar@{->}[dd]^(0.5){\vp_{H}}_{\simeq} & & 
M(H)(1)[2] \ar@{->}[dd]^(0.5){\vp_{H}(1)[2]}_{\simeq} \ar@{->}[rr]^(0.52){\pd_{\Aff(G), H} \circ M(s_{0})} & & 
M(G)[1] \ar@{-->}[dd]^(0.5){\xi_{G}[1]} \\
& & & & & & \\
\Sym(\Mof(G)) \ar@{->}[rr]_(0.49){g^{\bl}} & & 
\Sym(\Mof(H)) \ar@{->}[rr]_(0.43){h^{\bl}} & & 
\Sym(\Mof(H))(1)[2] \ar@{->}[rr]_(0.53){-f^{\bl}[1]} & & \Sym(\Mof(G))[1].} 
\end{align*} 
Shifting this diagram to the left, we may apply the triangulated category axiom labelled (TR3) in 
\cite[Definition 1.1.2]{Nee} to conclude that there exists, as indicated by the dotted arrow, a morphism:
\begin{align*} 
\xi_{G}: M(G) \lra \Sym(\Mof(G))
\end{align*} 
such that the diagram commutes. The isomorphy of $\vp_{H}$ and $\vp_{H}(1)[2]$ now imply that $\xi_{G}$ is an isomorphism, by Proposition \ref{tricatisompropfund}.
\end{prooff}

\begin{subappendices}

\section{Alternative description of $\Symn(\shfF)$}

We show in Proposition \ref{symnaltdescrip} of this section that $\Symn(\shfF)$ is the sheafification $\Symnp(\shfF)^{\dg}$ of a certain presheaf $\Symnp(\shfF)$.

\begin{notation}
Let $\shfF$ be a sheaf in the category $\cata$ and $n$ a non-negative integer. We denote by: 
\begin{itemize}
\item[\tn{(i)}] $\shfF^{\oxpn}$ the $n^{\tn{th}}$ tensor power presheaf defined by $\shfF^{\oxpn}(U) := \shfF(U)^{\ox n}$ for any section $U$. Thus $\shfF^{\ox n}$ is the sheafification of $\shfF^{\oxpn}$.
\item[\tn{(ii)}] $\Symn_{\tn{pre}}(\shfF)$ the presheaf defined by:
\begin{align*}
\Symnp(\shfF)(U) := \im \left( \frac{1}{n!} \sum_{\sigma \in \Sym(n)} \sigma_{\shfF(U)}: 
\shfF(U)^{\ox n} \lra \shfF(U)^{\ox n} \right),
\end{align*} 
for any section $U$. Note that the image is in the category of $\Q$-vector spaces. Equivalently, one could define:
\begin{align*}
\Symnp(\shfF) := \im \left( \frac{1}{n!} \sum_{\sigma \in \Sym(n)} \sigma_{\shfF, \tn{pre}}: \shfF^{\oxpn} \lra \shfF^{\oxpn} \right).
\end{align*} 
Here, the image is in the category of presheaves of $\Q$-vector spaces with transfers. 
\end{itemize}
Thus $\Symnp(\shfF)$ is the $n^{\tn{th}}$ symmetric power of $\shfF$ in the category of presheaves with transfers, and hence is equipped with a canonical projection: 
\begin{align*}
\pi_{\shfF}^{n, \tn{pre}}: \shfF^{\oxpn} \lra \Symnp(\shfF)
\end{align*}
satisfying the universal property described in Proposition \ref{symunprop}.
\end{notation}

\begin{remark} \label{sheafifysym} 
Let $\shfF$ be a sheaf in the category $\cata$, let $n$ be a non-negative integer, let $\sigma$ be an element of the symmetric group 
$\Sym(n)$, and denote by $\theta: \shfF^{\oxpn} \ra \shfF^{\ox n}$ the sheafification morphism. Denote by:
\begin{align*}
\sigma_{_{\shfF}, \tn{pre}}: \shfF^{\oxpn} \lra \shfF^{\oxpn} \qaq \sigma_{\shfF}: \shfF^{\ox n} \lra \shfF^{\ox n} 
\end{align*}
respectively, the rearrangement of factors endomorphisms of $\shfF^{\oxpn}$ and $\shfF^{\ox n}$ associated to the permutation $\sigma$. Explicitly $\sigma_{_{\shfF}, \tn{pre}}$ is given on a any section $U$ by:
\begin{align*}
\mathsmap{\sigma_{_{\shfF}, \tn{pre}}(U)}{\shfF(U)^{\ox n}}{\shfF(U)^{\ox n}}{x_{1} \ox \cdots \ox x_{n}}
{x_{\sigma(1)} \ox \cdots \ox x_{\sigma(n)},}
\end{align*}
and $\sigma:\shfF^{\ox n} \ra \shfF^{\ox n}$ is defined to be the unique morphism induced by the universal property of sheafification, rendering the diagram:
\begin{align} 
\xymatrix@C=14pt{   
\shfF^{\oxpn} \ar@{->}[rrr]^(0.5){\sigma_{\shfF, \tn{pre}}} \ar@{->}[dd]_(0.5){\theta} & & & 
\shfF^{\oxpn} \ar@{->}[dd]^(0.5){\theta}\\
& & \\
\shfF^{\ox n} \ar@{->}[rrr]_(0.5){\sigma} & & & \shfF^{\ox n}} \label{sigsquare}
\end{align} 
commutative. 
\end{remark}

\begin{construction} \label{abmorphconst}
Let $\shfF$ be a sheaf in the category $\cata$, let $n$ be a non-negative integer, and denote by $\theta: \shfF^{\oxpn} \ra \shfF^{\ox n}$ the sheafification morphism. We construct as follows morphisms:
\begin{align*}
a: \shfF^{\ox n} \lra \Symnp(\shfF) \qaq b: \shfF^{\oxpn} \lra \Symn(\shfF)
\end{align*}
such that $a \circ \sigma_{\shfF} = a$ and $b \circ \sigma_{\shfF, \tn{pre}} = b$, for all permutations $\Sym(n)$. Consider the diagram: 
\begin{align*} 
\xymatrix@C=14pt{   
\shfF^{\oxpn} \ar@{->}[rr]^(0.5){\theta} \ar@{->}[dd]^(0.5){\sigma_{\shfF, \tn{pre}}} 
\ar@/_4pc/@{->}[ddddr]_(0.5){\pi_{\shfF, \tn{pre}}^{n}} & & \shfF^{\ox n} \ar@{->}[dd]_(0.5){\sigma_{\shfF}} \\
& & \\
\shfF^{\oxpn} \ar@{->}[rr]^(0.5){\theta} \ar@{->}[ddr]^(0.5){\pi_{\shfF, \tn{pre}}^{n}} & & \shfF^{\ox n} \ar@{->}[ddl]^(0.5){a} \\
& & \\
& \Symnp(\shfF). &} 
\end{align*} 
We define $a$ to be the unique map, induced by the universal property of the sheafification morphism $\theta$, such that the lower triangle commutes. In Remark \ref{sheafifysym} we saw that the top rectangle commutes, the left-hand triangle commutes by permutation invariance of the canonical projection $\pi_{\shfF, \tn{pre}}^{n}$. Thus the whole diagram commutes, implying that both $a$ and 
$a \circ \sigma_{\shfF}$ satisfy $? \circ \theta = \pi_{\shfF, \tn{pre}}^{n}$, hence must be equal, again by the universal property of $\theta$. Define now $b$ to be the composite morphism:
\begin{align*}
b : \shfF^{\oxpn} \slra{\theta} \shfF^{\ox n} \slra{\pi^{n}_{\shfF}} \Symn(\shfF).
\end{align*}
The commutativity of the square (\ref{sigsquare}) in Remark \ref{sheafifysym} says that 
$\theta \circ \sigma_{\shfF, \tn{pre}} = \sigma_{\shfF} \circ \theta$, justifying the equality (ii) in the chain:
\begin{align*}
b \circ \sigma_{\shfF, \tn{pre}} & \lgeqt{(i)} \pi^{n}_{\shfF} \circ \theta \circ \sigma_{\shfF, \tn{pre}} \\
& \lgeqt{(ii)} \pi^{n}_{\shfF} \circ \sigma_{\shfF} \circ \theta \\
& \lgeqt{(iii)} \pi^{n}_{\shfF} \circ \theta \lgeqt{(iv)} b.
\end{align*}
Equalities (i) and (iv) are the definition of $b$, and equality (iii) follows from the permutation invariance of $\pi^{n}_{\shfF}$.
\end{construction}

\begin{proposition} \label{symnaltdescrip}
For any sheaf $\shfF$ in the category $\cata$ and non-negative integer $n$, the sheaf $\Symn(\shfF)$ is the sheafification $\Symn_{\tn{pre}}(\shfF)^{\dg}$of the presheaf $\Symn_{\tn{pre}}(\shfF)$.
\end{proposition}
\begin{prooff}
Preserving the notation of Construction \ref{abmorphconst}, consider the diagram:
\begin{align} 
\xymatrix@C=25pt{   
& & \shfF^{\ox n} \ar@{->}[dd]^(0.5){a} \ar@{->}[rr]^(0.5){\pi_{\shfF}^{n}} & & \Symn(\shfF) \ar@{->}[dd]^(0.5){\vp} \\
& & & & \\
\shfF^{\oxpn} \ar@{->}[ddrr]_(0.5){b} \ar@{->}[dd]_(0.5){\theta} 
\ar@{->}[uurr]^(0.5){\theta} \ar@{->}[rr]_(0.45){\pi_{\shfF, \tn{pre}}^{n}} & & 
\Symnp(\shfF) \ar@{->}[rr]^(0.5){\xi} \ar@{->}[dd]^(0.5){\psi^{\prm}} & & \Symnp(\shfF)^{\dg} \ar@{->}[ddll]^(0.5){\psi} \\
& & & & \\
\shfF^{\ox n} \ar@{->}[rr]_(0.45){\pi_{\shfF}^{n}} & & \Symn(\shfF), & &} \label{bigsheafdiag}
\end{align} 
in which: 
\begin{itemize}
\item[\tn{(i)}] $\xi: \Symnp(\shfF) \ra \Symnp(\shfF)^{\dg}$ is the sheafification morphism.
\item[\tn{(ii)}] The morphism $a$ is defined in Construction \ref{abmorphconst}, and renders the top-left triangle commutative.
\item[\tn{(iii)}] $\vp: \Symn(\shfF) \ra \Symnp(\shfF)^{\dg}$ is the unique morphism rendering the top-right square commutative, induced by the universal property of $\Symn(\shfF)$. We note here that the composite morphism $\xi \circ a$ is invariant under right composition of factor permutation morphisms $\sigma_{\shfF}: \shfF^{\ox n} \ra \shfF^{\ox n}$, seeing as $a$ has this property, as observed in 
Construction \ref{abmorphconst}.
\item[\tn{(iv)}] $b$ is the composite morphism $\pi_{\shfF}^{n}$ defined in Construction \ref{abmorphconst}.
\item[\tn{(v)}] $\psi^{\prm}: \Symnp(\shfF) \ra \Symn(\shfF)$ is the unique morphism rendering the top-right triangle in the lower-left square commutative, induced by the universal property of the presheaf $\Symnp(\shfF)$. Note we observed in Construction \ref{abmorphconst} that $b$ is invariant under right composition of the factor permutation morphisms\\ 
$\sigma_{\shfF, \tn{pre}}: \shfF^{\oxpn} \ra \shfF^{\oxpn}$.
\item[\tn{(vi)}] $\psi: \Symnp(\shfF)^{\dg} \ra \Symn(\shfF)$ is the unique morphism, induced by the universal property of the sheafification morphism $\xi$, rendering the bottom-right triangle commutative.
\end{itemize}
Thus the whole diagram commutes. In particular, both of the subtriangles:
\begin{align*} 
\xymatrix@C=14pt{ 
\shfF^{\oxpn} \ar@{->}[rr]^(0.5){\theta} \ar@{->}[ddrr]_(0.5){b} & & \shfF^{\ox n} \ar@{->}[dd]^(0.5){\psi^{\prm} \circ a} \\
& & \\
& & \Symn(\shfF)} \qaq
\xymatrix@C=14pt{ 
\shfF^{\oxpn} \ar@{->}[rr]^(0.5){\theta} \ar@{->}[ddrr]_(0.5){b} & & \shfF^{\ox n} \ar@{->}[dd]^(0.5){\pi_{\shfF}^{n}} \\
& & \\
& & \Symn(\shfF)}
\end{align*} 
of (\ref{bigsheafdiag}) commute, implying that $\psi^{\prm} \circ a = \pi_{\shfF}^{n}$, by the universal property of the sheafification morphism $\theta$. The right half of the diagram (\ref{bigsheafdiag}) can therefore be rewritten:
\begin{align*} 
\xymatrix@C=20pt{ 
\shfF^{\ox n} \ar@{->}[rr]^(0.45){\pi_{\shfF}^{n}} \ar@{->}[ddddrr]_(0.5){\pi_{\shfF}^{n}} & & \Symn(\shfF) \ar@{->}[dd]^(0.5){\vp} \\
& & \\
& & \Symnp(\shfF)^{\dg} \ar@{->}[dd]^(0.5){\psi} \\
& & \\
& & \Symn(\shfF),}
\end{align*} 
and so the composite morphism $\psi \circ \vp$ satisfies $? \circ \pi_{\shfF}^{n} = \pi_{\shfF}^{n}$. But since we also have obviously $1_{\Symn(\shfF)} \circ \pi_{\shfF}^{n} = \pi_{\shfF}^{n}$, it follows that $\psi \circ \vp = 1_{\Symn(\shfF)}$, by the universal property of $\Symn(\shfF)$. One can show similarly that $\vp \circ \psi = 1_{\Symnp(\shfF)^{\dg}}$. Hence $\vp$ and $\psi$ are mutually inverse isomorphisms.
\end{prooff}

\end{subappendices}

\newpage

\bibliographystyle{acm}	
\bibliography{/Users/enrightward/Desktop/DISS_NEWFONT/masterbib_NF.bib}

\end{document}